\documentclass[twoside, 10pt]{article}

\usepackage[hmargin= 0.8in, twosideshift=0in,height=9.3in]{geometry}
\usepackage{amssymb,amsthm, amsfonts,amsmath,mathrsfs, amsxtra, url, hyperref}
\usepackage{delarray,verbatim}
\usepackage{natbib}
\renewcommand{\cite}{\citet}

\usepackage{fancyhdr}
\fancypagestyle{plain}{ %
\fancyhf{} 
}

\newtheorem {thm}{Theorem}[section]
\newtheorem {prop}{Proposition}[section]
\newtheorem {lemm}{Lemma}[section]
\newtheorem {deff}{Definition}[section]
\newtheorem {cor}{Corollary}[section]

\newtheorem {rem}{Remark}[section]

\newtheorem {eg}{Example}[section]

\def\ba{\begin{array}}
\def\ea{\end{array}}
\def\beq{\begin{equation}}
\def\bes{\begin{equation*}}
\def\ees{\end{equation*}}
\def\bea{\begin{eqnarray}}
\def\eea{\end{eqnarray}}
\def\beas{\begin{eqnarray*}}
\def\eeas{\end{eqnarray*}}
\def\bi{\begin{itemize}}
\def\ei{\end{itemize}}

\def\a{\alpha}
\def\g{\gamma}
\def\d{\delta}
\def\e{\varepsilon}
\def\z{\zeta}
\def\k{\kappa}
\def\l{\lambda}

\def\si{\sigma}

\def\t{\tau}

\def\o{\omega}

\def\vf{\varphi}

\def\D{\Delta}
\def\G{\Gamma}
\def\L{\Lambda}
\def\O{\Omega}

\def\P{\Psi}

\def\Th{\Theta}

\def\bE{{\bf E}}
\def\bF{{\bf F}}

\def\cD{{\cal D}}
\def\cE{{\cal E}}
\def\cF{{\cal F}}
\def\cG{{\cal G}}
\def\cH{{\cal H}}
\def\cI{{\cal I}}
\def\cJ{{\cal J}}

\def\cM{{\cal M}}
\def\cN{{\cal N}}

\def\cP{{\cal P}}

\def\cS{{\cal S}}

\def\cU{{\cal U}}

\def\hB{\mathbb{B}}
\def\hC{\mathbb{C}}

\def\hN{\mathbb{N}}

\def\hR{\mathbb{R}}

\def\hZ{\mathbb{Z}}

\def\sB{\mathscr{B}}

\def\sD{\mathscr{D}}
\def\sE{\mathscr{E}}

\def\sG{\mathscr{G}}
\def\sH{\mathscr{H}}

\def\sP{\mathscr{P}}

\def\fS{\mathfrak{S}}

\def\fU{\mathfrak{U}}

\def\({\textnormal{(}}
\def\){\textnormal{)}}
\def\[{[\neg[}
\def\]{]\neg]}

\def\no{\noindent}

\def\ss{\smallskip}
\def\ms{\medskip}

\def\q{\quad}
\def\qq{\qquad}
\def\nts{\negthinspace}
\def\neg{\negthinspace}
\def\dneg{\neg \neg}
\def\tneg{\neg \neg \neg}

\def\fa{\,\forall \,}
\def\pa{\partial}

\def\es{\emptyset}

\def\ol{\overline}
\def\ul{\underline}
\def\ua{\mathop{\uparrow}}
\def\da{\mathop{\downarrow}}

\def\={=\nts \nts=\nts \nts=\nts \nts=}
\def\lan{\langle}
\def\ran{\rangle}
\def\({\textnormal{(}}
\def\){\textnormal{)}}
\def\wt{\widetilde}

\def\cd{\cdot}
\def\cds{\cdots}

\def\essinf{\mathop{\rm essinf}}
\def\esssup{\mathop{\rm esssup}}
\def\liminf{\mathop{\ul{\rm lim}}}
\def\limsup{\mathop{\ol{\rm lim}}}

\def\dtp{{\hbox{$dt \times dP$-a.s.}}}

\def\qed{\hfill $\Box$ \medskip}   
\def\dfnn{\stackrel{\triangle}{=}}

\def\b1{{\bf 1}}

\def\hb{\hbox}
\def\dis{\displaystyle}

\newcommand{\lmt}[1]{ \underset{#1}{\lim}}
\newcommand{\lmtu}[1]{ \underset{#1}{\lim} \neg \ua}

\begin{document}

\title{\bf   Optimal Stopping for Non-linear Expectations }

\author{
Erhan Bayraktar,\thanks{ \noindent Department of
  Mathematics, University of Michigan, Ann Arbor, MI 48109; email:
erhan@umich.edu. E. Bayraktar is supported in part by the National Science Foundation under an applied mathematics research grant and a Careergrant, DMS-0906257 and DMS-0955463, respectively, and in part by the Susan M. Smith Professorship. 
.}\q Song Yao\thanks{
\noindent Department of
  Mathematics, University of Michigan, Ann Arbor, MI 48109; email: songyao@umich.edu. } }

\date{ }

\maketitle

 \begin{abstract}

We develop a theory for solving continuous time optimal stopping
problems for non-linear expectations. Our motivation is to consider
problems in which the stopper uses risk measures to evaluate future
rewards.
\end{abstract}

 \ms   {\bf Keywords: }\: Nonlinear expectations,
  Optimal stopping, Snell envelope, Stability, $g$-expectations.


\section{Introduction}
We solve continuous time optimal stopping problems in which the
reward is evaluated using \emph{non-linear} expectations. Our
purpose is to use criteria other than
  the expected value to evaluate the present value of future rewards.
 Such criteria include \emph{risk measures}, which are not necessarily linear.
Given a filtered probability space $(\O,\cF, P, \bF=\{\cF_t\}_{t \in
[0,T]})$ satisfying the \emph{usual assumptions}, we define a
filtration-consistent non-linear expectation ($\bF$-expectation for
short)  with domain $\Lambda$ as a collection of operators
$\big\{\cE[\cd|\cF_t]: \L \mapsto \L_t\dfnn \L \cap L^0(\cF_t)
\big\}_{t \in [0,T]}$ satisfying ``Monotonicity",
``Time-Consistency", ``Zero-one Law" and ``Translation-Invariance".
This definition is similar to the one proposed in \cite{Pln}. A
notable example of an $\bF$-expectation is the so-called
\emph{g-expectation}, introduced by \cite{Peng-97}. A fairly large
class of \emph{convex risk measures} (see e.g.
\cite{Follmer_Schied_2004} for the definition of risk measures) are
$g$-expectations (see \cite{CHMP}, \cite{Pln}, \cite{MaYao_08} and
\cite{HMPY-07}).

 \ss We consider two optimal stopping problems. In the first one, the
stopper aims to find an optimal stopping time when there are
multiple priors and the {\it Nature} is in cooperation with the
stopper; i.e., the stopper finds an optimal stopping time that
attains
  \bea \label{eq:game-coop}
  Z(0) \dfnn \underset{(i,\rho) \in \cI \times \cS_{0,T}}{\sup}
 \; \cE_i[Y_{\rho}+H_{\rho}^i|\cF_0],
  \eea
in which $\sE =\{\cE_i \}_{i \in \cI}$ is a \emph{stable} class of
$\bF$-expectations, $\cS_{0,T}$ is the set of stopping times that
take value in $[0,T]$. The reward process $Y$ is a right-continuous
$\bF$-adapted process and for any $\nu \in \cS_{0,T}$, $Y_\nu$
belongs to $\L^\# \dfnn \{\xi \in \L\;|\; \xi \ge c, ~ a.s. \hb{ for
some }c \in \hR\}$, where $\Lambda$ is the common domain of the
elements in $\sE$. On the other hand, the \emph{model-dependent}
reward processes $\{H^i\}_{i \in \hN}$ is a family of
right-continuous adapted processes with $H^i_0=0$ that is
\emph{consistent} with $\sE$. We will express the solution of this
problem in terms of the $\sE$-\emph{upper Snell envelope} $Z^0$ of
$Y_t$, the smallest RCLL $ \bF$-adapted process dominating $Y$ such
that $Z^{i,0}\dfnn \{ Z^0_t+H_t^i \}_{t \in [0,T]}$
  is an $\wt{\cE}_i$-supermartingale for each $i \in \cI$.

 \ss The construction of the Snell envelope is not straightforward.
First, for any $i\in \cI $, the conditional expectation $\cE_i[\xi
|\cF_{\nu}]$, $\xi \in \L$ and $\nu \in  \cS_{0,T}$ may not be well
defined. However, we show that $t \to \cE_i[\xi |\cF_{t}]$ admits a
right-continuous modification $t \to \widetilde{\cE}_i[\xi|\cF_\cd]$
for any $\xi \in \L^\#$ and that $\widetilde{\cE}_i$ is itself an
$\bF$-expectation on $\L^\#$ such that $\widetilde{\cE}_i[\cd
|\cF_\nu]$ is well defined on $\L^\#$ for any $\nu \in \cS_{0,T}$.
In terms of $\widetilde{\cE}_i$ we have that \bea
\label{eq:game-coop-rep}
  Z(0)= \underset{(i,\rho) \in \cI \times \cS_{0,T}}{\sup}
 \; \widetilde{\cE}_i[Y_{\rho}+H_{\rho}^i|\cF_0].
  \eea
Finding a RCLL modification requires the development of an
upcrossing theorem. This theorem relies on the strict monotonicity
of $\cE_i$ and other mild hypotheses, one of which is equivalent to
having lower semi-continuity (i.e. Fatou's lemma). Thanks to the
right continuity of $t \to \widetilde{\cE}_i[\xi|\cF_t]$, we also
have an optional sampling theorem for right-continuous
$\widetilde{\cE}_i$-supermartingales. Another important tool in
finding an optimal stopping time, the dominated convergence theorem
is also developed under another mild assumption.

\ss The stability assumption we make on the family $\sE$ is another
essential ingredient in the construction of the Snell envelope. It
guarantees that the class $\sE$ is closed under \emph{pasting}:  for
any $i,j \in \cI$ and $\nu \in \cS_{0,T}$ there exists a $k \in \cI$
such that
$\widetilde{\cE}_k[\xi|\cF_{\sigma}]=\widetilde{\cE}_i\big[
\widetilde{\cE}_j [\xi|\cF_{\nu \vee \si}] \big|\cF_\si\big]$, for
any $\sigma \in \cS_{0,T}$. Under this assumption it can then be
seen, for example, that
 the collection of random variables $\left\{\widetilde{\cE}_i\left[X(\rho)
 +H_{\rho}^i-H_{\nu}^i\Big|\cF_\nu\right], \,( i,
\rho) \in \cI \times \cS_{\nu,T}\right\}$ is directed upwards. When
the constituents of $\sE$ are linear expectations, the notion of
stability of this collection is given by \cite[Definition
6.44]{Follmer_Schied_2004}, who showed that pasting two probability
measures equivalent to $P$ at a stopping time one will result in
another probability measure equivalent to $P$. Our result in
Proposition~\ref{tau_A}
 shows that we have the same pasting property for $\bF$-expectations.
 As we shall see, the stability  assumption is  crucial in showing that
 the Snell envelope is a supermartingale. This property of the Snell envelope is a generalization of
 \emph{time consistency}, i.e.,
\bea \label{eqn-uuu03}
 \underset{i  \in \cI  }{\esssup}\, \widetilde{\cE}_i [\xi|\cF_\nu]
 =\underset{i  \in \cI  }{\esssup}\, \widetilde{\cE}_i
 \left[\underset{i  \in \cI  }{\esssup}\, \widetilde{\cE}_i
 [\xi|\cF_\sigma]\bigg|\cF_\nu\right],\q a.s., \q \fa \nu, \si \in \cS_{0,T}
 \hb{ with }  \nu \le \si,~a.s.
 \eea
 \cite[Theorem 12]{Delbaen_2006} showed in the linear expectations case
 that the time consistency \eqref{eqn-uuu03} is equivalent to the
 stability.

 \ss When the reward $t \to Y_t+H_t^i$ is
``$\sE$-uniformly-left-continuous" and each non-linear expectation
in $\sE$ is convex, we can find an optimal stopping time
$\ol{\t}(0)$ for \eqref{eq:game-coop} in terms of the Snell
envelope. Then we can solve the problem \bea
\label{eq:single}
 \underset{\rho \in \cS_{0,T}}{\sup}
 \; \cE_i[Y_\rho+H^i_\rho |\cF_0],
\eea when $\cE_i[\cdot|\cF_t]$ has among other properties strict
monotonicity, lower semi-continuity, dominated convergence theorem
and the upcrossing lemma. Note that although, $\underset{i \in \cI
}{\esssup}\, \widetilde{\cE}_i [\cdot|\cF_t]$ has similar properties
to $\widetilde{\cE}_i[\cdot | \cF_t]$ (and that might lead one to
think that \eqref{eq:game-coop} can actually be considered as a
special case of \eqref{eq:single}), the former does not satisfy
strict monotonicity, the upcrossing lemma, and the dominated
convergence theorem. One motivation for considering optimal stopping
with multiple priors is to solve optimal stopping problems for
``non-linear expectations" which do not satisfy these properties.

 \ss We show that the collection of $g$-expectations with uniformly
Lipschitz generators satisfy the  uniform left continuity
assumption. Moreover, a $g$-expectation satisfies all the
assumptions we ask of each $\cE_i$ for the upcrossing theorem,
Fatou's lemma and the dominated convergence theorem to hold; and
pasting of $g$-expectations results in another $g$-expectation. As a
result the case of $g$-expectations presents a non-conventional
example in which we can determine an optimal stopping time for
\eqref{eq:game-coop}. In fact, in the $g$-expectation example we can
even find an optimal \emph{prior} $i^* \in \cI$, i.e.,
 \bea \label{eq:opt-cont-st}
  Z(0)= \cE_{i^*}[Y_{\ol{\t}(0)}+H_{\ol{\t}(0)}^{i^*}|\cF_0].
  \eea

 \ss In the second problem, the
\emph{stopper} tries to find a robust optimal stopping time that
attains
 \bea \label{eq:robust}
V(0) \dfnn \underset{\rho \in \cS_{0,T}}{\sup}
 \; \underset{i \in \cI}{\inf}\;\cE_i \big[Y_\rho \neg
 +\dneg H^i_\rho \big|\cF_0\big].
  \eea
 Under the ``$\sE$-uniform-right-continuity" assumption, we find an
optimal stopping time in terms of the $\sE$-\emph{lower Snell
envelope}. An immediate by-product is the following minimax theorem
\bea V(0)=
  \underset{i \in \cI}{\inf}  \; \underset{\rho \in \cS_{0,T}}{\sup} \cE_i\big[Y_\rho \neg
 +\dneg H^i_\rho \big|\cF_0\big].
\eea

 \ss Our paper was inspired by \cite{Kara_Zam_2006} and
\cite{Kara_Zam_2008}, which developed a martingale approach to
solving \eqref{eq:game-coop} and \eqref{eq:robust}, when $\sE$ is a
class of linear expectations. In particular, \cite{Kara_Zam_2006}
considered the \emph{controller-stopper} problem
 \bea \sup_{\rho \in
\cS_{0,T}}\sup_{U \in
\fU}\bE^{u}\left[g\big(X(\rho)\big)+\int_0^{\rho}h(s,X,U_s)ds\right],
 \eea
  where $X(t)=x+\int_0^{t}f(s,X,U_s)ds+\int_0^{t} \sigma(s,X)
dW^{U}_s$. In this problem, the stability condition is automatically
satisfied. Here, $g$ and $h$ are assumed to be bounded measurable
functions. Our results on $g$-expectations extend the results of
\cite{Kara_Zam_2006}  from bounded rewards to rewards satisfying
linear growth. \cite{Delbaen_2006}, \cite{Kara_Zam_2005} also
considered \eqref{eq:game-coop} when the $\cE_i$'s are linear
expectations. The latter paper made a \emph{convexity} assumption on
the collection of equivalent probability measures instead of a
stability assumption. On the other hand, the discrete time version
of the robust optimization problem was analyzed by
\cite{Follmer_Schied_2004}. Also see \cite[Sections 5.2 and
5.3]{CDK-2006}.


 \ss  The rest of the paper is organized as follows: In Section~\ref{sec:notation}
 we will introduce some notations that will be used throughout the paper.
 In Section~\ref{ch_2}, we define what we mean by an $\bF$-expectation $\cE$, propose some basic
hypotheses on $\cE$ and discuss their implications such as Fatou's
lemma, dominated convergence theorem and upcrossing lemma. We show
that $t \to \cE[\cdot|\cF_t]$ admits a right-continuous modification
which is also an $\bF$-expectation and satisfies Fatou's lemma and
the dominated convergence theorem. This step is essential since
$\cE[\cdot| \cF_{\nu}],~\nu \in \cS_{0,T}$ may not be well defined.
We also show that the optional sampling theorem holds. The results
in Section~\ref{ch_2} will be the backbone of our analysis in the
later sections.

 \ss    In Section~\ref{ch_3} we introduce the stable class of
$\bF$-expectations and review the properties of essential extremum.
In Section~\ref{co_game} we solve \eqref{eq:game-coop-rep}, and find
an optimal stopping time in terms of the $\sE$-\emph{upper Snell
envelope}. On the other hand, in Section~\ref{nonco_game} we solve
the robust optimization problem \eqref{eq:robust} in terms of the
$\sE$-\emph{lower Snell envelope}. In Section~\ref{section_remark},
we give some interpretations and remarks on our results in the
previous sections. In Section~\ref{ch_app} we consider the case when
$\sE$ is a certain collection of $g$-expectations. We see that in
this framework, our assumptions on each $\cE_i$, the stability
condition and the uniform left/right continuity conditions are
naturally satisfied. We also determine an optimal prior $i^* \in
\cI$ satisfying \eqref{eq:opt-cont-st}. Moreover, we show how the
controller and stopper problem of \cite{Kara_Zam_2006} fits into our
g-expectations framework. This lets us extend their result from
bounded rewards to rewards satisfying linear growth. In this
section, we also solve optimal stopping problem for quadratic
$g$-expectations. The proofs of our results are presented in
Section~\ref{sec:Proofs}.

\subsection{Notation}\label{sec:notation}
Throughout this paper, we fix a finite time horizon $T>0$ and
consider a complete probability space $(\O,\cF, P)$ equipped with a
right continuous filtration $\bF \dfnn \{\cF_t\}_{t \in [0,T]}$, not
necessarily Brownian one, such that $\cF_0 $ is generated by all
$P$-null sets in $\cF$ (in fact, $\cF_0$ collects all measurable
sets with probability $0$ or $1$).
 Let $\cS_{0,T}$ be the collection of all $\bF$-stopping times $\nu$ such that $0 \le \nu \le T$,
a.s. For any $ \nu, \si \in \cS_{0,T}$ with $\nu \le \si$, a.s., we
define $\cS_{\nu,\si} \dfnn \{ \rho \in \cS_{0,T}\, |\; \nu \le \rho
\le \si, ~a.s.\}$ and let $\cS^F_{\nu,\si}$ denote all finite-valued
stopping times in $\cS_{\nu,\si}$.
 We let $\cD=\{k2^{-n} \;|\; k \in \hZ,\, n \in \hN\}$ denote the set of all
dyadic rational numbers and set $\cD_T \dfnn \big([0,T)\cap \cD
\big) \cup \{T\} $. For any $t \in [0,T]$ and $n \in \hN$, we also
define
 \bea \label{dyadic}
 q^+_n(t) \dfnn \frac{\lceil 2^nt \rceil}{2^n}\land T.
 \eea
It is clear that  $ q^+_n(t)  \in  \cD_T$.

 \ss    In what follows we let $\cF'$ be a generic sub-$\si$-field of $\cF $ and let
$\hB$ be a generic Banach space with norm $|\cd|_{\hB}$. The
following spaces of functions will be used in the sequel.

 \ss \no (1) For $0\le p\le\infty$, we define
 \bi
 \item $L^p(\cF';\hB)$ to be the space of all $\hB$-valued,
$\cF'$-measurable random variables $\xi$ such that
$E(|\xi|^p_{\hB})<\infty$. In particular, if $p=0$,
  $L^0(\cF';\hB)$ stands for the space of all $\hB$-valued,
$\cF'$-measurable random variables; and if $p=\infty$,
$L^\infty(\cF';\hB)$ denotes the space of all $\hB$-valued,
$\cF'$-measurable random variables $\xi$ with $\|\xi\|_\infty \dfnn
\underset{\o \in \O}{\esssup}\,|\xi(\o)|_{\hB}<\infty$.

\item  $L^p_\bF([0,T];\hB)$ to be the space of all $\hB$-valued,
$\bF$-adapted processes $X$ such that $ 
E \neg\int_0^T \neg |X_t|^p_{\hB} dt <\infty$. In particular, if
$p=0$,   $L^0_\bF([0,T];\hB)$ stands for the space of all
$\hB$-valued, $\bF$-adapted processes; and if $p=\infty$,
$L^\infty_\bF([0,T];\hB)$ denotes the space of all $\hB$-valued,
$\bF$-adapted processes $X$ 
with $\|X\|_\infty \dfnn \underset{(t,\o) \in [0,T] \times \O}
{\esssup}\, |X_t(\o)|_{\hB}<\infty$.

\item $\hC^p_\bF([0,T];\hB) \dfnn \{X \in L^p_\bF([0,T];\hB) :\, \hb{$X$ has continuous
paths}\}$.

\item $\cH^p_\bF([0,T];\hB) \dfnn \{X \in L^p_\bF([0,T];\hB) :\, \hb{$X$ is predictably
measurable}\}$.
 \ei

 \ss \no (2) For $p \ge 1$, we define a Banach space
  \beas
  M^p_\bF([0,T];\hB)=\bigg\{X \in \cH^0_\bF([0,T];\hB) :\,
  \|X \|_{M^p} \dfnn \bigg\{E\Big[\big(\int_0^T |X_s|^2_{\hB}ds \big)^{p/2}
  \Big]\bigg\}^{ 1/p} < \infty\bigg\},
  \eeas
 and denote $ M_\bF([0,T];\hB)   \dfnn    \underset{p \ge
 1}{\cap}M^p_\bF([0,T];\hB)$.

 \ss \no (3) We further define
\beas
  L^{\neg e}(\cF';\hB) &\dfnn&  \Big\{\xi \in L^0(\cF';\hB)
  :\, E\big[e^{\l |\xi|_{\hB} }\big]<\infty \hb{ for all } \l > 0 \Big\},     \\
 \hC^e_\bF([0,T];\hB)  &\dfnn&   \Big\{ X \in \hC^0_\bF([0,T];\hB) :\,
    E\Big[\exp\big\{\l \underset{t \in [0,T]}{\sup}|X_t|_{\hB}
\big\}\Big]<\infty \hb{ for all } \l > 0 \Big\}.
  \eeas

 \if{0}
 For any $p>0$, we denote $\cM^p(\hR^d)$ to be the space of
all $\hR^d$-valued predictable processes $X$ such that
 \bea
 \label{Mp}
 \|X \|_{\cM^p} \dfnn \Big(E\big(\int_0^T |X_s|^2ds \big)^{p/2}\Big)^{1 \land 1/p}
  < \infty.
 \eea
We note that for $p \geq 1$, $\cM^p(\hR^d)$ 
is a Banach space with the norm  $\|\cd\|_{\cM^p}$, 
and for $p \in (0,1)$, $\cM^p(\hR^d)$ 
is a complete metric space with the distance defined through
(\ref{Mp}).
 \fi

 \no If $d=1$, we shall drop $\hB=\hR$ from the above notations
\big(e.g., $L^p_\bF([0,T])=L^p_\bF([0,T];\hR)$,
$L^p({\cF_T})=L^p({\cF_T};\hR)$\big). In this paper, all
$\bF$-adapted processes are supposed to be real-valued unless
specifying otherwise.

 \section{$\bF$-expectations and Their Properties} \label{ch_2}
\setcounter{equation}{0}

 \ms

 We will define non-linear expectations on subspaces of
 $L^0(\cF_T)$ satisfying certain algebraic properties, which are
 listed in the definition below.

\ms
\begin{deff}
Let $\sD_T$ denote the collection of all non-empty subsets $\L$ of
$L^0(\cF_T)$ satisfying:
 \bi
 \item[\(D1\)] $0,1 \in \L $; 
\item[\(D2\)] $\L$ is closed under addition and under multiplication with
indicator random variables. 
 Namely, for any $\xi,\eta \in \L$ and
 $A \in \cF_T$, both $\xi+\eta$ and $\b1_A \xi $ belong to $\L$;
 \item[\(D3\)] $\L$ is positively solid: For any $\xi ,\eta \in L^0(\cF_T)$ with $0 \le \xi \le \eta $, a.s.,
 if $\eta \in \L$, then $\xi \in \L$ as well.
 \ei
\end{deff}

\begin{rem}\label{domain}
\no (1)  Each $\L \in \sD_T$ is also closed under maximization
``\,$\vee$" and under minimization
 ``\,$\land$": In fact, for any $\xi,\eta \in \L$, since the set $\{\xi > \eta\} \in \cF_T$,
 (D2) implies that $ \xi \vee \eta =\xi \b1_{\{\xi > \eta\}}+\eta \b1_{\{\xi \le
\eta\}} \in \L $. Similarly, $ \xi \land \eta \in \L$;

\ms \no (2) For each $\L \in \sD_T$, (D1)-(D3) imply that $c \in \L$
for any $c \ge 0$;

 \ms  \no (3)
$\sD_T$ is closed under intersections:
 If $\{\L_i\}_{i \in \cI}$ is a subset of $\sD_T$, then $ \dis \underset{i \in I}{\cap} \L_i \in \sD_T$;
  $\sD_T$ is closed under unions of increasing sequences: If
$\{\L_n\}_{n \in \hN} \subset \sD_T$ such that $\L_n \subset
\L_{n+1}$ for any $n \in \hN$, then $ \underset{n \in \hN}{\cup}
\L_n \in \sD_T$;

\ms \no (4) It is clear that $L^p(\cF_T) \in \sD_T$ for all $0 \le p
\le \infty$.

\end{rem}

 \begin{deff} An $\bF$-consistent non-linear expectation \($\bF$-expectation for short\)
 is a pair \($\cE$, $\L$\) in which
 $\L \in \sD_T$ and $\cE$ denotes a family of operators
$\big\{\cE[\cd|\cF_t]: \L \mapsto \L_t\dfnn \L \cap L^0(\cF_t)
\big\}_{t \in [0,T]}$ satisfying the following hypothesis for any $\xi,\eta \in \L$ and $t
\in [0,T]$:
  \bi
 \item[ {\bf (A1)} ] {\it ``Monotonicity (positively strict)":}\;  $ \cE [\xi |\cF_t] \le  \cE [\eta |\cF_t] $,
$a.s.$ if $\xi \le \eta$, a.s.; Moreover, if $0 \le \xi \le \eta$
a.s. and $ \cE [\xi|\cF_0] = \cE [\eta|\cF_0]$, then $\xi=\eta$,
a.s.;

 \item[ {\bf (A2)} ] {\it ``Time Consistency":}\; $\cE\big[\cE[\xi|\cF_t]\big|\cF_s \big]=\cE[\xi|\cF_s]$,
a.s. for any $0 \le s\le t \le T$;

\item[ {\bf (A3)} ] {\it ``Zero-one Law":}\;
$\cE [\b1_A \xi |\cF_t]=\b1_A \cE [\xi |\cF_t]$, $a.s.$ for any $A
\in \cF_t$;

\item[ {\bf (A4)} ] {\it ``Translation Invariance":}\;
$ \cE [\xi+\eta |\cF_t]=\cE [\xi |\cF_t]+\eta$, a.s. if $\eta \in
\L_t$.
 \ei
\end{deff}

\ss We denote the domain $\L$ by $Dom(\cE)$ and define
 \beas
 Dom_\nu(\cE) \dfnn Dom(\cE)\cap L^0(\cF_\nu),\q \fa  \nu \in
 \cS_{0,T}.
 \eeas
For any $\xi, \eta \in Dom(\cE)$ with $\xi=\eta$, a.s., (A1) implies
that $\cE[\xi|\cF_t]=\cE[\eta|\cF_t]$, a.s. for any $t \in [0,T]$,
which shows that the $\bF$-expectation ($\cE$, $Dom(\cE)$) is
well-defined. Moreover, since $Dom_0(\cE)  = Dom(\cE) \cap
L^0(\cF_0) \subset L^0(\cF_0)=\hR$, $\cE[\cd|\cF_0]$ is a
real-valued function on $Dom(\cE)$. In the rest of the paper, we
will substitute $\cE[\cd]$ for $\cE[\cd|\cF_0]$.

 \begin{rem}
Our definition of $\bF$-expectations is similar to that of
$\cF^X_t$-consistent non-linear expectations introduced in
\cite[page 4]{Pln}.
 \end{rem}

\begin{eg} \label{3eg}
The following pairs satisfy (A1)-(A4); thus they are
$\bF$-expectations:

\ms \no (1)  $\big(\{E[\cd|\cF_t]\}_{t \in [0,T]}, L^1(\cF_T)\big)$:
the linear expectation $E$ is a special $\bF$-expectation with
domain $L^1(\cF_T)$;

 \ms \no (2) $\big(\{\cE_g[\cd|\cF_t]\}_{t \in [0,T]}, L^2(\cF_T)\big)$: the $g$-expectation
 with generator $g(t,z)$ Lipschitz in $z$ (see \cite{Peng-97}, \cite{CHMP}
 or Subsection \ref{subsection_gexp} of the present
 paper);

 \ms \no (3) $\big(\{\cE_g[\cd|\cF_t]\}_{t \in [0,T]}, L^{\neg
e}(\cF_T) \big)$: the
 $g$-expectation with generator $g(t,z)$ having quadratic growth in
 $z$ (see Subsection \ref{subs_g2_exp} of this paper).
\end{eg}

\ss  $\bF$-expectations can alternatively be introduced in a
more classical way:

\begin{prop} \label{prop_Fexp}
Let $\cE^o: \L \mapsto \hR$ be a mapping on some $\L \in \sD_T$
satisfying:

\ss \no \(a1\)  \;  For any $\xi,\eta \in \L$ with $\xi \le \eta$,
$a.s.$, we have $\cE^o[\xi] \le \cE^o[\eta]$. Moreover, if
$\cE^o[\xi] = \cE^o[\eta]$, then $\xi=\eta$, $a.s.$;

\ss \no \(a2\) \; For any $\xi \in \L$ and $t \in [0,T]$, there
exists a unique random variable $\xi_t \in \L_t $ such that
 $\cE^o[\b1_A  \xi+\g]=\cE^o\big[\b1_A  \xi_t+\g \big]$ holds
 for any $A \in \cF_t$ and $\g \in \L_t $.

\no Then $\{\cE^o[\xi|\cF_t]\dfnn \xi_t,\, \xi \in \L\}_{t \in [0,T]
}$ defines an $\bF$-expectation with domain $\L$.
\end{prop}

\begin{rem}
For a mapping $\cE^o$ on some $\L \in \sD_T$ satisfying (a1) and
(a2), the implied operator $\cE^o[\cd |\cF_0]$ is also from $\L$ to
$\hR$, which, however, may not equal to $\cE^o $. In fact, one can
only deduce that $\cE^o[\xi]=\cE^o\big[\cE^o[\xi|\cF_0]\big]$ for any
$ \xi \in \L$.
\end{rem}

  From now on, when we say an $\bF$-expectation $\cE$, we
will refer to the pair $\big(\cE, Dom(\cE)\big)$.
 Besides (A1)-(A4), the $\bF$-expectation $\cE$ has the following properties:

\begin{prop} \label{3addition}
For any $\xi,\eta \in Dom(\cE)$ and $t \in [0,T]$, we have

\ss \no (1) \; {\it ``Local Property":}\; $\cE[\b1_A \xi+
\b1_{A^c}\eta|\cF_t]=\b1_A\cE[ \xi|\cF_t]+ \b1_{A^c}\cE[ \eta|\cF_t]
$, a.s. for any $A \in \cF_t$;

\ss \no (2)\; {\it ``Constant-Preserving":}\; $\cE [\xi
|\cF_t]=\xi$, a.s. if $\,\xi \in Dom_t(\cE)$;

\ss \no (3)\; {\it ``Comparison":}\; Let $\xi, \eta \in L^0(\cF_\nu)
$ for some $\nu \in \cS_{0,T}$. If $\,\eta \ge c$, a.s. for some $c
\in \hR$, then $\xi \le (\hb{or } =)\; \eta$, a.s. if and only if
$\,\cE[\b1_A \xi] \le (\hb{or } =)\; \cE[ \b1_A \eta]$ for all $A
\in \cF_\nu$.
\end{prop}

\ss The following two subsets of $Dom(\cE)$ will be of interest:
 \bea \label{eq:dms}
 \q   Dom^+(\cE) \dfnn \{\xi \in
Dom(\cE)   :\, \xi \ge 0, \; a.s. \},~\;  Dom^\#(\cE)\dfnn \{\xi \in
  Dom(\cE)  :\, \xi \ge c, \;a.s. \hb{~for some }c=c(\xi) \in \hR\}
  .\q
 \eea

\begin{rem}
The restrictions of $\cE$ on $Dom^+(\cE)$ and on $Dom^\#(\cE)$,
namely $\big(\cE, Dom^+(\cE) \big)$ and $\big(\cE, Dom^\#(\cE)
\big)$ respectively, are both $\bF$-expectations: To see this, first
note that $Dom^+(\cE)$ and $Dom^\#(\cE)$ both belong to $\sD_T$. For
any $t \in [0,T]$, (A1) and Proposition \ref{3addition} (2) imply
that for any $\xi \in Dom^\#(\cE)$
 \beas
     \cE[\xi|\cF_t] \ge \cE\big[c(\xi)\big|\cF_t\big]=c(\xi), \q a.s., \q
     \hb{thus} \q \cE[\xi|\cF_t] \in Dom^\#(\cE),
 \eeas
 which shows that $\cE[\cd|\cF_t]$ maps $Dom^\#(\cE)$ into $Dom^\#(\cE) \cap L^0(\cF_t)$.
 Then it is easy to check that the restriction of $\cE=\big\{\cE[\cd|\cF_t]\big\}_{t\in[0,T]}$ on $Dom^\#(\cE)$
 satisfies (A1) to (A4), thus it is an $\bF$-expectation. Similarly, $\big(\cE, Dom^+(\cE)
 \big)$ is also an $\bF$-expectation.
 \if{0}
 Moreover, the restrictions of
 $\cE$ on $Dom^-(\cE) \dfnn \{\xi \in
Dom(\cE) \;|\; \xi \le 0, \; a.s. \}$ and on $
Dom_\natural(\cE)\dfnn \{\xi \in
  Dom(\cE)\;|\; \xi \le c, \;a.s. \hb{~for some }c=c(\xi) \in \hR\}$ are $\bF$-expectations as well.
 \fi

\ss We should remark that restricting $\cE$ on any subset $\L'$ of
$Dom(\cE)$, with $\L' \in \sD_T$, may not result in an
$\bF$-expectation, i.e. $(\cE, \Lambda')$ may not be an
$\bF$-expectation.
\end{rem}

 \begin{deff}
(1)  An $\bF$-adapted process $X=\{X_t\}_{t \in [0,T]}$ is called an
  ``$\cE$-process" if $ X_t \in Dom(\cE)$ for
any $t \in [0,T]$;

\ss \no (2) An $\cE$-process $X$ is said to be an
$\cE$-supermartingale
 (resp.
$\cE$-martingale, $\cE$-submartingale) if for any $0\le s < t \le
T$, $\,\cE[X_t|\cF_s] \le (\hb{resp.}\,=,\; \ge)~X_s$, a.s.
 \end{deff}

  Given a $\nu \in \cS^F_{0,T}$ taking values in a finite set
$\{t_1 < \cds< t_n\}$, if $X$ is an $\cE$-process, (D2) implies that
$ X_\nu = \sum^n_{i=1} \b1_{\{\nu =t_i\}} X_{t_i} \in Dom(\cE)$,
thus $X_\nu \in Dom_\nu(\cE)$. 
 Since $\big\{X^\xi_t\dfnn \cE[\xi|\cF_t]\big\}_{t \in [0,T]}$ is an
$\cE$-process for any $\xi \in Dom(\cE)$, we can define an operator
$\cE[\cd|\cF_\nu]$ from $Dom(\cE)$ to $Dom_\nu(\cE)$ by
 \beas
  \cE[\xi|\cF_\nu] \dfnn X^\xi_\nu , \q \hb{for any }~ \xi \in
  Dom(\cE),
 \eeas
 which allows us to state a basic Optional Sampling Theorem for $\cE$.

\begin{prop} \label{op_sa} (Optional Sampling Theorem)
Let $X$ be an $\cE$-supermartingale  (resp. $\cE$-martingale,
$\cE$-submartingale). Then for any $\nu, \si \in \cS^F_{0,T}$, $ \cE
[X_\nu |\cF_\si] \le (\hb{resp.}\,=,\; \ge)\, X_{\nu \land \si}$,
a.s.
\end{prop}
 In particular, 
  applying Proposition \ref{op_sa} to each
$\cE$-martingale $\{\cE[\xi|\cF_t]\}_{t \in
 [0,T]}$, in which $\xi \in Dom(\cE)$, yields the following result.
 \begin{cor}\label{cor_os}
 For any $\xi  \in Dom(\cE)$ and $\nu,\si \in \cS^F_{0,T}$,
 $\cE\big[\cE[\xi|\cF_\nu]\big|\cF_\si \big]=\cE[\xi|\cF_{\nu \land
 \si}]$, a.s.
 \end{cor}

\begin{rem} Corollary \ref{cor_os} extends the ``Time-Consistency" (A2) to the
case of finite-valued stopping times.
 \end{rem}

 \ms $\cE[\cd|\cF_\nu]$ inherits other properties of $\cE[\cd|\cF_t]$ as well:

\begin{prop} \label{properties_2}
For any $\xi, \eta  \in Dom(\cE)$ and $\nu \in \cS^F_{0,T}$, it
holds that
  \bi
 \item[(1)] {\it ``Monotonicity (positively strict)":}\;  $ \cE [\xi |\cF_\nu] \le  \cE [\eta |\cF_\nu] $,
$a.s.$ if $   \xi \le \eta$, $a.s.$; Moreover, if $0 \le \xi \le
\eta$, a.s. and $ \cE [\xi |\cF_\si] = \cE [\eta |\cF_\si]$, $a.s.$
for some $\si \in \cS^F_{0,T}$, then $\xi=\eta$, $a.s.$;

 \item[(2)] {\it ``Zero-one Law":}\;
$\cE [\b1_A \xi |\cF_\nu]=\b1_A \cE [\xi |\cF_\nu]$, a.s. for any $A
\in \cF_\nu$;

 \item[(3)] {\it ``Translation Invariance":}\; 
$ \cE [\xi+\eta |\cF_\nu]=\cE [\xi |\cF_\nu]+\eta$, a.s. if $ \eta
\in Dom_\nu(\cE) $;

 \item[(4)] {\it ``Local Property":}\; $\cE[\b1_A \xi+ \b1_{A^c}\eta|\cF_\nu]
  =\b1_A\cE[ \xi|\cF_\nu]+ \b1_{A^c}\cE[ \eta|\cF_\nu] $, $a.s.$ for any $A \in \cF_\nu$;

 \item[(5)] {\it ``Constant-Preserving":}\;  $\cE [\xi |\cF_\nu]=\xi$, $a.s.$, \hb{if} $\,\xi \in
  Dom_\nu(\cE) $.
 \ei
 \end{prop}

We make the following basic hypotheses on the
$\bF$-expectation $\cE$. These hypotheses will be essential in developing Fatou's lemma,
the Dominated Convergence Theorem and the Upcrossing Theorem.

\ms \no {\bf Hypotheses}
 \bi
\item[{\bf (H0)}] For any $A \in \cF_T$ with $P(A)>0$,
we have $\underset{n\to \infty}{\lim} \cE[n\b1_A]=\infty$;

\item[{\bf (H1)}] 
For any $\xi \in Dom^+(\cE)$ and any $\{A_n\}_{n \in \hN} \subset
 \cF_T$ with $\underset{n \to \infty}{\lim} \neg \ua \b1_{A_n}=1 $, a.s.,
we have $ \underset{n \to \infty}{\lim}  \neg \ua \cE[\b1_{A_n}\xi ]
= \cE[\xi]$;

\item[{\bf (H2)}] 
For any $\xi, \eta \in Dom^+(\cE)$ and any $\{A_n\}_{n \in \hN}
 \subset \cF_T$ with $\underset{n \to \infty}{\lim} \neg \da \b1_{A_n}=0 $, a.s.,
we have $\underset{n \to \infty}{\lim} \neg \da \cE[
\xi+\b1_{A_n}\eta] = \cE[\xi]$.
 \ei

\begin{rem}
The linear expectation $E$ on $L^1(\cF_T)$ clearly satisfies
(H0)-(H2). We will show that Lipschitz and quadratic
$g$-expectations also satisfy (H0)-(H2) in Propositions \ref{gexp}
and \ref{g2exp} respectively.
\end{rem}

\ss The $\bF$-expectation $\cE$ satisfies the following Fatou's
Lemma and Dominated Convergence Theorem.

\begin{thm}\label{fatou} (Fatou's Lemma) (H1) is equivalent to the lower semi-continuity of $\cE$:
If a sequence $\{\xi_n\}_{n \in \hN} \subset Dom^+(\cE)$ converges
a.s. to some $\xi \in Dom^+(\cE)$, then for any $\nu \in
\cS^F_{0,T}$, we have
 \bea  \label{eqn-b01}
 \cE[\xi|\cF_\nu] \le \underset{n \to \infty}{\liminf}
 \cE[\xi_n|\cF_\nu], \q a.s.,
 \eea
where the right hand side of (\ref{eqn-b01}) could be equal to infinity with
 non-zero probability.
\end{thm}

\begin{rem} \label{rem_fatou}
  In the case of the linear expectation $E$,
a converse to (\ref{eqn-b01}) holds: For any non-negative sequence
$\{\xi_n\}_{n \in \hN} \subset L^1(\cF_T)$ that converges a.s. to
some $\xi \in L^0(\cF_T)$, if $\underset{n \to \infty}{\liminf}
 E[\xi_n]<\infty$, then $ \xi  \in  L^1(\cF_T)$.
 However, this statement may not be the case for an arbitrary
 $\bF$-expectation. That is, $\underset{n \to \infty}{\liminf}
 \cE[\xi_n]<\infty$ may not imply that $\xi \in Dom^+(\cE)$ given that
 $  \{\xi_n\}_{n \in \hN} \subset Dom^+(\cE) $ is a sequence convergent a.s. to some $\xi \in L^0(\cF_T)$.
 \(See Example \ref{counter_eg_fatou} for a counterexample in the
case of
 a Lipschitz $g$-expectation.\)
 \end{rem}

 \begin{thm} \label{DCT} (Dominated Convergence Theorem)
 Assume (H1) and (H2) hold. Let $\{\xi_n\}_{n \in \hN}$ be a sequence in $
Dom^+(\cE)$ that converges a.s. If there is an $\eta \in Dom^+(\cE)$
such that $\xi_n \le \eta$ a.s. for any $n \in \hN$, then the limit
$\xi$ of
$\{\xi_n\}_{n \in \hN}$ belongs to 
$ Dom^+(\cE)$, and for any $\nu \in \cS^F_{0,T}$, we have
 \beas  
  \underset{n \to \infty}{\lim} \cE[\xi_n|\cF_\nu]=\cE[\xi|\cF_\nu], \q
 a.s.
 \eeas
 \end{thm}

\ss Next, we will derive an Upcrossing Theorem for
$\cE$-supermartingales, which is crucial in obtaining an RCLL
(right-continuous, with limits from the left) modification for the
process $\dis \left\{\cE[\xi|\cF_t]\right\}_{t \in [0,T]}$ as long
as $\xi \in Dom(\cE)$ is bounded from below. Obtaining a right
continuous modification is crucial, since otherwise the conditional
expectation $\cE[\xi|\cF_{\nu}]$ may not be well defined for any
$\nu \in \cS_{0,T}$.

\ss Let us first recall what the ``number of upcrossings" is:  Given
a real-valued process $\{X_t\}_{t \in [0,T]}$ and two real numbers
$a<b$, for any finite subset $F$ of $[0,T]$, we can define the
``number of upcrossings" $U_F(a,b; X(\o))$ of the interval $[a,b]$
by the sample path $\{X_t(\o)\}_{t \in F} $ as follows:
 Set $\nu_0=-1$, and for any $j=1,2, \cds$ we recursively define
  \beas
 \nu_{2j-1}(\o) &\dfnn& \min\{t \in F :\, t > \nu_{2j-2}(\o), X_t(\o) < a\} \land T  \in \cS^F_{0,T}, \\
  \nu_{2j}(\o) &\dfnn& \min\{t \in F :\, t > \nu_{2j-1}(\o), X_t(\o) > b \} \land T \in \cS^F_{0,T},
  \eeas
with the convention that $\min \es = \infty$. Then $U_F(a,b; X(\o))$
is defined to be the largest integer $j$ for which $\nu_{2j}(\o)<T$.
If $I \subset [0,T]$ is not a finite set, we define
 \beas 
 U_I\left(a,b;X(\o)\right) \dfnn \sup\{U_F(a,b;
X(\o)) :\,  F  \hb{ is a finite subset of } I  \}.
 \eeas

 \ms It will be convenient to introduce a subcollection of $\sD_T$
  \beas
  \wt{\sD}_T
 \dfnn \left\{ \L \in  \sD_T  :\, \hR \subset \L
 \right\}.
  \eeas
  Clearly, $\wt{\sD}_T$ contains all $L^p(\cF_T)$, $0 \le p \le
  \infty$. In particular, $L^\infty(\cF_T) $ is the smallest element
of $\wt{\sD}_T$ in the following sense:

  \begin{lemm} \label{lemm_wtsd}
 For each $\L \in  \wt{\sD}_T$, $L^\infty(\cF_T) \subset \L$.
 \end{lemm}

 \ss \no {\bf Proof:} For any $\xi \in L^\infty(\cF_T)$, we have $-
\|\xi\|_\infty, 2\|\xi\|_\infty \in \hR \subset
 \L$. Since $0\le \xi+ \|\xi\|_\infty
\le 2\|\xi\|_\infty$, a.s., (D3) implies that $\xi+ \|\xi\|_\infty
\in \L$. Then we can deduce from (D2) that $\xi=(\xi+
\|\xi\|_\infty)+(-\|\xi\|_\infty) \in \L$. \qed

 \ss For any $\bF$-adapted process $X$, we define its   right-limit
process by
 \beas 
     X^+_t \dfnn \underset{n \to
\infty}{\liminf} X_{q^+_n(t)}, \qq \hb{ for any } t \in [0,T],
 \eeas
 where   $q^+_n(t)$ is defined in (\ref{dyadic}). Since the filtration $\bF$ is right-continuous,
 we see that   $X^+$ is an $\bF$-adapted process.

\ms It is now the time to present our Upcrossing Theorem for
$\cE$-supermartingales.

    \begin{thm} \label{upcrossing} (Upcrossing Theorem)
Assume that (H0), (H1) hold and that $Dom(\cE) \in \wt{\sD}_T$. For
any $\cE$-supermartingale $X$, we assume either that $X_T \ge c$,
a.s. for some $c \in \hR$ or that the operator $\cE[\cd]$ is
concave: For any $\xi, \eta \in Dom(\cE)$
 \bea \label{def_cE_concave}
 \cE[\l\xi+(1-\l)\eta ]\ge \l\cE[\xi ] +(1-\l) \cE[\eta ], \q \fa \l \in
(0,1).
 \eea
Then for any two real numbers $a<b$, it holds that $ P\big(
U_{\cD_T} (a, b; X )<\infty\big)=1$. Thus we have
 \bea \label{limits0}
  P \Big(     X^+_t = \underset{n \to \infty}{\lim} X_{q^+_n(t)}
\hb{ for any } t \in [0,T] \Big)= 1.
 \eea
As a result, $X^+$ is an RCLL process.
 \end{thm}

In the rest of this section, we assume that the $\bF$-expectation
$\cE$ satisfies (H0)-(H2) and  that $Dom(\cE) \in \wt{\sD}_T$. The
following proposition will play a fundamental role throughout this
paper.

 \begin{prop} \label{prop_RCLL}
Let $X$ be a non-negative $\cE$-supermartingale.
  (1) Assume
either that $ \underset{t \in \cD_T}{\esssup}\, X_t \in Dom^+(\cE) $
or that for any sequence $\{\xi_n\}_{n \in \hN} \subset Dom^+(\cE)$
convergent a.s. to some $\xi \in L^0(\cF_T)$,
 \bea \label{ass_fatou}
 \hb{$\underset{n \to \infty}{\liminf} \cE[\xi_n ]<\infty~$ implies $~\xi
\in Dom^+(\cE)$. }
 \eea
 Then for any $\nu \in \cS_{0,T}$,   $X^+_\nu $   belongs to $
 Dom^+(\cE)$;

\ms \no (2) If $X^+_t \in Dom^+(\cE)$ for any $t \in [0,T]$, then
 $~ X^+$ is an RCLL $\cE$-supermartingale such that for any $t \in
 [0,T]$,
  $ 
  X^+_t  \le  X_t$, a.s.;

 \ms \no (3) Moreover, if the function $t \mapsto \cE[X_t]$ from $[0,T]$ to $\hR$ is right
continuous, then $ X^+$ is an RCLL modification of $X$. Conversely,
if $X$ has a right-continuous modification, then the function $t
\mapsto \cE[X_t]$ is right continuous.
 \end{prop}

  \ss Now we add one more hypothesis to the $\bF$-expectation $\cE$:
 \bi
  \item[{\bf (H3)}] For any $\xi \in
Dom^+(\cE)$ and $\nu \in \cS_{0,T}$, $X^{\xi,+}_\nu   \in
  Dom^+(\cE)$.
 \ei
In light of Proposition \ref{prop_RCLL} (1), (H3) holds if $
\underset{t \in \cD_T}{\esssup}\, \cE[\xi|\cF_t] \in Dom^+(\cE) $ or
if $\cE$ satisfies (\ref{ass_fatou}).

 \ss For each $\xi \in Dom^\#(\cE)$, we define $\xi'\dfnn \xi-c(\xi) \in
 Dom^+(\cE)$. Clearly $X^{\xi'} \dfnn \big\{ \cE[\xi'|\cF_t]\big\}_{t \in [0,T]}$ is
 a non-negative $\cE$-martingale. By (A2),
$\cE\big[X^{\xi'}_t\big]=\cE\big[\cE[\xi'|\cF_t]\big]=\cE[\xi']
 $ for any $t \in [0,T]$, which means that $t \mapsto \cE\big[X^{\xi'}_t\big]$ is continuous function
 on $[0,T]$. Thanks to Proposition \ref{prop_RCLL} (2) and (H3), the process $X^{\xi',+}_t
\dfnn \underset{n \to \infty}{\liminf} X^{\xi'}_{q^+_n(t)}$, $t \in
[0,T]$ is an RCLL modification of $X^{\xi'}$. Then for any $\nu \in
\cS_{0,T}$, we define
  \bea \label{cE_tau}
 \wt{\cE}[\xi|\cF_\nu] \dfnn  X^{\xi',+}_\nu \neg +c(\xi)
  \eea
as the conditional $\bF$-expectation of $\xi$ at the stopping time
 $\nu \in \cS_{0,T}$. Since we have assumed $Dom(\cE) \in \wt{\sD}_T$,
Lemma \ref{lemm_wtsd}, (H3), (D2)
 as well as the non-negativity of $ X^{\xi',+}_\nu$ imply that
    \bea \label{cE_tau1}
 \wt{\cE}[\xi|\cF_\nu] \in Dom^\#(\cE),
  \eea
 which shows that $\wt{\cE}[\cd|\cF_\nu]$ is an
operator from $  Dom^\#(\cE)$ to $ Dom^\#_\nu(\cE)\dfnn Dom^\#(\cE)
\cap L^0(\cF_\nu) $. In fact, $ \big\{\wt{\cE}[\cd|\cF_t]\big\}_{t
\in [0,T]} $ defines a $\bF$-expectation on $Dom^\#(\cE)$, as the
next result shows.

 \begin{prop} \label{prop_tilde_E}
 For any $\xi \in
 Dom^\#(\cE)$, $ \wt{\cE}[\xi|\cF_\cd] $ is an RCLL modification of $
\cE[\xi |\cF_\cd] $.
 $\big\{\wt{\cE}[\cd |\cF_t]\big\}_{t \in [0,T]}$ is an
 $\bF$-expectation with domain $Dom(\wt{\cE})=Dom^\#(\cE) \in \wt{\sD}_T$ and satisfying (H0)-(H2); thus
 all preceding results are applicable to $\wt{\cE}$.  \end{prop}

 \ss \no {\bf Proof:} As $Dom(\cE) \in
\wt{\sD}_T$ is assumed, we see that $Dom^\#(\cE)$ also belongs to
$\wt{\sD}_T$. Fix $\xi  \in Dom^\#(\cE)$. Since $X^{\xi',+}$ is an
RCLL modification of $X^{\xi'}$, (A4) implies that for any $t \in
[0,T]$
 \bea \label{eqn-axa01}
 \wt{\cE}[\xi|\cF_t]= X^{\xi',+}_t \neg +c(\xi)=
 \cE[\xi'|\cF_t]+c(\xi)=\cE[\xi'+c(\xi)|\cF_t]= \cE[\xi |\cF_t], \q
 a.s.
 \eea
Thus $ \wt{\cE}[\xi|\cF_\cd] $ is actually an RCLL modification of $
\cE[\xi |\cF_\cd] $. Then it is easy to show that the pair
$\big(\wt{\cE}, Dom^\#(\cE)\big)$ satisfies (A1)-(A4) and (H0)-(H2);
thus it is an $\bF$-expectation.  \qed

\ss We restate Proposition \ref{prop_RCLL} with respect to
$\wt{\cE}$ for future use.

 \begin{cor} \label{cor_RCLL}
Let $X$ be an $\wt{\cE}$-supermartingale such that $\underset{t \in
 [0,T]}{\essinf}\, X_t \ge c$, a.s. for some
$c \in \hR$. 

\ms \no  (1) If  $ \underset{t \in \cD_T}{\esssup}\, X_t \in
Dom^\#(\cE) $ or if (\ref{ass_fatou}) holds, then  
$X^+_\nu $ belongs to $ Dom^\#(\cE)$ for any $\nu \in \cS_{0,T}$;

\ms \no (2) If $X^+_t \in Dom^\#(\cE)$ for any $t \in [0,T]$, then
$~ X^+$ is an RCLL $\wt{\cE}$-supermartingale such that for any $t
\in [0,T]$, $ X^+_t  \le X_t$, a.s.

\ms \no (3) Moreover, if the function $t \mapsto \wt{\cE}[X_t]$ from
$[0,T]$ to $\hR$ is right continuous, then $ X^+$ is an RCLL
modification of $X$. Conversely, if $X$ has a right-continuous
modification, then the function $t \mapsto \wt{\cE}[X_t]$ is right
continuous.
 \end{cor}

 The next result is the
Optional Sampling Theorem of $\wt{\cE}$ for the stopping times in
$\cS_{0,T}$.

\begin{thm} \label{op_sa2} (Optional Sampling Theorem 2) 
Let $X$ be a right-continuous $\wt{\cE}$-supermartingale (resp.
$\wt{\cE}$-martingale, $\wt{\cE}$-submartingale) such that
 $\underset{t \in \cD_T}{\essinf}\, X_t \ge c$, a.s. for some $c \in \hR$.
If $X_\nu \in Dom^\#(\cE)$ for any $\nu \in \cS_{0,T}$, then for any
$\nu, \si \in \cS_{0, T}$, we have
 \beas
 \wt{\cE}[X_\nu |\cF_\si] \le (\hb{resp.}\,=,\; \ge)\, X_{\nu \land \si}, \q
 a.s.
 \eeas
\end{thm}

\no Using the Optional Sampling Theorem, we are able to extend
Corollary \ref{cor_os} and Proposition \ref{properties_2} to the
operators $\wt{\cE}[\cd|\cF_\nu]$, $\nu \in \cS_{0,T}$.

 \begin{cor}\label{cor_os2}
 For any $\, \xi  \in Dom^\#(\cE)$ and $\nu,\si \in \cS_{0,T}$, we have
  \bea   \label{eq:cor_os2}
   \wt{\cE}\big[\wt{\cE}[\xi|\cF_\nu]\big|\cF_\si \big]=\wt{\cE}[\xi|\cF_{\nu \land
 \si}], \q a.s.
  \eea
 \end{cor}

 \ss \no {\bf Proof:} Since $\big(\wt{\cE}, Dom^\#(\cE)\big)$ is an
 $\bF$-expectation by Proposition \ref{prop_tilde_E}, for any $\xi \in  Dom^\#(\cE)$, (A2) implies
 that the RCLL process $\wt{X}^\xi \dfnn \big\{ \wt{\cE}[\xi|\cF_t]\big\}_{t \in [0,T]}$
 is an $\wt{\cE}$-martingale. For any $t \in [0,T]$,  \eqref{eqn-axa01} and Proposition \ref{3addition} (2) show that
 \beas
   \wt{X}^\xi_t=\wt{\cE}[\xi|\cF_t]\ge \wt{\cE}[c(\xi)|\cF_t]=\cE [c(\xi)|\cF_t]= c(\xi) , \q
   a.s.,
 \eeas
 which implies that $ \underset{t \in
 [0,T]}{\essinf}\,\wt{X}^\xi_t \ge c(\xi) $, a.s.
 Then \eqref{cE_tau1} and Theorem \ref{op_sa2} give rise to
 \eqref{eq:cor_os2}.  \qed

\begin{prop} \label{properties_3}
For any $\xi, \eta  \in Dom^\#(\cE)$ and $\nu \in \cS_{0,T}$, it
holds that
  \bi
 \item[(1)] {\it ``Strict Monotonicity":}\;  $ \wt{\cE} [\xi |\cF_\nu] \le  \wt{\cE} [\eta |\cF_\nu] $,
$a.s.$ if $   \xi \le \eta$, $a.s.$; Moreover, if $ \wt{\cE} [\xi
|\cF_\si] = \wt{\cE} [\eta |\cF_\si]$, $a.s.$ for some $\si \in
\cS_{0,T}$, then $\xi=\eta$, $a.s.$;

 \item[(2)] {\it ``Zero-one Law":}\;
$\wt{\cE} [\b1_A \xi |\cF_\nu]=\b1_A \wt{\cE} [\xi |\cF_\nu]$, a.s.
for any $A \in \cF_\nu$;

 \item[(3)] {\it ``Translation Invariance":}\; 
$ \wt{\cE} [\xi+\eta |\cF_\nu]=\wt{\cE} [\xi |\cF_\nu]+\eta$, a.s.
if $ \eta \in  Dom^\#_\nu(\cE) $;

 \item[(4)] {\it ``Local Property":}\; $\wt{\cE}[\b1_A \xi+ \b1_{A^c}\eta|\cF_\nu]
  =\b1_A\wt{\cE}[ \xi|\cF_\nu]+ \b1_{A^c}\wt{\cE}[ \eta|\cF_\nu] $, $a.s.$
  for any $A \in \cF_\nu$;

 \item[(5)] {\it ``Constant-Preserving":}\;  $\wt{\cE} [\xi |\cF_\nu]=\xi$, $a.s.$, if $\,\xi \in
  Dom^\#_\nu(\cE) $.
 \ei

\end{prop}

\begin{rem} Corollary \ref{cor_os2}, Proposition \ref{properties_3} (2) and (\ref{eqn-axa01}) imply that
for any $\xi \in Dom^\#(\cE)$ and $\nu \in \cS_{0,T}$,
 \bea \label{eqn-g01}
  \cE[\b1_A \xi]=\wt{\cE}[\b1_A \xi] =\wt{\cE}\big[\wt{\cE}[ \b1_A
 \xi|\cF_\nu ]\big]=\wt{\cE}\big[\b1_A \wt{\cE}[  \xi|\cF_\nu
 ]\big]= \cE\big[\b1_A \wt{\cE}[  \xi|\cF_\nu
 ]\big], \q \fa   A \in \cF_\nu.
 \eea
 In light of Proposition \ref{3addition} (3), $\wt{\cE}[
 \xi|\cF_\nu ] $ is the unique element (up to a $P$-null set) in $Dom^\#_\nu(\cE)$ that makes (\ref{eqn-g01})
 hold. Therefore, we see that the random variable $\wt{\cE}[
 \xi|\cF_\nu ]$ defined by \eqref{cE_tau} is exactly the conditional $\bF$-expectation
 of $\xi$ at the stopping time $\nu$ in the classical sense.
  \end{rem}

 \ms In light of Corollary \ref{cor_os2} and Proposition \ref{properties_3}, we can generalize Fatou's
Lemma (Theorem \ref{fatou}) and the Dominated Convergence Theorem
(Theorem \ref{DCT}) to the conditional $\bF$-expectation
$\wt{\cE}[\cd|\cF_\nu]$, $\nu \in \cS_{0,T}$.

\begin{prop}\label{fatou2} (Fatou's Lemma 2) 
Let $\{\xi_n\}_{n \in \hN}$ be a sequence in $ Dom^\#(\cE)$ that
converges a.s. to some $\xi \in Dom^\#(\cE)$ and satisfies
$\underset{n \in \hN}{\essinf}\, \xi_n \ge c$, a.s. for some $c \in
\hR$, then for any $\nu \in \cS_{0,T}$, we have
 \bea  \label{eqn-d01}
 \wt{\cE}[\xi|\cF_\nu] \le \underset{n \to \infty}{\liminf}
 \wt{\cE}[\xi_n|\cF_\nu], \q a.s.,
 \eea
 where the right hand side of (\ref{eqn-d01}) could be equal to infinity with non-zero probability.
\end{prop}

 \begin{prop} \label{DCT2} (Dominated Convergence Theorem 2) 
Let $\{\xi_n\}_{n \in \hN}$ be a sequence in $ Dom^\#(\cE)$ that
converges a.s. and that satisfies $\underset{n \in \hN}{\essinf}\,
\xi_n \ge c$, a.s. for some $c \in \hR$. If there is an $\eta \in
Dom^\#(\cE)$ such that $\xi_n \le \eta $ a.s. for any $n \in \hN$,
then the limit $\xi$ of
$\{\xi_n\}_{n \in \hN}$ belongs to 
$ Dom^\#(\cE)$ and for any $\nu \in \cS_{0,T}$, we have
 \bea  \label{eqn-d02}
    \underset{n \to \infty}{\lim}
 \wt{\cE}[\xi_n|\cF_\nu]=\wt{\cE}[\xi|\cF_\nu], \q a.s.
 \eea
 \end{prop}

 \ss \no {\bf Proof of Propositions \ref{fatou2} and \ref{DCT2}:}
 In the proofs of Theorem
\ref{fatou} and Theorem \ref{DCT}, we only need to replace
$\{\xi_n\}_{n \in \hN}$ and $\cE[\cd|\cF_t]$ by $\{\xi_n-c\}_{n \in
\hN}$ and $\wt{\cE}[\cd|\cF_\nu]$ respectively. Instead of (A1),
(A3) and (A4), we apply Proposition \ref{properties_3} (1)-(3).
Moreover, since (A2) is only used on $Dom^+(\cE) $ in the proofs of
Theorem \ref{fatou} and Theorem \ref{DCT}, we can substitute
Corollary \ref{cor_os2} for it. Eventually, a simple application of
Proposition \ref{properties_3} (3) yields (\ref{eqn-d01}) and
(\ref{eqn-d02}). \qed

\section{Collections of $\bF$-Expectations}
\label{ch_3}

\setcounter{equation}{0}

\ms

In this section, we will show that \emph{pasting} of two
$\bF$-expectations at a given stopping time is itself an
$\bF$-expectation.  Moreover, pasting preserves (H1) and (H2). We
will then introduce the concept of a \emph{stable} class of
$\bF$-expectations, which are collections closed under pasting. We
will solve the optimal stopping problems introduced in
\eqref{eq:game-coop} and \eqref{eq:robust} over this class of
$\bF$-expectations. Before we show the pasting property of
$\bF$-expectations, we introduce the concept of convexity for an
$\bF$-expectation and give one of the consequences of having
convexity:
\begin{deff}
An $\bF$-expectation $\cE$ is called ``positively-convex" if for any
$\xi, \eta \in Dom^+(\cE)$, $\l \in (0,1)$ and $t  \in [0,T]$
 \beas
 \cE[\l\xi+(1-\l)\eta|\cF_t]\le \l\cE[\xi|\cF_t] +(1-\l) \cE[\eta|\cF_t], \q
 a.s.
 \eeas
 \end{deff}

\begin{lemm}\label{lem_pconvex}
  Any positively-convex $\bF$-expectation satisfies (H0). Moreover,
  an $\bF$-expectation $\cE$ is positively-convex if and only if
  the implied $\bF$-expectation $\big(\wt{\cE}, Dom^\#(\cE)\big)$ is convex, i.e., for any
$\xi, \eta \in Dom^\#(\cE)$, $\l \in (0,1)$ and $t  \in [0,T]$
 \bea \label{eqn_convex}
 \wt{\cE}[\l\xi+(1-\l)\eta|\cF_t]\le \l\wt{\cE}[\xi|\cF_t] +(1-\l) \wt{\cE}[\eta|\cF_t], \q
 a.s.
 \eea
 \end{lemm}

\begin{prop} \label{tau_A}
Let $ \cE_i, \cE_j$ be two $\bF$-expectations with the same domain
$\L \in \wt{\sD}_T$ and satisfying (H1)-(H3). For any $\nu \in
\cS_{0,T}$, we define the pasting of $\cE_i, \cE_j$ at the stopping
time $\nu$ to be the following RCLL $\bF$-adapted process
 \bea \label{tau_ij}
 \cE^\nu_{i,j}[\xi |\cF_t] \dfnn \b1_{\{\nu\le
t\}}  \wt{\cE}_j [\xi|\cF_t]  + \b1_{\{\nu> t\}}\wt{\cE}_i\big[
\wt{\cE}_j [\xi|\cF_\nu]  \big|\cF_t\big], \q    \fa t \in [0,T]
 \eea
for any $\xi \in \L^\#=\{\xi \in \L :\,  \xi \ge c, ~ a.s.\hb{ for
some }c=c(\xi) \in \hR\}$. Then $\cE^{\nu }_{i,j}$ is an
$\bF$-expectation
 with domain $\L^\# \in \wt{\sD}_T$ and satisfying (H1) and (H2). Moreover,
 if $\cE_i$ and $\cE_j$ are both positively-convex, $\cE^\nu_{i,j}$ is convex in the sense of
 \eqref{eqn_convex}.
\end{prop}

\ss In particular, for any $\si \in \cS_{0,T}$, applying Proposition
\ref{properties_3} (4) and (5), we obtain
 \bea
 \cE^\nu_{i,j}[\xi |\cF_\si]&\dneg =& \dneg \b1_{\{\nu\le
\si\}}  \wt{\cE}_j [\xi|\cF_\si]  + \b1_{\{\nu>
\si\}}\wt{\cE}_i\big[ \wt{\cE}_j [\xi|\cF_\nu]\big|\cF_\si\big]
=\b1_{\{\nu\le \si\}} \wt{\cE}_i\big[\wt{\cE}_j
[\xi|\cF_\si]\big|\cF_\si\big]  +
 \b1_{\{\nu> \si\}} \wt{\cE}_i\big[ \wt{\cE}_j [\xi|\cF_\nu]
\big|\cF_\si\big] \nonumber \hspace{0.5cm} \\
&\dneg =&  \dneg  \wt{\cE}_i\big[\b1_{\{\nu\le \si\}} \wt{\cE}_j
[\xi|\cF_\si]+ \b1_{\{\nu> \si\}} \wt{\cE}_j [\xi|\cF_\nu]
\big|\cF_\si\big]= \wt{\cE}_i\big[  \wt{\cE}_j [\xi|\cF_{\nu \vee
\si}] \big|\cF_\si\big] , \q a.s., \label{tau_ij2}
 \eea
 where we used the fact that $\{\nu > \si \} \in \cF_{\nu \land \si}
$ thanks to \cite[Lemma 1.2.16]{Kara_Shr_BMSC}.

\begin{rem}
 Pasting may not preserve (H0). From now on, we will replace the (H0) assumption by the positive
 convexity,
 which implies the former and which is an invariant property under
 pasting thanks to the previous two results. Positive convexity is also important in
constructing an optimal stopping time of \eqref{eq:game-coop} (see
Theorem \ref{SN_exist}).
\end{rem}

 \ms All of the ingredients are in place to introduce what we
mean by a stable class of $\bF$-expectations.
 As we will see in Lemma \ref{lem_02}, stability assures that
the essential supremum or infimum over the class can be approximated
by an increasing or decreasing sequence in the class.

\begin{deff} \label{def_stable_class}
 A class $\sE =\{\cE_i \}_{i \in \cI}$ of $\bF$-expectations is said to be
 ``stable"
 if 

\ss \no (1)  All $\cE_i,\, i \in \cI$ are positively-convex
$\bF$-expectations with the same domain
 $\L  \in \wt{\sD}_T$ and they satisfy \(H1\)-\(H3\);

 \ss \no (2)  $\sE$ is closed under pasting: namely, for any $i, j \in \cI$, $\nu \in \cS_{0,T}$,
there exists a $k= k(i,j, \nu ) \in \cI$ such that $\cE^{\nu}_{i,j}$
 coincides with $\wt{\cE}_k$ on $ \L^\# $. 

 \ss We shall denote $Dom(\sE) \dfnn \L^\#$, thus $Dom(\sE)= Dom^\#(\cE_i)
\in \wt{\sD}_T$ for any $i \in \cI$. Moreover, if $\sE' =\{\cE_i
\}_{i \in \cI'}$ satisfies (2) for some non-empty subset $\cI'$ of
 $\cI$, then we call $\sE'$ a stable subclass of $\sE$, clearly $Dom(\sE')=Dom(\sE)$.
\end{deff}

\begin{rem}
  The  notion of ``pasting" for linear expectations was given by
  \cite[Definition 6.41]{Follmer_Schied_2004}. The counterpart of Proposition~\ref{tau_A} for the linear expectations,
  which states that pasting two probability measures equivalent to $P$
  results in another probability measure equivalent to $P$, is given by \cite[Lemma 6.43]{Follmer_Schied_2004}.
  Note that in the case of linear expectations, (H1), (H2) and the convexity are trivially
  preserved because pasting in that case gives us a linear expectation.
On the other hand, the notion of stability for linear expectations
was given by \cite[Definition 6.44]{Follmer_Schied_2004}. The
stability is also referred to as ``fork convexity" in stochastic
control theory, ``m-stability" in stochastic analysis or
``rectangularity" in decision theory (see the introduction of
\cite{Delbaen_2006} and \cite{Bion_2009} for details).
\end{rem}

 \begin{eg} \label{eg_stable_class}
 (1) Let $\cP$ be the set of all probability measures equivalent
 to $P$, then $\sE_{\cP} \dfnn \{ E_Q \}_{Q \in \cP} $ is a stable
 class of linear expectations; see \cite[Proposition
 6.45]{Follmer_Schied_2004}.

\ss \no (2) Consider a collection $\fU$ of admissible control
processes. For any $U \in \fU$, let $P^U$ be the equivalent
probability measure defined via \cite[(5)]{Kara_Zam_2008}\;(or
\cite[(2.5)]{Kara_Zam_2006}), then $\sE_{\fU} \dfnn \{ E_{P^U} \}_{U
\in \fU} $ is a stable class of linear expectations; see Subsection
\ref{subsection_karatzas} of the present paper.

 \ss \no (3) For any $M>0$, a family $\sE_M$ of convex Lipschitz $g$-expectations
 with Lipschitz coefficient  $K_g \le M$ is an example
of stable class of non-linear expectations; see Subsection
\ref{subsection_gexp} of this paper.
 \end{eg}

The following Lemma gives us a tool for checking whether a random
variable is inside the domain $Dom(\sE)$ of a stable class $\sE$.

\begin{lemm} \label{lem_dom_sharp}
Given a stable class $\sE$ of $\bF$-expectations, a random variable
$\xi$ belongs to $Dom(\sE)$ if and only if $c \le \xi \le \eta$,
 a.s. for some $c \in \hR$ and $\eta \in Dom(\sE)$. 
\end{lemm}

 \ss \no {\bf Proof:} Consider a random variable $\xi$. If $\xi \in Dom(\sE)$,
 since $Dom(\sE)=Dom^\#(\cE_i) $ for any $i \in \cI$,
 we know that there exists a $c=c(\xi) \in \hR$ such that $\xi \ge c(\xi)$, a.s.

 \ss On the other hand,
 if $c \le \xi \le \eta$, a.s. for some $c \in \hR$ and $\eta \in
 Dom(\sE)$, it follows that $0 \le \xi-c \le \eta-c$, a.s. Since $Dom(\sE) \in
 \wt{\sD}_T$, we see that $-c, c \in \hR \subset Dom(\sE)$. Then (D2)
 shows that $\eta-c \in Dom(\sE)$ and thus (D3) implies that $\xi-c \in Dom(\sE)$, which further leads to that
 $\xi=(\xi-c)+c \in Dom(\sE)$ thanks to (D2) again. \qed
 \if{0}
 \ss Moreover, for any $\xi \in L^\infty(\cF_T)$, since
 $-\|\xi\|_\infty \le \xi \le \|\xi\|_\infty$, a.s. and since $\|\xi\|_\infty \in \hR \subset
 Dom(\sE)$, we see from the above argument that $\xi \in Dom(\sE)$.
 Hence $L^\infty(\cF_T) \subset Dom(\sE)$, proving the lemma. \qed
 \fi

 \ss We end this section by reviewing some basic properties of the essential
 supremum and essential infimum (for their definitions,
see e.g. \cite[Proposition VI-\b1-1]{Neveu_1975}, or \cite[Theorem
A.32]{Follmer_Schied_2004}).

 \begin{lemm}\label{lem_ess}
Let $\{\xi_j\}_{j \in \cJ}$ and $\{\eta_j\}_{j \in \cJ}$ be two
families of 
random variables of $L^0(\cF)$  with the
same index set $\cJ$.

\ss \no  (1)  If $\xi_j \le (=)\; \eta_j$, a.s. for any $j \in \cJ$,
 then $\underset{j \in \cJ}{\esssup}\, \xi_j \le (=)\; \underset{j \in
\cJ}{\esssup}\, \eta_j$, a.s.

\ss \no  (2)  For any $A \in \cF$, it holds a.s. that
 $\;\underset{j \in \cJ}{\esssup}\, \big( \b1_A \xi_j + \b1_{A^c} \eta_j
\big) = \b1_A \underset{j \in \cJ}{\esssup}\, \xi_j + \b1_{A^c}
\underset{j \in \cJ}{\esssup}\, \eta_j$; In particular, $\underset{j
\in \cJ}{\esssup}\, \big( \b1_A \xi_j  \big) =  \b1_A \underset{j
\in \cJ}{\esssup}\, \xi_j $, a.s.

\ss \no (3) For any random variable $\g \in L^0(\cF)$ and any
 $\a>0$, we have $\underset{j \in \cJ}{\esssup}\,  (\a \xi_j + \g ) = \a \, \underset{j
\in \cJ}{\esssup}\, \xi_j + \g $, a.s.

\ss \no Moreover, (1)-(3) hold when we replace $\underset{j \in
\cJ}{\esssup}\,$  by $\underset{j \in \cJ}{\essinf}\,$.
\end{lemm}

 \section{Optimal Stopping with Multiple Priors} \label{co_game}

\setcounter{equation}{0}

\ss

 In this section, we will solve an optimal stopping problem in
which the objective of the stopper is to determine an optimal
stopping time $\tau^*$ that satisfies
  \bea \label{eq:defn-coopgm}
  \underset{(i,\rho) \in \cI \times \cS_{0,T}}{\sup}
 \; \cE_i[Y_{\rho}+H_{\rho}^i ]
 = \underset{i \in \cI}{\sup}\; \cE_i[Y_{\t^*}+H_{\t^*}^i ],
  \eea
 where $\sE= \{\cE_i \}_{i \in \cI}$ is a stable class of
$\bF$-expectations, $Y$ is a primary reward process and $H^i$ is a
\emph{model}-dependent cumulative reward process. (We will outline
the assumptions on the reward processes below.) To find an optimal
stopping time, we shall build a so-called ``$\sE$-upper Snell
envelope", which we will denote by $Z^0$, of the reward process $
Y$. 
Namely, $Z^0$ is the smallest RCLL $\bF$-adapted process dominating
$Y$ such that $Z^0 + H^i$ is an $\wt{\cE}_i$-supermartingale for any
$i \in \cI$. We will show under certain assumptions that the first
time $Z^0$ meets $Y$ is an optimal stopping time for
\eqref{eq:defn-coopgm}.

 \ss We start by making some assumptions on the reward processes: 
 Let $\sE= \{\cE_i \}_{i \in \cI}$ be
a stable class of $\bF$-expectations accompanied by a family $\sH
\dfnn \{H^i \}_{i \in \cI}$ of right-continuous $\bF$-adapted
processes that satisfies

 \ss \no (S1)  For any $ i \in  \cI  $, $H^i_0 =0$, a.s. and
 \bea
    H^i_{\nu,\rho} \dfnn H^i_\rho-H^i_\nu \in Dom(\sE) ,\q
    \fa \nu, \rho \in \cS_{0, T}\hb{ with $\nu \le \rho$, a.s.}   \label{eqn-cxc01}
  \eea

   \no \qq   Moreover, if no member of $\sE$ satisfies (\ref{ass_fatou}), we
   assume that there exists a $j \in \cI$ such that
   \bea
     \z^j \dfnn \underset{s,t \in \cD_T;
s<t}{\esssup}\,H^j_{s,t} \in Dom(\sE). \label{ass_zi}
  \eea

    \no (S2) There exists a $ C_H <0$ such that for any $i \in \cI$,
$ \underset{s,t \in \cD_T; s<t}{\essinf}\, H^i_{s,t} \ge C_H $, a.s.

  \ss \no (S3)  For any $\nu \in \cS_{0,T}$ and $i,j \in
 \cI$, it holds for any $0 \le s < t \le T $ that
  $  H^k_{s,t}  =  H^i_{ \nu  \land s, \nu  \land t }
  +  H^j_{\nu  \vee s, \nu  \vee t }$, a.s.,
  where $k=k\big(i,j, \nu \big) \in \cI$ is the index defined in Definition \ref{def_stable_class} (2).

\begin{rem} \label{rem_H_process}

 \ss \no (1) For any $i \in \cI$, (S2) and the right-continuity of $H^i$
imply that except on a null set $N(i)$
  \bea \label{H_LB}
    H^i_{s,t} \ge C_H, \q \hb{for any }\; 0 \le s < t \le T,  \q
   \hb{thus} \q    H^i_{\nu,\rho} \ge C_H  , \q \fa \nu, \rho \in \cS_{0, T}\hb{ with $\nu \le \rho$, a.s.}
  \eea

 \ss \no (2) If \eqref{ass_zi} is assumed for some $j \in \cI$, we can deduce from the right-continuity of $H^j$ and \eqref{H_LB} that
 except on a null set $N $
 \beas
      C_H \le   H^j_{s,t} \le \z^j, ~\; \hb{for any }\; 0 \le s < t \le T,  ~\;
   \hb{thus}  ~\;    C_H \le   H^j_{\nu,\rho} \le \z^j , ~\; \fa \nu, \rho \in \cS_{0, T}\hb{ with $\nu \le \rho$, a.s.}
 \eeas
  Then Lemma \ref{lem_dom_sharp} implies that \eqref{eqn-cxc01} holds for
  $j$. Hence we see that \eqref{ass_zi} is a stronger condition than
  \eqref{eqn-cxc01}.

 \ss \no (3) Since $H^i$, $H^j$ and $H^k$ are all right-continuous
processes,   (S3) is equivalent to the   statement that    a.s.
      \bea   \label{h_tau_A}
 H^k_{s,t}  =  H^i_{ \nu  \land s, \nu  \land t }
  +  H^j_{\nu  \vee s, \nu  \vee t }  ,  \q \fa \; 0 \le s < t \le T .
  \eea
  \end{rem}

 \ss  Now we give an example of $\sH$.

  \begin{lemm} \label{lem_H_eg}
  Let $\{ h^i  \}_{i \in
\cI}$ be a family of progressive processes satisfying the following
assumptions:

  \ss \no (h1)  For any $ i \in  \cI  $ and $\nu, \rho \in \cS_{0, T}$
with $\nu \le \rho$, a.s., $\int_\nu^\rho h^i_t \, dt   \in Dom(\sE)
$.
 Moreover, if no member of $\sE$ satisfies

   \no \qq    (\ref{ass_fatou}),
 we assume that there exists a $j \in \cI$ such that  $  \int_0^T |h^j_t| \, dt   \in Dom(\sE)$.

  \ss \no (h2)  There exists a $\,c < 0 $ such that for any $i \in
\cI$, $h^i_t \ge c$, $\dtp$

  \ss \no (h3)  For any $\nu \in \cS_{0,T}$ and $i, j \in \cI$, it
holds for any $t \in [0,T] $ that
  $   h^k_t   = \b1_{\{\nu \le t \}} h^j_t+ \b1_{\{\nu > t \}} h^i_t $,
   $\dtp$, where

   \no \qq  $k=k\big(i,j, \nu \big) \in \cI$ is the index defined in Definition
 \ref{def_stable_class} (2).

  \ss \no   Then $\big\{H^i_t \dfnn \int_0^t h^i_s ds ,\, t \in
 [0,T]\big\}_{i \in \cI}$ is a family of right-continuous $\bF$-adapted
processes satisfying (S1)-(S3).
 \end{lemm}

 \noindent \textbf{Standing assumptions on $Y$ in this section.}

 \ss Let $Y$ be a right-continuous $\bF$-adapted process that
 satisfies

 \ss \no (Y1)  For any $\nu \in \cS_{0,T}$, $Y_\nu \in Dom(\sE)$.

 \ss \no  (Y2)  $ \underset{(i, \rho ) \in \cI \times \cS_{0,T} }
 {\sup} \cE_i\left[Y^i_\rho \right] <  \infty$, where $Y^i \dfnn \{Y_t+ H^i_t \}_{t \in [0,T]}$.
 Moreover, if no member of $\sE$ satisfies
 (\ref{ass_fatou}), then we

 \ss \no \qq   assume that
  \bea
    \qq \z_Y \dfnn \underset{(i, \rho, t) \in \cI \times \cS_{0,T} \times \cD_T}
 {\esssup}\, \wt{\cE}_i[Y^i_\rho |\cF_t]  \in  Dom(\sE).
 \label{eqn-m48}
  \eea

  \no (Y3)   $\underset{t  \in \cD_T }{\essinf}\, Y_t  \ge
C_Y$, a.s. for some $ C_Y <0$.

\ss
\begin{rem} \label{rem_fatou2}
(1) For any $i \in \cI$,  (A4) and (\ref{eqn-axa01}) imply that
$\cE_i$ satisfies (\ref{ass_fatou}) if and only if $\wt{\cE}_i$
 satisfies the following: Let  $\{\xi_n\}_{n \in \hN} \subset
Dom(\sE)$ be a sequence converging a.s. to some $\xi \in
L^0(\cF_T)$. If $\underset{n \in \hN}{\inf} \xi_n  \ge c$, a.s. for
some $c \in \hR$, then  $\underset{n \to \infty}{\liminf}
\wt{\cE}_i[\xi_n ]<\infty $ implies $ \xi \in Dom (\sE)$. The proof
of this equivalence is similar to that of Corollary \ref{cor_RCLL}.

\ss \no (2) It is clear that \eqref{eqn-m48} implies $ \underset{(i,
\rho ) \in \cI \times \cS_{0,T} }
 {\sup} \cE_i\left[Y^i_\rho \right] <  \infty$.

\ss \no (3) In light of (Y3) and the right-continuity of $Y$, it
holds except on a null set $N$ that
 \bea \label{Y_LB}
   Y_t  \ge C_Y, \q \fa t \in [0,T], \q \hb{thus} \q Y_\nu \ge
   C_Y,
   \q \fa \nu \in \cS_{0,T}.
  \eea
  Then for any $i \in \cI$, Remark \ref{rem_H_process} (1) implies that except on a
  null set $\wt{N}(i)$
   \bea \label{Y_LB2}
   Y^i_\nu =Y_\nu+H^i_\nu \ge C_* \dfnn C_Y+C_H, \q \fa \nu \in \cS_{0,T}.
   \eea
\end{rem}

    The following lemma states that the supremum or infimum over a
stable class of $\bF$-expectations can be approached by an
increasing or decreasing sequence in the class.

 \begin{lemm}  \label{lem_02}
 Let $\nu \in \cS_{0,T}$ and $\cU$ be a non-empty subset of
 $\cS_{\nu,T}$ such that
  \beas 
    \rho_1 \b1_A +  \rho_2 \b1_{A^c} \in \cU, \q \fa~~ \rho_1, \rho_2
    \in \cU, ~~~\fa A \in \cF_\nu.
  \eeas
Let $\{X(\rho) \}_{\rho \in \cU} \subset Dom(\sE) $ be a family of
random variables, indexed by $\rho$, such that for any $ \nu, \si
\in \cU,~
   \b1_{\{\nu=\si\}}X(\nu)=\b1_{\{\nu=\si\}}X(\si)$, a.s.,
  then for any stable subclass $\sE' =\{\cE_i
\}_{i \in \cI'}$ of $\sE$, there exist two sequences $\{(i_n,
\rho_n)\}_{n \in \hN}$ and $\{(i'_n, \rho'_n)\}_{n \in \hN}$ in
$\cI' \times \cU$ such that
   \bea
  \underset{(i, \rho) \in \cI' \times \cU}{\esssup}\,   \wt{\cE}_i
  \big[X(\rho)+H^i_{\nu,\rho} \big|\cF_\nu\big]
  &\tneg=&  \tneg \underset{n \to \infty}{\lim} \dneg \ua
    \wt{\cE}_{i_n} \big[X(\rho_n)+H^{i_n}_{\nu, \rho_n}
 \big|\cF_\nu\big]  ,\q  a.s.,  \label{eqn-k70} \\
  \underset{(i, \rho) \in \cI' \times \cU}{\essinf}\,   \wt{\cE}_i
  \big[X(\rho)+H^i_{\nu,\rho} \big|\cF_\nu\big]
  &\tneg =& \tneg \underset{n \to \infty}{\lim} \dneg \da
     \wt{\cE}_{i'_n} \big[X(\rho'_n)+H^{i'_n}_{\nu, \rho'_n} \big|\cF_\nu\big],\q
  a.s.  \label{eqn-k71}
  \eea
 \end{lemm}

\ms 
 For any $\nu \in \cS_{0,T}$ and $i \in \cI$, let us define
  \beas
  Z(\nu) \dfnn  \underset{(i, \rho) \in \cI \times \cS_{\nu,T} }
 {\esssup}\, \wt{\cE}_i[Y_\rho+ H^i_{\nu,\rho}  |\cF_\nu]\in \cF_\nu
 \q \hb{and} \q Z^i(\nu) \dfnn Z(\nu)+H^i_\nu.   
 \eeas
 Clearly, taking $\rho=\nu$ above yields that
  \bea  \label{eqn-q104x}
Y_\nu \le   Z(\nu),  \q  a.s.
 \eea

 \ss The following two Lemmas give the bounds on $Z(\nu)$,
$Z^i(\nu)$, $i \in \cI$, and show that they all belong to
 $Dom(\sE)$.

\begin{lemm} \label{Z_bound}
For any $\nu \in \cS_{0,T}$ and $i \in \cI$
  \bea \label{Z_LB}
     Z(\nu) \ge  C_*   \q \hb{and} \q    Z^i(\nu)  \ge   C_Y+2C_H  ,  \q a.s.
  \eea
Moreover, if no member of $\sE$ satisfies (\ref{ass_fatou}), then we
further have
 \bea \label{Z_UB}
     Z(\nu) \le  \z_Y -C_H   \q \hb{and} \q    Z^i(\nu)  \le   \z_Y -C_H +H^i_\nu   ,  \q
     a.s.,
  \eea
  where $\z_Y -C_H$ and $ \z_Y -C_H + H^i_\nu$ both belong to $Dom(\sE)$.
  \end{lemm}

\begin{lemm} \label{lem_Z_nu}
For any $\nu \in \cS_{0,T}$ and $i \in \cI$, both $Z(\nu)$ and
$Z^i(\nu)$ belong to $Dom(\sE)$.
\end{lemm}

 \ms 
 In the next two propositions, we will see
that the $\bF$-adapted process $\{Z(t)\}_{t \in [0,T]}$ has an RCLL
modification $Z^0$, and that both $\big\{Z^i(t)\big\}_{t \in [0,T]}$
and $Z^{i,0} \dfnn \big\{Z^0_t+H^i_t\big\}_{t \in
 [0,T]}$ are $\wt{\cE}_i$-supermartingales for any $i \in \cI$.

\begin{prop} \label{prop_01}
 For any $\nu, \si \in \cS_{0,T}$ and $ \g \in \cS_{\nu,T}$, we have
  \bea
  Z(\nu) &=& Z(\si), \q  a.s. \hb{ on ~~} \{ \nu=\si \}, \label{D_1} \\
 \underset{i \in \cI}{\esssup}\, \wt{\cE}_i [ Z(\g) +  H^i_{\nu,\g} |\cF_\nu]
 &=& \tneg \underset{(i, \rho) \in \cI \times \cS_{\g,T}
  }{\esssup}\,\wt{\cE}_i[Y_\rho+ H^i_{\nu,\rho} |\cF_\nu] \le Z(\nu), \q a.s.   \label{D_2}
  \eea
 \end{prop}


\begin{prop} \label{Z_RCLL}
Given $i \in \cI$, for any $\nu ,\rho \in \cS_{0,T}$ with $\nu \le
\rho$, a.s., we have
 \bea \label{eqn-xxx01}
 \wt{\cE}_i[Z^i(\rho) |\cF_{\nu}] \le  Z^i(\nu), \q a.s.
 \eea
In particular, $\big\{Z^i(t)\big\}_{t \in [0,T]}$ is an
$\wt{\cE}_i$-supermartingale. Moreover, the process
$\big\{Z(t)\big\}_{t \in [0,T]}$ admits an RCLL modification $Z^0$.
The process $Z^{i,0} \dfnn \big\{Z^0_t+H^i_t\big\}_{t \in
 [0,T]}$ is also an $\wt{\cE}_i$-supermartingale.
\end{prop}

We call $Z^0$ the ``$\sE$-upper Snell envelope" of the reward
process $Y$. From \eqref{eqn-q104x} and their right-continuity, we
see that $Z^0$ dominates $Y $ in the following sense:

 \begin{deff}
  We say that process $X$ ``dominates" process
  $ X'$ if $P\big(X_t \ge X'_t, ~\fa t \in [0,T]\big)=1$.
 \end{deff}

 \begin{rem} \label{rem_domin}
 \ss \no (1) If $X $ dominates $X'$, then $X_\nu \ge X'_\nu$,
 a.s. for any $\nu \in \cS_{0,T}$.

 \ss \no (2) Let $X$ and $ X'$ be two right-continuous $\bF$-adapted processes.
 If $P(X_t \ge X'_t )=1$ holds for all $t$
 in a countable dense subset of $[0,T]$, then $X $ dominates $X'$.
 \end{rem}

 The following Proposition demonstrates that $Z^0$ is the smallest RCLL $\bF$-adapted process
 dominating $Y$ such that $Z^{i,0}$ is an $\wt{\cE}_i$-supermartingale for any $i \in
 \cI$.

 \begin{prop} \label{Z0}
 The process $Z^0$ dominates the process $Y$. Moreover,
 for any $\nu \in \cS_{0,T}$ and $i \in \cI$, we have $ Z^0_\nu,~ Z^{i,0}_\nu \in
 Dom(\sE)$ and
 \bea
  Z^0_\nu=Z(\nu),     \q Z^{i,0}_\nu =Z^i(\nu) , \q  a.s.  \label{eqn-h42}
  \eea
 Furthermore, if $X$ is another RCLL $\bF$-adapted process dominating $Y$ such that $ X^i \dfnn \{X_t +
H^i_t\}_{t \in [0,T]}$ is an $\wt{\cE}_i$-supermartingale for any $i
\in
 \cI$, then $X$ also dominates $Z^0$.
 \end{prop}

As a consequence of Proposition~\ref{Z0} and \eqref{Z_LB}, we have
for any $\nu \in \cS_{0,T}$ and $i \in \cI$ that
 \bea
    Z^0_\nu  \ge  C_* ,  \q   Z^{i,0}_\nu
   \ge C_Y+2C_H, \q a.s.  \label{Z0_LB}
 \eea

 \if{0}
  \begin{rem}
  If $\{X_t\}_{t \in [0,T)}$ is a c\`adl\`ag $\cE$-supermartingale satisfying
  \beas
   P\big\{X_t \ge Y_t\big\}=1, \q \fa t \in [0,T),
  \eeas
  then it is still true that $X$ dominates $Z^0$.  
  \end{rem}

 Our main concern in this section is to find an optimal stopping
time for the cooperative game.
\begin{deff} \label{def_Emax}
  A stopping time $ \t^* \in  \cS_{0,T}$ is called an ``$\sE$-optimal" stopping time for the cooperative game if
 \beas
  \underset{i \in \cI}{\sup} \wt{\cE}_i [Y^i_{\t^*} ]=Z(0)
 =\underset{(i,\rho) \in \cI \times \cS_{0,T}}{\sup} \dneg \wt{\cE}_i [Y^i_\rho ].
 \eeas
\end{deff}
\fi

  In what follows, we first give \emph{approximately} optimal
stopping times. This family of stopping times will be crucial in
finding an optimal stopping time for \eqref{eq:defn-coopgm}.

 \begin{deff}
  For any $\d \in (0,1)$ and $\nu \in \cS_{0,T}$, we define 
 \beas
 &&\t_\d(\nu)  \dfnn  \inf\big\{ t \in [\nu ,T] :\,  Y_t \ge \d Z^0_t +(1-\d)(C_Y+2C_H)  \big\}\land T \in
 \cS_{\nu,T} \\
 \hb{and} &&
   J_\d(\nu) \dfnn \underset{i \in \cI}{\esssup}\,\wt{\cE}_i[Z^0_{\t_\d(\nu)}
  +  H^i_{\nu, \t_\d(\nu)} |\cF_\nu].
 \eeas
 \end{deff}

 \begin{rem} \label{delta_t}
(1) For any $\d \in (0,1)$ and $\nu \in \cS_{0,T}$, the
right-continuity of $Y$ and $Z^0$ implies that $\{\t_\d(t)\}_{t\in
[0,T]}$
 is also a right-continuous process. Moreover, since $Z^0_T =Z(T)=Y_T$, a.s., we can deduce from (Y3) that
 $Y_T > \d Z^0_T +(1-\d)(C_Y+2C_H)$. Then the right-continuity of processes $Y$ and $Z^0$ implies that
 \beas  
 Y_{\t_\d(\nu) } \ge \d Z^0_{\t_\d(\nu)} +(1-\d)(C_Y+2C_H)  , \q a.s.
 \eeas
 (2) For any $\nu \in \cS_{0,T}$,  we can deduce from (\ref{eqn-h42}) and (\ref{D_2}) that
  \bea \label{eqn-exe03}
   J_\d(\nu)=\underset{i \in
 \cI}{\esssup}\, \wt{\cE}_i\big[Z^0_{\t_\d(\nu)}+ H^i_{\nu,\t_\d(\nu)} \big|\cF_\nu\big]
    =\underset{i \in
 \cI}{\esssup}\,\wt{\cE}_i\big[Z(\t_\d(\nu))+ H^i_{\nu,\t_\d(\nu)} \big|\cF_\nu\big]
  \le Z(\nu)= Z^0_\nu, \q a.s.
  \eea
 \end{rem}

    The following two results show that for any  $\d \in (0,1)$, the
family $\{J_\d(\nu)\}_{\nu \in \cS_{0,T}}$ possesses similar
properties to  $\{Z(\nu)\}_{\nu \in \cS_{0,T}}$.
\begin{lemm} \label{lem_Jd}
 For any $\d \in (0,1)$ and $\nu \in \cS_{0,T}$,
 we have $J_\d(\nu) \in Dom(\sE)$. And for any $\si \in \cS_{0,T}$,
 $J_\d(\nu)= J_\d(\si)$, a.s. on $\{ \nu=\si \}$.
\end{lemm}

\begin{prop} \label{prop_Jd}
Given $\d \in (0,1)$, the followings statements hold:

 \ss \no  (1)  For any $i \in \cI$,
$\{J^i_\d(t) \dfnn J_\d(t) + H^i_t \}_{t \in [0,T]}$ is an
$\wt{\cE}_i $-supermartingale;

  \ss \no (2)  $\{J_\d(t) \}_{t \in
[0,T]}$ admits an RCLL modification $ J^{\d,0}$ such that the
process $J^{\d,i,0} \dfnn \{J^{\d,0}_t + H^i_t \}_{t \in [0,T]} $ is
an

 \no \qq   $\wt{\cE}_i $-supermartingale for each $i \in \cI$;

 \ss \no   (3)    For any $\nu \in \cS_{0,T}$, $J^{\d,0}_\nu \in Dom(\sE)$
 and $J^{\d,0}_\nu=J_\d(\nu)$, a.s.

 \end{prop}

 \ss Fix $\nu \in \cS_{0,T}$. The right continuity of $Z^0$ and \eqref{Z0_LB} imply that
  the stopping times $ \t_\d(\nu)$ are increasing in
$\d$. Therefore, we can define the limiting stopping time
 \bea \label{nu_star}
 \ol{\t}(\nu) \dfnn \underset{\d \nearrow 1}{\lim} \t_\d(\nu).
 \eea

  To show that $\ol{\t}(0) \in \cS_{0,T}$ is an optimal stopping time for \eqref{eq:defn-coopgm},
  we need the family of processes
$\{Y^i\}_{ i \in \cI}$ to be uniformly continuous from the left over
the stable class $\sE$.

 \begin{deff}
 The family  $\{Y^i\}_{ i \in \cI}$ is called ``$\sE$-uniformly-left-continuous"
 if for any  $   \nu, \rho \in \cS_{0,T}$ with $\nu \le \rho$, a.s. and for
 any sequence $\{\rho_n\}_{n \in \hN}\subset \cS_{\nu, T}$ increasing a.s. to
 $\rho$, we can find a subsequence $\{n_k\}_{k \in \hN}$ of $\hN$
  such that
  \bea \label{eqn-q170}
  \underset{k \to \infty}{\lim} \underset{i \in \cI}{\esssup}\,
   \Big|\wt{\cE}_i\big[\hb{$\frac{n_k}{n_k-1}$} Y_{\rho_{n_k}}+H^i_{\rho_{n_k}} \big|\cF_\nu\big]-
\wt{\cE}_i\big[ Y^i_\rho \big|\cF_\nu\big]\Big|=0, \q a.s.
  \eea
 \end{deff}

 The next theorem shows that $\ol{\t}(\nu)$ is not only
the first time when process $Z^0$ meets the process $Y$ after $\nu$,
but it is also an optimal stopping time after $\nu$. The assumption
that the elements of the stable class $\sE$ are convex (see
\eqref{eqn_convex}) becomes crucial in the proof of this result.

 \begin{thm} \label{SN_exist}
 Assume that 
 $\{Y^i\}_{ i \in \cI}$ is ``$\sE$-uniformly-left-continuous".
 Then for each $\nu \in \cS_{0,T}$, the stopping time $\ol{\t}(\nu)$
defined by (\ref{nu_star}) satisfies
 \bea  \label{eqn-a013}
  Z(\nu)=\underset{i \in \cI}{\esssup}\, \wt{\cE}_i\neg  \big[ Y_{\ol{\t}(\nu) }
 \neg + \neg H^i_{\nu, \ol{\t}(\nu) } \neg \big|\cF_\nu\big]
  \neg = \underset{i \in \cI}{\esssup}\, \wt{\cE}_i \neg \big[  Z(\ol{\t}(\nu))
 \neg + \neg  H^i_{\nu, \ol{\t}(\nu)}  \neg \big|\cF_\nu\big]
   \neg = \underset{i \in \cI}{\esssup}\, \wt{\cE}_i \big[  Z(\rho)
 \neg + \neg   H^i_{\nu,\rho}\big|\cF_\nu\big]
, ~\; a.s.
 \eea
for any $\rho \in \cS_{\nu, \ol{\t}(\nu)} $ and $  \ol{\t}(\nu) =
\t_1(\nu) \dfnn  \inf\big\{t \in [\nu,T] :\, Z^0_t=Y_t\big\}$, a.s.
 \end{thm}

  \ss  In particular, taking $\nu=0$ in \eqref{eqn-a013}, one can deduce
from \eqref{eqn-axa01} that $ \ol{\t}(0)  = \inf\big\{t \in [0,T]
:\, Z^0_t=Y_t\big\}$ satisfies
 \beas
  \underset{(i,\rho) \in \cI \times \cS_{0,T}}{\sup}
 \; \cE_i[Y_{\rho}+H_{\rho}^i ]=\underset{(i,\rho) \in \cI \times \cS_{0,T}}{\sup}
 \;  \wt{\cE_i}[Y_{\rho}+H_{\rho}^i ]
 = Z(0)=\underset{i \in \cI}{\sup} \, \wt{\cE}_i\neg \big[  Y_{\ol{\t}(0) }+H^i_{\ol{\t}(0) }  \big]
  =\underset{i \in \cI}{\sup} \, \cE_i\neg \big[ Y_{\ol{\t}(0) }+H^i_{\ol{\t}(0) }
  \big].
 \eeas
 Therefore, we see that $\ol{\t}(0)$, the first time the Snell envelope $Z^0$ meets the process $Y$ after time
$t=0$, is an optimal stopping time for \eqref{eq:defn-coopgm}.

 \section{Robust Optimal Stopping
  }\label{nonco_game}


\setcounter{equation}{0}

\ms

In this section we analyze the {\it robust} optimal stopping problem
in which the stopper
 aims to find an optimal stopping time $\tau_*$
  that satisfies
  \bea \label{eq:defn-non-coopgm}
  \underset{ \rho  \in  \cS_{0,T}}{\sup}\, \underset{ i \in \cI}{\inf}
 \; \cE_i[Y^i_{\rho}  ]
 = \underset{ i \in \cI}{\inf} \; \cE_i[Y^i_{\t_*} ],
  \eea
 where $\sE=\{\cE_i \}_{i \in \cI}$ is a stable class of
$\bF$-expectations and $Y^i=Y+H^i$, $i \in \cI$ is the
model-dependent reward process introduced in
\eqref{eq:defn-non-coopgm}.
 (We will modify the assumptions on the reward processes shortly). In order to
find an optimal stopping time we construct the lower and the upper
values of the optimal stopping problem at any stopping time $\nu \in
\cS_{0,T}$, i.e.,
 \beas 
 \ul{V}(\nu) \dfnn \underset{\rho \in \cS_{\nu,T}}{\esssup}\,
 \Big(  \underset{i \in \cI}{\essinf}\,\wt{\cE}_i \big[Y_\rho \neg
 +H^i_{\nu,\rho}|\cF_\nu \big] \Big)
 ,\q  \ol{V}(\nu) \dfnn \underset{i \in
 \cI}{\essinf}\, \Big(\, \underset{\rho \in \cS_{\nu,T}}{\esssup}\,
  \wt{\cE}_i \big[Y_\rho \neg +H^i_{\nu,\rho}|\cF_\nu \big] \Big).
 \eeas
With the so-called ``$\sE$-uniform-right-continuity" condition on
$\{Y^i\}_{i \in \cI}$, we can show that at any $\nu \in \cS_{0,T}$,
 $\ul{V}(\nu)$ and $\ol{V}(\nu)$ coincide with each
other (see Theorem~\ref{V_process}). We denote the common value, the
\emph{value} of the robust optimal stopping problem, as $V(\nu)$ at
$\nu$. We will show that up to a stopping time  $\ul{\t}(0)$ (see
Lemma~\ref{ul_nu}), at which we have $V(\ul{\t}(0))=Y_{\ul{\t}(0)}$,
a.s., the stopped value process $\big\{ V\big( \ul{\t}(0) \land
t\big) \big\}_{t \in [0,T]} $ admits an RCLL modification $V^0$. The
main result in this section, Theorem~\ref{V_RC}, shows that the
first time $V^0$ meets $Y$ is an optimal stopping time for
\eqref{eq:defn-non-coopgm}.

\ss \noindent{\textbf{Standing assumptions on $\sH$ and $Y$ in this
section.}} Let us introduce
 \beas
 R^i(\nu) \dfnn
   \underset{\rho \in \cS_{\nu,T}}{\esssup}\,
  \wt{\cE}_i \big[Y_\rho\neg +H^i_{\nu,\rho}|\cF_\nu], \q  \hb{for any } i \in \cI \hb{ and }\nu \in \cS_{0,T}.
 \eeas
To adapt the results we obtained for the family $\{Z(\nu)\}_{\nu \in
\cS_{0,T}}$ to each family $\{R^i(\nu)\}_{\nu \in \cS_{0,T}}$, $i
\in \cI$, we assume that $\sH = \{H^i \}_{i \in \cI}$ is a family of
right-continuous $\bF$-adapted processes satisfying (S2), (S3) and,

 \bi
\item[(S1')] For any $i \in
 \cI$, $H^i_0 =0$, a.s. and \eqref{eqn-cxc01} holds. If $\cE_i$ does not
 satisfy (\ref{ass_fatou}), then we assume that $  \z^i = \underset{s,t \in \cD_T; s<t}{\esssup}\,H^i_{s,t}
  \in Dom(\sE) $.
 \ei
On the other hand, we assume that $Y$ is a right-continuous
$\bF$-adapted process that satisfies (Y1), (Y3) and
 \bi
\item[(Y2')]
   For any $i \in
 \cI$, $\underset{  \rho   \in \cS_{0,T} }
 {\sup} \wt{\cE}_i\left[Y^i_\rho \right] <  \infty$. If $\cE_i$ does not
 satisfy (\ref{ass_fatou}), then we assume that $ \underset{( \rho, t) \in \cS_{0,T} \times \cD_T}
 {\esssup}\, \wt{\cE}_i[Y^i_\rho |\cF_t]  \in  Dom(\sE)$.
  \ei
 We also assume that for any $i \in \cI$, $Y^i$ is ``quasi-left-continuous"
under $\wt{\cE}_i$: for any $\nu, \rho \in \cS_{0,T}$ with
$\nu \le \rho$, a.s. and for any sequence $\{\rho_n\}_{n \in
\hN}\subset \cS_{\nu, T}$ increasing a.s. to $\rho$, we can find a
subsequence $\big\{n_k=n^{(i)}_k \big\}_{k \in \hN }$ of $\hN$ such that
  \bea      \label{eqn-q180}
  \underset{n \to \infty}{\lim}  \wt{\cE}_i\big[\hb{$\frac{n_k}{n_k-1}$} Y_{\rho_{n_k}}
  +H^i_{\rho_{n_k}} \big|\cF_\nu\big]
    =\wt{\cE}_i\big[ Y^i_\rho \big|\cF_\nu\big], \q a.s.
  \eea

  \begin{rem} \label{rem_nonco_YH}
  (S1') and (Y2') are analogous to (S1) and (Y2) respectively
  while the quasi-left-continuity \eqref{eqn-q180} is the counterpart
  of the $\sE$-uniform-left-continuity \eqref{eqn-q170}. 
  It is obvious that
  (S1') implies (S1) and that \eqref{eqn-q170} gives rise to
  \eqref{eqn-q180}. Moreover, \eqref{eqn-m48} implies (Y2'): In
  fact, for any $i \in \cI$, one can deduce from \eqref{Y_LB2} that
   \beas
     C_* \le \underset{( \rho, t) \in \cS_{0,T} \times \cD_T}
 {\esssup}\, \wt{\cE}_i[Y^i_\rho |\cF_t]  \le  \underset{(i, \rho, t) \in \cI \times \cS_{0,T} \times \cD_T}
 {\esssup}\, \wt{\cE}_i[Y^i_\rho |\cF_t], \q a.s.
   \eeas
   Then Lemma \ref{lem_dom_sharp} implies that $\underset{( \rho, t) \in \cS_{0,T} \times \cD_T}
 {\esssup}\, \wt{\cE}_i[Y^i_\rho |\cF_t]  \in  Dom(\sE)$, and it
 follows that $\underset{  \rho   \in \cS_{0,T} }
 {\sup} \wt{\cE}_i\left[Y^i_\rho \right] <  \infty$. \qed
 \end{rem}

  \ss Fix $i \in \cI$.
  Applying Lemma \ref{lem_Z_nu}, \eqref{Y_LB}, \eqref{H_LB}, \eqref{D_2}, Proposition \ref{Z_RCLL},
  Proposition \ref{Z0} and Theorem \ref{SN_exist}
   to the family $\{R^i(\nu)\}_{\nu \in \cS_{0,T}}$, we obtain
 \begin{prop} \label{prop_Ri}
\no (1) For any $\nu \in \cS_{0,T}$,  $R^i(\nu)$ 
 belongs to $Dom(\sE)$ and satisfies
  \bea \label{eqn-q104}
 C_Y \le Y_\nu \le \underset{ \rho \in  \cS_{\nu,T}}{\esssup}\, \wt{\cE}_i [Y_\rho+
    H^i_{\nu,\rho}|\cF_\nu]= R^i(\nu), ~\;\; a.s. , \q \hb{thus} ~\;\;
    C_* \le Y^i_\nu, \q a.s. 
 \eea

\no  (2) For any $\nu, \si \in \cS_{0,T}$ and $ \g \in \cS_{\nu,T}$,
we have
  \bea
  R^i(\nu) &=& R^i(\si), \q  a.s. \hb{ on ~~} \{ \nu=\si \}, \label{eqn-k777} \\
 \wt{\cE}_i [ R^i(\g) +  H^i_{\nu,\g}|\cF_\nu]
 &=& \tneg \underset{ \rho  \in \cS_{\g,T}}{\esssup}\,
  \wt{\cE}_i[Y_\rho +  H^i_{\nu,\rho}|\cF_\nu] \le R^i(\nu), \q a.s.   \label{eqn-k26}
  \eea

\no  (3)
 The process $\big\{R^i(t)\big\}_{t \in [0,T]}$ admits an RCLL
 modification $R^{i,0}$, called ``$\cE_i$ Snell
envelope", such that
 $\{   R^{i,0}_t+H^i_t \}_{t \in [0,T]}$ is an
 $\wt{\cE}_i$-supermartingale and that for any $\nu \in \cS_{0,T}$
 \bea
 R^{i,0}_\nu = R^i(\nu), \q a.s.  \label{eqn-p04}
 \eea

\no  (4)
For any $\nu \in \cS_{0,T}$, $ \t^i(\nu) \dfnn \inf\{t \in [\nu,T]
:\, R^{i,0}_t=Y_t\}$
is an optimal stopping time with respect to $\sE^i$ after time
$\nu$, i.e., for any $\g \in \cS_{\nu, \t^i(\nu)}$,
 \bea
   R^i(\nu)  =   \wt{\cE}_i \big[ Y_{\t^i(\nu) }+
    H^i_{\nu, \t^i(\nu) }  |\cF_\nu]   =   \wt{\cE}_i \big[  R^i(\t^i(\nu))+
   H^i_{\nu, \t^i(\nu)}  |\cF_\nu]
   =   \wt{\cE}_i \big[  R^i (\g) +H^i_{\nu,\g}|\cF_\nu], \q
   a.s.     \label{eqn-k25} 
 \eea
\end{prop}

\begin{cor}
For any $\nu \in \cS_{0,T}$, both $ \ul{V}(\nu) $ and $ \ol{V}(\nu)$
belong to $Dom(\sE)$.
\end{cor}

\ss \no {\bf Proof:} Fix $(l, \rho)  \in \cI \times \cS_{\nu,T} $,
for any $ i \in \cI$, \eqref{Y_LB}, \eqref{H_LB} and Proposition
\ref{properties_3} (5) imply that
 \beas
  \wt{\cE}_i \big[Y_\rho \neg
 +H^i_{\nu,\rho}|\cF_\nu] \ge \wt{\cE}_i \big[C_Y+C_H |\cF_\nu]  = C_* , \q  a.s.
 \eeas
Taking the essential infimum over $i \in \cI$ on the left-hand-side
yields that
 \beas 
 \q  \;\;  C_* \le   \underset{i \in \cI}{\essinf}\,\wt{\cE}_i \big[Y_\rho \neg
 +H^i_{\nu,\rho}|\cF_\nu] \le \neg \underset{\rho \in \cS_{\nu,T}}{\esssup}\,
 \Big(  \underset{i \in \cI}{\essinf}\,\wt{\cE}_i \big[Y_\rho \neg
 +H^i_{\nu,\rho}|\cF_\nu] \Big) \neg = \ul{V}(\nu) \le \ol{V}(\nu) =\underset{i \in \cI}{\essinf}\,
  R^i(\nu) \le   R^l(\nu)  , ~\;
  a.s.
  \eeas
   Since $R^l(\nu) \in Dom(\sE)$ by Proposition \ref{prop_Ri} (1),
 a simple application of Lemma \ref{lem_dom_sharp} proves the
corollary. \qed

\ss As we will see in the next lemma since the stable class $\sE$ is
closed under pasting \big(see Definition \ref{def_stable_class}
(2)\big), $\ol{V}(\nu)$ can be approximated by a decreasing sequence
that belongs to the family $\{R^i(\nu)\}_{i \in \cI}$.

\begin{lemm} \label{essinf_lim}
 For any $\nu \in \cS_{0,T}$, there exists a sequence $\{i_n\}_{n \in \hN} \subset \cI$ such that
  \bea \label{eqn-p05}
  \ol{V}(\nu) =\underset{i \in \cI}{\essinf}\,
  R^i(\nu) = \underset{n \to \infty}{\lim} \dneg \da R^{i_n}(\nu), \q a.s.
  \eea
\end{lemm}

  Thanks again to the stability of $\sE$ under pasting, the infimum of
the family $\{\t^i(\nu)\}_{i \in \cI}$ of optimal stopping times can
be approached by a decreasing sequence in the family. As a result
the infimum is also a stopping time.

\begin{lemm} \label{ul_nu}
For any $ \nu \in \cS_{0,T}$, there exists a sequence $\{i_n\}_{n
\in \hN} \subset \cI$ such that
  \beas
  \ul{\t}(\nu)\dfnn \underset{i \in \cI}{\essinf}\, \t^i(\nu) = \underset{n \to \infty}{\lim} \dneg \da \t^{i_n}(\nu),
  ~\;  a.s., ~ \hb{ thus }\; \ul{\t}(\nu)  \in \cS_{\nu,T}.
  \eeas
\end{lemm}

The family of stopping times $\{\ul{\t}(\nu)\}_{\nu \in \cS_{0,T}}$
will play a critical role in this section. The next lemma basically
shows that if $\wt{\cE}_j$ and $\wt{\cE}_k$ behave the same after
some stopping time $\nu$, then $R^{j,0}$ and $R^{k,0}$ are identical
after $\nu$:

\begin{lemm} \label{lem_trancate}
For any $i,j \in \cI$ and $\nu \in \cS_{0,T}$, let $k=k(i,j,\nu) \in
\cI$ as in Definition \ref{def_stable_class}.
For any $\si \in \cS_{\nu,T}$,
 we have $ R^{k,0}_\si= R^k(\si)=R^j(\si)=R^{j,0}_\si$, a.s.
\end{lemm}

 \no {\bf Proof:} For any $\rho \in
\cS_{\si,T}$, 
 applying Proposition \ref{properties_3} (5) to $\wt{\cE}_i $,
 we can deduce from (\ref{h_tau_A}) and (\ref{tau_ij2}) that
 \beas
  \wt{\cE}_{k}\big[ Y_\rho +H^{k}_{\si,\rho} \big|\cF_{\si}\big]
& =& \wt{\cE}_{k}\big[ Y_\rho +H^{j}_{\si,\rho} \big|\cF_{\si}\big]
 =  \cE^{\nu}_{i, j}\big[ Y_\rho +H^{j}_{\si,\rho}  \big|\cF_{\si}\big]
 =\wt{\cE}_i\big[  \wt{\cE}_{j} \big[Y_\rho+H^{j}_{\si,\rho}\big|\cF_{\nu \vee \si}\big] \big|\cF_{\si}\big]\\
 & =& \wt{\cE}_i\big[  \wt{\cE}_{j}
\big[Y_\rho+H^{j}_{\si,\rho}\big|\cF_{\si}\big] \big|\cF_{\si}\big]
 =  \wt{\cE}_{j} \big[Y_\rho+H^{j}_{\si,\rho}\big|\cF_{\si}\big]   , \q a.s.
 \eeas
 Then (\ref{eqn-p04}) implies that
 \beas
  R^{k,0}_{\si}= R^k (\si )
  = \underset{\rho \in \cS_{\si, T}}{\esssup}\,\wt{\cE}_{k}
  \big[Y_\rho+ H^k_{\si,\rho} \big|\cF_{\si}\big]
  =  \underset{\rho \in \cS_{\si,T}}{\esssup}\,
  \wt{\cE}_{j} \big[Y_\rho+H^{j}_{\si,\rho}\big|\cF_{\si}\big]
 = R^{j}(\si)=R^{j,0}_{\si}, \q a.s.,
 \eeas
which proves the lemma. \qed

We now introduce the notion of the uniform right continuity of the
family $\{Y^i\}_{i \in \cI}$ over the stable class $\sE$. With this
assumption on the reward processes, we can show that at any $\nu \in
\cS_{0,T}$, $\ul{V}(\nu)=\ol{V}(\nu)$, a.s., thus the robust optimal
stopping problem has a value $V(\nu)$ at $\nu$.

 \begin{deff}
 The family  $\{Y^i\}_{ i \in \cI}$ is called ``$\sE$-uniformly-right-continuous"
 if for any $\nu \in \cS_{0,T}$ and for any sequence
  $\{\nu_n\}_{n \in \hN } \subset \cS_{\nu,T}$ decreasing a.s. to $\nu
  $, we can find a subsequence of $\{\nu_n\}_{n \in \hN }$ (we still
  denote it by $\{\nu_n\}_{n \in \hN }$) such that
  $        \underset{n \to \infty}{\lim}\underset{i \in \cI}{\esssup}\,
  \Big|\wt{\cE}_i[Y^i_{\nu_n}|\cF_\nu]-Y^i_\nu \Big|=0$,   a.s.

 \end{deff}

\begin{thm} \label{V_process}
Suppose that $\{Y^i\}_{ i \in \cI}$ is
``$\sE$-uniformly-right-continuous". Then for any $\nu \in
\cS_{0,T}$, we have
 \bea \label{eqn-k51}
   \ul{V}(\nu) = \underset{i \in \cI}{\essinf}\, \wt{\cE}_i  \big[ Y_{\ul{\t}(\nu)}
   +\neg H^i_{\nu,\ul{\t}(\nu)}  \big|\cF_\nu \big] = \ol{V}(\nu)
    \ge Y_\nu, \q a.s.
  \eea
We will denote the common value by $V(\nu)(=\ul{V}(\nu)=
\ol{V}(\nu))$. Observe that $\ul{\t}(0)$ is optimal for the robust
optimal stopping problem in \eqref{eq:defn-non-coopgm}.
\end{thm}

\noindent \textbf{Standing assumption on $Y$ for the rest of this
section.}
 We assume that the family of processes $\{Y^i\}_{ i
\in \cI}$ is ``$\sE$-uniformly-right-continuous".

 \begin{prop} \label{V_Y_meet}
For any $\nu \in \cS_{0,T}$, we have
$V(\ul{\t}(\nu))=Y_{\ul{\t}(\nu)}$, a.s.
 \end{prop}

       Note that $\ul{\t}(\nu)$ may not be the first time after $\nu$
when the value of the robust optimal stopping problem is equal to
the primary reward. Actually, since the process $\{V(t)\}_{t \in
[0,T]}$ is not necessarily right-continuous, $\inf\{t \in [\nu, T]
\;|\; V(t)=Y_t\} $ may not even be a stopping time. We will address
this issue in the next two results.

 \begin{prop} \label{R_sub}
Given $ i \in \cI$, for any $\nu, \rho \in \cS_{0,T}$ with $  \nu
\le \rho$, a.s., we have
 \bea \label{eqn-m31}
  \underset{k \in \cI}{\essinf}\,\wt{\cE}_k[V^k(\rho)|\cF_\nu] \le
  V^i(\nu), \q a.s.,
 \eea
where $V^i(\nu) \dfnn V(\nu)+H^i_\nu \in Dom(\sE)$. Moreover if
$\rho \le \ul{\t}(\nu)$, a.s., then
 \bea
   &&\q \wt{\cE}_i[V^i(\rho)|\cF_\nu]  \ge  V^i(\nu), \q a.s. \hspace{1cm}  \label{eqn-l20}
 \eea
 In particular, the stopped process $\big\{V^i\big( \ul{\t}(0) \land t\big)\big\}_{t
\in [0,T]}$ is an $\wt{\cE}_i$-submartingale.
\end{prop}

   \ss  Now we show that the stopped value process
   $\big\{ V\big(\ul{\t}(0) \land t\big) \big\}_{t \in [0,T]} $ admits an RCLL
    modification $V^0$. As a result, the first time when the process $V^0$ meets the process $Y$ after time
$t=0$ is an optimal stopping time of the robust optimal stopping
problem.

    \begin{thm} \label{V_RC}
 Assume that for some $i' \in \cI$, $\z^{i'} = \underset{s,t \in \cD_T; s<t}{\esssup}\,H^{i'}_{s,t}
  \in Dom(\sE) $ and that there exists a concave
 $\bF$-expectation $\cE'$ (not necessarily in $\sE$) satisfying (H0) and (H1) 
 such that
  \bea \label{ass_odd}
  Dom(\cE')\supset \{-\xi: \xi \in
 Dom(\sE)\} \; \hb{and for every
$\wt{\cE}_{i'}$-submartingale $X$,\; $-X$ is an
$\cE'$-supermartingale.}
 \eea
 We also assume that for any $\rho \in \cS_{0,T}$, there exists a
 $j=j(\rho) \in \cI$ such that $ \underset{  t \in \cD_T} {\esssup}\, \wt{\cE}_j[Y^j_\rho |\cF_t] \in
 Dom(\sE)$.

 \noindent (1) Then the stopped value process
   $\big\{ V\big(\ul{\t}(0) \land t\big) \big\}_{t \in [0,T]} $ admits an RCLL
    modification $V^0$ (called ``$\sE$-lower Snell
envelope" of $\,Y$) such that for any $\nu \in \cS_{0,T}$
    \bea  \label{V0}
    V^0_\nu = V(\ul{\t}(0) \land \nu ), \q  a.s.
    \eea
  (2) Consequently,
    \bea  \label{eqn-m15}
  \t_V  \dfnn \inf\{ t \in [0, T] :\, V^0_t=Y_t \}
   \eea
 is  a stopping time. In fact, it is an optimal stopping time of \eqref{eq:defn-non-coopgm}.
 \if{0}
    i.e. \beas \underset{i\in \cI}{\inf} \Big( \underset{ \rho \in
   \cS_{0,T}}{\sup} \cE_i \big[Y^i_\rho   \big]\Big)=
    \underset{ \rho \in \cS_{0,T}}{\sup}
   \Big(\underset{i\in
   \cI}{\inf}\cE_i \big[Y^i_\rho   \big]\Big)
   = \underset{i \in \cI}{\inf}  \cE_i   \big[ Y^i_{\t_V}   \big] .   \eeas
   \fi
 \end{thm}

  \section{Remarks on Sections \ref{co_game} \& \ref{nonco_game} }
 \label{section_remark}

 \setcounter{equation}{0}

\ms \no  {\bf \large Remark 1.}

\ms

Let $\sE= \{\cE_i \}_{i \in \cI}$ be a stable class of
$\bF$-expectations. For any $\xi \in Dom(\sE)$ and $\nu \in
\cS_{0,T}$, we define
 \beas
  \ol{\sE}[\xi |\cF_\nu] \dfnn \underset{i  \in \cI }{\esssup}\, \wt{\cE}_i
  [\xi|\cF_\nu]  \q   \hb{ and }  \q
   \ul{\sE}[\xi |\cF_\nu] \dfnn \underset{i  \in \cI }{\essinf}\, \wt{\cE}_i
   [\xi|\cF_\nu]
  \eeas
 as the maximal and minimal expectation of $\xi$ over $\sE$ at the stopping time
 $\nu$. It is worth pointing out that $\ol{\sE}$ is not an $\cF$-expectation on $Dom(\sE)$ simply because
$\ol{\sE}[\xi|\cF_t]$ may not belong to $Dom(\sE)$ for some $\xi \in
Dom(\sE)$ and $t \in [0,T]$. On the other hand, we will see in
Example \ref{eg_olsE} that neither $\ol{\sE}$ nor $ \ul{\sE}$
satisfy strict monotonicity. Moreover, as we shall see in the same
example, $\ol{\sE}$ does not satisfy (H2) while $\ul{\sE}$ does not
satisfy (H1); thus we do not have a dominated convergence theorem
for either $\ol{\sE}$ or $\ul{\sE}$. Note also that $\ul{\sE}$ may
not even be convex.

 \ss Our results in Sections~\ref{co_game} and \ref{nonco_game}
 can be interpreted as a first step in extending the results for the optimal stopping problem
$ \underset{ \rho  \in \cS_{0,T}}{\sup} \, \cE_i[Y_{\rho}]$, in
which $\cE_i$ ($i \in \cI$) is an $\bF-$expectation
 satisfying positive convexity and the assumptions (H1)-(H3), to
optimal stopping problems for other non-linear expectations, such as
$\ol{\sE}$ and $\ul{\sE}$, which may fail to satisfy these
assumptions.

 \begin{eg} \label{eg_olsE}
 Consider a probability space $\big([0,\infty), \sB[0,\infty),
 \bF=\{\mathcal{F}_t\}_{t \geq 0}, P \big)$ be a filtered probability space in which $P$ is defined by
 \beas
 P(A) \dfnn  \int_A e^{-x} dx, \q \fa A \in \sB[0,\infty).
 \eeas
 We assume that the filtration $\bF$ satisfies the usual hypothesis.
Let $\cP$ denote the set of all probability measures equivalent
 to $P$. For any $n \in \hN$, we define a $P_n \in \cP$ by
  \beas
 P_n(A) \dfnn n \int_A e^{-nx} dx, \q \fa A \in \sB[0,\infty).
 \eeas
As discussed in Example \ref{eg_stable_class},
 $ \sE=\{ E_Q \}_{Q \in \cP} $ is a stable class.  For any $h >0$, one can deduce that
  \beas
 ~\;\;\; 1 = \neg \underset{Q \in \cP}{\sup} E_Q [1 ] \neg  \ge \neg
 \ol{\sE}\big[\b1_{[0,h]}\big]
  \neg =\neg  \underset{Q \in \cP}{\sup} E_Q\big[\b1_{[0,h]}\big]
    \neg\ge    \underset{n \in
   \hN}{\sup}E_{P_n}\big[\b1_{[0,h]}\big] \neg = \underset{n \in
   \hN}{\sup}\, P_n[0,h]   =   \underset{n \in
   \hN}{\lim} \big(1-e^{-nh}\big)\neg =1,
 \eeas
 where we used the fact that $\wt{E}_Q= E_Q$ for any $Q \in \cP$
  since $E_Q[\xi|\cF_\cd]$ is an RCLL process for any $\xi \in  L^1\big([0,\infty), \\ \sB[0,\infty), P
 \big) $. Hence, we have
 \beas 
   \ol{\sE}\big[\b1_{[0,h]}\big]=1, \q \fa h>0,
  \eeas
 which implies that $\ol{\sE}$ does not satisfy strict monotonicity.

 \ms Moreover, $\ol{\sE}$ does not satisfy (H2). For $\xi=0$, $\eta=1$ and
 $A_n=[0, \frac{1}{n}]$, $n \in \hN$, it follows that
 \beas
 \underset{n \to \infty  }{\lim}\neg \da \ol{\sE}
  [\xi+\b1_{A_n}\eta]=\underset{n \to \infty  }{\lim}\ol{\sE}
  [ \b1_{[0, \frac{1}{n}]} ]=1 \ne 0= \underset{Q \in \cP}{\sup}
  \wt{E}_Q[0]=\ol{\sE}[0]=\ol{\sE}[\xi].
 \eeas

 On the other hand, it is simple to see that $\ul{\sE}[\b1_{[h, \infty)}]=0
 $ for any $h>0$, which means that $\ul{\sE}$ does not satisfy strict monotonicity either.
 Furthermore, $\ul{\sE}$ does not satisfy (H1). For $\xi=1$ and $A_n=[\frac{1}{n},\infty)$, $n \in \hN$, we have that
 \beas
 \hspace{3.4cm} \underset{n \to \infty  }{\lim}\neg \ua \ul{\sE}
  [ \b1_{A_n}\xi]=\underset{n \to \infty  }{\lim}\ul{\sE}
  [ \b1_{[\frac{1}{n},\infty)}]=0 \ne 1 = \underset{Q \in \cP}{\inf}
  \wt{E}_Q[1]=\ul{\sE}[1]=\ul{\sE}[\xi]. \hspace{3.3cm} \hb{ \qed }
 \eeas
 \end{eg}

  Although it does not satisfy strict monotonicity, $\ul{\sE}$ is
almost an $\bF$-expectations on $Dom(\sE)$ as we will see next.

\begin{prop} \label{prop_olsE}
For any $t \in  [0,T]$,  $ \ul{\sE}[\cd|\cF_t] $ is an operator from
$Dom(\sE)$ to $Dom_t(\sE) \dfnn Dom(\sE)\cap L^0(\cF_t)$. Moreover,
the family of operators $ \big\{ \ul{\sE}[\cd|\cF_t] \big\}_{t \in
[0,T]}$ satisfies (A2)-(A4) as well as
 \bea \label{eq: mono_ulsE}
   \ul{\sE}[\xi|\cF_t] \le  \ul{\sE}[\eta|\cF_t], \q a.s. \q \hb{for
   any }\xi,\;\eta \in Dom(\sE) \hb{ with } \xi \le \eta,~a.s.
 \eea
\end{prop}

\ms \no {\bf \large Remark 2.}

\ms We have found that the first time $\ol{\t}(0)$ when the Snell
envelope $Z^0$ meets the process $Y$ is an optimal stopping time for
\eqref{eq:defn-coopgm} while the first time $\t_V$ when the process
$V^0$ meets the process $Y$ is an optimal stopping time for
\eqref{eq:defn-non-coopgm}. It is natural to ask whether
$\ol{\t}(0)$ (resp. $\t_V$) is the minimal optimal stopping time of
\eqref{eq:defn-coopgm} (resp. \eqref{eq:defn-non-coopgm}).
 This answer is affirmative when $\sE$ is a singleton. Let $\cE$ be a positively-convex $\bF$-expectation satisfying
(H1)-(H3) and let $Y$ be a right-continuous $\bF$-adapted process
satisfying (Y1), (Y3) and the following
 \beas
   \underset{ \rho \in  \cS_{0,T} }
 {\sup} \cE\left[Y_\rho \right] <  \infty; \hb{ if $\sE$ does not
 satisfy (\ref{ass_fatou}), then }
    \underset{( \rho, t) \in  \cS_{0,T} \times \cD_T}
 {\esssup}\, \wt{\cE}[Y_\rho |\cF_t]  \in  Dom^\#(\cE).
 \eeas
(Note that we have here merged the cumulative reward process $H$
into the primary reward process $Y$.) If $\t \in \cS_{0,T}$ is an
optimal stopping time for \eqref{eq:defn-coopgm}, i.e. $\underset{
\rho  \in  \cS_{0,T}}{\sup}
 \; \cE [Y_{\rho}  ] =  \cE [Y_{\t} ] $,
 Proposition \ref{Z_RCLL} and \eqref{eqn-h42} imply that
 \beas
  \underset{ \rho  \in  \cS_{0,T}}{\sup}
 \; \cE [Y_{\rho}  ]=  \underset{ \rho  \in  \cS_{0,T}}{\sup}
 \; \wt{\cE} [Y_{\rho}  ] = Z(0) \ge \wt{\cE} [Z(\t)] =\wt{\cE} \big[Z^0_\t\big]=\cE \big[Z^0_\t\big] \ge \cE
[Y_\t]= \underset{ \rho  \in  \cS_{0,T}}{\sup}
 \; \cE [Y_{\rho}  ] ,
 \eeas
thus $\cE [Z^0_\t] = \cE [Y_\t] $. The second part of (A1) then
implies that $Z^0_\t=Y_\t$, a.s. Hence $\ol{\t}(0) \le \t$, a.s.,
which means that $\ol{\t}(0)$ is the minimal stopping time for
\eqref{eq:defn-coopgm}.

\ss However, this is not the case in general. Let $\sE=\{\cE_i\}_{i
\in \cI}$ be a stable class of $\bF$-expectations and let $Y$ be a
right-continuous $\bF$-adapted process satisfying (Y1)-(Y3). We take
 $H^i \equiv 0$ for any $i \in \cI$. If $\t \in \cS_{0,T}$ is an
optimal stopping time for \eqref{eq:defn-coopgm}, i.e. $\underset{
(i,\rho)  \in  \cI \times \cS_{0,T}}{\sup}
 \; \cE_i [Y_{\rho}  ] = \underset{
 i   \in  \cI  }{\sup} \, \cE_i [Y_{\t} ]$,
\eqref{D_2} and \eqref{eqn-h42} then imply that
 \beas
  \underset{ (i, \rho)  \in \cI \times \cS_{0,T}}{\sup}
 \; \cE_i [Y_{\rho}  ]&= & \underset{ (i, \rho)  \in \cI \times \cS_{0,T}}{\sup}
 \; \wt{\cE}_i [Y_{\rho}  ] = Z(0) \ge \underset{  i   \in \cI }{\sup}\, \wt{\cE}_i [Z(\t)]
 = \underset{  i   \in \cI }{\sup}\,  \cE_i [Z(\t)] \\
 &=& \underset{  i   \in \cI }{\sup}\, \cE_i \big[Z^0_\t\big]
\ge \underset{  i   \in \cI }{\sup}\, \cE_i [Y_\t]= \underset{
(i,\rho) \in  \cI \times  \cS_{0,T}}{\sup}
 \; \cE_i [Y_{\rho}  ] ,
 \eeas
thus $ \ol{\sE}\big[Z^0_\t\big] = \underset{  i   \in \cI }{\sup}\,
\cE_i \big[Z^0_\t\big] = \underset{  i   \in \cI }{\sup}\, \cE_i
[Y_\t] =\ol{\sE}[Y_\t] $. However, this may not imply that $Z^0_\t=
Y_\t$, a.s. since $\ol{\sE}$ does not satisfy strict monotonicity as
we have seen in Example \ref{eg_olsE}.

\ms Now we further assume that $Y$ satisfies (Y2'), if $\t' \in
\cS_{0,T}$ is an optimal stopping time for
\eqref{eq:defn-non-coopgm}, i.e. $ \underset{\rho \in
\cS_{0,T}}{\sup}\underset{  i  \in \cI  }{\inf}
 \; \cE_i [Y_{\rho}  ]$ $= \underset{
 i   \in  \cI  }{\inf} \, \cE_i [Y_{\t'} ]$, \eqref{eqn-m31} and
Theorem \ref{V_process} imply that
 \beas
 \underset{\rho
\in \cS_{0,T}}{\sup}\underset{  i  \in \cI  }{\inf}
 \; \cE_i [Y_{\rho}  ]&= & \underset{\rho
\in \cS_{0,T}}{\sup}\underset{  i  \in \cI  }{\inf}
 \; \wt{\cE}_i [Y_{\rho}  ] =\ul{V}(0) = V(0) \ge \underset{  i   \in \cI }{\inf}\, \wt{\cE}_i [V(\t')]
 = \underset{  i   \in \cI }{\inf}\,  \cE_i [V(\t')] \\
&\ge& \underset{  i   \in \cI }{\inf}\, \cE_i [Y_{\t'}]=
\underset{\rho \in \cS_{0,T}}{\sup}\,\underset{  i  \in \cI  }{\inf}
 \; \cE_i [Y_{\rho}  ],
 \eeas
 thus $\ul{\sE}[V(\t')] = \underset{  i   \in \cI }{\inf}\,  \cE_i [V(\t')]
 = \underset{  i   \in \cI }{\inf}\, \cE_i [Y_{\t'}] =\ul{\sE}[Y_{\t'}] $. However,
this may not imply that $V(\t')= Y_{\t'}$, a.s. since $\ul{\sE}$
does not satisfy strict monotonicity, which we have also seen in
Example \ref{eg_olsE}. \big(If $V(\t')$ were a.s. equal to
$Y_{\t'}$, for any $i \in \cI$, applying \eqref{D_1} to singleton
$\{\cE_i\}$, we would deduce from \eqref{V0} and Lemma \ref{lem_ess}
that
 \beas
 V^0_{\t' \land \t_V} &=& V( \t' \land \t_V) = \ol{V}(\t' \land \t_V) = \underset{i \in \cI}{\essinf}\,R^i(\t' \land
\t_V)=\underset{i \in \cI}{\essinf}\, \big( \b1_{\{\t' \le \t_V\}}
R^i( \t' )+ \b1_{\{\t' >  \t_V\}} R^i(\t_V) \big)\\
&=& \b1_{\{\t' \le \t_V\}} \underset{i \in \cI}{\essinf}\, R^i( \t'
)+ \b1_{\{\t' > \t_V\}}\underset{i \in \cI}{\essinf}\, R^i(\t_V)=
\b1_{\{\t' \le \t_V\}} \ol{V}( \t' )+ \b1_{\{\t' > \t_V\}}\ol{V}
 (\t_V)\\
& = & \b1_{\{\t' \le \t_V\}} V( \t' )+ \b1_{\{\t' > \t_V\}}V(\t_V)
=\b1_{\{\t' \le \t_V\}} V( \t' )+ \b1_{\{\t' > \t_V\}}V^0_{\t_V} \\
&=& \b1_{\{\t' \le \t_V\}} Y_{\t'}+ \b1_{\{\t' > \t_V\}} Y_{\t_V} =
Y_{\t' \land \t_V}, \q a.s.,
 \eeas
which would further imply that $\t_V = \t' \land \t_V $, a.s., thus
$\t_V \le \t' $, a.s.\big)

 \section{Applications} \label{ch_app}

\setcounter{equation}{0}

\ms

In this section, we take a $d$-dimensional Brownian motion $B$ on
the probability space $(\O,\cF, P)$ and consider the Brownian
filtration generated by it:
 \bea
   \label{defn_BM}
  \bF= \Big\{\cF_t \dfnn \si\Big( \si\big(B_s; s\in [0,t]\big) \cup \cN \Big) \Big\}_{t \in
  [0,T]}, \hb{ where $\cN$ collects all $P$-null sets in $\cF$.}
  \eea
  We also let $\sP$ denote the predictable $\si$-algebra with respect
to $\bF$.

\subsection{Lipschitz $g$-Expectations}  \label{subsection_gexp}

\ss

   Suppose that a ``generator" function $g=g(t,\o, z): [0,T] \times \O
\times \hR^d \mapsto \hR$
 satisfies
 \bea
 \label{g-cond}
\tneg\tneg \left\{ \ba{l}
(i)\phantom{i}  ~ g(t, \o,  0)= 0, ~\;  \dtp  \\
(ii) ~ g \mbox{ is Lipschitz in } z \mbox{ 
for some } K_g>0: \mbox{ it holds $\dtp$ that} \\
\qq \qq |g(t, \o, z_1)\neg-\neg g(t, \o, z_2)|\le K_g |z_1\neg -\neg
z_2|,
 \q  \fa   z_1, z_2 \in  \hR^d.
 \ea \right.
 \eea
For any $\xi \in L^2(\cF_T)$, it is well known from \cite{PP-90}
that the backward stochastic differential equation (BSDE)
 \bea \label{BSDE}
  \G_t =\xi+\int_t^T g(s, \Th_s)ds-\int_t^T  \Th_s dB_s, \qq t \in [0,T]
 \eea
admits a unique solution $\big(\G^{\xi, g}, \Th^{\xi, g}\big) \in
\hC^2_\bF([0,T])\times \cH^2_\bF([0,T];\hR^d)$ (for convenience, we
will denote (\ref{BSDE}) by BSDE$(\xi, g)$ in the sequel),
 based on which \cite{Peng-97} introduced the so-called
``$g$-expectation" of $\xi$ by
 \bea   \label{def_g_exp}
 \cE_g[\xi|\cF_t] \dfnn \G^{\xi, g}_t , \q t \in [0,T].
 \eea

 \ss To show that the $g$-expectation $\cE_g$ is
an $\bF$-expectation with domain $Dom(\cE_g)=L^2(\cF_T)$, we first
note that $L^2(\cF_T) \in \wt{\sD}_T$. 
The (strict) Comparison Theorem for
BSDEs (see e.g. \cite[Theorem 35.3]{Peng-97}) then shows that (A1)
holds for the family of operators $\big\{\cE_g[ \cd |\cF_t]:
L^2(\cF_T) \mapsto L^2(\cF_t)\big\}_{t \in [0,T]}$, while the
uniqueness of the solution $(\G^{\xi, g},\Th^{\xi, g})$ to the
BSDE$(\xi, g)$ implies that the family $\big\{\cE_g[ \cd |\cF_t]
\big\}_{t \in [0,T]}$ satisfies (A2)-(A4) (see e.g. \cite[Lemma
36.6]{Peng-97} and \cite[Lemma 2.1]{CHMP}). Therefore, $\cE_g$ is an
$\bF$-expectation with domain $Dom(\cE_g)=L^2(\cF_T)$.

 \ss Moreover, the generator $g$ characterizes $\cE_g$ in the following ways:

\ss \no (1) Theorem 3.2 of \cite{Jiang_2008} (see also Proposition 10 of \cite{Eman}) shows that
   $\cE_g[ \cd |\cF_t] $ is a convex (resp. concave) operator on
$L^2(\cF_T)$ for any $t \in[0,T]$ if and only if  the generator $g$ is {\it convex} (resp. {\it concave}) 
in $z$, i.e., it holds $\dtp$ that
 \bea \label{def-g-convex}
  g(t, \l z_1 + (1 - \l)z_2)\neg \le \; (\hb{resp.}\, \ge ) ~   
 \l g(t,  z_1)+(1-\l)g(t, z_2), \q  \fa \l \in (0,1), ~ \fa z_1 ,z_2 \in \hR^d .
 \eea

 \no (2) Let $\tilde{g}$ be another generator satisfying
(\ref{g-cond}). If it holds $\dtp$ that
 \beas
 g(t,z) \ge  \tilde{g}(t,z), \q  \fa  z \in \hR^d,
 \eeas
  the Comparison Theorem for BSDEs (see e.g. \cite{EPQ-97} or \cite[Proposition 5.1]{Pln}) shows that for any
$\xi \in L^2(\cF_T)$ and $t \in [0,T]$
 \bea \label{eqn-y11}
\cE_g[\xi|\cF_t]\ge\cE_{\tilde{g}}[\xi|\cF_t], \q a.s.
 \eea
 Thanks to Theorem 4.1 of \cite{BCHMP}, the reverse statement also holds given that almost surely, the mapping $t  \to g(t,z)   $ is continuous for any $z \in \hR^d$.

  \ss  \no (3)  $   g^-(t, \o, z) \dfnn -g(t, \o, -z)$, $ (t, \o, z)
 \in [0,T]  \times \O \times \hR^d$
  also satisfies (\ref{g-cond}). Its corresponding $g$-expectation $\cE_{g^-}$ relates to
  $\cE_g$ in that for any $ \xi \in
    L^2(\cF_T)$ and $t \in [0,T]$
    \bea \label{eqn-g-neg}
    \cE_{g^-}[\xi|\cF_t]= -\cE_g[-\xi|\cF_t],\q  a.s.
    \eea
 \big(In fact, multiplying both sides of BSDE$(-\xi,g)$ by $-1$, we see
 that the pair $\big(-\G^{-\xi, g}, -\Th^{-\xi, g}\big) $ solves the
 BSDE$(\xi, g^-)$.\big)

 \ms To show that the
$g$-expectation $\cE_g$ satisfies (H0)-(H3), we need two basic
inequalities it satisfies.  

\begin{lemm} \label{g-inequ}
Let $g$ be a generator satisfying (\ref{g-cond}).

\ss \no (1) For any $\xi \in L^2(\cF_T)$, we have
 \beas
  \Big\| \underset{t \in [0,T]}{\sup} \big|\cE_g[\xi|\cF_t] \big|
  \Big\|_{L^2(\cF_T)} +\left\|\Th^{\xi,g} \right\|_{L^2_\bF([0,T]; \hR^d)}
     \le C e^{(K_g+K^2_g)T}  \|\xi\|_{L^2(\cF_T)},
 \eeas
where $C$ is a universal constant independent of $\xi$ and $g$.

\ms \no (2) For any $\mu \ge K_g$ and $\xi, \eta \in L^2(\cF_T)$, it
holds a.s. that
 \beas
  \big| \cE_g[\xi|\cF_t]- \cE_g[\eta|\cF_t] \big|
  \le  \cE_{g_\mu}\big[|\xi-\eta|\big|\cF_t\big], \q   \fa t \in [0,T],
 \eeas
 where the generator $g_\mu$ is defined by $g_\mu(z) \dfnn \mu |z|$, $  z \in
\hR^d$.
\end{lemm}

\ss \no {\bf Proof:} 
 A simple application of
\cite[Proposition 2.2]{BCHMP} yields (1).
On the other hand, (2) is a mere generalization of \cite[Proposition
3.7, inequality (60)]{Pln} by taking into account the continuity of
processes $\cE_g[\xi|\cF_\cd] $ and $\cE_{g_\mu}[\xi|\cF_\cd] $ for
any $\xi \in L^2(\cF_T)$. \qed

\begin{prop}  \label{gexp}
Let $g$ be a generator satisfying (\ref{g-cond}). Then
$\cE_g$ 
 satisfies (H0)-(H3).
\end{prop}

 \begin{rem} \label{rem_gexp_H3}
  Since $\cE_g[ \xi |\cF_\cd]$ is a continuous process
for any $\xi \in L^2(\cF_T)$, we see from \eqref{cE_tau} that
$\wt{\cE}_g[\cd |\cF_\nu]$ is just a restriction of $\cE_g[\cd
|\cF_\nu]$ to $L^{2,\#}(\cF_T)
 \dfnn \{ \xi \in L^2(\cF_T): \xi \ge c, \;a.s. \hb{~for some }c = c(\xi) \in \hR
 \}$ for any $\nu \in \cS_{0,T}$.
  \end{rem}

\ss Thanks to Proposition~\ref{gexp}, all results on
 $\bF$-expectations $\cE$ and $\wt{\cE}$ in Section \ref{ch_2} are
 applicable to $g$-expectations. In the following example we deliver the promise we made
 in Remark~\ref{rem_fatou}. This example indicates that for
 some $g$-expectations, $\underset{n \to \infty}{\liminf}
 \cE_g[\xi_n]< \infty $ is not a sufficient condition for $ \underset{n \to \infty}{\lim} \xi_n
 \neg \in \neg Dom^+(\cE_g)=L^{2,+}(\cF_T)\dfnn \neg \{\xi \in L^2(\cF_T):  \xi \ge 0, ~a.s.\} $
  given that $\{\xi_n\}_{n \in \hN}$ is an a.s. convergent sequence in $  Dom^+(\cE_g)
  $.

\begin{eg} \label{counter_eg_fatou}
Consider a probability space $([0,1], \sB[0,1],\l)$, where $\l$ is
the Lebesgue measure on $[0,1]$. We define a generator
$\tilde{g}(z)\dfnn -|z|$, $z \in \hR^d$. For any $n \in \hN$, it is
clear that the random
 variable $\big\{\xi_n(\o)\dfnn \o^{-\frac{1}{2}+\frac{1}{n+2}}\big\}_{\o \in [0,1]}
  \in  L^{2,+}(\cF_T)=Dom^+(\tilde{g})$. Proposition \ref{3addition} (2)
  then implies that
  \beas
 0= \cE_{\tilde{g}}[0] \le \cE_{\tilde{g}}[\xi_n] = \G^{\xi_n \neg,\hspace{0.3mm} \tilde{g}}_0= \xi_n-\int_0^T
 |\Th^{\xi_n \neg,\hspace{0.3mm} \tilde{g}}_s|ds-\int_0^T
  \Th^{\xi_n \neg,\hspace{0.3mm} \tilde{g}}_s dB_s \le  \xi_n-\int_0^T
  \Th^{\xi_n \neg,\hspace{0.3mm} \tilde{g}}_s dB_s.
  \eeas
Taking the expected value of the above inequality yields that
 \bea \label{eqn-bxb01}
  0 \le \cE_{\tilde{g}}[\xi_n] \le  E\big[ \xi_n-\hb{$\int_0^T
  \Th^{\xi_n \neg,\hspace{0.3mm} \tilde{g}}_s dB_s$}\big]=E[ \xi_n]=\int_0^1
  \o^{-\frac{1}{2}+\frac{1}{n+2}}
  d\o=\frac{1}{\frac{1}{2}+\frac{1}{n+2}}<2.
 \eea
 Since $\{\xi_n\}_{n \in \hN}$ is an increasing sequence, we can deduce from (A1) and (\ref{eqn-bxb01}) that
  $\,0 \le \underset{n \to \infty}{\lim} \neg \ua \cE_{\tilde{g}}[\xi_n]\le
  2$. However, $\underset{n \to \infty}{\lim} \neg \ua \xi_n =\big\{\o^{-\frac{1}{2}}\big\}_{\o \in
  [0,1]}$ does not belong to $L^{2,+}(\cF_T)=Dom^+(\tilde{g})$. \qed
 \end{eg}

 Similar to Proposition \ref{tau_A}, pasting two $g$-expectations at any
 stopping time generates another $g$-expectation.

\begin{prop} \label{tau_A_g}
Let $g_1$, $g_2$ be two generators satisfying (\ref{g-cond}) with
Lipschitz coefficients $K_1$ and $K_2$ respectively. For any $\nu
\in \cS_{0,T}$, we define the pasting of $ \cE_{g_1}, \cE_{g_2}$ at
$\nu$ to be the following continuous $\bF$-adapted process
 \bea \label{g_pasting}
  \cE^{\nu }_{g_1,g_2} [\xi  |\cF_t ] \dfnn
\b1_{\{\nu\le t\}} \cE_{g_2}[\xi |\cF_t]  + \b1_{\{\nu>
t\}}\cE_{g_1}\big[  \cE_{g_2} [\xi |\cF_\nu]  \big|\cF_t\big], \q
\fa t \in [0,T]
 \eea
 for any $\xi \in L^2(\cF_T)$. Then $\cE^{\nu}_{g_1,g_2}$ is exactly
 the $g$-expectation $\cE_{g^\nu}$ with
 \bea   \label{def_g_nu}
 g^{\nu}(t,\o, z)  \dfnn \b1_{\{  \nu (\o) \le t \} }  g_2(t,\o,z)+ \b1_{\{ \nu (\o)>t \}}
    g_1(t,\o,z)  , \q   (t,\o,z) \in [0,T]\times \O \times \hR^d,
 \eea
 which is a generator satisfying (\ref{g-cond}) with the Lipschitz
coefficient $K_1 \vee K_2$.
 \end{prop}

\ms Fix $M>0$, we denote by $\sG_M$ the collection of all
 convex generators $g$ satisfying
 (\ref{g-cond}) with Lipschitz coefficient $K_g \le M$.
  Proposition \ref{tau_A_g} shows that the family of convex $g$-expectations $ \sE_M \dfnn \{\cE_g\}_{g \in \sG_M }$
  is closed under the pasting \eqref{g_pasting}. To wit, $\sE_M$ is a stable
  class of $g$-expectations in the sense of Definition
 \ref{def_stable_class}. In what follows we let $\sG'$ be a non-empty subset of $\sG_M$ such
that $\sE' \dfnn \{\cE_g\}_{g \in \sG'}$ is closed under pasting.
  \if{0}
   Correspondingly, $\{\wt{\cE}_g\}_{g \in
  \sG'}$ is closed under the general pasting \eqref{tau_ij},
  thus we see from Proposition \ref{gexp} 
  that $\sE'$ is a stable class in the sense of Definition
 \ref{def_stable_class}. In light of
 Remark \ref{rem_gexp_H3}, when applying the optimal stopping theory
 of $\bF$-expectations (as presented in Section \ref{co_game} and
 \ref{nonco_game}) to the stable class $\sE'$, we can replace
 $\wt{\cE}_g$ by $\cE_g$ for each $g \in \sG'$. However, such
extension does not mean that we can take away (S2) and (Y3), the
lower-bound requirements on reward processes $Y$ and $\sH$
respectively.
In fact, as we see in the proof of Proposition \ref{Z_RCLL}, the
deduction of the RCLL modification $Z^0$ of the process
$\big\{Z(t)\big\}_{t\in [0,T]}$ totally relies on applying Corollary
\ref{cor_RCLL} to certain process $\big\{Z^j(t)\big\}_{t\in [0,T]}$.
This requires that the process $\big\{Z^j(t)\big\}_{t\in [0,T]}$ is
bounded from below by some $c \in \hR$, which is in turn assured by
(S2) and (Y3) \big(see \eqref{Z_LB}\big). {\bf \large (compare
Remark ); (H3) becomes futile, which is used to assume
 that $\wt{\cE}[\xi|\cF_\nu] \in Dom^\#(\cE)$ for any $\xi \in
Dom^\#(\cE)$ and $\nu \in \cS_{0,T}$.}
 \fi
Now we make the following assumptions on the reward processes:

 \ms \no {\bf Standing assumptions on the reward processes in this subsection.}
 Let $Y$ be a continuous $\bF$-adapted process with
 \bea \label{ass_Y_gexp}
 \z'_Y  \dfnn \Big(\underset{t  \in \cD_T }{\esssup}\, Y_t \Big)^+ \in L^2(\cF_T) ~\;\;
 \eea
 and satisfying (Y3). Moreover, for any $g \in \sG'$, we suppose that the
  model-dependent cumulative reward process is in the form of
 \beas
 H^g_t \dfnn \int_0^t h^g_s ds , \q \fa t \in [0,T],
 \eeas
 where $\{ h^g_t, t \in [0,T] \}_{g \in
\sG'}$ is a family of predictable processes satisfying the following
assumptions:

 \ss \no ($\tilde{h}$1) There exists a $\,c' < 0 $ such that for any
$g \in \sG'$, $h^g_t \ge c'$, $\dtp$

 \ss \no ($\tilde{h}$2) 
 The random variable $ \o \mapsto \int_0^T   h'(t,\o) \, dt $
 belongs to $L^2(\cF_T)$ with $ h'(t,\o) \dfnn \Big(\underset{g \in
\sG'  }{\esssup}\,  h^g_t(\o) \Big)^+ $   \big(the essential
supremum is taken
with respect to the product measure space $([0,T]\times \O,\, \sP,
\l \times P )$, 
  where $\l$ denotes the Lebesgue measure on $[0,T]$\big).

\ss \no ($\tilde{h}$3) For any $\nu \in \cS_{0,T}$ and $g_1, g_2 \in
 \sG'$, it holds for any $0 \le s < t \le T $ that
  \beas
   h^{g^\nu}_t   = \b1_{\{\nu \le t \}} h^{g_2}_t+ \b1_{\{\nu > t \}} h^{g_1}_t    ,  \q
   \dtp,
  \eeas
  where $g^\nu$ is defined in \eqref{def_g_nu}.

\ms Then the triple $\big(\sE', \sH' \dfnn \{H^g \}_{g \in \sG'} ,
Y\big)$ satisfies all the conditions stated in Section \ref{co_game}
and \ref{nonco_game}. Thus we can carry out the optimal stopping
theory developed for $\bF$-expectations to $(\sE', \sH', Y)$ as we
will see next.

 \begin{thm} \label{prop_g_result}
The stable class $\sE'$ satisfies \eqref{ass_odd},  the family of
processes $\sH'$ satisfies (S1') (thus (S1) see Remark
\ref{rem_nonco_YH}), (S2) and (S3), while the process $Y$ satisfies
(Y1), \eqref{eqn-m48} \big(thus (Y2'), again by Remark
\ref{rem_nonco_YH}\big) and (Y3). Moreover, the family of processes
$\big\{Y^g_t \dfnn Y_t+H^g_t,~ t\in [0,T]\big\}_{g \in \sG'}$ is
both ``$\sE'$-uniformly-left-continuous" \big(thus satisfies
\eqref{eqn-q180}, see also Remark \ref{rem_nonco_YH}\big) and
``$\sE'$-uniformly-right-continuous".
  \end{thm}

 \subsection{Existence of an Optimal Prior in \eqref{eq:defn-coopgm} for $g$-Expectations}

 \ss

  For certain collections of $g$-expectations, we can even determine an
optimal generator $g*$ in the following sense:
 \beas
  \cE_{g*}[Y^{g*}_{\ol{\t}(0)}  ]=\underset{g \in \cG}{\sup} \, \cE_g  \big[ Y^g_{\ol{\t}(0) } \big]
  = \underset{(g,\rho) \in \cG \times \cS_{0,T}}{\sup} \; \cE_g[Y^g_{\rho}        ],
 \eeas
where the optimal stopping time $\ol{\t}(0)$ is defined as in
Theorem \ref{SN_exist}.

 \ms Let $S$ be a separable metric space with metric $|\cd|_S$ such that $S$ is a
countable union of non-empty compact subsets. We denote by $\fS$ the
Borel $\si$-algebra of $S$ and take $\cH^0_\bF([0,T]; S)$ as the
space
  of admissible control strategies. For any $U \in \cH^0_\bF([0,T]; S)$,
  we define the generator
   \bea \label{defn_gU}
   g_{{}_U}(t, \o, z) \dfnn g^o(t, \o, z, U_t(\o)),
 \eea
  where the function $g^o(t, \o, z, u): [0,T]
\times \O \times \hR^d \times S \mapsto \hR$ satisfies:

  \ss \no ($g^o1$)  $g^o$ is $\sP \otimes \sB(\hR^d) \otimes \fS/
 \sB(\hR)$-measurable.

 \ss \no ($g^o2$)  It holds $\dtp$ that
   $   g^o(t, \o,  0, u)= 0$ for any $       u \in S$.

 \ss \no ($g^o3$)  $g^o$ is Lipschitz in $z$: For some $K_o>0$, it holds $\dtp$ that
 \beas
 |g^o(t, \o, z_1, u)-  g^o(t, \o, z_2, u)|\le K_o |z_1 - z_2|,
 \q  \fa   z_1, z_2 \in  \hR^d, ~\fa u \in S.
 \eeas

  \no ($g^o4$)  $g^o$ is convex in $z$: It holds $\dtp$ that
  \beas
  g^o(t, \o, \l z_1 \neg + \neg(1 \neg- \neg \l)z_2, u)\le \l g^o(t, \o, z_1, u)\neg
  + \neg (1\neg -\neg \l) g^o(t, \o, z_2,
  u), ~\;    \fa \l \in (0,1),\; \fa  z_1, z_2 \in  \hR^d, \; \fa u \in
  S. \q
  \eeas

 \ss Now fix a non-empty subset $\fU$ of $\cH^0_\bF([0,T]; S)$ that preserves ``pasting",
 i.e.,
  for any $\nu \in \cS_{0,T}$ and $U^1,U^2 \in \fU$,
 \bea      \label{def_U_nu}
   U^\nu_t(\o) \dfnn  \b1_{\{  \nu (\o) \le t \} }  U^2_t(\o)+ \b1_{\{
\nu (\o)>t \}} U^1_t(\o) , \q   (t,\o) \in [0,T]\times \O,
 \eea
 also belongs to $\fU$. Then it is easy to check that $\{\cE_{g_{{}_U}} \}_{ U \in \fU} \subset \sE_{K_o} $
forms a stable class of $g$-expectations.

 \ms Let $Y$ still be a continuous $\bF$-adapted process
 satisfying \eqref{ass_Y_gexp} and (Y3). For any $U \in \fU$,
 assume that the model dependent reward process has a density which is given by
 \beas
 h^U_t(\o) \dfnn h(t, \o, U_t(\o)), \q (t,\o) \in [0,T] \times \O,
 \eeas
where $h(t, \o, u): [0,T] \times \O \times S \mapsto \hR $ is a $\sP
\otimes \fS /\sB(\hR)$-measurable function satisfying the following
assumptions:

 \ss \no ($\hat{h}$1) For some $c < 0 $, it holds $\dtp$ that
 $  h(t, \o, u) \ge c$ for any  $u \in S$.

 \ss \no  ($\hat{h}$2) The random variable $\o \mapsto \int_0^T
\hat{h}(t,\o) dt$ belongs to
 $L^2(\cF_T)$ with $ \hat{h}(t,\o) \dfnn \Big(\underset{U \in \fU}{\esssup}\, h^U_t(\o)\Big)^+ $
 \big(the essential supremum is
taken
with respect to the product measure space $([0,T]\times \O,\, \sP, 
 \l \times P )$, where $\l$ denotes the Lebesgue
measure on $[0,T]$\big).

  \ms   It is easy to see that $\{ h^U_t, t \in [0,T] \}_{U \in \fU}$ is
 a family of predictable processes satisfying ($\tilde{h}$1)-($\tilde{h}$3).
 Hence, we can apply the optimal stopping theory
developed for $\bF$-expectations to the triple
$\big(\{\cE_{g_{{}_U}} \}_{ U \in \fU}, \{ h^U \}_{U \in \fU}, Y
\big)$ thanks to Theorem \ref{prop_g_result}.
  Now let us construct a so-called {\it Hamiltonian} function
 \beas
  H(t, \o, z, u) \dfnn g^o(t,\o, z, u)+h(t,\o,u), \q (t, \o, z, u) \in
  [0,T] \times \O \times  \hR^d \times S.
 \eeas
   We assume that for any $(t, \o, z) \in [0,T] \times \O \times
 \hR^d$, there exists a $u=u^*(t, \o, z) \in S$ such that
  \bea   \label{eqn-ttt71}
   \underset{u \in S}{\sup} H(t, \o, z, u) = H\big(t, \o, z, u^*(t, \o, z) \big).
  \eea
 (This is valid, for example, when the metric space $S$ is compact
 and the mapping $u \mapsto H(t, \o, z, u)$ is continuous.)
 Then it can be shown (see \cite[Lemma 1]{Benes_1970} or
\cite[Lemma 16.34]{Elliott_1982}) that the mapping $u^*: [0,T]
\times \O \times \hR^d \mapsto S$ can be selected to be $\sP \otimes
\sB(\hR^d)/\fS$-measurable.

\ms The following theorem is the main result of this subsection.

\begin{thm} \label{prop_optim_control}
There exists a $\,U^* \in \fU$ such that
    $
  \underset{(U, \rho) \in \fU \times \cS_{0,T}}{\sup}
 \ \cE_{g_{{}_U}}\big[Y^U_\rho  \big]
 =\cE_{g_{{}_{U^*}}} \big[Y^{U^*}_{\ol{\t}(0)} \big]
 $,
  where the stopping time $\ol{\t}(0)$ is as in Theorem \ref{SN_exist}.
\end{thm}

\subsection{The Cooperative Game of Karatzas and Zamfirescu [2006]
Revisited} \label{subsection_karatzas}


\ms

  In this subsection, we apply the results of the last subsection
  to extend the results of
  \cite{Kara_Zam_2006}. Let us first recall their setting:

 \ms \no $\bullet$
  Consider the canonical space
 $(\O,\cF)= \big(C([0,T]; \hR^d), \sB(C([0,T]; \hR^d)) \big)$ endowed with Wiener measure $P$, under which the coordinate
 mapping process $B(t,\o)=\o(t)$, $t \in [0,T]$ becomes a standard $d$-dimensional Brownian motion.
 We still take the filtration $\bF$ generated by the
  Brownian motion $B$ (see \eqref{defn_BM})
  and let $\sP$ denote the predictable $\si$-algebra with respect to $\bF$.

 \ms \no $\bullet$ It is well-known (see e.g. \cite[Theorem
14.6]{Elliott_1982})
 that given a $x \in \hR^d$, there exists a pathwise
unique, strong solution $X(\cd)$ of the stochastic equation
 \beas
  X(t)=x+\int_0^t \si(s,X) dB_s , \q   t \in [0,T],
 \eeas
   where the diffusion term $\si(t,\o)$ is a
   $ \hR^{d \times d}$-valued predictable process satisfying:

\ss \no ($\si$1) $\int_0^T |\si(t, \vec{0})|^2 dt < \infty $ and
$\si(t, \o)$ is nonsingular for any $(t,\o) \in [0,T]\times \O$.

\ss \no ($\si$2) There exists a $K>0$ such that for any $ \o,\wt{\o}
\in \O$ and $t \in [0,T]$
 \bea \label{eqn-ttt25}
 \|\si^{-1}(t,\o)\| \le
K ~\; \hb{and} ~ \; \big|\si_{ij}(t,\o)-\si_{ij}(t,\wt{\o})
 \big|\le K \|\o-\wt{\o}\|^*_t, \q \fa  1\le i,j\le n
 \eea
 with $\|\o \|^*_t\dfnn \underset{s \in [0,t]}{\sup} |\o(s)|$.

 \ms \no $\bullet$ Let $f(t, \o,  u): [0,T]
\times \O  \times S \mapsto \hR^d$ is a $\sP \otimes \fS/
\sB(\hR^d)$-measurable function such that

 \ss \no (\,i) For any $u \in S$, the
mapping $(t,\o) \mapsto f(t, \o,
 u)$ is predictable (i.e. $\sP$-measurable);

  \ss \no (ii) With the same $K$ as in
 \eqref{eqn-ttt25},
  \bea \label{eqn-ttt41}
     |f(t, \o,  u)| \le K \big(1+\|\o\|^*_t\big), \q \fa  (t, \o,
 u) \in [0,T] \times \O  \times S .
  \eea
   For any $U \in \wt{\fU} \dfnn \cH^0_\bF([0,T];
S)$, \cite[page 166]{Kara_Zam_2006} shows that
 \beas
  \frac{dP_U}{dP}  \dfnn \exp\bigg\{ \int_0^T  \lan \si^{-1}(t,X)f(t,X,U_t), dB_t \ran
 -\frac12 \int_0^T \big|\si^{-1}(t,X)f(t,X,U_t)\big|^2 dt \bigg\}
 \eeas
 defines a probability measure $P_U$.

 \ms The objective of \cite{Kara_Zam_2006} is to find an optimal stopping time $\t^* \in \cS_{0,T}$ and an optimal
 control strategy $U^* \in \wt{\fU}$ that maximizes the expected reward
 \beas
 E_U\bigg[\vf(X\big(\rho)\big)  + \int_0^\rho h (s, X, U_s) ds \bigg]
 \eeas
 over $(\rho, U) \in \cS_{0,T} \times \wt{\fU}$. Here $\vf: \hR^d \mapsto
  \hR$ is a bounded continuous function, and $h(t, \o, u): [0,T] \times \O  \times S \mapsto \hR$
 is a $\sP  \otimes \fS/ \sB(\hR)$-measurable function such that $|h(t, \o, u) | \le K$ for any $(t, \o,
 u) \in [0,T] \times \O  \times S $ \big(with the same $K$ that appears in
 \eqref{eqn-ttt25}\big).

 \ms Corollary 8 of \cite{Kara_Zam_2006} shows that under \eqref{eqn-ttt41},  the process
  \beas
      \wt{Z}(t) \dfnn  \underset{(U, \rho) \in
\wt{\fU} \times \cS_{t,T} }
 {\esssup}\, E_U\bigg[ \vf\big(X(\rho)\big) + \int_t^\rho h (s, X, U_s) ds  \bigg|\cF_t\bigg] , \q  t \in
 [0,T]
  \eeas
  admits an RCLL modification $\wt{Z}^0$, and that the first time
  processes $\wt{Z}^0$ and $\big\{\vf\big(X(t)\big)\big\}_{t\in
  [0,T]}$ meet with each other, i.e. $\ol{\t}(0)\dfnn   \inf\big\{t \in [0,T]
\;|\; \wt{Z}^0_t=\vf\big(X(t)\big)\big\} $, is an optimal stopping
time. That is,
 \bea \label{eqn-ttt43}
 \underset{(U, \rho) \in \wt{\fU} \times \cS_{0,T} }
 {\sup} E_U\bigg[ \vf\big(X(\rho)\big) + \int_0^\rho h (s, X, U_s) ds  \bigg]
 =\underset{ U  \in \wt{\fU}    }
 {\sup} E_U\bigg[ \vf\big(X(\ol{\t}(0))\big) + \int_0^{\ol{\t}(0)} h (s, X, U_s) ds
 \bigg].
 \eea

 \ss  Moreover, if for any $(t, \o, z) \in [0,T]
\times \O \times
 \hR^d$, there is a $u^*(t, \o, z) \in S$ which is  $\sP \otimes
\sB(\hR^d)/\fS$-measurable such that
  \bea   \label{eqn-ttt45}
   \underset{u \in S}{\sup} \, \wt{H}(t, \o, z, u) = \wt{H}\big(t, \o, z, u^*(t, \o, z) \big)
  \eea
with $\wt{H}(t, \o, z, u) \dfnn \big\lan \si^{-1}(t,\o)f(t,\o, u),z
\big\ran+ h(t, \o, u)$, $ (t, \o, z, u) \in [0,T] \times \O \times
\hR^d \times S$, then there further exists an optimal control
strategy $U^* \in \wt{\fU}$ (see \cite[Section 8]{Kara_Zam_2006})
such that
 \bea \label{eqn-ttt47}
 \underset{(U, \rho) \in \wt{\fU} \times \cS_{0,T} }
 {\sup} E_U\bigg[ \vf\big(X(\rho)\big) + \int_0^\rho h (s, X, U_s) ds  \bigg]
 =  E_{U^*}\bigg[ \vf\big(X(\ol{\t}(0))\big) + \int_0^{\ol{\t}(0)} h (s, X, U^*_s) ds \bigg].
 \eea

 In the main result of this subsection, we will show that the assumption of \cite{Kara_Zam_2006}
  that $\vf$ and $h$ are bounded from above by constants can be relaxed
 and replaced by linear-growth
 conditions. This comes, however, at the cost of strengthening the assumption stated in \eqref{eqn-ttt41}.
 \begin{prop} \label{prop_Karatzas}
With the same $K$ as in \eqref{eqn-ttt25}, we assume that
 \bea \label{eqn-ttt55}
  -K \le \vf(x) \le K|x|,  \q \fa  x \in \hR^d
  \eea
   and that for a.e. $t \in [0,T]$
 \bea \label{eqn-ttt51}
 |f(t, \o,  u)| \le K \q \hb{and} \q  -K \le   h(t, \o,u) \le K\|\o\|^*_T
  ,   \q  \fa (\o, u) \in \O \times S.
 \eea
Then the process $\big\{ \wt{Z}(t)\big\}_{t \in [0,T]} $ has a RCLL
modification $\wt{Z}^0$, and the first time $\ol{\t}(0)$ when
 the process $\wt{Z}^0$ meets the process $\big\{\vf\big(X(t)\big)\big\}_{t\in
  [0,T]}$ is an optimal stopping time; i.e., it  satisfies \eqref{eqn-ttt43}.
 Moreover, if there exists a measurable mapping $u^*: [0,T]
\times \O \times \hR^d \mapsto S$ satisfying \eqref{eqn-ttt45},
 then there exists an optimal control
strategy $U^* \in \wt{\fU}$ such that \eqref{eqn-ttt47} holds.
 \end{prop}

\subsection{Quadratic $g$-Expectations} \label{subs_g2_exp}

\ms

 Now we consider a \emph{quadratic} generator $\hat{g}=\hat{g}(t,\o, z): [0,T] \times \O
\times \hR^d \mapsto \hR$
 that satisfies
  \bea
  \label{g2-cond}
  \tneg\tneg \left\{ \ba{l}
  (i)\phantom{ii}  ~ ~  \hat{g}(t, \o,  0)= 0, ~\; \dtp  \\
  (ii)\phantom{i}~ \hb{For some $\k>0$,  it holds $\dtp$ that}  \\
   \qq \qq   \left|\frac{\pa \hat{g}}{\pa z}(t, \o, z)\right| \le \k(1+|z|),  \q  \fa   z  \in  \hR^d  . \\
 (iii)  ~  \hat{g} \hb{ is convex in $z$ in the sense of
 \eqref{def-g-convex}.}
 \ea \right.
 \eea
Note that under (ii), (i) is equivalent to the following statement:
It holds $\dtp$ that
  \bea \label{eqn-ttt03}
     | \hat{g} (t, \o, z)| \le \k\big(|z|+\hb{$\frac12$}|z|^2\big),  \q \fa   z  \in  \hR^d.
  \eea
 In fact, it is clear that \eqref{eqn-ttt03} implies (i). Conversely,
 for $\dtp$ $(t,\o) \in [0,T]\times \O$, one can deduce that for any
 $z \in \hR^d$,  $ |  \hat{g} (t, \o, z)| =\big|  \hat{g} (t, \o, z)-\hat{g} (t, \o, 0)\big|
 = \big| \int_0^1   \frac{\pa \hat{g}}{\pa z}(t, \l z) z   d\l \big|
 \le  \k \int_0^1  (1+\l|z|)  |z|   d\l
 =\k\big(|z|+\hb{$\frac12$}|z|^2\big) $.

\ms For any $  \xi \in L^{\neg e}(\cF_T)$, \cite[Corollary 6]{BH-07}
\big(where we take $f=g$, thus $\a(t) \equiv \frac{\k}{2}$ and
$(\beta, \g)=(0, 2\k )$\big) shows that the quadratic BSDE$( \xi,
\hat{g})$
 admits a unique solution $\big(\G^{\xi, \hat{g}}, \Th^{\xi, \hat{g}}\big) \in
 \hC^e_\bF([0,T]) \times  M_\bF([0,T];\hR^d)$.
 Hence we can correspondingly define the ``quadratic" $g$-expectation of $\xi$ by
  \beas
 \cE_{\hat{g}}[\xi|\cF_t] \dfnn \G^{\xi, {\hat{g}}}_t, \q t \in
 [0,T].
 \eeas

 \ss To show that the quadratic $g$-expectation $\cE_{\hat{g}}$ is an $\bF$-expectation with
domain $Dom(\cE_{\hat{g}})=L^{\neg e}(\cF_T)$, we first note that
$L^{\neg e}(\cF_T) \in \wt{\sD}_T$ \big(Clearly, $L^{\neg e}(\cF_T)$
satisfies (D1) and (D3) and $\hR \subset L^{\neg e}(\cF_T)$. For any
$\xi, \eta \in L^{\neg e}(\cF_T) $, $A \in \cF_T$ and $\l > 0$, we
have
 $ E\big[e^{\l |\b1_A\xi| }\big] \le  E\big[e^{\l | \xi| }\big]
 <\infty$ and $ E\big[e^{\l | \xi+\eta| }\big] \le  E\big[e^{\l | \xi | }e^{\l | \eta | }\big]
 \le \frac12 E\big[ e^{2\l | \xi | }\big]+\frac12 E\big[ e^{2\l | \eta |
 }\big] <\infty$, thus (D2) also holds for $L^{\neg e}(\cF_T)$\big).
Similar to the Lipschitz $g$-expectation case, the uniqueness of the
solution $(\G^{\xi, \hat{g}},\Th^{\xi, \hat{g}})$ to the quadratic
BSDE$(\xi, \hat{g})$ implies that the family of operators
$\big\{\cE_{\hat{g}}[ \cd |\cF_t]: L^{\neg e}(\cF_T) \mapsto L^{\neg
e}(\cF_t)\big\}_{t \in [0,T]} $ satisfies (A2)-(A4) (cf.
\,\cite[Lemma 36.6]{Peng-97} and \cite[Lemma 2.1]{CHMP}), while a
comparison theorem for quadratic BSDEs (see e.g. \cite[Theorem
5]{BH-07}) and the following proposition show that (A1) also holds
for the family $\big\{\cE_{\hat{g}}[ \cd |\cF_t] \big\}_{t \in
[0,T]}$.

\begin{prop} \label{thm_strict_com}
Let $\hat{g}$ be a quadratic generator satisfying \eqref{g2-cond}.
For any $\xi^1, \xi^2 \in L^{\neg e}(\cF_T)$, if $\xi^1 \geq \xi^2$,
a.s., then it holds a.s. that
 \bea \label{scomp1}
 \G^{\xi_1, \hat{g}}_t \ge  \G^{\xi_2, \hat{g}}_t,\q \fa t \in [0,T].
 \eea
 Moreover, if $\,\G^{\xi_1, \hat{g}}_\nu =  \G^{\xi_2, \hat{g}}_\nu
  $, a.s. for some $\nu \in \cS_{0,T}$, then
   \bea \label{scomp2}
   \xi^1 = \xi^2, \q a.s.
   \eea
\end{prop}

\ss Therefore, the quadratic $g$-expectation $\cE_{\hat{g}}$ is an
$\bF$-expectation with domain $Dom(\cE_{\hat{g}})=L^{\neg
e}(\cF_T)$. Similar to the Lipschitz $g$-expectation case, the
convexity \eqref{g2-cond}(iii) of the quadratic generator $\hat{g}$
as well as Theorem 5 of \cite{BH-07} determine that $\cE_{\hat{g}}[
\cd |\cF_t] $ is a convex operator on $L^{\neg e}(\cF_T)$ for any $t
\in[0,T]$. Hence, $\cE_{\hat{g}}$ satisfies (H0) thanks to Lemma
\ref{lem_pconvex}. To see $\cE_{\hat{g}}$ also 
satisfying (H1)-(H3), we need the following stability result.

 \begin{lemm} \label{lem_g2_stable}
 If $\xi_n \to \xi$, a.s. and $\,E\big[e^{\l |\xi|
 }\big]+\underset{n \in \hN}{\sup}E\big[e^{\l |\xi_n|
 }\big]< \infty $ for any $\l > 0$, then
 \bea \label{g2_stable}
 \underset{n \to \infty}{\lim}  E\bigg[   \underset{t \in [0,T]}{\sup}
 \Big|\cE_{\hat{g}}[\xi_n|\cF_t]-\cE_{\hat{g}}[\xi|\cF_t]\Big|   \bigg]=0 .
 \eea
 \end{lemm}

 \ss \no {\bf Proof: } Taking $f_n\equiv g$ and $f=g$ in Proposition 7 of \cite{BH-07} yields that
 \beas
 \underset{n \to \infty}{\lim}  E\bigg[\exp \Big\{ p \underset{t \in [0,T]}{\sup}
 \Big|\cE_{\hat{g}}[\xi_n|\cF_t]-\cE_{\hat{g}}[\xi|\cF_t]\Big|  \Big\}\bigg]=0, \q \fa p \ge 1.
 \eeas
 Then \eqref{g2_stable} follows since
  $
  E\Big[   \underset{t \in [0,T]}{\sup}
 \Big|\cE_{\hat{g}}[\xi_n|\cF_t]\neg-\neg \cE_{\hat{g}}[\xi|\cF_t]\Big|
 \Big] \neg \le \neg E\Big[\neg \exp \Big\{  \underset{t \in [0,T]}{\sup}
 \Big|\cE_{\hat{g}}[\xi_n|\cF_t] \neg-\neg \cE_{\hat{g}}[\xi|\cF_t]\Big|  \Big\} \Big]
 $ for any $n \in \hN$.  \qed

\begin{prop}  \label{g2exp}
Let $\hat{g}$ be a quadratic generator satisfying \eqref{g2-cond}.
Then the quadratic
$g$-expectation $\cE_{\hat{g}}$ 
 satisfies (H0)-(H3).
\end{prop}

Similar to Remark \ref{rem_gexp_H3}, since $\cE_{\hat{g}}[ \xi
|\cF_\cd]$ is a continuous process for any $\xi \in L^{\neg
e}(\cF_T)$, we see from \eqref{cE_tau} that $\wt{\cE}_{\hat{g}}[\cd
|\cF_\nu]$ is just a restriction of $\cE_{\hat{g}}[\cd |\cF_\nu]$ to
$L^{\neg e,\#}(\cF_T)
 \dfnn \{ \xi \in L^{\neg e}(\cF_T): \xi \ge c, \;a.s. \hb{~for some }c \in \hR
 \}$ for any $\nu \in \cS_{0,T}$. Therefore, all results on $\bF$-expectations
 $\cE$ and $\wt{\cE}$ in Section \ref{ch_2} work for quadratic $g$-expectations.

 \ms The next result, which shows the existence of an optimal stopping time for a quadratic $g$-expectation,
 is the main result of this subsection.

\begin{thm} \label{prop_g2exp_result}
Let $\hat{g}$ be a quadratic generator satisfying \eqref{g2-cond}.
For any right-continuous $\bF$-adapted process $Y$ with $ \hat{\z}_Y
\dfnn \Big(\underset{t \in \cD_T}{\esssup}\, Y_t \Big)^+ \neg \in
L^{\neg e}(\cF_T)$
 and satisfying (Y3), we have

  \hspace{6.3cm} $  \underset{  \rho \in  \cS_{0,T}}{\sup}  \, \cE_{\hat{g}}[Y_{\rho} ] =   \cE_{\hat{g}}[Y_{\ol{\t}(0)}
 ]$,

 \no where $\ol{\t}(0)$ is as in Theorem \ref{SN_exist}. 
  \end{thm}

\section{Proofs} \label{sec:Proofs}

\setcounter{equation}{0}

\subsection{Proofs of Section \ref{ch_2}}

\ss \no {\bf Proof of Proposition \ref{prop_Fexp}:}
For any $\xi \in \L$ and $t \in [0,T]$, let us define $
\cE^o[\xi|\cF_t]\dfnn \xi_t$. We will check that the system $
\big\{\cE^o[\xi |\cF_t] ,\, \xi \in \L \big\}_{t \in [0,T] }$
satisfies (A1)-(A4); thus it is an $\bF$-expectation with domain $\L$.

 \ms
\no 1) For any $ \eta \in \L $ with $  \xi \le \eta$, $a.s.$, we set
$A \dfnn \big\{\cE^o [\xi |\cF_t] >  \cE^o [\eta |\cF_t] \big\}\in
\cF_t$, thus $\b1_A \cE^o [\xi |\cF_t] \ge \b1_A \cE^o [\eta |\cF_t] $.
It follows from (a1) and (a2) that
 \beas
 \cE^o\big[\b1_A \cE^o [\xi |\cF_t]\big] \ge \cE^o\big[ \b1_A
 \cE^o [\eta |\cF_t]\big]= \cE^o [ \b1_A \eta ]\ge \cE^o [ \b1_A \xi ]
 =\cE^o\big[ \b1_A \cE^o [\xi |\cF_t]\big],
 \eeas
which shows that $\cE^o\big[\b1_A \cE^o [\xi |\cF_t]\big]
 =\cE^o\big[ \b1_A \cE^o [\eta |\cF_t]\big]$.
 Then the ``strict monotonicity" of (a1) further implies that
 $\b1_A \cE^o [\xi |\cF_t]= \b1_A \cE^o [\eta |\cF_t] $, a.s.,
  thus $P(A)=0$, i.e., $\cE^o [\xi |\cF_t] \le  \cE^o [\eta |\cF_t]$, a.s.

\ms Moreover, if $0 \le \xi \le \eta$, a.s. and $ \cE^o [\xi|\cF_0] =
\cE^o [\eta |\cF_0]$, applying (a2) with $A=\O$ and $\g=0$, we obtain
 \beas
 \cE^o[\xi] =\cE^o\big[ \cE^o [\xi |\cF_0]\big]
   =\cE^o\big[ \cE^o [\eta |\cF_0 ]\big]  =\cE^o [ \eta].
 \eeas
 Then the strict monotonicity of (a1) implies that $\xi
=\eta$, a.s., proving (A1).

\ms \no 2) Let $0 \le s\le t \le T$, for any $A \in \cF_s \subset
\cF_t$ and $\g \in \L_s \subset \L_t$, one can deduce that
 \beas
\cE^o\Big[\b1_A \cE^o\big[\cE^o[\xi|\cF_t]\big|\cF_s \big] +\g
\Big]=\cE^o\big[\b1_A \cE^o[\xi|\cF_t]+\g \big]=\cE^o\big[\b1_A  \xi
+\g \big].
 \eeas
 Since $\cE^o\big[\cE^o[\xi|\cF_t]\big|\cF_s \big] \in \cF_s $, (a2)
 implies that $ \cE^o[\xi|\cF_s ]= \xi_s=\cE^o\big[\cE^o[\xi|\cF_t]\big|\cF_s \big]$, proving (A2).

\ms \no 3)  Fix $A \in \cF_t $, for any $\tilde{A} \in \cF_t $ and
$\g \in \L_t$, we have
 \beas
 \cE^o\big[\b1_{\tilde{A}} \big(\b1_A \cE^o[\xi|\cF_t] \big) +\g \big]
  = \cE^o\big[\b1_{\tilde{A}\cap A} \cE^o[\xi|\cF_t] +\g \big]
  =\cE^o\big[\b1_{\tilde{A}\cap A} \xi  +\g \big]
  =\cE^o\big[\b1_{\tilde{A}}(\b1_A  \xi )  +\g \big].
\eeas Since $\b1_A \cE^o[\xi|\cF_t] \in \cF_t$, (a2) implies that
$\cE^o[ \b1_A  \xi|\cF_t] =\b1_A \cE^o[\xi|\cF_t] $, proving (A3).

\ms \no 4) For any $A \in \cF_t$ and $\eta,\, \g \in \L_t$, (D2)
implies that  $\b1_A\eta +\g \in \L_t $, thus we have
 \beas
 \cE^o\big[\b1_A \big( \cE^o[\xi|\cF_t]+\eta \big) +\g \big]=  \cE^o\big[\b1_A  \cE^o[\xi|\cF_t]
 +(\b1_A\eta +\g) \big]=   \cE^o\big[\b1_A \xi  +(\b1_A\eta +\g) \big]
 =  \cE^o\big[\b1_A (\xi  +\eta) +\g \big].
 \eeas
Then it follows from (a2) that
 $ \cE^o[\xi  +\eta|\cF_t]= \cE^o[\xi|\cF_t]+\eta$,
 proving (A4).   \qed

\ss \no {\bf Proof of Proposition \ref{3addition}:} (1) For any $A
\in \cF_t $, using (A3) twice, we obtain
 \beas
  \cE[\b1_A \xi \neg +\neg \b1_{A^c}\eta|\cF_t]
 &\tneg  =& \tneg   \b1_A \cE[\b1_A \xi \neg+\neg \b1_{A^c}\eta|\cF_t]
 \neg+\neg \b1_{A^c}\cE[\b1_A \xi \neg + \neg \b1_{A^c}\eta|\cF_t]\neg
 = \neg \cE[\b1_A( \b1_A \xi \neg+\neg \b1_{A^c}\eta) |\cF_t]
  \neg+\neg \cE[\b1_{A^c}(\b1_A \xi \neg+\neg \b1_{A^c}\eta)|\cF_t]\\
 & \tneg = & \tneg
   \cE[  \b1_A \xi  |\cF_t]+ \cE[  \b1_{A^c}\eta |\cF_t]
   = \b1_A\cE[ \xi|\cF_t]+ \b1_{A^c}\cE[ \eta|\cF_t], \q a.s.
\eeas
 (2) Applying (A3) with a
null set $A$ and $\xi=0$, we obtain
 $\cE[ 0|\cF_t]=  \cE[ \b1_A 0 |\cF_t] =\b1_A \cE[0|\cF_t]=0$, a.s.
 If $\xi \in Dom_t(\cE)$, (A4) implies that
 $ \cE[ \xi |\cF_t]= \cE[0 +\xi |\cF_t]= \cE[0  |\cF_t]+\xi=\xi$,
 a.s.

\ms \no (3) If $\xi \le \eta$, a.s., (A1) directly implies that for
any $A \in \cF_\nu$, $\cE[\b1_A \xi] \le  \cE[ \b1_A \eta]$. On the
other hand, suppose that $\cE[\b1_A \xi] \le \cE[ \b1_A \eta]$ for
any $A \in \cF_\nu$. We set $\tilde{A} \dfnn \{\xi
> \eta \} \in \cF_\nu$, thus $\b1_{\tilde{A}} \xi \ge
\b1_{\tilde{A}} \eta \ge c \land 0$, a.s. Using (A1) we see that
$\cE[\b1_{\tilde{A}} \xi] \ge \cE[ \b1_{\tilde{A}} \eta]$; hence
$\cE[\b1_{\tilde{A}} \xi] = \cE[ \b1_{\tilde{A}} \eta]$. Then (A4)
implies that
 \beas
 \cE[\b1_{\tilde{A}}\xi -c \land 0 ] =\cE[\b1_{\tilde{A}} \xi] -c \land 0= \cE[
 \b1_{\tilde{A}} \eta] -c \land 0 = \cE[ \b1_{\tilde{A}}\eta -c \land 0 ].
 \eeas
   Applying the second part of
(A1), we obtain that $\b1_{\tilde{A}} \xi-c \land 0= \b1_{\tilde{A}}
\eta -c \land 0$, a.s., which implies that $P(\tilde{A})=0$, i.e.
$\xi \le \eta$, a.s. \qed

 \ss \no {\bf Proof of Proposition \ref{op_sa}:}
We shall only consider the $\cE$-supermartingale case, as the other
cases can be deduced similarly. We first show that for any $s \in
[0,T]$ and $\nu \in \cS^F_{0,T}$
 \bea \label{eqn-f05}
 \cE [X_\nu |\cF_s ]\le
X_{\nu \land s}, \q a.s.
 \eea
 To see this, we note that since
$\{\nu\le s \} \in \cF_s$, (A3) and Proposition \ref{3addition} (2)
imply that
 \bea \label{eqn-e011}
  \cE [X_\nu  |\cF_s ]
  &= &  \b1_{\{\nu > s \}} \cE [X_\nu  |\cF_s ]+\b1_{\{\nu\le s \}}  \cE [X_\nu  |\cF_s ]
 = \cE\big[\b1_{\{\nu > s \}} X_{\nu \vee s} \big|\cF_s\big]
 +  \cE\big[\b1_{\{\nu\le s \}} X_{\nu \land s}\big|\cF_s\big] \nonumber \\
 &=& \b1_{\{\nu > s \}}\cE\big[ X_{\nu \vee s} \big|\cF_s\big]
 +  \b1_{\{\nu\le s \}} \cE\big[ X_{\nu \land s}\big|\cF_s\big]
 = \b1_{\{\nu > s \}}\cE\big[ X_{\nu \vee s} \big|\cF_s\big]
 +  \b1_{\{\nu\le s \}}  X_{\nu \land s}, \q a.s. \hspace{1cm}
 \eea
Suppose that $\nu_s \dfnn \nu \vee s$ takes values in a finite
subset $\{t_1 < \cds< t_n\}$ of $[s,T]$. Then (A4) implies that
 \beas
\cE [X_{\nu_s}  |\cF_{t_{n-1}}]=\cE\big[
 \b1_{\{{\nu_s}=t_n\}} X_{t_n} \big|\cF_{t_{n-1}}\big]+\sum^{n-1}_{i=1}
 \b1_{\{{\nu_s}=t_i\}} X_{t_i},\q a.s.
 \eeas
 Since $\{{\nu_s}=t_n\}=\{{\nu_s} > t_{n-1}\} \in \cF_{t_{n-1}}$, (A3) shows that
 \beas
  \cE\big[ \b1_{\{{\nu_s}=t_n\}} X_{t_n} \big|\cF_{t_{n-1}}\big]
   =\b1_{\{{\nu_s}=t_n\}}\cE [ X_{t_n} |\cF_{t_{n-1}} ]
   \le 
    \b1_{\{{\nu_s}=t_n\}} X_{t_{n-1} }, \q a.s.
 \eeas
 Thus it holds a.s. that $\dis \cE [X_{\nu_s}  |\cF_{t_{n-1}}]
 \le 
 \b1_{\{{\nu_s}> t_{n-2} \}}
 X_{t_{n-1}} +\sum^{n-2}_{i=1} \b1_{\{{\nu_s}=t_i\}} X_{t_i}$.
  Applying $\cE[\cd |\cF_{t_{n-2}} ] $ on both sides, we can
 further deduce from (A2)-(A4) that
 \beas
 \cE\big[X_{\nu_s} \big|\cF_{t_{n-2}}\big]&=& \cE\Big[\cE\big[X_{\nu_s}
\big|\cF_{t_{n-1}}\big]\Big|\cF_{t_{n-2}}\Big]
 \le 
   \b1_{\{{\nu_s}> t_{n-2} \}}\cE[X_{t_{n-1}}|\cF_{t_{n-2}}] +\sum^{n-2}_{i=1}
 \b1_{\{{\nu_s}=t_i\}} X_{t_i} \hspace{1.5cm}  \\
 &\le& \tneg 
 \b1_{\{{\nu_s}> t_{n-2} \}}X_{t_{n-2}}  +\sum^{n-2}_{i=1}
 \b1_{\{{\nu_s}=t_i\}} X_{t_i}
 =\b1_{\{{\nu_s}> t_{n-3} \}}X_{t_{n-2}}  +\sum^{n-3}_{i=1}
 \b1_{\{{\nu_s}=t_i\}} X_{t_i} ,\q a.s.
  \eeas
  Inductively, it follows that $\cE [X_{\nu_s}  |\cF_{t_1} ]
   \le 
   X_{t_1} $, a.s. Applying (A2) once again, we obtain
   \beas
\cE\big[X_{\nu_s}\big|\cF_s\big]=\cE\Big[\cE\big[X_{\nu_s}\big|\cF_{t_1}\big]\Big|\cF_s\Big]
  \le 
   \cE [ X_{t_1} |\cF_s ]\le 
   X_s, \q  a.s.,
   \eeas
which together with (\ref{eqn-e011}) implies that
 \beas
 \cE [X_\nu  |\cF_s] \le 
 \b1_{\{\nu > s \}}X_s
 +  \b1_{\{\nu\le s \}}  X_{\nu \land s}= X_{\nu \land s}, \q a.s., \q
 \hb{proving (\ref{eqn-f05})}.
 \eeas
 Let $\si \in \cS^F_{0,T}$ taking values in a
finite set $\{s_1 < \cds< s_m\}$, then
 \beas
 \hspace{3cm}  \cE [X_\nu  |\cF_\si  ]&\tneg=& \tneg \sum^m_{j=1} \b1_{\{\si=s_j\}}
 \cE [X_\nu |\cF_{s_j} ] \le 
 \sum^m_{j=1} \b1_{\{\si=s_j\}} X_{\nu \land s_j}
  = X_{\nu \land \si}, \q a.s. \hspace{3cm} \hb{\qed}
 \eeas

\ss \no {\bf Proof of Proposition \ref{properties_2}:} Given $\xi
\in Dom(\cE)$, we let $\nu \in \cS^F_{0,T}$ take values in a finite
set $\{t_1 < \cds< t_n\}$. 

\ms \no 1) For any $ \eta \in Dom(\cE) $ with $  \xi \le \eta$,
$a.s.$, (A1) implies that
 \beas
  \cE[\xi|\cF_\nu]=\sum^n_{i=1}
\b1_{\{\nu=t_i\}}\cE[\xi|\cF_{t_i}]\le \sum^n_{i=1}
\b1_{\{\nu=t_i\}}\cE[\eta|\cF_{t_i}]= \cE[\eta|\cF_\nu],\q a.s.
  \eeas
  Moreover, if $0 \le \xi \le
\eta$, a.s. and $\cE[\xi|\cF_\si]=  \cE[\eta|\cF_\si]$, a.s. for
some $\si \in \cS^F_{0,T}$, we can apply Corollary \ref{cor_os} to
obtain
 \beas
  \cE[\xi]=\cE\big[\cE[\xi|\cF_\si]\big]=\cE\big[\cE[\eta|\cF_\si]\big]=\cE[\eta].
 \eeas
 The second part of (A1) then implies that $ \xi= \eta $, a.s.,
proving (1).

 \ms \no  2) For any $A \in \cF_\nu$, it is clear that $A \cap
\{\nu=t_i \} \in
  \cF_{t_i}$ for each $i \in \{1,\cds,n\}$. Hence we can deduce from (A3) that
 \beas
 \cE [\b1_A \xi |\cF_\nu]&=&\sum^n_{i=1}
\b1_{\{\nu=t_i\}} \cE [\b1_A \xi |\cF_{t_i}] =\sum^n_{i=1}  \cE
[\b1_{\{\nu=t_i\}\cap A} \xi
|\cF_{t_i}]=\sum^n_{i=1}\b1_{\{\nu=t_i\}\cap A}  \cE
[ \xi |\cF_{t_i}]\\
 &=&  \b1_A \sum^n_{i=1} \b1_{\{\nu=t_i\}}  \cE [ \xi
|\cF_{t_i}]=\b1_A   \cE [ \xi |\cF_\nu],\q a.s., \q  \hb{proving
(2).}
 \eeas

 \no 3) For any $\eta \in Dom_\nu(\cE) $, since $ \b1_{\{\nu=t_i \}}
\eta \in
  Dom_{t_i}(\cE)$ for each $i \in \{1,\cds,n\}$, (A3) and (A4) imply that
 \beas
  \cE [  \xi+\eta |\cF_\nu]&=&\sum^n_{i=1} \b1_{\{\nu=t_i\}} \cE [\xi+\eta
 |\cF_{t_i}]=\sum^n_{i=1}  \cE[\b1_{\{\nu=t_i\}} \xi+\b1_{\{\nu=t_i\}}\eta
 |\cF_{t_i}]  =  \sum^n_{i=1} \Big( \cE[\b1_{\{\nu=t_i\}} \xi
 |\cF_{t_i}]+\b1_{\{\nu=t_i\}}\eta\Big)\\
 &=& \sum^n_{i=1}  \b1_{\{\nu=t_i\}} \cE[ \xi
 |\cF_{t_i}]+\eta =  \cE[ \xi |\cF_\nu]+\eta , \q a.s., \q  \hb{proving (3).}
 \eeas

 The proof of (4) and (5) is similar to that of Proposition
 \ref{3addition} (1) and (2) by applying the just obtained ``Zero-one Law" and ``Translation Invariance".
  \qed

 \ss \no {\bf Proof of Theorem \ref{fatou}:} (H1) is an easy
 consequence of the lower semi-continuity \eqref{eqn-b01}. In fact, for any $\xi \in Dom^+(\cE)$ and any $\{A_n\}_{n \in \hN} \subset
 \cF_T$ with $\underset{n \to \infty}{\lim} \neg \ua \b1_{A_n}=1
 $ a.s., $ \{ \b1_{A_n}\xi \}_{n \in \hN}$ is an increasing sequence
 converging to $\xi$. Then applying the lower semi-continuity with $\nu=0$ and using (A1), we obtain
  $
  \cE[\xi ] \le \underset{n \to \infty}{\lim } \neg \ua \cE[\b1_{A_n}\xi ]\le \cE[\xi
  ]$; so (H1) follows.

 On the other hand, to show that (H1) implies the lower semi-continuity,
 we first extend (H1) as follows: For any $\xi \in Dom^+(\cE)$ and any $\{A_n\}_{n
\in \hN} \subset \cF_T$ with $\underset{n \to \infty}{\lim} \neg \ua
\b1_{A_n}=1 $, a.s., it holds for any $t \in [0,T]$ that
 \bea \label{ass_H1'}
  \underset{n \to \infty}{\lim}\neg \ua \cE[\b1_{A_n}\xi|\cF_t]=\cE[\xi
|\cF_t] , \q  a.s.
 \eea
 In fact, by (A1), it holds a.s. that
$\big\{\cE[\b1_{A_n}\xi|\cF_t]\big\}_{n \in \hN}$ is an increasing
sequence bounded from above by $\cE[\xi|\cF_t]$. Hence, $
\underset{n \to \infty}{\lim} \neg \ua \cE[\b1_{A_n}\xi|\cF_t] \le
\cE[\xi |\cF_t ] $, a.s. Assuming that $ \underset{n \to
\infty}{\lim} \neg \ua \cE[\b1_{A_n}\xi|\cF_t] < \cE[\xi |\cF_t ] $
with a positive probability, we can find an $\e
>0$ such that the set $A_\e = \big\{\underset{n \to \infty}{\lim} \neg \ua
\cE[\b1_{A_n}\xi|\cF_t] \le \cE[\xi |\cF_t ]-\e \big\} \in \cF_t $
still has positive probability. Hence for any $n \in \hN$, we have
 \beas
  \b1_{A_\e} \cE[\b1_{A_n}\xi|\cF_t] \le \b1_{A_\e} \underset{n \to
  \infty}{\lim} \neg \ua
\cE[\b1_{A_n}\xi|\cF_t] \le \b1_{A_\e} \big(\cE[\xi |\cF_t
]-\e\big), \q a.s.
 \eeas
Then (A1)-(A4) imply that
 \beas
   \cE[\b1_{A_\e}  \b1_{A_n}\xi]+  \e & = &  \cE[\b1_{A_\e}  \b1_{A_n}\xi+  \e ]
   = \cE\big[\cE[\b1_{A_\e}  \b1_{A_n}\xi+  \e |\cF_t]\big]
  =\cE\big[\b1_{A_\e} \cE[\b1_{A_n}\xi |\cF_t] +  \e  \big]  \\
  &\le&  \cE\big[ \b1_{A_\e}  \cE[\xi |\cF_t] + \e \b1_{A^c_\e}\big]
 =\cE\big[ \cE[\b1_{A_\e}  \xi + \e \b1_{A^c_\e} |\cF_t ] \big]
 =\cE [ \b1_{A_\e} \xi + \e \b1_{A^c_\e} ].
  \eeas
Using (A4), (H1) and (A1), we obtain
 \beas
 \cE[\b1_{A_\e} \xi +\e ] = \cE[\b1_{A_\e} \xi  ]+\e= \underset{n \to
 \infty}{\lim} \neg \ua
\cE[\b1_{A_n} \b1_{A_\e} \xi  ]+\e  \le \cE [ \b1_{A_\e} \xi + \e
\b1_{A^c_\e}   ]  \le \cE [ \b1_{A_\e} \xi + \e   ],
 \eeas
thus $\cE[\b1_{A_\e} \xi +\e ] =\cE [ \b1_{A_\e} \xi + \e
\b1_{A^c_\e}   ] $. Then the second part of (A1) implies that $
\b1_{A_\e} \xi +\e   =  \b1_{A_\e} \xi + \e \b1_{A^c_\e} $, a.s.,
which can hold only if $P(A_\e)=0$. This results in a contradiction.
Thus $\underset{n \to \infty}{\lim} \neg \ua \cE[\b1_{A_n}\xi|\cF_t]
= \cE[\xi |\cF_t ]$, a.s., proving (\ref{ass_H1'}).

\ss Next, we show that (\ref{eqn-b01}) holds for each deterministic
stopping time $\nu=t \in [0,T]$. For any $j, n \in \hN$, we define
$\dis A^j_n \dfnn \cap^\infty_{k =n} \{|\xi-\xi_k|<1/j \} \in
\cF_T$. (A1) and (A4) imply that for any $k \ge n$
 \beas
  \cE[\b1_{A^j_n} \xi |\cF_t] \le \cE[\b1_{\{ |\xi-\xi_k|<1/j \}} \xi |\cF_t]
  \le \cE[\xi_k +1/j|\cF_t] = \cE[\xi_k |\cF_t]+ 1/j  , \q a.s.
 \eeas
Hence, except on a null set $N^j_n $, the above inequality holds for
any $k \ge n$. As $k \to \infty$, it holds on $\big(N^j_n \big)^c$
that
 \beas
 \cE[\b1_{A^j_n}\xi|\cF_t]  \le \underset{k \to \infty}{\liminf}\cE[\xi_k|\cF_t]+ 1/j  .
 \eeas
(Here it is not necessary that $\underset{k \to
\infty}{\liminf}\cE[\xi_k|\cF_t]< \infty$, a.s.) Since 
$\xi_n  \to \xi $, a.s. as $n \to \infty$, it is clear that
$\underset{n \to \infty}{\lim} \dneg \ua \b1_{A^j_n} =1$, a.s. Then
(\ref{ass_H1'}) implies that $\dis \cE[\xi|\cF_t] = \underset{n \to
\infty}{\lim} \neg \ua \cE[\b1_{A^j_n}\xi|\cF_t]$ holds except on a
null set $ N^j_0$. Let $N^j= \cup^\infty_{n=0}N^j_n$. It then holds
on $\big( N^j \big)^c$ that
 \beas
\cE[\xi|\cF_t] =  \underset{n \to \infty}{\lim } \neg \ua
\cE[\b1_{A^j_n}\xi|\cF_t] \le \underset{k \to \infty}{\liminf}
\cE[\xi_k|\cF_t]+ 1/j .
 \eeas
 As $j \to \infty$, it holds except on the null set $\cup^\infty_{j=1}N^j
 $ that
 \bea \label{eqn-e03}
  \cE[\xi|\cF_t]  \le \underset{n \to
\infty}{\liminf}\cE[\xi_n|\cF_t].
 \eea

 \ms Let $\nu \in \cS^F_{0,T}$ taking values in a finite set $\{t_1 < \cds<
t_n\}$. 
Then we can deduce from (\ref{eqn-e03}) that
 \bea
 \cE[\xi|\cF_\nu] = \sum^n_{i=1}\b1_{\{\nu=t_i\}}\cE[\xi|\cF_{t_i}]
 \le \sum^n_{i=1} \b1_{\{\nu=t_i\}} \underset{n \to \infty}{\liminf}
 \cE[\xi_n|\cF_{t_i}]
 = \neg \underset{n \to \infty}{\liminf}
 \sum^n_{i=1} \b1_{\{\nu=t_i\}} \cE[\xi_n|\cF_{t_i}]= \underset{n \to \infty}{\liminf}
\cE[\xi_n|\cF_\nu],~\; a.s., \label{eqn-e04}
 \eea
which completes the proof. \qed

\ss \no {\bf Proof of Theorem \ref{DCT}:} We first show an extension
of (H2): For any $\xi, \eta \in Dom^+(\cE)$ and any $\{A_n\}_{n \in
\hN} \subset \cF_T$ with $\underset{n \to \infty}{\lim} \neg \da
\b1_{A_n}=0 $, a.s., it holds a.s. that
 \bea \label{ass_H2'}
 \underset{n \to \infty}{\lim} \neg \da \cE[ \xi+\b1_{A_n}\eta|\cF_t] =
 \cE[\xi|\cF_t], \q  a.s.
 \eea
 In fact, by (A1), it holds
a.s. that $\big\{\cE[\xi+\b1_{A_n}\eta|\cF_t]\big\}_{n \in \hN}$ is
a decreasing sequence bounded from below by $\cE[\xi|\cF_t]$. Hence,
$ \underset{n \to \infty}{\lim} \neg \da
\cE[\xi+\b1_{A_n}\eta|\cF_t] \ge \cE[\xi |\cF_t ] $, a.s. Assume
that $ \underset{n \to \infty}{\lim} \neg \da
\cE[\xi+\b1_{A_n}\eta|\cF_t] > \cE[\xi |\cF_t ] $ with a positive
probability, then we can find an $\e
>0$ such that the set $A'_\e = \big\{\underset{n \to \infty}{\lim} \neg \da
\cE[\xi+\b1_{A_n}\eta|\cF_t] \ge \cE[\xi |\cF_t ]+\e \big\} \in
\cF_t$ still has positive probability. For any $n \in \hN$, (A4)
implies that
 \beas
 \b1_{A'_\e}\cE[\xi+\b1_{A_n}\eta|\cF_t] \ge \b1_{A'_\e} \underset{n \to \infty}{\lim} \neg \da
\cE[\xi+\b1_{A_n}\eta|\cF_t] \ge \b1_{A'_\e} \big(\cE[\xi |\cF_t
]+\e\big) = \b1_{A'_\e} \cE[\xi+\e |\cF_t ] , \q a.s.
 \eeas
  Applying (A1)-(A3), we obtain
   \beas
\cE[\b1_{A'_\e}\xi + \b1_{A_n} \b1_{A'_\e} \eta ]&=&
\cE\big[\cE[\b1_{A'_\e}\xi + \b1_{A_n} \b1_{A'_\e} \eta |\cF_t]\big]
= \cE\big[\b1_{A'_\e}\cE[\xi + \b1_{A_n} \eta |\cF_t]\big] \ge
\cE\big[\b1_{A'_\e} \cE[\xi+\e |\cF_t ]\big] \\
&=& \cE\big[ \cE[\b1_{A'_\e} (\xi+\e) |\cF_t ]\big]=\cE\big[
\b1_{A'_\e} (\xi+\e) \big].
 \eeas
 Thanks to (H2) we further have
 \beas
 \cE[\b1_{A'_\e} \xi ]= \underset{n \to \infty}{\lim} \neg \da
\cE[\b1_{A'_\e} \xi+ \b1_{A_n} \b1_{A'_\e} \eta ]\ge \cE\big[
\b1_{A'_\e} (\xi+\e)  \big] \ge \cE\big[ \b1_{A'_\e} \xi  \big],
 \eeas
thus $ \cE[\b1_{A'_\e} \xi ]=\cE[\b1_{A'_\e} (\xi+\e)  \big]$. Then
the second part of (A1) implies that $P(A'_\e)=0$, which yields a
contradiction. Therefore, $ \underset{n \to \infty}{\lim} \neg \da
\cE[\xi+\b1_{A_n}\eta|\cF_t] = \cE[\xi |\cF_t ]$, a.s., proving
(\ref{ass_H2'}).

 \ms Since the sequence $\{\xi_n\}_{n \in \hN}$ is bounded above by $\eta$, it holds a.s.
that $ \xi = \underset{n \to \infty}{\lim} \xi_n \le \eta $, thus
(D3) implies that $ \xi \in Dom(\cE)$.
Then Fatou's Lemma (Theorem \ref{fatou}) implies that 
for any $\nu \in \cS^F_{0,T}$,
 \bea \label{eqn-e05}
 \cE[\xi|\cF_\nu] \le \underset{n \to
\infty}{\liminf} \cE[\xi_n|\cF_\nu], \q a.s.
 \eea
On the other hand, we first fix $t \in [0,T]$. For any $j, n \in
\hN$, define $\dis A^j_n \dfnn \cap^\infty_{k =n} \{ |\xi-\xi_k|<1/j
\} \in \cF_T$. Then one can deduce that for any $k \ge n$
  \beas
  \cE[\xi_k|\cF_t] \le \cE[\b1_{A^j_n}(\xi+1/j)+\b1_{(A^j_n)^c}\eta|\cF_t]
  \le  \cE[\xi+1/j+\b1_{(A^j_n)^c}(\eta-\xi)|\cF_t] , \q  a.s.
 \eeas
 Hence, except on a null set $N^j_n$, the above inequality holds for any $k \ge n$.
 As $k \to \infty$, it holds on $\big(N^j_n\big)^c$ that
 \beas
 \underset{k \to \infty}{\limsup} \cE[\xi_k|\cF_t] \le   \cE[\xi+1/j+\b1_{(A^j_n)^c}(\eta-\xi)|\cF_t] .
 \eeas
Since $\xi \in L^0(\cF_T)$ and $\xi_n  \to \xi $, a.s. as $n \to
\infty$, it is clear that $\underset{n \to \infty}{\lim} \dneg \ua
\b1_{A^j_n} =1$, a.s.
 Then (\ref{ass_H2'}) and (A4) imply that except on a null
set $ N^j_0$, we have
 \beas
   \underset{n \to \infty}{\lim } \neg \da \cE[\xi+1/j+\b1_{(A^j_n)^c}(\eta-\xi)|\cF_t]
     = \cE[ \xi +1/j |\cF_t]=\cE[ \xi  |\cF_t]+1/j.
  \eeas
Let $N^j= \cup^\infty_{n=0}N^j_n$, thus it holds on $\big( N^j
\big)^c$ that
 \beas
   \underset{k \to \infty}{\limsup} \cE[\xi_k|\cF_t] \le \cE[ \xi  |\cF_t]+1/j.
  \eeas
As $j \to \infty$, it holds except on the null set
$\cup^\infty_{j=1}N^j
 $ that
 $  \underset{n \to \infty}{\limsup} \cE[\xi_n|\cF_t] \le \cE[ \xi
  |\cF_t]. $
Then for any $\nu \in \cS^F_{0,T}$, using an argument similar to
(\ref{eqn-e04}), we can deduce that
 \beas
 \underset{n \to \infty}{\limsup}
\cE[\xi_n|\cF_\nu] \le \cE[\xi|\cF_\nu] , \q a.s.,
 \eeas
 which together with (\ref{eqn-e05}) proves the theorem.   \qed

\ss \no {\bf Proof of Theorem \ref{upcrossing}:} Let $F =\{t_1<t_2<
\cds <t_d\}$ be any finite subset of $\cD_T$. For $j=1,\cds,d$, we
define $A_j=\{\nu_j<T\}\in \cF_{\nu_j}$,
clearly, $A_j \supset A_{j+1}$. Let $\dis d'=\Big\lfloor\,
\frac{d}{2}\,\Big\rfloor$,
 one can deduce that $\dis U_F\big(a, b; X \big)=\sum^{d'}_{j=1} \b1_{A_{2j}}$ and that
 \beas
 \b1_{\cup^{d'}_{j=1}(A_{2j-1} \backslash A_{2j})}   (X_T-a)
 \ge \b1_{\cup^{d'}_{j=1}(A_{2j-1} \backslash A_{2j})} \b1_{\{X_T<a\}}  (X_T-a)
  \ge \b1_{\{X_T<a\}} (X_T-a)=-(a-X_T)^+.
 \eeas
  Since $X_T \in Dom(\cE)$ and $L^\infty(\cF_T) \subset
Dom(\cE)$ (by Lemma \ref{lemm_wtsd}), we can deduce from (D2) that
  \beas
 (b-a) U_F\big(a, b; X \big)-(a-X_T)^+  =\sum^{d'}_{j=1} \b1_{A_{2j}}(b-a)  +\b1_{\{X_T<a\}} (X_T-a)
\in Dom(\cE).
  \eeas
 Then Proposition \ref{properties_2} (1)-(3) and Proposition \ref{op_sa}
imply that
 \beas
&& \hspace{-1cm} \cE\big[ (b-a) U_F\big(a, b; X
  \big)-(a-X_T)^+\big|\cF_{\nu_{2d'}}\big]
\le (b-a) \sum^{d'}_{j=1} \b1_{A_{2j}}
+ \cE\big[\b1_{\cup^{d'}_{j=1} (A_{2j-1} \backslash A_{2j})} (X_T-a)\big|\cF_{\nu_{2d'}}\big] \nonumber \\
 &=&\tneg (b-a) \sum^{d'}_{j=1} \b1_{A_{2j}}
 \neg + \b1_{\cup^{d'}_{j=1} (A_{2j-1} \backslash A_{2j})}\big( \cE[X_T |\cF_{\nu_{2d'}}]-a\big)
 \le  (b-a) \sum^{d'}_{j=1} \b1_{A_{2j}}
 \neg+ \b1_{\cup^{d'}_{j=1} (A_{2j-1} \backslash A_{2j})} \big(  X_{\nu_{2d'}}\neg
 -a\big), \q
 a.s.
 \eeas
  Applying $\cE[ \cd
|\cF_{\nu_{2d'-1}}] $ to the above inequality, using Proposition
\ref{properties_2} (1)-(3) and Proposition \ref{op_sa} once again,
we obtain
 \beas
\hspace{2cm} && \hspace{-2cm} \cE\big[ (b-a) U_F\big(a, b; X \big)-(a-X_T)^+\big|\cF_{\nu_{2d'-1}}\big] \nonumber \\
&\le &  \cE\Big[  (b-a) \sum^{d'-1}_{j=1} \b1_{A_{2j}}
 \neg+  \big( \b1_{A_{2d'-1}}
 +   \b1_{\cup^{d'-1}_{j=1} (A_{2j-1} \backslash A_{2j})}\big)
 \big(  X_{\nu_{2d'}} -a\big)   \Big|\cF_{\nu_{2d'-1}}\Big] \nonumber \\
 &= &    (b-a) \sum^{d'-1}_{j=1} \b1_{A_{2j}}+\cE\Big[
\big(\b1_{A_{2d'-1}} +   \b1_{\cup^{d'-1}_{j=1} (A_{2j-1} \backslash
A_{2j})}\big)
\big(  X_{\nu_{2d'}} -a\big) \Big|\cF_{\nu_{2d'-1}}\Big] \nonumber \\
&=&  (b-a) \sum^{d'-1}_{j=1} \b1_{A_{2j}}+\big(\b1_{A_{2d'-1}} +
\b1_{\cup^{d'-1}_{j=1} (A_{2j-1} \backslash A_{2j})}\big)
\big( \cE  [   X_{\nu_{2d'}}    |\cF_{\nu_{2d'-1}} ]-a\big) \nonumber \\
&\le &  (b-a) \sum^{d'-1}_{j=1} \b1_{A_{2j}}+\big( \b1_{A_{2d'-1}}
+ \b1_{\cup^{d'-1}_{j=1} (A_{2j-1} \backslash A_{2j})} \big) (  X_{\nu_{2d'-1}}-a ) \nonumber \\
&\le&  (b-a) \sum^{d'-1}_{j=1} \b1_{A_{2j}}+ \b1_{\cup^{d'-1}_{j=1}
(A_{2j-1} \backslash A_{2j})}  (  X_{\nu_{2d'-1}} -a ) , \q a.s.,
 \eeas
where we used the fact that $X_{\nu_{2d' }} >b  $ on $A_{2d' }$ in
the first inequality and the fact that $X_{\nu_{2d'-1}} <a $ on
$A_{2d'-1}$ in the last inequality. Similarly, applying $\cE[ \cd
|\cF_{\nu_{2d'-2}}] $ to the above inequality yields that
 \beas
    \cE\big[ (b-a) U_F\big(a, b; X
\big)-(a-X_T)^+\big|\cF_{\nu_{2d'-2}}\big]
 \le  (b-a) \sum^{d'-1}_{j=1} \b1_{A_{2j}}+ \b1_{\cup^{d'-1}_{j=1} (A_{2j-1}
\backslash A_{2j})}  (  X_{\nu_{2d'-2}} -a ) , \q a.s.
 \eeas
 Iteratively applying
$\cE[\cd|\cF_{\nu_{2d'-3}}]$, $\cE[\cd|\cF_{\nu_{2d'-4}}]$ and so
on, we eventually obtain that
 \bea \label{eqn-vvv01}
\cE\big[ (b-a) U_F\big(a, b; X \big)-(a-X_T)^+ \big] \le 0.
 \eea

\ss  We assume first that $X_T \ge c$,
a.s. for some $c \in \hR$. Since $(a-X_T)^+ \le \neg |a|+|c|$, it directly follows from 
(A4) that
 \bea \label{eqn-e06}
  0 \ge \cE\big[ (b-a) U_F\big(a, b; X \big)-(a-X_T)^+ \big]\ge
     \cE\big[ (b-a) U_F\big(a, b; X \big) \big]-(|a|+|c|).
 \eea
Let $\{F_n \}_{n \in \hN}$ be an increasing sequence of finite
subsets of $\cD_T$ with $\dis \cup_{n \in \hN} F_n =\cD_T $, thus $
\underset{n \to \infty}{\lim} \neg \ua U_{F_n} (a, b; X ) =
U_{\cD_T} (a, b; X )$. Fix $M \in \hN$, we see that
 \bea \label{eqn-vvv03}
 \underset{n \to \infty}{\lim} \neg \ua \b1_{\{U_{F_n} (a, b; X ) >M\}} = \b1_{\cup_n \{U_{F_n} (a, b; X ) >M\}}
 =   \b1_{  \{U_{\cD_T} (a, b; X ) >M\}}.
 \eea
For any $n \in \hN$, we know from (\ref{eqn-e06}) that $\cE[(b-a)M
\b1_{ \{U_{F_n} (a, b; X ) >M\}} ] \le \cE[(b-a) U_{F_n} (a, b; X )
] \le |a|+|c|$, thus Fatou's Lemma (Theorem \ref{fatou}) implies
that
 \bea  \label{eqn-vvv05}
 \cE\big[(b-a)M  \b1_{  \{U_{\cD_T} (a, b; X ) =\infty \}} \big]
 &\le&
\cE\big[(b-a)M  \b1_{  \{U_{\cD_T} (a, b; X ) >M\}} \big] \nonumber \\
       &\le&
      \underset{n \to \infty}{\lim} \neg \ua \cE\big[(b-a)M  \b1_{  \{U_{F_n} (a, b; X ) >M\}}\big]
\le |a|+|c|.
 \eea

 \ss On the other hand, if $\cE[\cd]$ is concave, then we can deduce from \eqref{eqn-vvv01} that
 \beas
  0 \ge \cE\big[ (b-a) U_F\big(a, b; X \big)-(a-X_T)^+ \big]\ge
   \frac12 \cE\big[ 2(b-a) U_F\big(a, b; X \big) \big]+\frac12 \cE\big[  -2(a-X_T)^+
\big].
 \eeas
Mimicking the arguments in \eqref{eqn-vvv03} and \eqref{eqn-vvv05},
we obtain that
 \beas
 \cE\big[(b-a)2M  \b1_{  \{U_{\cD_T} (a, b; X ) =\infty \}} \big]
  \le -  \cE\big[  -2(a-X_T)^+ \big].
 \eeas
 where $ -2(a-X_T)^+ =   \b1_{\{X_T<a\}} 2 (X_T-a) \in Dom(\cE)$ thanks to
(D2). Also note that (A1) and Proposition \ref{properties_2} (5)
imply that $\cE\big[  -2(a-X_T)^+ \big]\le \cE[0]=0 $.

 \ms Using (H0) in both cases above yields that $P( U_{\cD_T} (a, b; X )
=\infty )=0$, i.e., $U_{\cD_T} (a, b; X )< \infty$, a.s. Then a
classical argument (see e.g. \cite[Proposition
1.3.14]{Kara_Shr_BMSC})
 shows that
 \beas
 P \Big( \hb{both }\underset{s \nearrow t,\,  s \in \cD_T}{\lim} X_s \hb{
and }  \underset{s \searrow t,\, s \in \cD_T}{\lim} X_s \hb{ exist
for any } t \in [0,T] \Big)= 1.
 \eeas
This completes the proof. \qed

\ss \no {\bf Proof of Proposition \ref{prop_RCLL}:} We can deduce
from (\ref{limits0}) that except on a null set $N$
 \bea
   &&   X^+_t = \underset{n \to \infty}{\lim}
   X_{q^+_n(t)} \le \underset{s \in \cD_T}{\esssup}\, X_s \q \hb{for any } t \in
   [0,T], \hspace{1cm} \label{eqn-j10} \\
 \hb{thus}  & &  X^+_\nu = \underset{n \to \infty}{\lim}
   X_{q^+_n(\nu)} \le \underset{s \in \cD_T}{\esssup}\, X_s \q \hb{for any } \nu \in
   \cS_{0,T}. \label{eqn-j110}
 \eea

\no  \emph{Proof of (1)}: Case I. For any $\nu \in \cS_{0,T}$, if $
\underset{s \in \cD_T}{\esssup}\, X_s \in Dom^+(\cE)$, (D3) and
(\ref{eqn-j110}) directly imply that   $ X^+_\nu$
belongs to $Dom(\cE)$.

 \ss Case II. Assume that $\cE$ satisfies
(\ref{ass_fatou}). For any $n \in \hN$, since $X$ is an
$\cE$-supermartingale and since $  q^+_n(\nu)  \in
\cS^F_{0,T}$, Corollary \ref{cor_os} and Proposition \ref{op_sa}
imply that
 \beas
 &&\cE\big[ X_{q^+_n(\nu)}\big] =\cE \Big[\cE\big[ X_{q^+_n(\nu)}\big|\cF_{q^+_{n+1}(\nu)}\big] \Big]
 \le  \cE\big[ X_{q^+_{n+1}(\nu)}
 \big] \le X_0  .
  \eeas
Hence, $\big\{\cE\big[X_{q^+_n(\nu)}\big]\big\}_{n \in \hN}$ is an
increasing non-negative sequence  that  is  bounded from
above by $X_0 \in [0, \infty)$. (\ref{ass_fatou}) and
(\ref{eqn-j110}) then imply that   $X^+_\nu$ belongs to $Dom(\cE)$, proving statement (1).

\ms \noindent \emph{Proof of  (2)}: Now suppose that $X^+_t \in
Dom^+(\cE)$ for any $t \in [0,T]$. First, we shall show that for $t
\in [0,T]$ and $A \in \cF_t$
 \bea \label{eqn-xx01}
  \cE\big[ \b1_A    X^+_t  \big] =
  \underset{n \to \infty }{\lim} \cE\big[  \b1_A X_{q^+_n(t)}
  \big].
 \eea
Since the distribution function $x \mapsto P\{X^+_t \le x\}$
jumps up at most on a countable subset $S$ of $[0,\infty)$, we can
find a sequence $\{K_j\}^\infty_{j=1} \subset [0, \infty) \backslash
S$ increasing to $\infty$. 
 Fix $m, j \in \hN$,
(A1)-(A3) imply that for any $n \ge m$
 \beas 
    \cE\big[\b1_A \b1_{\{X_{q^+_n(t)}< K_j\}}( X_{q^+_n(t)} \land K_j )\big]
   &=&\cE\big[ \b1_A \b1_{\{X_{q^+_n(t)}< K_j\}} X_{q^+_n(t)} \big]
   \ge 
   ~\cE\Big[ \b1_A \b1_{\{X_{q^+_n(t)}< K_j\}} \cE\big[ X_{q^+_m(t)} \big|\cF_{q^+_n(t)}\big]\Big]  \nonumber \\
   &=&\cE\Big[  \cE\big[ \b1_A \b1_{\{X_{q^+_n(t)}< K_j\}}  X_{q^+_m(t)} \big|\cF_{q^+_n(t)}\big]\Big]
    =\cE\big[  \b1_A \b1_{\{X_{q^+_n(t)}< K_j\}}  X_{q^+_m(t)}  \big].
 \eeas
 Since $K_j \notin S$, $P\{X^+_t={K_j}\}=0$, one can easily deduce
 from (\ref{eqn-j10}) that  $\underset{n \to \infty}{\lim} \b1_{\{X_{q^+_n(t)}< K_j\}}=\b1_{\{X^+_t<
 K_j\}}$, a.s.\;\big(In fact, for almost every $\o \in \{X^+_t< K_j\}$ (resp. $\{X^+_t> K_j\}$), there exists an
$N(\o)\in \hN$ such that $X_{q^+_n(t)}< (\hb{resp.}\, >)\, K_j$ for
any $n \ge N(\o)$,
  which means $\dis\underset{n \to \infty}{\lim} \b1_{\{X_{q^+_n(t)}< K_j\}}(\o)
  =1 (\hb{resp.} \; 0) =  \b1_{\{X^+_t< K_j\}}(\o)$\big).
Applying the Dominated Convergence Theorem (Theorem \ref{DCT})
twice, we obtain \beas
   \cE\big[ \b1_A \b1_{\{X^+_t< K_j\}}   X^+_t  \big]
  &=& \cE\big[ \b1_A \b1_{\{X^+_t< K_j\}} ( X^+_t \land K_j ) \big]
  =\underset{n \to \infty}{\lim}\cE\big[ \b1_A \b1_{\{X_{q^+_n(t)}< K_j\}}( X_{q^+_n(t)} \land K_j)
  \big]\\
 &  \ge & 
 \underset{n \to \infty}{\lim}
 \cE\big[  \b1_A \b1_{\{X_{q^+_n(t)}< K_j\}}  X_{q^+_m(t)}  \big]
  =\cE\big[  \b1_A \b1_{\{X^+_t< K_j\}}  X_{q^+_m(t)}  \big].
 \eeas
Since $ \underset{j \to \infty}{\lim} \dneg \ua \b1_{\{X^+_t< K_j\}}
= 1$, a.s., the Dominated Convergence Theorem again implies that
 \beas
  \cE\big[ \b1_A    X^+_t  \big]= \underset{j \to \infty}{\lim}
    \cE\big[ \b1_A \b1_{\{X^+_t< K_j\}}   X^+_t  \big]
    \ge 
    \underset{j \to \infty}{\lim}  \cE\big[  \b1_A \b1_{\{X^+_t< K_j\}}  X_{q^+_m(t)}  \big]
    =\cE\big[  \b1_A    X_{q^+_m(t)}  \big],
 \eeas
which leads to that $\cE\big[ \b1_A    X^+_t  \big] \ge
  \underset{m \to \infty }{\limsup} \cE\big[  \b1_A X_{q^+_m(t)}
  \big]$. 
  Fatou's Lemma (Theorem \ref{fatou}) gives the reverse inequality, thus proving (\ref{eqn-xx01}).
Since $X$ is an
$\cE$-supermartingale, 
using (\ref{eqn-xx01}), (A2) and (A3), we obtain
 \beas
\cE\big[ \b1_A X^+_t \big]=\underset{n \to \infty }{\lim } \cE\big[
\b1_A X_{q^+_n(t)} \big] =\underset{n \to \infty }{\lim } \cE\big[
\cE[ \b1_A X_{q^+_n(t)} |\cF_t] \big] = \underset{n \to \infty
}{\lim } \cE\Big[ \b1_A \cE\big[ X_{q^+_n(t)} \big|\cF_t\big] \Big]
\le
\cE \big[ \b1_A X_t \big]
 \eeas
for any $A \in
 \cF_t$, which further implies that
$ X^+_t  \le  
X_t$, a.s. thanks to Proposition \ref{3addition} (3).

\ms Next, we show that $X^+$ is an
 $\cE$-supermartingale: 
 For any $0\le s < t \le T$,
it is clear that $ q^+_n(s) \le q^+_n(t) $ for any $n
\in \hN$. 
 For any $A \in \cF_s$, (A3) and Corollary \ref{cor_os} imply that for any $n \in \hN$
 \beas
\cE\big[\b1_A X_{q^+_n(s)}\big] \ge 
\cE\Big[\b1_A \cE\big[
X_{q^+_n(t)}\big|\cF_{q^+_n(s)}\big]\Big]=\cE\Big[ \cE\big[ \b1_A
X_{q^+_n(t)}\big|\cF_{q^+_n(s)}\big]\Big]= \cE\big[ \b1_A
X_{q^+_n(t)}\big].
 \eeas
  As $n \to \infty$, (\ref{eqn-xx01}),  (A2) and (A3) imply that
 \beas
 \cE\big[\b1_A  X^+_s \big]= \underset{n \to \infty }{\lim} \cE\big[
\b1_A X_{q^+_n(s)} \big]
 \ge 
 \underset{n \to \infty }{\lim} \cE\big[ \b1_A X_{q^+_n(t)} \big]=\cE\big[\b1_A  X^+_t \big]
 =\cE\big[\cE[\b1_A  X^+_t |\cF_s] \big]
 = \cE\big[\b1_A \cE[ X^+_t|\cF_s] \big].
 \eeas
  Then Proposition \ref{3addition} (3) implies that
 $X^+_s \ge  
 \cE[ X^+_t|\cF_s]$, a.s., thus $\{X^+_t\}_{t \in [0,T]}$ is
 an RCLL $\cE$-supermartingale.

 \ms \noindent \emph{Proof of (3):} If $t \mapsto \cE[X_t]$ is right continuous, for any
$t \in [0,T]$, (\ref{eqn-xx01}) implies that
 \beas
 \cE [ X^+_t  ] = \underset{n \to \infty }{\lim} \cE\big[
X_{q^+_n(t)} \big]
 = \cE [  X_t   ].
  \eeas
 Then the second part of (A1) imply that $ X^+_t=X_t$, a.s., which means that
 $ X^+$ is an RCLL modification of $X$. On the other hand, if
 $ \widetilde{X}$ is a right-continuous modification
 of $X$,  we see from (\ref{limits0}) that except on a null set $\tilde{N}$
 \beas
   X^+_t  =  \underset{n \to \infty }{\lim}  X_{q^+_n(t)}, \q   \widetilde{X}_t
    =  \underset{n \to \infty }{\lim}  \widetilde{X}_{q^+_n(t)}, \q  \widetilde{X}_t
    =X_t,  \q \hb{and} \q  \widetilde{X}_{q^+_n(t)}=X_{q^+_n(t)} \hb{  for any } n \in \hN.
 \eeas
 Putting them together, it holds on $\tilde{N}^c$ that
 \bea \label{eqn-xx03}
  X^+_t  =  \underset{n \to \infty }{\lim}  X_{q^+_n(t)}
  = \underset{n \to \infty }{\lim}  \widetilde{X}_{q^+_n(t)}=\widetilde{X}_t
  =X_t.
 \eea
  Since $X$ is an $\cE$-supermartingale,  (A2) implies that
 for any $0 \le t_1< t_2 \le T$, $ \cE[X_{t_1}]\ge \cE\big[ \cE[X_{t_2}|\cF_{t_1}] \big]=
   \cE[X_{t_2}]$, which shows that the function $t \mapsto \cE[X_t]$ is decreasing.
Then (\ref{eqn-xx01})  and (\ref{eqn-xx03}) imply that for any $t
\in [0,T]$
 \beas
  \cE[X_t] \ge  \underset{s \da t}{\lim}\, \cE[X_s]= \underset{n \to \infty}{\lim} \cE\big[X_{q^+_n(t)}\big]
    = \cE [ X^+_t   ]= \cE [ X_t ],
 \eeas
thus $\underset{s \da t}{\lim} \,\cE[X_s]=\cE [ X_t ]$, i.e., the
function $t \mapsto \cE[X_t]$ is right continuous. \qed

\ss \no {\bf Proof of Corollary \ref{cor_RCLL} :} Since $\underset{t
\in [0,T]}{\essinf}\, X_t \ge c$, a.s., we can deduce from (A4) that
$X^c \dfnn \{ X_t-c\}_{t \in [0,T]} $ is a non-negative
$\cE$-supermartingale. If $ \underset{t \in \cD_T}{\esssup}\, X_t
\in Dom^\#(\cE) $ \big ((D2) implies that $ \underset{t \in
\cD_T}{\esssup}\, X_t \in Dom^\#(\cE) $ is equivalent to $
\underset{t \in \cD_T}{\esssup}\, X^c_t \in Dom^+(\cE)$\big) or if
(\ref{ass_fatou}) holds, Proposition \ref{prop_RCLL} (1) shows that
for any $\nu \in \cS_{0,T}$,   $(X^c)^+_\nu $
belongs to $ Dom^+(\cE)$. Because
 \bea \label{eqn-axa02}
    (X^c)^+_t=X^+_t -c, \q \fa t
 \in [0,T],
 \eea
(D2) and the non-negativity of   $(X^c)^+$ imply that
 \beas
  X^+_\nu=
(X^c)^+_\nu+c \in Dom^\#(\cE).
 \eeas
On the other hand, if $X^+_t \in Dom^\#(\cE)$ for any $t \in [0,T]$,
(D2) implies that the non-negative random variable $(X^c)^+_t =X^+_t
-c $ belongs to $Dom^+(\cE)$. Hence, Proposition \ref{prop_RCLL} (2)
show that  $ (X^c)^+$ is an RCLL $\cE$-supermartingale such that for
any $t \in [0,T]$, $  (X^c)^+_t \le  X^c_t$, a.s. Then
(\ref{eqn-axa02}), (\ref{eqn-axa01}) and (A4) imply that $X^+$ is an
RCLL $\wt{\cE}$-supermartingale such that for any $t \in [0,T]$, $
X^+_t  \le  X_t$, a.s.  Moreover, if $t \mapsto \wt{\cE}[X_t]$ is a
right-continuous function (which is equivalent to the right
continuity of $t \mapsto \cE[X^c_t]$), then we know from Proposition
\ref{prop_RCLL} (2) that for any $t \in [0,T]$, $ (X^c)^+_t =
X^c_t$, a.s., or equivalently, $X^+_t = X_t$, a.s. Conversely, if
$X$ has a right-continuous modification, so does $X^c$, then
Proposition \ref{prop_RCLL} (2) once again shows that $t \mapsto
\cE[X^c_t]$ is right continuous, which is equivalent to the right continuity of $t \mapsto \wt{\cE}[X_t]$.
This completes the proof. \qed

 \ss \no {\bf Proof of Theorem \ref{op_sa2}:} We shall only consider the $\wt{\cE}$-supermartingale case,
 as the other cases can be deduced easily by similar arguments. Fix $t \in
 [0,T]$, we let $\{\nu^t_n\}_{n \in \hN}$ be a decreasing sequence in
$\cS^F_{t,T}$ such that $ \underset{n \to \infty}{\lim} \nu^t_n=\nu
\vee t    $. Since $\underset{t \in \cD_T}{\essinf}\, X_t \ge c$,
a.s., it holds a.s. that $X_t \ge c$ for each $t \in \cD_T$. The
right-continuity of the process $X$ then implies that except on a
null set $N$, $X_t \ge c$ for any $t \in [0,T]$. Thus
 we see from (A4) that $X^c \dfnn \{ X_t-c\}_{t \in
[0,T]} $ is a non-negative $\cE$-supermartingale. For any $n \in
\hN$ and $A \in \cF_t \subset \cF_{\nu \vee t}$, (A2), (A3) and
Proposition \ref{op_sa} imply that
  \bea  \label{eqn-xx05}
   \cE  [ \b1_A X^c_{\nu^t_n}] = \cE \big[   \cE  [\b1_A X^c_{\nu^t_n} |\cF_t] \big]
   =  \cE \big[ \b1_A  \cE  [X^c_{\nu^t_n} |\cF_t] \big] \le  \cE  [ \b1_A X^c_t  ].
  \eea
     We also have that
 $ 
  \cE\big[ \b1_A  X^c_{\nu \vee t}  \big]
   =  \underset{n \to \infty }{\lim }  \cE \big[  \b1_A X^c_{\nu^t_n}
  \big].
 $ The proof is similar to that of (\ref{eqn-xx01}).
\big(We only need to replace $X^+_t$ by $X^c_{\nu \vee t}$ and
$X_{q^+_n(t)}$ by $X^c_{\nu^t_n}$ in the proof of (\ref{eqn-xx01})
\big). As $n \to \infty$ in (\ref{eqn-xx05}), (A2) and (A3) imply
that
 \beas
    \cE [ \b1_A X^c_t ] \ge  \underset{n \to \infty }{\lim }  \cE\big[  \b1_A X^c_{\nu^t_n}  \big]
  =  \cE  [ \b1_A  X^c_{\nu \vee t}  ]= \cE \big[  \cE [\b1_A  X^c_{\nu \vee t}  |\cF_t] \big]
  =  \cE \big[ \b1_A  \cE [ X^c_{\nu \vee t} |\cF_t]\big].
 \eeas
Applying Proposition \ref{3addition} (3), we obtain that $
 \cE [ X^c_{\nu \vee t} |\cF_t] \le X^c_t$, a.s.
Then (A4) and (\ref{eqn-axa01}) imply that
 \beas
  \wt{\cE}[ X_{\nu \vee t}
|\cF_t]= \cE [ X_{\nu \vee t} |\cF_t]= \cE [ X^c_{\nu \vee t}+c
|\cF_t] =\cE [ X^c_{\nu \vee t} |\cF_t]+c \le X^c_t+c= X_t , \q a.s.
 \eeas
 Since $\{\nu\le t \} \in \cF_t$, we can deduce from (A3) and (A4) that
 \beas
  \wt{\cE} [X_\nu  |\cF_t ]
   &= &     \wt{\cE}\big[\b1_{\{\nu > t \}} X_{\nu \vee t}
 +  \b1_{\{\nu\le t \}} X_{\nu \land t}\big|\cF_t\big]
 =\b1_{\{\nu > t \}}\wt{\cE}\big[ X_{\nu \vee t} \big|\cF_t\big]
 +  \b1_{\{\nu\le t \}} X_{\nu \land t} \\
 &\le&  \b1_{\{\nu > t \}}X_t
 +  \b1_{\{\nu\le t \}} X_{\nu \land t}= X_{\nu \land t} \q a.s.
 \eeas
Hence, we can find a null set $\wt{N}$ such that except on
$\wt{N}^c$
 \beas
 \wt{\cE} [X_\nu |\cF_t ]\le 
 X_{\nu \land t}, \hb{~for any $t \in  \cD_T$ and the paths of $ \wt{\cE}[ X_\nu |\cF_\cd]$ and $ X_{\nu \land \cd}$
  are all RCLL.}
 \eeas
 As a result, on $\wt{N}^c$
 \beas
 \hspace{2.8cm} \wt{\cE} [X_\nu |\cF_t ]\le 
 X_{\nu \land t}, \q \fa t \in
 [0,T], \q \hb{thus} \q
 \wt{\cE} [X_\nu |\cF_\si ]\le 
 X_{\nu \land \si}, \q \fa \si \in \cS_{0,T}. \hspace{2.8cm} \hb{\qed}
 \eeas

 \ss \no {\bf Proof of Proposition \ref{properties_3}:}
  1) If $  \xi \le \eta$, a.s., by (A1), it holds except on a null set $N$ that
 \beas
  \wt{\cE}[\xi|\cF_t] \le \wt{\cE}[\eta|\cF_t], \hb{ for any $ t \in
  \cD_T$ and that the paths of $ \wt{\cE}[ \xi |\cF_\cd]$ and $\wt{\cE}[ \eta
|\cF_\cd]$ are all RCLL,}
 \eeas
 which implies that on $N^c$
 \beas
 \wt{\cE}[\xi|\cF_t] \le \wt{\cE}[\eta|\cF_t], \q \fa  t \in
 [0,T], \q \hb{thus }  \q \wt{\cE}[\xi|\cF_\nu] \le
 \wt{\cE}[\eta|\cF_\nu].
 \eeas
 Moreover, if $\wt{\cE}[\xi|\cF_\si]=
\wt{\cE}[\eta|\cF_\si]$, a.s. for some $\si \in \cS_{0,T}$, we can
apply (\ref{eqn-axa01}) and Corollary \ref{cor_os2} to get
 \beas
 \cE[\xi]=\wt{\cE}[\xi]
=\wt{\cE}\big[\wt{\cE}[\xi|\cF_\si]\big]=
\wt{\cE}\big[\wt{\cE}[\eta|\cF_\si]\big]=\wt{\cE}[\eta]=\cE[\eta].
 \eeas
 Then (A4) implies that
 $ \cE[\xi-c(\xi)]=\cE[\xi]-c(\xi)
=\cE[\eta]-c(\xi)=\cE[\eta-c(\xi)]$. Clearly, $0 \le \xi-c(\xi) \le
\eta -c(\xi) $, a.s. The second part of (A1) then implies that
$\xi-c(\xi) = \eta -c(\xi)$, a.s., i.e. $\xi= \eta $, a.s., proving
(1).

\ms \no 2) For any $A \in \cF_\nu$ and $\eta \in Dom^\#_\nu(\cE)$,
we let $\{\nu_n\}_{n \in \hN}$ be a decreasing sequence in
$\cS^F_{0,T}$ such that $ \underset{n \to \infty}{\lim} \neg \da
\nu_n=\nu$, a.s. For any $n \in \hN$, since $A \in \cF_{\nu_n}$ and
$\eta \in Dom^\#_{\nu_n}(\cE)$, Proposition \ref{properties_2} (2)
and (3) imply that
 \bea \label{eqn-f10}
  \wt{\cE} [\b1_A \xi |\cF_{\nu_n}] =\b1_A \wt{\cE} [\xi
|\cF_{\nu_n}], \q \hb{and} \q \wt{\cE} [  \xi+\eta
|\cF_{\nu_n}]=\wt{\cE} [  \xi |\cF_{\nu_n}]+\eta  , \q  a.s.
 \eea
 Then we can find a null set $N'$ such that except on $N'$
 \beas
 \hb{ (\ref{eqn-f10}) holds for any $n \in \hN$ and the paths
 of $ \wt{\cE} [\b1_A \xi |\cF_\cd]$, $\wt{\cE} [\xi |\cF_\cd]$ and
$\wt{\cE} [\xi+\eta |\cF_\cd]$ are all RCLL.}
 \eeas
 As $n \to \infty$, it holds on $(N')^c$ that
 \beas
 &&\wt{\cE} [\b1_A \xi |\cF_\nu] = \underset{n \to \infty}{\lim } \wt{\cE} [\b1_A \xi |\cF_{\nu_n}]
 =\underset{n \to \infty}{\lim } \b1_A \wt{\cE} [\xi
|\cF_{\nu_n}] =\b1_A \wt{\cE} [\xi |\cF_\nu], \\
\hb{and that} && \wt{\cE} [ \xi+\eta |\cF_\nu] = \underset{n \to
\infty}{\lim } \wt{\cE} [ \xi+\eta |\cF_{\nu_n}] = \underset{n \to
\infty}{\lim }\wt{\cE} [ \xi |\cF_{\nu_n}]+\eta   =\wt{\cE} [ \xi
|\cF_\nu]+\eta, \hspace{2.5cm}
 \eeas
proving (2) and (3). Proofs of (4) and (5) are similar to those of
Proposition \ref{3addition} (1) and (2). The proofs can be carried out by applying the just
obtained ``Zero-one Law" and ``Translation Invariance".  \qed

\subsection{Proofs of Section \ref{ch_3}}

\ss \no {\bf Proof of Lemma \ref{lem_pconvex}:} (1) Let $\cE$ be a
positively-convex $\bF$-expectation. For any $A \in \cF_T$ and $n
\in \hN$, (D1) and (D2) imply that $\b1_A, n\b1_A \in Dom(\cE)$.
Then the positive-convexity of $\cE$ and Proposition \ref{3addition}
(2) show that
  \bea \label{eqn-hxh01}
  \cE[\b1_A]=\cE\Big[ \frac{1}{n} \big(n \b1_A\big)\Big]\le \frac{1}{n}\cE[n
  \b1_A]
 +\big(1-\frac{1}{n}\big)\cE[0]= \frac{1}{n}\cE [n \b1_A]+\big(1-\frac{1}{n}\big)\cd 0
 = \frac{1}{n}\cE [n \b1_A].
  \eea
 Since $P(A)>0$, one can deduce from the second part of (A1) that $\cE [ \b1_A
  ]>0$. Letting $n \to \infty$ in \eqref{eqn-hxh01} yields that
 \beas
  \underset{n \to \infty}{\lim}\cE[n \b1_A] \ge \underset{n \to
  \infty}{\lim} n \cE[\b1_A]= \infty,
 \eeas
 thus $\cE$ satisfies (H0). Moreover, for any
$\xi, \eta \in Dom^\#(\cE)$, $\l \in (0,1)$ and $t  \in [0,T]$,
 we can deduce from \eqref{eqn-axa01}, (A4) and the positive-convexity of $\cE$ that
 \beas
 \wt{\cE}[\l\xi+(1-\l)\eta|\cF_t] &=& \cE[\l\xi+(1-\l)\eta|\cF_t]=
 \cE[\l\big(\xi-c(\xi)\big)+(1-\l)\big(\eta-c(\eta)\big)|\cF_t]+\l
 c(\xi)+ (1-\l)c(\eta)\\ &\le&
 \l \cE[  \xi-c(\xi) |\cF_t]+\l
 c(\xi)+(1-\l)\cE[ \eta-c(\eta) |\cF_t]+ (1-\l)c(\eta)\\
 &=& \l \cE[  \xi  |\cF_t] +(1-\l)\cE[ \eta  |\cF_t]
=    \l\wt{\cE}[\xi|\cF_t] +(1-\l)\wt{\cE}[\eta|\cF_t], \q a.s.,
 \eeas
 which shows that $\wt{\cE}$ is convex in the sense of
 \eqref{eqn_convex}. On the other hand, if $\wt{\cE}$ satisfies \eqref{eqn_convex},
 since $Dom^+(\cE) \subset Dom^\#(\cE)$, one can easily deduce from
 \eqref{eqn-axa01} that $\cE$ is positively-convex.  \qed

 \ss \no {\bf Proof of Proposition \ref{tau_A}:}
   We first check that
  $\cE^{\nu}_{i,j}$ satisfies (A1)-(A4). For this purpose, let $\xi, \eta \in \L^\#$ and $t \in
  [0,T]$.

 \ss \no 1) If $  \xi \le \eta$, a.s.,
 applying Proposition \ref{properties_3} (1) to $\wt{\cE}_j$ yields that
 $ \wt{\cE}_j[\xi|\cF_{\nu \vee t}] \le \wt{\cE}_j[\eta|\cF_{\nu \vee t}]$,
 a.s. Then (A1) of $\wt{\cE}_i$ and (\ref{tau_ij2}) imply that
  \beas
   \cE^{\nu}_{i,j}\big[\xi \big|\cF_t\big]
   = \wt{\cE}_i\big[ \wt{\cE}_j[\xi|\cF_{\nu \vee t}] \big|\cF_t \big]
    \le \wt{\cE}_i\big[\wt{\cE}_j[\eta|\cF_{\nu \vee t}] \big|\cF_t
    \big]=\cE^{\nu}_{i,j}\big[\eta \big|\cF_t\big],\q a.s.
  \eeas
 Moreover, if $0 \le \xi \le \eta$
a.s. and $ \cE^{\nu}_{i,j} [\xi   ] = \cE^{\nu }_{i,j} [\eta
]$\;\big(\,i.e. $ \wt{\cE}_i\big[ \wt{\cE}_j[\xi|\cF_{\nu}]   \big]
= \wt{\cE}_i\big[\wt{\cE}_j[\eta|\cF_{\nu}]
    \big] $ by (\ref{tau_ij2})\,\big), the second part of (A1) implies that $ \wt{\cE}_j[\xi|\cF_{\nu}]
    = \wt{\cE}_j[\eta|\cF_{\nu}]$, a.s. Further applying the
second part of Proposition \ref{properties_3} (1), we obtain $\xi
=\eta$, a.s., proving (A1) for $\cE^{\nu}_{i,j}$.

\ss \no 2) Next, we let $0 \le s\le t \le T$ and set $\Xi_t \dfnn
\cE^{\nu}_{i,j}\big[\xi \big|\cF_t\big]$. Applying Proposition
\ref{properties_3} (2) to $\wt{\cE}_i$ and $\wt{\cE}_j$, we obtain
 \beas 
\cE^{\nu }_{i,j}\big[ \Xi_t \big|\cF_s\big]
 \neg = \neg \b1_{\{\nu\le s\}}  \wt{\cE}_j \big[\Xi_t\big|\cF_s\big]\neg
 + \neg \b1_{\{\nu> s\}}\wt{\cE}_i\big[ \wt{\cE}_j
\big[\Xi_t\big|\cF_\nu\big]  \big|\cF_s\big] \neg =\neg
 \wt{\cE}_j \big[\b1_{\{\nu\le s\}} \Xi_t\big|\cF_s\big] \neg +\neg
\wt{\cE}_i\big[ \wt{\cE}_j \big[\b1_{\{\nu>
s\}}\Xi_t\big|\cF_\nu\big] \big|\cF_s\big], \q a.s.,
 \eeas
where we used the fact that $\{\nu > s \} \in \cF_{\nu \land s} $
thanks to \cite[Lemma 1.2.16]{Kara_Shr_BMSC}. Then (A3) and (A2)
imply that
 \bea \label{eqn-bxb02}
 \wt{\cE}_j \big[\b1_{\{\nu\le s\} }
\Xi_t\big|\cF_s\big]=\wt{\cE}_j \big[\b1_{\{\nu\le s\}  }
\wt{\cE}_j\big[\xi\big|\cF_t\big]\big|\cF_s\big]=\b1_{\{\nu\le s\}
 }\wt{\cE}_j \big[ \wt{\cE}_j\big[\xi\big|\cF_t\big]\big|\cF_s\big]
=\b1_{\{\nu\le s\}  } \wt{\cE}_j\big[\xi\big|\cF_s\big], \q a.s.
 \eea
On the other hand, we can deduce from (\ref{tau_ij}) that
 \beas
  \b1_{\{\nu > s \} } \Xi_t
   =    \b1_{\{s <\nu\le t\}  } \wt{\cE}_j \big[\xi\big|\cF_t\big] +
\b1_{\{\nu> t\}  }\wt{\cE}_i\big[  \wt{\cE}_j
 [\xi |\cF_\nu ]  \big|\cF_t \big]
  =     \b1_{\{s <\nu\le t\} } \wt{\cE}_j \big[\xi\big|\cF_t\big] +
\b1_{\{\nu> t\}  }\wt{\cE}_i\big[ \wt{\cE}_j  [\xi |\cF_\nu ]
\big|\cF_{\nu \land t}\big],\q a.s.
 \eeas
Since both $\{s <\nu\le t\}
 =\{ \nu >s\}\cap \{ \nu >t\}^c $ and $\{\nu > t \}$
 belong to $\cF_{\nu \land t}$,  Proposition
\ref{properties_3} (3) and (2) as well as Corollary \ref{cor_os2}
imply that
 \beas
\q && \hspace{-1.6cm} \wt{\cE}_j \big[\b1_{\{\nu > s \}  } \Xi_t
\big|\cF_\nu \big]
  =    \wt{\cE}_j \big[ \b1_{\{s <\nu\le t\} }\wt{\cE}_j
 [\xi |\cF_t ] \big|\cF_\nu \big] +  \b1_{\{\nu> t\}
}\wt{\cE}_i\big[ \wt{\cE}_j  [\xi |\cF_\nu ]
\big|\cF_{\nu \land t}\big]\\
&=&  \b1_{\{s <\nu\le t\} } \wt{\cE}_j \big[ \wt{\cE}_j
 [\xi |\cF_t ] \big|\cF_\nu \big] +  \b1_{\{\nu> t\}
}\wt{\cE}_i\big[ \wt{\cE}_j  [\xi |\cF_\nu ] \big|\cF_t \big]=
\b1_{\{s <\nu\le t\}  }    \wt{\cE}_j  [\xi |\cF_{\nu \land t} ] +
\wt{\cE}_i\big[  \b1_{\{\nu> t\} } \wt{\cE}_j  [\xi
|\cF_\nu ] \big|\cF_t\big] \\
&=&    \wt{\cE}_i\big[ \b1_{\{s <\nu\le t\}  }    \wt{\cE}_j  [\xi
|\cF_{\nu \land t} ]+ \b1_{\{\nu> t\} } \wt{\cE}_j [\xi |\cF_\nu ]
\big|\cF_t\big]=    \wt{\cE}_i\big[ \b1_{\{s <\nu \}  } \wt{\cE}_j
[\xi |\cF_\nu ] \big|\cF_t\big]  , \q a.s.
 \eeas
Taking $\wt{\cE}_i\big[\cd \big|\cF_s\big] $ of both sides as well
as using (A2) and (A3) of $\wt{\cE}_i$, we obtain
 \beas
  \q \wt{\cE}_i\big[ \wt{\cE}_j  [\b1_{\{\nu > s \}  } \Xi_t  |\cF_\nu  ] \big|\cF_s\big]
  \neg = \neg  \wt{\cE}_i\Big[  \wt{\cE}_i\big[ \b1_{\{s <\nu \}  } \wt{\cE}_j [\xi |\cF_\nu
] \big|\cF_t\big] \Big|\cF_s\Big] \neg  = \neg \wt{\cE}_i\big[
\b1_{\{s <\nu \}  } \wt{\cE}_j [\xi |\cF_\nu ] \big|\cF_s\big] \neg
=\neg \b1_{\{\nu> s\} } \wt{\cE}_i\big[  \wt{\cE}_j [\xi |\cF_\nu ]
\big|\cF_s\big], \q a.s.,
 \eeas
 which together with (\ref{eqn-bxb02}) yields that
 \beas
   \cE^{\nu }_{i,j}\big[\cE^{\nu }_{i,j} [\xi
 |\cF_t ] \big|\cF_s\big] =\b1_{\{\nu\le s\}  }
\wt{\cE}_j\big[\xi\big|\cF_s\big] + \b1_{\{\nu> s\} }
\wt{\cE}_i\big[ \wt{\cE}_j [\xi |\cF_\nu ] \big|\cF_s\big]= \cE^{\nu
}_{i,j}\big[\xi \big|\cF_s\big], ~\; a.s., ~\; \hb{proving (A2) for
} \cE^{\nu}_{i,j} .
 \eeas

 \ss \no 3) For any $A \in \cF_t$, using (\ref{tau_ij2}), (A3) of $\wt{\cE}_i$ as well as applying
 Proposition \ref{properties_3} (2) to $\wt{\cE}_j$, we obtain
 \beas
 \cE^{\nu }_{i,j}\big[\b1_{A}\xi \big|\cF_t\big]
  = \wt{\cE}_i \big[ \b1_{A} \wt{\cE}_j[\xi|\cF_{\nu \vee t}] \big|\cF_t \big]
  = \b1_{A} \wt{\cE}_i \big[ \wt{\cE}_j[ \xi|\cF_{\nu \vee t}] \big|\cF_t \big]
  =\b1_{A} \cE^{\nu }_{i,j}\big[ \xi \big|\cF_t\big], ~\; a.s., ~\; \hb{proving (A3) for }
\cE^{\nu}_{i,j} .
 \eeas
Similarly, we can show that (A4) holds for $\cE^{\nu}_{i,j}$ as
well. Therefore, $\cE^{\nu }_{i,j}$ is an $\bF$-expectation
 with domain $ \L^\#$. Since $\L \in \wt{\sD}_T$, i.e. $ \hR \subset \L$, it follows easily that
 $ \hR \subset \L^\#$, which shows that $\L^\# \in \wt{\sD}_T$.

 \ss \no 4) Now we show that $\cE^{\nu}_{i,j}$ satisfies (H1) and (H2):
For any $\xi \in \L^+$ and any $\{A_n\}_{n \in \hN} \subset \cF_T$
with $\underset{n \to \infty}{\lim} \neg \ua \b1_{A_n}=1 $, a.s.,
 the Dominated Convergence Theorem (Proposition \ref{DCT2}) implies that
$ \underset{n \to \infty}{\lim} \dneg \ua \wt{\cE}_j [\b1_{A_n} \xi
|\cF_\nu  ]=\wt{\cE}_j [\xi |\cF_\nu  ]$, a.s. Furthermore, using
(\ref{tau_ij2}) and applying the Dominated Convergence Theorem to
$\wt{\cE}_i$ yield that
 \beas
 \underset{n \to \infty}{\lim} \dneg \ua  \cE^{\nu}_{i,j}\big[\b1_{A_n}\xi \big]
  = \underset{n \to \infty}{\lim}
\dneg \ua \wt{\cE}_i \big[ \wt{\cE}_j\big[\b1_{A_n} \xi\big|\cF_\nu
\big] \big] = \wt{\cE}_i \big[ \wt{\cE}_j [\xi |\cF_\nu  ] \big]
=\cE^{\nu}_{i,j}\big[ \xi \big], \q \hb{proving (H1) for
}\cE^{\nu}_{i,j}.
 \eeas
 With a similar argument, we can show that $\cE^{\nu}_{i,j}$ also satisfies (H2).

\ss \no 5) If both $\cE_i$ and $\cE_j$ are positively-convex, so are
$\wt{\cE}_i$ and $\wt{\cE}_j$ thanks to (\ref{eqn-axa01}). To see
 that $\cE^{\nu }_{i,j}$ is convex
in the sense of \eqref{eqn_convex}, we fix $\xi, \eta \in \L^\#
 $, $\l \in (0,1)$ and $t \in
 [0,T]$. For any $s \in [0,T]$, we have
  \beas
  \wt{\cE}_j[\l\xi+(1-\l)\eta|\cF_s]\le \l \wt{\cE}_j[\xi|\cF_s]
  +(1-\l) \wt{\cE}_j[\eta|\cF_s], \q a.s.
   \eeas
Since $\wt{\cE}_j[\l\xi+(1-\l)\eta|\cF_\cd]$,
 $\wt{\cE}_j[\xi|\cF_\cd]$ and $\wt{\cE}_j[\eta|\cF_\cd]$ are all RCLL
 processes, it holds except on a null set $N$ that
  \beas
  \wt{\cE}_j[\l\xi+(1-\l)\eta|\cF_s]&\le& \l \wt{\cE}_j[\xi|\cF_s] +(1-\l)
 \wt{\cE}_j[\eta|\cF_s], \q \fa s \in [0,T], \\
 \hb{thus}  \qq \wt{\cE}_j[\l\xi+(1-\l)\eta|\cF_{\nu \vee t} ] &\le& \l \wt{\cE}_j[\xi|\cF_{\nu \vee t}] +(1-\l)
 \wt{\cE}_j[\eta|\cF_{\nu \vee t}].
  \eeas
 Then (\ref{tau_ij2}) implies that
 \beas
\hspace{1.6cm} && \hspace{-2cm} \cE^{\nu }_{i,j} [ \l\xi+(1-\l)\eta
|\cF_t ]
  =  \wt{\cE}_i \big[ \wt{\cE}_j[\l\xi+(1-\l)\eta|\cF_{\nu \vee t}] \big|\cF_t \big]
 \le  \wt{\cE}_i \big[ \l \wt{\cE}_j[\xi|\cF_{\nu \vee t}] +(1-\l)
 \wt{\cE}_j[\eta|\cF_{\nu \vee t}] \big|\cF_t \big], \\
  &\le & \l \wt{\cE}_i \big[  \wt{\cE}_j[\xi|\cF_{\nu \vee t}]  \big|\cF_t \big]
+(1-\l) \wt{\cE}_i \big[ \wt{\cE}_j[\eta|\cF_{\nu \vee t}]
\big|\cF_t \big]=  \l \cE^{\nu }_{i,j} [  \xi  |\cF_t
]+(1-\l)\cE^{\nu }_{i,j} [ \eta |\cF_t ] , \q a.s. \hspace{1.6cm}
\hb{\qed}
  \eeas

 \subsection{Proofs of Section \ref{co_game}}

\ss \no {\bf Proof of Lemma \ref{lem_H_eg}:} For any $i \in \cI$,
it is clear that
 $H^i_0=0$ and that \eqref{eqn-cxc01} directly follows from  (h1).
 For any $s,t \in \cD_T$
with $s<t$, we can deduce from (h2) that
 \bea \label{eqn-wxw100}
  H^i_{s,t}=\int_s^t h^i_r
dr \ge c\int_s^t  ds \ge cT  , \q  a.s.,
 \eea
 which implies that
$\underset{s,t \in \cD_T; s<t}{\essinf}\,H^i_{s,t} \ge cT$, a.s.
 Thus (S2) holds with $ C_H = cT$.

 \ss If no member of $\sE$ satisfies (\ref{ass_fatou}),
  then $  \int_0^T |h^j_t| \, dt   \in Dom(\sE)$ for some $j \in \cI$ is assumed.
 For any $s,t \in \cD_T$
with $s<t$, we can deduce from \eqref{eqn-wxw100} and (h2) that
 \beas
   C_H   \le   H^j_{s,t}  \le  \int_s^t
   | h^j_r |  dr \le \int_0^T  | h^j_r |  dr, \q a.s.,
 \eeas
 which implies that
  $
   C_H  \le   \underset{s,t \in \cD_T;
s<t}{\esssup}\,H^j_{s,t} \le  \int_0^T  | h^j_r |  dr$ a.s.  Then
Lemma \ref{lem_dom_sharp} shows that
 $\underset{s,t \in \cD_T;
s<t}{\esssup}\,H^j_{s,t} \in Dom(\sE)$, i.e. \eqref{ass_zi}.
 Moreover, we can derive (S3) directly from (h3). \qed

\ss \no {\bf Proof of Lemma \ref{lem_02}:} For any $i,j \in \cI'$
and $\rho_1, \rho_2 \in \cU$, we consider the event
 \beas
 A \dfnn \Big\{ \wt{\cE}_i  \big[X(\rho_1)+H^i_{\nu,\rho_1} \big|\cF_\nu \big]
 \le \wt{\cE}_j  \big[X(\rho_2) +   H^j_{\nu,\rho_2} \big|\cF_\nu \big]
 \Big\} \in \cF_\nu,
 \eeas
 and define stopping times $\rho \dfnn  \rho_2 \b1_A +  \rho_1
 \b1_{A^c} \in \cU$ and $\nu(A) \dfnn \nu \b1_A  +T \b1_{A^c}  \in \cS_{\nu,T}$.
 Since $\sE' =\left\{\cE_i \right\}_{i \in \cI'}$ is a stable
subclass of $\sE$,
 Definition \ref{def_stable_class} assures the existence of $k=k\big(i,j,\nu(A)\big) \in \cI'$
 such that $\wt{\cE}_k = \cE^{\nu(\neg
 A)}_{i,j}$. Applying Proposition
\ref{properties_3} (5) to $\wt{\cE}_j$ and Proposition
\ref{properties_3} (3) \& (2) to $\wt{\cE}_i$, we can deduce from
(\ref{tau_ij2}) that for any $\xi \in Dom(\sE) $
 \bea
 \wt{\cE}_k[\xi |\cF_\nu]&=&   \cE^{\nu(\neg A)}_{i,j}[\xi |\cF_\nu]
  =  \wt{\cE}_i\big[  \wt{\cE}_j [\xi|\cF_{\nu(A) \vee
\nu}] \big|\cF_\nu\big]= \wt{\cE}_i\big[ \b1_A \wt{\cE}_j
[\xi|\cF_{\nu}] + \b1_{A^c} \wt{\cE}_j
[\xi|\cF_T]\big|\cF_\nu\big] \nonumber \\
  &=&  \wt{\cE}_i\big[ \b1_A \wt{\cE}_j
[\xi|\cF_{\nu }] + \b1_{A^c}  \xi \big|\cF_\nu\big] =\b1_A
\wt{\cE}_j [\xi|\cF_{\nu }] +\b1_{A^c} \wt{\cE}_i [    \xi  |\cF_\nu
 ], \q a.s. \label{eqn-fxf01}
 \eea
 Moreover, (\ref{h_tau_A}) implies that
  \beas
      H^k_{\nu, \rho}
 \neg =\neg   H^i_{\nu(A) \land \nu, \nu(A)  \land \rho }
 \neg+\neg    H^j_{\nu(A)  \vee \nu, \nu(A)  \vee \rho }
 \neg=\neg  \b1_{A^c} H^i_{\nu, \rho_1}    \neg +\neg  \b1_A H^j_{\nu,\rho_2} , \q a.s.
  \eeas
Then applying Proposition \ref{properties_3} (2) to $\wt{\cE}_i$ and
$\wt{\cE}_j$, we see from \eqref{eqn-fxf01} that
 \beas
  \wt{\cE}_k\big[X(\rho)+H^k_{\nu,\rho}  \big|\cF_\nu\big]
  &=& 
       \b1_A\wt{\cE}_j\big[X(\rho) +H^k_{\nu,\rho} \big|\cF_\nu\big]
  + \b1_{A^c} \wt{\cE}_i\big[X(\rho) +H^k_{\nu,\rho} \big|\cF_\nu\big]\\
 &=&\wt{\cE}_j \big[\b1_A X(\rho_2) + \b1_A  H^j_{\nu,\rho_2}  \big|\cF_\nu\big]
 + \wt{\cE}_i\big[\b1_{A^c} X(\rho_1)+\b1_{A^c}  H^i_{\nu,\rho_1}  \big|\cF_\nu\big] \\
 &=& \b1_A   \wt{\cE}_j\big[ X(\rho_2)+H^j_{\nu,\rho_2}
 \big|\cF_\nu\big]
 + \b1_{A^c}   \wt{\cE}_i\big[ X(\rho_1)+H^i_{\nu,\rho_1} \big|\cF_\nu\big] \\
 &=&  \wt{\cE}_i\big[ X(\rho_1)+H^i_{\nu,\rho_1} \big|\cF_\nu\big]   \vee
  \wt{\cE}_j\big[ X(\rho_2)+H^j_{\nu,\rho_2} \big|\cF_\nu\big]  ,\q a.s.
 \eeas
Similarly, taking $\rho' \dfnn  \rho_1 \b1_A +  \rho_2 \b1_{A^c}$
and $k'=k\big(i,j,\nu(A^c)\big) $, we obtain
 \beas
 \wt{\cE}_{k'} \big[X(\rho')+H^{k'}_{\nu,\rho'}  \big|\cF_\nu\big]
 =  \wt{\cE}_i\big[ X(\rho_1)+H^i_{\nu,\rho}\big|\cF_\nu\big]
 \land   \wt{\cE}_j\big[ X(\rho_2)+H^j_{\nu,\rho}
 \big|\cF_\nu\big] , \q a.s.
 \eeas
 Hence, the family $\Big\{\wt{\cE}_i\big[X(\rho)+H^i_{\nu, \rho}    \big|\cF_\nu \big]
 \Big\}_{( i, \rho) \in \cI' \times \cU}$ is closed under pairwise maximization
and pairwise minimization. Thanks to \cite[Proposition
VI-\b1-1]{Neveu_1975},
we can find two sequences $\left\{(i_n, \rho_n)\right\}_{n \in \hN}$
and $\left\{(i'_n, \rho'_n)\right\}_{n \in \hN}$ in $\cI' \times
\cU$ such that (\ref{eqn-k70}) and (\ref{eqn-k71}) hold.
  \qed

\ss \no {\bf Proof of Lemma \ref{Z_bound}:} We fix $\nu \in
\cS_{0,T}$. For any $(i, \rho) \in \cI \times
 \cS_{\nu,T}$, \eqref{Y_LB}, \eqref{H_LB} and Proposition \ref{properties_3} (5) show that
 $   \wt{\cE}_i\big[Y_\rho+H^i_{\nu,\rho}
 \big|\cF_\nu\big] \ge \wt{\cE}_i [ C_* |\cF_\nu ] = C_*$ , a.s.
Taking the essential supremum over  $(i, \rho) \in \cI \times
 \cS_{\nu,T}$ gives
  \beas
    Z(\nu)= \underset{(i, \rho) \in \cI \times \cS_{\nu,T} }
 {\esssup}\, \wt{\cE}_i\big[Y_\rho+H^i_{\nu,\rho} \big|\cF_\nu\big] \ge C_*  , \q a.s.
 \eeas
 Then for any $i \in \cI$, \eqref{H_LB} implies that $Z^i(\nu)=Z(\nu)+H^i_\nu \ge C_*+
 C_H=C_Y+2C_H$, a.s.

\ms  If no member of $\sE$ satisfies (\ref{ass_fatou}) (thus
 (\ref{eqn-m48}) is assumed), then for any $(i, \rho) \in \cI \times
\cS_{\nu,T}$, it holds a.s. that
 \beas
   \wt{\cE}_i\big[Y^i_\rho \big|\cF_t\big] \le \z_Y, \q \fa   t \in
  \cD_T.
  \eeas
  Since $\wt{\cE}_i\big[Y^i_\rho \big|\cF_\cd\big]$ is an RCLL
  process, it holds except on a null set $N=N(i, \rho)$ that
 \beas
 \wt{\cE}_i\big[Y^i_\rho \big|\cF_t\big] \le \z_Y, \q \fa t \in
 [0,T], \q \hb{thus} \q \wt{\cE}_i\big[Y^i_\rho \big|\cF_\nu\big] \le
 \z_Y.
 \eeas
  Moreover, Proposition \ref{properties_3} (3) and (\ref{H_LB}) imply that
   \beas
\z_Y \ge \wt{\cE}_i\big[Y^i_\rho \big|\cF_\nu\big]=
\wt{\cE}_i\big[Y_\rho+H^i_{\nu,\rho}
 \big|\cF_\nu\big]+H^i_\nu \ge \wt{\cE}_i\big[Y_\rho+H^i_{\nu,\rho}
 \big|\cF_\nu\big]+ C_H, \q a.s.
   \eeas
Taking essential supremum over  $(i, \rho) \in \cI \times
 \cS_{\nu,T}$ yields that
  \beas
    Z(\nu)= \underset{(i, \rho) \in \cI \times \cS_{\nu,T} }
 {\esssup}\, \wt{\cE}_i\big[Y_\rho+H^i_{\nu,\rho} \big|\cF_\nu\big]  \le  \z_Y -C_H, \q a.s.
 \eeas
 where $\z_Y -C_H \in Dom(\sE)$ thanks to (\ref{eqn-m48}) and
(D2). Hence, for any $i \in \cI$, we have $ Z^i(\nu) =Z(\nu)+H^i_\nu
\le \z_Y-C_H+H^i_\nu$, a.s. And \eqref{eqn-cxc01} together with (D2)
imply that $\z_Y -C_H+ H^i_\nu \in Dom(\sE)$. \qed

\ss \no {\bf Proof of Lemma \ref{lem_Z_nu}:} If no member of $\sE$
satisfies (\ref{ass_fatou}),
then we see from Lemma \ref{Z_bound} that
 \beas
  C_* \le  Z(\nu)  \le  \z_Y -C_H, \q a.s.,
 \eeas
and that $\z_Y -C_H \in Dom(\sE)$. Hence $Z(\nu) \in
Dom(\sE)$ thanks to Lemma \ref{lem_dom_sharp}. 

 \ms On the other hand, if $\cE_j$ satisfies (\ref{ass_fatou}) for some
$j \in \cI$,
 letting $(X, \cI', \cU)=(Y, \cI,
\cS_{\nu,T})$ in Lemma \ref{lem_02}, we can find a sequence
$\left\{(i_n, \rho_n)\right\}_{n \in \hN}$ in $\cI \times
\cS_{\nu,T}$ such that
 \beas 
  Z(\nu)=\underset{(i, \rho) \in \cI \times \cS_{\nu,T}}{\esssup}\,
   \wt{\cE}_i \big[Y_\rho+H^i_{\nu, \rho}   \big|\cF_\nu
  \big]
  = \underset{n \to \infty}{\lim} \dneg \ua
   \wt{\cE}_{i_n}\big[Y_{\rho_n}+H^{i_n}_{\nu, \rho_n}
   \big|\cF_\nu\big], \q a.s.
  \eeas
For any $n \in \hN$, it follows from Definition
\ref{def_stable_class} that there exists $k_n=k(j,i_n, \nu) \in \cI$
such that
 $  \wt{\cE}_{k_n}=\cE^\nu_{j,i_n}$.
Applying Proposition \ref{properties_3} (3) to $\wt{\cE}_{k_n}$, we
 can deduce from \eqref{H_LB}, (\ref{tau_ij2}) and (\ref{h_tau_A}) that
 \beas
   \wt{\cE}_{k_n}\big[Y^{k_n}_{\rho_n}
 \big]-C_H &=&  \wt{\cE}_{k_n}\big[Y_{\rho_n}+ H^{k_n}_{\rho_n}-C_H
 \big] =\wt{\cE}_{k_n}\big[Y_{\rho_n}+H^{k_n}_{\nu, \rho_n}+ H^{k_n}_\nu-C_H
 \big] \ge \wt{\cE}_{k_n}\big[Y_{\rho_n}+ H^{k_n}_{\nu,\rho_n} \big]  \\
   &=&  \cE^\nu_{j,i_n}\big[Y_{\rho_n}+ H^{k_n}_{\nu, \rho_n} \big]= \wt{\cE}_j\big[   \wt{\cE}_{i_n}\neg
    \big[Y_{\rho_n}+ H^{k_n}_{\nu,\rho_n}  \big|\cF_\nu
    \big]\big]  =   \wt{\cE}_j\Big[\wt{\cE}_{i_n}\big[Y_{\rho_n}+H^{i_n}_{\nu,\rho_n}
       \big|\cF_\nu\big]\Big],
  \eeas
  which together with (Y2) shows that
   \beas
     \underset{n \to \infty}{\liminf}  \cE_j
     \Big[ \wt{\cE}_{i_n}\big[Y_{\rho_n}+H^{i_n}_{\nu,\rho_n}\big|\cF_\nu\big]\Big]
    \le  \underset{(i, \rho) \in \cI \times \cS_{0,T}}{\sup} \wt{\cE}_i \big[Y^i_\rho
 \big]-C_H   < \infty.
   \eeas
 For any $n \in \hN$, (\ref{Y_LB}),
 (\ref{H_LB}) and Proposition \ref{properties_3} (5) imply that
 \beas
 \wt{\cE}_{i_n}\big[Y_{\rho_n}+H^{i_n}_{\nu,\rho_n}\big|\cF_\nu\big]
\ge   \wt{\cE}_{i_n} [C_* |\cF_\nu ]= C_* , \q a.s.
 \eeas
Therefore, we can deduce from Remark \ref{rem_fatou2} (1) that
 \beas
    Z(\nu)  = \underset{n \to \infty}{\lim} \dneg \ua
   \wt{\cE}_{i_n}\big[Y_{\rho_n}+H^{i_n}_{\nu,\rho_n}\big|\cF_\nu\big]
    \in  Dom(\sE).
  \eeas
For any $i \in \cI$, (\ref{eqn-cxc01}) and (D2) imply that $Z^i(\nu)
= Z(\nu)+ H^i_\nu \in Dom(\sE)$. \qed

\ss \no {\bf Proof of Proposition \ref{prop_01}:} To see
(\ref{D_1}), we first note that the event $A \dfnn \left\{\nu=\si
\right\}$ belong to $\cF_{\nu \land \si} $ thanks to \cite[Lemma
1.2.16]{Kara_Shr_BMSC}. For any $i \in \cI$ and $\rho \in \cS_{\nu,
T}$, we define $  \rho(A)  \dfnn  \rho \b1_A+ T\b1_{A^c} $,
 which clearly belongs to $\cS_{\si, T}$.  Proposition \ref{properties_3} (2) and
 (3) then imply that
  \beas
   \b1_A  \wt{\cE}_i\big[Y_{\rho }+H^i_{\nu, \rho} \big|\cF_\nu\big]
  &=&\b1_A \Big(\wt{\cE}_i\big[Y_{\rho
}+H^i_\rho
  \big|\cF_\nu\big]-H^i_\nu \Big)
 =  \b1_A \Big(\wt{\cE}_i\big[Y_{\rho }+H^i_\rho
  \big|\cF_\si \big]-H^i_\si \Big)
  = \b1_A  \wt{\cE}_i\big[Y_{\rho }+H^i_{\si, \rho}
  \big|\cF_\si \big]\\
  &=&   \wt{\cE}_i\big[\b1_A \big(Y_{ \rho(A) }+H^i_{\si, \rho(A) } \big)  \big|\cF_\si \big]
  =\b1_A \wt{\cE}_i\big[  Y_{ \rho(A) }+H^i_{\si, \rho(A) }    \big|\cF_\si  \big] \\
  &\le&     \b1_A \underset{(i, \g) \in \cI \times \cS_{\si,T}
  }{\esssup}\, \wt{\cE}_i\big[  Y_\g+H^i_{\si,\g} \big|\cF_\si
 \big]  =\b1_A  Z(\si), \q a.s.
 \eeas
Taking the essential supremum of the left-hand-side over $(i, \rho)
\in
 \cI \times \cS_{\nu,
T}$ and applying Lemma \ref{lem_ess} (2), we obtain
 \beas
   \b1_A Z(\nu) =\b1_A  \underset{(i, \rho) \in \cI \times \cS_{\nu, T} }{\esssup}\,
   \wt{\cE}_i\big[Y_{\rho }+H^i_{\nu,\rho}\big|\cF_\nu\big]
  =   \underset{(i, \rho) \in \cI \times \cS_{\nu, T}
  }{\esssup}\, \Big( \b1_A \wt{\cE}_i\big[Y_{\rho }+H^i_{\nu,\rho}\big|\cF_\nu\big] \Big)
    \le  \b1_A  Z(\si),  \q a.s.
 \eeas
 Reversing the roles of $\nu$ and $\si$, we obtain (\ref{D_1}).

 \ms  As to (\ref{D_2}), since $\cS_{\g, T} \subset \cS_{ \nu, T}$,
it is clear that
 \beas
  \underset{(i, \rho) \in \cI \times \cS_{\g, T}
  }{\esssup}\,   \wt{\cE}_i\big[Y_\rho+H^i_{\nu,\rho}\big|\cF_\nu\big]
   \le \underset{(i, \rho) \in \cI \times \cS_{\nu, T}
  }{\esssup}\,   \wt{\cE}_i\big[Y_\rho+H^i_{\nu,\rho}\big|\cF_\nu\big] = Z(\nu), \q a.s.
 \eeas
  Letting $(X, \nu, \cI', \cU)=(Y, \g, \cI, \cS_{\g, T})$ in Lemma \ref{lem_02}, we can find a sequence
 $\big\{ (i_n, \rho_n )\big\}_{n \in \hN}$ in $ \cI \times \cS_{\g, T}$ such that
 \beas
Z(\g)  =\underset{(i, \rho) \in \cI \times \cS_{\g,T}}{\esssup}\,
   \wt{\cE}_i \big[Y_\rho+H^i_{\g, \rho}   \big|\cF_\g  \big]
     = \underset{n \to \infty}{\lim} \dneg \ua
    \wt{\cE}_{i_n}\big[Y_{\rho_n} +H^{i_n}_{\g, \rho_n} \big|\cF_\g\big] , \q a.s.
 \eeas
Now fix $ j \in \cI$. For any $n \in \hN$, it follows from
Definition \ref{def_stable_class} that there exists
  a $k_n =k(j,i_n,\g) \in \cI$ such that  $\wt{\cE}_{k_n}=\cE^\g_{j,i_n}$.
Applying Proposition \ref{properties_3} (3) to $\wt{\cE}_{i_n}$, we
can deduce from (\ref{tau_ij2}), (\ref{h_tau_A}) that
 \bea
 && \hspace{-1.2cm} \underset{(i, \rho) \in \cI \times \cS_{\g, T}
  }{\esssup}\,   \wt{\cE}_i\big[Y_\rho+ H^i_{\nu,\rho} \big|\cF_\nu
  \big]
   \ge \wt{\cE}_{k_n}\big[Y_{\rho_n}+ H^{k_n}_{\nu,\rho_n} \big|\cF_\nu \big]
   = \cE^\g_{j,i_n} \big[Y_{\rho_n}+ H^{k_n}_{\nu,\rho_n}  \big|\cF_\nu  \big]
   = \wt{\cE}_j\big[ \wt{\cE}_{i_n} \big[Y_{\rho_n}
  + H^{k_n}_{\nu, \rho_n}\big|\cF_\g \big] \big|\cF_\nu\big] \nonumber \\
  &  =&    \wt{\cE}_j\big[ \wt{\cE}_{i_n} \big[Y_{\rho_n}
  + H^{k_n}_{\g,\rho_n} \big|\cF_\g \big] + H^{k_n}_{\nu, \g} \big|\cF_\nu\big]
    =  \wt{\cE}_j\Big[ \wt{\cE}_{i_n} \big[Y_{\rho_n}
  + H^{i_n}_{\g, \rho_n}  \big|\cF_\g \big] +H^j_{\nu,\g} \big|\cF_\nu\Big], \q
  a.s. \label{eqn-cxc03}
  \eea
  For any $n \in \hN$, Proposition \ref{properties_3} (5), (\ref{Y_LB}) and (\ref{H_LB}) show that
  \beas
 C_Y+2C_H =  \wt{\cE}_{i_n} \big[C_* \big|\cF_\g \big]+ C_H \le   \wt{\cE}_{i_n} \big[Y_{\rho_n}
  + H^{i_n}_{\g, \rho_n} \big|\cF_\g \big] +H^j_{\nu,\g}  \le Z(\g) +H^j_{\nu,\g} , \q
  a.s.,
  \eeas
  where $ Z(\g) +H^j_{\nu,\g}  \in Dom(\sE)$ thanks to Lemma \ref{lem_Z_nu}, \eqref{eqn-cxc01} and (D2).
  Then the Dominated Convergence Theorem (Proposition \ref{DCT2}) and (\ref{eqn-cxc03}) imply that
\beas
   \wt{\cE}_j \big[ Z(\g) +H^j_{\nu, \g} \big|\cF_\nu \big] =  \underset{n \to \infty}{\lim}
    \wt{\cE}_j\Big[ \wt{\cE}_{i_n}\neg \big[Y_{\rho_n}
  \neg+\neg H^{i_n}_{\g, \rho_n}  \big|\cF_\g \big]\neg + \neg H^j_{\nu,\g} \big|\cF_\nu\Big]
   \le  \underset{(i, \rho) \in \cI \times \cS_{\g, T}
  }{\esssup}\, \dneg \wt{\cE}_i\big[Y_\rho\neg + \neg H^i_{\nu, \rho} \big|\cF_\nu\big]
     , \q a.s.
 \eeas
Taking the essential supremum of the left-hand-side over $j \in
\cI$, we obtain
 \bea \label{eqn-q102}
\underset{j \in \cI}{\esssup}\, \wt{\cE}_j \big[ Z(\g) + \neg
H^j_{\nu, \g} \big|\cF_\nu \big]\le \underset{(i, \rho) \in \cI
\times \cS_{\g, T}
  }{\esssup}\, \dneg \wt{\cE}_i\big[Y_\rho+H^i_{\nu,\rho}\big|\cF_\nu\big]
     , \q a.s.
 \eea
  On the other hand, for any $i \in \cI$ and $\rho \in \cS_{\g,
 T}$, 
applying Corollary \ref{cor_os2} and Proposition \ref{properties_3}
(3), we obtain
 \beas
  \wt{\cE}_i \big[ Y_\rho+H^i_{\nu,\rho} \big|\cF_\nu \big]
 & = & \wt{\cE}_i\big[\wt{\cE}_i\big[ Y_\rho+H^i_{\nu,\rho} \big|\cF_\g\big]\big|\cF_\nu \big]
  = \wt{\cE}_i\big[\wt{\cE}_i\big[ Y_\rho+ H^i_{\g, \rho} \big|\cF_\g\big]+H^i_{\nu,\g}\big|\cF_\nu \big]\\
  &\le& \wt{\cE}_i\big[Z(\g)+H^i_{\nu,\g}\big|\cF_\nu\big] \le
  \underset{i \in \cI}{\esssup}\, \wt{\cE}_i \big[ Z(\g) + \neg H^i_{\nu,\g} \big|\cF_\nu \big], \q a.s.
 \eeas
Taking the essential supremum of the left-hand-side over $(i, \rho)
\in \cI \times \cS_{\g,
 T}$ yields that
  \beas
 \underset{(i, \rho) \in \cI \times \cS_{\g, T}
  }{\esssup}\, \wt{\cE}_i \big[ Y_\rho +H^i_{\nu,\rho} \big|\cF_\nu \big]
   \le  \underset{i \in \cI}{\esssup}\, \wt{\cE}_i \big[ Z(\g) + \neg H^i_{\nu,\g} \big|\cF_\nu \big], \q
   a.s.,
 \eeas
which together with (\ref{eqn-q102}) proves (\ref{D_2}). \qed

\ss \no {\bf Proof of Proposition \ref{Z_RCLL}:} For any $i \in \cI$
and $\nu ,\g \in \cS_{0,T}$ with $\nu \le \g$, a.s., Proposition
\ref{properties_3} (3), (\ref{D_2}) imply that
 \beas 
 \wt{\cE}_i\big[Z^i(\rho) \big|\cF_{\nu}\big] = \wt{\cE}_i\big[Z(\rho)+\neg
 H^i_{\nu,\rho} \big|\cF_{\nu}\big]+\neg H^i_\nu
 \le \underset{i \in \cI}{\esssup}\, \wt{\cE}_i\big[Z(\rho)+\neg
 H^i_{\nu,\rho} \big|\cF_{\nu}\big]+\neg H^i_\nu 
  \le   Z(\nu)+\neg H^i_\nu =  Z^i(\nu), \q a.s.,
 \eeas
 which implies that $\left\{Z^i(t)\right\}_{t \in [0,T]}$ is an
$\wt{\cE}_i$-supermartingale.
  Proposition \ref{prop_tilde_E}, Theorem
\ref{upcrossing} and \eqref{Z_LB} then show that \\$ \Big\{
Z^{i,+}_t \dfnn \underset{n \to \infty}{\liminf}
Z^i\big(q^+_n(t)\big)\Big\}_{t \in [0,T]} $
  defines an RCLL process. Moreover, (\ref{Z_LB}) implies that
  \bea \label{eqn-cxc11}
  \underset{t \in [0,T]}{\essinf}\, Z^i(t) \ge C_Y+2C_H, \q  a.s.
  \eea

 If $\cE_j$ satisfies
(\ref{ass_fatou}) for some $j \in \cI$,
 Corollary \ref{cor_RCLL} and (\ref{eqn-cxc11}) imply that
 \bea
 Z^{j,+}_\nu \in Dom^\#(\cE_j)=Dom(\sE), \q \fa \nu \in \cS_{0,T}, \hspace{3.5cm} \label{eqn-lxl01}\\
 \hb{and that $Z^{j,+}$ is an RCLL $ \wt{\cE}_j$-supermartingale
 such that for any $t \in [0,T]$, $Z^{j,+}_t \le Z^j(t)$,\; a.s.} \label{eqn-lxl02}
 \eea
 Otherwise, if no member of $\sE$ satisfies (\ref{ass_fatou}),
 we suppose that
\eqref{ass_zi} holds for some $j \in \cI$.
  Then Lemma \ref{Z_bound} and \eqref{ass_zi} imply that for any $t \in \cD_T$,
  \beas
  C_Y+2C_H \le  Z^j(t)= Z(t)+H^j_t \le \z_Y - C_H+ \z^j, \q
  a.s.
  \eeas
 Taking essential supremum of $Z^j(t)$ over $t \in \cD_T$ yields that
  \beas
  C_Y+2C_H \le \underset{t \in \cD_T}{\esssup}\, Z^j(t) \le \z_Y - C_H+ \z^j, \q
  a.s.,
  \eeas
where $ \z_Y - C_H+ \z^j \in  Dom(\sE)$ thanks to \eqref{eqn-m48},
\eqref{ass_zi} and (D2). Hence Lemma \ref{lem_dom_sharp} implies
that $\underset{t \in \cD_T}{\esssup}\, Z^j(t) \in
Dom(\sE)=Dom^\#(\cE_j)$. Applying Corollary \ref{cor_RCLL} and
(\ref{eqn-cxc11}) again yields \eqref{eqn-lxl01} and
\eqref{eqn-lxl02}.

 \ms To see that $Z^{j,+}$ is a modification of
 $\left\{Z^j(t)\right\}_{t \in  [0,T]}$, it suffices to show that for any
$t \in [0,T]$, $ Z^{j,+}_t \ge Z^j(t) $, a.s. Fix $t \in [0,T]$. For
any $(i, \nu) \in \cI \times \cS_{t, T}$, Definition
\ref{def_stable_class} assures that there exists a $k=k(j,i,t) \in
\cI$ such that $\wt{\cE}_k=\cE^t_{j,i}$. (S1) and (\ref{h_tau_A})
imply that
  \bea \label{eqn-q108}
    H^k_t= H^k_{0,t}  = H^j_{0,t} =H^j_t, \q \hb{and} \q
    H^k_{t,\nu} = H^i_{t,\nu}, \q a.s.
  \eea
 For any  $n \in \hN$, we set
$t_n \dfnn q^+_n(t)$ and define $  \nu_n  \dfnn (\nu+ 2^{-n} ) \land
T \in \cS_{t, T}$.  Let $m \ge n$, it is clear that $t_m \le t_n \le
\nu_n$, a.s. Then Proposition \ref{properties_3} (3) implies that
  \beas
      \wt{\cE}_k\big[Y^k_{\nu_n}\big|\cF_{t_m}\big] 
   =  \wt{\cE}_k\big[Y_{\nu_n}+ H^k_{t_m, \nu_n}  \big|\cF_{t_m}\big]
    +   H^{k}_{t_m} \le
     Z(t_m) + H^{k}_{t_m} =Z^j(t_m) + H^{k}_{t_m} -H^{j}_{t_m}, \q
     a.s.
  \eeas
  As $m \to \infty$, \eqref{eqn-q108} as well as the right-continuity of the
  processes
  $\wt{\cE}_k\big[Y^k_{\nu_n} \big|\cF_\cd\big]$, $H^k$ and $H^j$ imply that
 \beas
  \wt{\cE}_k\big[Y^k_{\nu_n}\big|\cF_t\big]
 =\underset{m \to \infty}{\lim}\wt{\cE}_k\big[Y^k_{\nu_n}\big|\cF_{t_m}\big]
 \le \underset{m \to \infty}{\liminf} Z^j({t_m})
  + H^k_t -H^{j}_t = \underset{m \to \infty}{\liminf} Z^j({t_m})= Z^{j,+}_t, \q a.s.
 \eeas
Since $ \underset{n \to \infty}{\lim} \dneg \da \nu_n =\nu$, a.s.,
 the right continuity of the process $Y^k$ implies that
 $ Y^k_{\nu_n}$ converges a.s. to $Y^k_\nu$, which belongs to $Dom(\sE)$ due to assumption (Y1) and \eqref{eqn-cxc01}.
 Then (\ref{Y_LB2}) and Fatou's
Lemma (Theorem \ref{fatou}) 
imply that
 \beas
 \wt{\cE}_k\big[Y^k_\nu\big|\cF_t\big]  \le \underset{n \to \infty}{\liminf}
\wt{\cE}_k\big[Y^k_{\nu_n}\big|\cF_t\big]
   \le Z^{j,+}_t, \q a.s.
 \eeas
Applying Proposition \ref{properties_3} (5) and (3) to $\wt{\cE}_j$
and $\wt{\cE}_i$ respectively, we can deduce from (\ref{tau_ij2})
and (\ref{eqn-q108}) that
 \bea \label{eqn-q122}
  Z^{j,+}_t \ge \wt{\cE}_k\big[Y^k_\nu\big|\cF_t\big]=
 \cE^t_{j,i}\big[Y^k_\nu \big|\cF_t\big]=\wt{\cE}_j\big[ \wt{\cE}_i \big[Y^k_\nu\big|\cF_t \big]
 \big|\cF_t\big]=\wt{\cE}_i \big[Y^k_\nu\big|\cF_t \big]=\wt{\cE}_i \big[Y_\nu+H^i_{t,\nu} \big|\cF_t \big]+H^j_t
 , \q a.s.
 \eea
Letting $(i, \nu)$ run throughout $\cI \times \cS_{t,
  T}$ yields that
   \beas
   Z^{j,+}_t   \ge  \underset{(i, \nu) \in \cI \times \cS_{t,
  T}}{\esssup}\, \wt{\cE}_i \big[Y_\nu+H^i_{t,\nu} \big|\cF_t \big]  +H^j_t
  = Z(t)+H^j_t =Z^j(t) , \q a.s.,
    \eeas
 which implies that $Z^{j,+} $ is an RCLL modification of $\left\{Z^j(t)\right\}_{t \in
  [0,T]}$. Correspondingly,
   $   Z^0 \dfnn \big\{ Z^{j,+}_t -H^j_t \big\}_{t\in [0,T]} $
 is an RCLL modification of $\left\{Z(t)\right\}_{t \in
  [0,T]}$. Moreover, for any $ i \in \cI$,
  $ Z^{i,0} \dfnn \big\{ Z^0_t+H^i_t\big\}_{t \in [0,T]}$ defines an RCLL modification of
  $\big\{Z^i(t)\big\}_{t \in [0,T]}$, thus it is an
  $\wt{\cE}_i$-supermartingale.  \qed

 \ss \no {\bf Proof of Proposition \ref{Z0}:} For any $t \in [0,T]$, we know from (\ref{eqn-q104x}) and
Proposition \ref{Z_RCLL} that $ Y_t \le Z(t) = Z^0_t$, a.s. Since
the processes $Y$ and $Z^0$ are both right continuous, it follows
from Remark \ref{rem_domin} (2) that $Z^0$ dominates $Y$.

\ms  If $\nu \in \cS^F_{0,T}$ takes values
 in a finite set $\{t_1 < \cds< t_n\}$, for any $\a \in \{1 \cds n\}$,
we can deduce from (\ref{D_1}) that
  \beas
   \b1_{ \{\nu=t_\a  \}}Z(\nu)= \b1_{ \{\nu=t_\a  \}}Z(t_\a)
   =\b1_{ \{\nu=t_\a  \}}Z^0_{t_\a}=\b1_{ \{\nu=t_\a  \}}Z^0_\nu,  \q  a.s.
  \eeas
 Summing the above expression over $\a$, we obtain
  \bea \label{eqn-q124}
   Z^0_\nu=Z(\nu), \q  a.s.
   \eea
For general stopping time $\nu \in \cS_{0,T}$, we let $\left\{\nu_n\right\}_{n \in \hN}$ 
 be a decreasing sequence in $\cS^F_{0,T}$ such that $\underset{n \to
 \infty}{\lim} \neg \da
\nu_n=\nu$, a.s. Thus for any $i \in \cI$,
the right-continuity of the process $Z^{i,0}$ shows that
 \bea \label{eqn-dxd01}
 Z^{i,0}_\nu = \underset{n \to \infty}{\lim} Z^{i,0}_{\nu_n}, \q
 a.s.
 \eea
 For any $n
\in \hN$, (\ref{eqn-q124}) and (\ref{Z_LB}) imply that
   \bea   \label{eqn-q134}
  Z^{i,0}_{\nu_n}=  Z^i(\nu_n) \ge C_Y+2C_H, \q a.s.
   \eea
 If $\cE_j$ satisfies
(\ref{ass_fatou}) for some $j \in \cI$, we can deduce from
(\ref{eqn-xxx01}) and (Y2) that
 \beas  
     \wt{\cE}_j\big[ Z^{j,0}_{\nu_n} \big]=\wt{\cE}_j\big[ Z^j(\nu_n)
    \big] \le 
   Z^j(0)=Z(0)= \underset{(i, \rho ) \in \cI \times \cS_{0,T} }
 {\sup} \wt{\cE}_i\big[Y_\rho+ H^i_\rho \big]<\infty,
 \eeas
thus $\underset{n \to \infty}{\liminf} \wt{\cE}_j\big[
Z^{j,0}_{\nu_n} \big]<\infty $. Then Remark \ref{rem_fatou2} (1)
implies that $ Z^{j,0}_\nu \in Dom(\sE)$.

\ms On the other hand, if no member of $\sE$ satisfies
(\ref{ass_fatou}),
 we suppose that \eqref{ass_zi} holds for some $j \in
\cI$. In light of Proposition \ref{Z_RCLL} and Lemma \ref{Z_bound},
it holds a.s. that
 \beas
 C_Y + 2 C_H \le  Z^{j,0}_t  =Z^0_t +H^j_t  =Z(t)+H^j_t  \le \z_Y-C_H +\z^j  , \q \fa t \in
 \cD_T,
 \eeas
where $ \z_Y-C_H +\z^j \in Dom(\sE)$ thanks to \eqref{eqn-m48},
\eqref{ass_zi} and (D2). Since $Z^{j,0}$ is an RCLL process, it
holds except on a null set $N$ that
 \bea \label{eqn-q734}
     C_Y + 2 C_H \le  Z^{j,0}_t     \le \z_Y-C_H +\z^j
   , \q \fa t \in [0,T], \q
  \hb{thus} \q   C_Y + 2 C_H \le  Z^{j,0}_\nu   \le \z_Y-C_H +\z^j.
 \eea
 Lemma
\ref{lem_dom_sharp} then implies that $Z^{j,0}_\nu \in Dom(\sE)$. We
have seen in both cases that $Z^{j,0}_\nu \in Dom(\sE)$ for some $j
\in \cI$.

\ms Since $Z^{j,0}$ is an RCLL $\wt{\cE}_j$-supermartingale by
Proposition \ref{Z_RCLL}, (\ref{eqn-q134}) and the Optional Sampling
Theorem (Theorem \ref{op_sa2}) imply that
$\wt{\cE}_j\big[Z^{j,0}_{\nu_n} \big|\cF_{\nu_{n+1}}\big] \le
Z^{j,0}_{\nu_{n+1}}$, a.s. for any $n \in \hN$. Applying Corollary
\ref{cor_os2} and Theorem \ref{op_sa2} once again, we obtain
  \bea   \label{eqn-lxl05}
 \wt{\cE}_j\big[Z^{j,0}_{\nu_n} \big|\cF_\nu\big]
 = \wt{\cE}_j\big[ \wt{\cE}_j\big[Z^{j,0}_{\nu_n} \big|\cF_{\nu_{n+1}}\big] \big|\cF_\nu \big] \le
 \wt{\cE}_j\big[Z^{j,0}_{\nu_{n+1}}\big|\cF_\nu\big]\le Z^{j,0}_\nu,
 \q a.s.,
  \eea
 which implies that $ \underset{n \to \infty}{\lim }
 \neg \ua \wt{\cE}_j\big[Z^{j,0}_{\nu_n} \big|\cF_\nu\big]\le Z^{j,0}_\nu$, a.s.
On the other hand, using \eqref{eqn-dxd01} and (\ref{eqn-q134}), we
can deduce from Proposition \ref{properties_3} (5) and Fatou's Lemma
(Theorem \ref{fatou}) that
 \beas
  Z^{j,0}_\nu =\wt{\cE}_j\big[ Z^{j,0}_\nu  \big|\cF_\nu\big]
  \le \underset{n \to \infty}{\lim }
 \neg \ua \wt{\cE}_j\big[ Z^{j,0}_{\nu_n} \big|\cF_\nu\big]\le
  Z^{j,0}_\nu, \q a.s.
 \eeas
 Then \eqref{eqn-q124} and \eqref{eqn-xxx01} imply that
 \bea \label{eqn-k101}
  Z^{j,0}_\nu
  = \underset{n \to \infty}{\lim }
 \neg \ua \wt{\cE}_j\big[ Z^{j,0}_{\nu_n} \big|\cF_\nu\big]
  =\underset{n \to \infty}{\lim }
 \neg \ua  \wt{\cE}_j\big[ Z^j(\nu_n) \big|\cF_\nu\big] \le
  Z^j(\nu),~\;\;
  a.s.,  ~\;\; \hb{thus } ~\; Z^0_\nu  \le  Z(\nu) ~\;\;a.s. 
 \eea

\ms On the other hand, for any $(i, \rho) \in \cI \times \cS_{\nu,
T}$ and $n \in \hN$, we define
 $\rho_n \dfnn \rho \vee \nu_n   \in \cS_{\nu_n,T}$.
  Proposition \ref{properties_3} (3) implies that
  \beas
 \wt{\cE}_i\big[  Y^i_{\rho_n} \big|\cF_{\nu_n}\big]
 = \wt{\cE}_i\big[ Y_{\rho_n}+H^i_{\nu_n,\rho_n}\big|\cF_{\nu_n}\big]+H^i_{\nu_n}
 \le Z(\nu_n) +H^i_{\nu_n} =  Z^i(\nu_n) , \q a.s.
 \eeas
 Taking $\wt{\cE}_i\big[\cd \big|\cF_{\nu}\big]$ on both sides, we see from Corollary \ref{cor_os2} that
 \beas
 \wt{\cE}_i\big[  Y^i_{\rho_n}
  \big|\cF_{\nu}\big] =  \wt{\cE}_i\Big[\wt{\cE}_i\big[  Y^i_{\rho_n}
  \big|\cF_{\nu_n}\big]\Big|\cF_\nu\Big]
\le   \wt{\cE}_i\big[ Z^i(\nu_n)\big|\cF_\nu\big] , \q a.s.
 \eeas
 It is easy to see that $\underset{n \to \infty}{\lim} \dneg \da \rho_n=\rho$, a.s.
  Using the right continuity of processes $Y$ and $H^i$, we
   can deduce from (\ref{Y_LB2}), Fatou's Lemma (Proposition
   \ref{fatou2}) and \eqref{eqn-k101} that
 \beas 
 \wt{\cE}_i \big[  Y^i_\rho  \big|\cF_\nu \big]
   \le  \underset{n \to \infty}{\liminf} \wt{\cE}_i \big[   Y^i_{\rho_n}  \big|\cF_\nu \big]
      \le  \underset{n \to \infty}{\lim } \dneg \ua \wt{\cE}_i\big[ Z^i(\nu_n) \big|\cF_\nu
            \big] =Z^{i,0}_\nu, \q a.s.
 \eeas
 Then Proposition \ref{properties_3} (3) again implies that
 \beas
 \wt{\cE}_i\big[  Y_\rho + H^i_{\nu,\rho} \big|\cF_\nu\big]
 = \wt{\cE}_i\big[  Y^i_\rho \big|\cF_\nu\big]- H^i_\nu
 \le Z^{i,0}_\nu- H^i_\nu=Z^0_\nu  ,\q a.s.
 \eeas
Taking the essential supremum over  $(i, \rho) \in \cI \times
\cS_{\nu, T}$ yields that $Z(\nu) \le Z^0_\nu$, a.s., which in
conjunction with
 (\ref{eqn-k101}) shows that $Z^0_\nu= Z(\nu)$, a.s., thus $Z^0_\nu
\in Dom(\sE)$ by Lemma \ref{lem_Z_nu}. Moreover, for any $i \in
\cI$, we have
 \beas
 Z^{i,0}_\nu=Z^0_\nu+H^i_\nu=Z(\nu)+H^i_\nu=Z^i(\nu), \q a.s.,
 \eeas
 thus $Z^{i,0}_\nu
\in Dom(\sE)$ thanks to Lemma \ref{lem_Z_nu} once again.
(\ref{eqn-h42}) is proved.

\ms  Now let $X $ be another RCLL $\bF$-adapted process dominating
$Y$ such that
 $X^i \dfnn \left\{X_t+H^i_t\right\}_{t \in
 [0,T]}$ is an $\wt{\cE}_i$-supermartingale for any $i \in \cI$.
  We fix $t \in [0,T]$. For any $i \in \cI$ and $\nu \in \cS_{t, T} $, we let
  $\{\nu_n\}_{n \in \hN}$ be a decreasing sequence in $\cS^F_{t,
 T}$ such that $\underset{n \to
 \infty}{\lim} \neg \da
\nu_n=\nu$, a.s. For any $n \in \hN$, since $X^i$ dominates $Y^i$,
Remark \ref{rem_domin} (1) shows that $X^i_{\nu_n} \ge Y^i_{\nu_n}$,
a.s.
 \if{0}
Moreover, (\ref{Y_LB})
  implies that
 \beas
 \underset{t \in \cD_T}{\essinf}\, X^i_t
  \ge  \underset{t \in \cD_T}{\essinf}\, Y^i_t \ge C_*, \q a.s.
 \eeas
  \fi
Then (A4), Proposition \ref{prop_tilde_E} and the Optional Sampling
Theorem (Theorem \ref{op_sa2}) imply that
 \beas
  \wt{\cE}_i\big[  Y_{\nu_n} + H^i_{t,\nu_n} \big|\cF_t\big]
   = \wt{\cE}_i\big[  Y^i_{\nu_n} \big|\cF_t\big]- H^i_t
     \le \wt{\cE}_i\big[ X^i_{\nu_n} \big|\cF_t\big]- H^i_t \le X^i_t- H^i_t= X_t, \q a.s.
 \eeas
 The right-continuity of the processes $Y$ and $H^i$ shows that $ Y_{\nu} + H^i_{t,\nu}
 =\underset{n \to \infty}{\lim} \big( Y_{\nu_n} + H^i_{t,\nu_n} \big)$, a.s.,
 thus it follows from \eqref{Y_LB}, \eqref{H_LB} and Fatou's Lemma
 (Proposition \ref{fatou2}) that
  \beas
 \wt{\cE}_i\big[ Y_{\nu} + H^i_{t,\nu} \big|\cF_t\big] \le \underset{n \to \infty}{\liminf}
 \wt{\cE}_i\big[  Y_{\nu_n} + H^i_{t,\nu_n} \big|\cF_t\big] \le X_t,
 \q a.s.
  \eeas
Taking the essential supremum of the left-hand-side over $(i, \nu)
\in \cI \times \cS_{t, T}$, we can deduce from Proposition
\ref{Z_RCLL} that
 \beas
 Z^0_t= Z(t)= \underset{(i,\nu) \in \cI \times \cS_{t, T} }{\esssup}\,
  \wt{\cE}_i\big[  Y_\nu + H^i_{t,\nu} \big|\cF_t\big]
  \le X_t , \q a.s.
 \eeas
 Since both $Z^0$ and $X$ are RCLL
processes, Remark \ref{rem_domin} (2) once again shows that $X$
dominates $Z^0$.  \qed

  \ss \no {\bf Proof of Lemma \ref{lem_Jd}:}
   For any $i \in \cI$, (\ref{Z0_LB}), (\ref{H_LB})
 as well as Proposition \ref{properties_3} (5) imply that
 \beas
  \wt{\cE}_i\big[Z^0_{\t_\d(\nu)}+ H^i_{\nu,\t_\d(\nu)} \big|\cF_\nu\big]
 \ge \wt{\cE}_i \big[C_*+ C_H   \big|\cF_\nu\big]= C_Y + 2 C_H, \q a.s.
 \eeas
 Taking the essential supremum of the left-hand-side over $i \in \cI$, we can deduce from \eqref{eqn-exe03} that
  \bea \label{eqn-wxw01}
 C_Y+ 2 C_H \le  \underset{i \in
 \cI}{\esssup}\,\wt{\cE}_i\big[Z^0_{\t_\d(\nu)}+ H^i_{\nu,\t_\d(\nu)}
 \big|\cF_\nu\big]=J_\d(\nu)
 \le Z(\nu), \;\; a.s. \q
  \eea
Then Lemma \ref{lem_dom_sharp} imply that $J_\d(\nu) \in Dom(\sE)$.
Let $\si$ be another stopping time in $\cS_{0,T}$. In light of
(\ref{eqn-h42}) and (\ref{D_1}), we see that
 \bea \label{eqn-j15}
 \b1_{\left\{\t_\d(\nu)=\t_\d(\si)\right\}}Z^0_{\t_\d(\nu)}=\b1_{\left\{\t_\d(\nu)
 =\t_\d(\si)\right\}}Z(\t_\d(\nu)) = \b1_{\left\{\t_\d(\nu)=\t_\d(\si)\right\}}
Z(\t_\d(\si))=\b1_{\left\{\t_\d(\nu)=\t_\d(\si)\right\}}
Z^0_{\t_\d(\si)}, \q a.s.
 \eea
It is clear that $\left\{\nu =\si\right\} \subset \left\{
\t_\d(\nu)=\t_\d(\si) \right\}$. Thus multiplying $\b1_{\left\{\nu
=\si\right\}}$ to both sides of (\ref{eqn-j15}) gives that
 \beas
 \b1_{ \left\{\nu =\si\right\} }Z^0_{\t_\d(\nu)} =\b1_{\left\{\nu
=\si\right\}} Z^0_{\t_\d(\si)}, \q a.s.
 \eeas
For any $i \in \cI$, applying Proposition \ref{properties_3} (2) and
recalling how $\wt{\cE}_i [\cd |\cF_\nu ]$ and $\wt{\cE}_i [\cd
|\cF_\si ]$ are defined in (\ref{cE_tau}), we obtain
 \beas
  \b1_{ \left\{\nu =\si\right\} } \wt{\cE}_i\big[  Z^0_{\t_\d(\nu)}+ H^i_{\nu,\t_\d(\nu)}  \big|\cF_\nu\big]
  &=& \b1_{ \left\{\nu =\si\right\} } \wt{\cE}_i\big[ Z^0_{\t_\d(\nu)}+ H^i_{\nu,\t_\d(\nu)} \big|\cF_\si  \big]
  =\wt{\cE}_i\big[ \b1_{ \left\{\nu =\si\right\} } Z^0_{\t_\d(\si)}
  +  \b1_{ \left\{\nu =\si\right\} } H^i_{\si,\t_\d(\si)} \big|\cF_\si \big] \\
  &=& \b1_{ \left\{\nu =\si\right\} } \wt{\cE}_i\big[  Z^0_{\t_\d(\si)}
  +H^i_{\si,\t_\d(\si)}  \big|\cF_\si  \big], \q a.s.,
 \eeas
 where we use the fact that $\left\{\nu=\si \right\} \in \cF_{\nu \land \si}$ thanks to \cite[Lemma
1.2.16]{Kara_Shr_BMSC}. Taking the essential supremum of both sides
over $i \in \cI$, Lemma \ref{lem_ess} (2) implies that
 $$
   \b1_{ \left\{\nu =\si\right\} }J_\d(\nu)
  = \underset{i \in
 \cI}{\esssup}\, \b1_{ \left\{\nu =\si\right\} } \, \wt{\cE}_i\big[
Z^0_{\t_\d(\nu)}+ H^i_{\nu,\t_\d(\nu)} \neg \big|\cF_\nu\big]
   \neg =   \underset{i \in
 \cI}{\esssup}\, \b1_{ \left\{\nu =\si\right\} }
 \wt{\cE}_i\big[  Z^0_{\t_\d(\si)}
  +H^i_{\si,\t_\d(\si)}\neg  \big|\cF_\si  \big]
  \neg = \neg \b1_{ \left\{\nu =\si\right\}}J_\d(\si), ~ a.s.,
 $$
which proves the lemma. \qed

 \ms \no  { \bf Proof of Proposition \ref{prop_Jd}:}

 \ss \no \emph{Proof of 1.}
We fix $ i \in \cI$\, and $\nu, \rho \in \cS_{0,T}$ with $ \nu \le
\rho$, a.s.  Taking $( \nu, \cI',\cU)=\big(  \rho , \cI, \{
\t_\d(\rho) \}\big)$ and $X\big(\t_\d(\rho)\big)=Z^0_{\t_\d(\rho)}$
in Lemma \ref{lem_02}, we can find a sequence $\{j_n\}^\infty_{n=1}$
in $\cI$ such that
  \beas
  J_\d(\rho)=\underset{j \in
 \cI}{\esssup}\,\wt{\cE}_j[Z^0_{\t_\d(\rho)}+ H^j_{\rho, \t_\d(\rho)}  |\cF_\rho]
  = \underset{n \to \infty}{\lim} \dneg \ua
  \wt{\cE}_{j_n}[ Z^0_{\t_\d(\rho)}+ H^{j_n}_{\rho, \t_\d(\rho)} |\cF_\rho], \q  a.s.
  \eeas
  For any $n \in \hN$, it follows from Definition \ref{def_stable_class} that there exists
  a $k_n =k(i,j_n, \rho)\in \cI$ such that
$\wt{\cE}_{k_n}=\cE^\rho_{i,j_n} $.
 \if{0}
  \bea \label{eqn-q110}
H^{k_n}_{\t_\d(\rho)}- H^{k_n}_{\rho}=H^{j_n}_{
\t_\d(\rho)}-H^{j_n}_{\rho } \q \hb{and}\q
H^{k_n}_{\rho}-H^{k_n}_{\nu} = H^i_{\rho}-H^i_{\nu}, \q a.s.
 \eea
 \fi
 Applying Proposition \ref{properties_3} (3) to $\wt{\cE}_{j_n}$, we
can deduce from (\ref{tau_ij2}) and (\ref{h_tau_A}) that
 \bea \label{eqn-j01}
   \wt{\cE}_{k_n}\big[Z^0_{\t_\d(\rho)}+H^{k_n}_{\nu,\t_\d(\rho)} \big|\cF_\nu \big]
  &=&  \cE^\rho_{i,j_n}\big[Z^0_{\t_\d(\rho)}+H^{k_n}_{\nu,\t_\d(\rho)} \big|\cF_\nu  \big]
  =  \wt{\cE}_i\big[ \wt{\cE}_{j_n} \big[Z^0_{\t_\d(\rho)}
 +H^{k_n}_{\nu,\t_\d(\rho)}\big|\cF_\rho \big] \big|\cF_\nu  \big]\nonumber \\
 &=& \wt{\cE}_i\big[ \wt{\cE}_{j_n} \big[Z^0_{\t_\d(\rho)}
 +H^{j_n}_{\rho,\t_\d(\rho)}\big|\cF_\rho \big] +H^i_{\nu,\rho}\big|\cF_\nu  \big] , \q a.s.
 \eea
Since $ \nu \le \rho$, a.s., we see that $\t_\d(\nu) \le
\t_\d(\rho)$, a.s. Due to (\ref{eqn-h42}) and (\ref{D_2}), we have
that
 \beas
    \wt{\cE}_{k_n}\big[Z^0_{\t_\d(\rho)}+H^{k_n}_{\t_\d(\nu),\t_\d(\rho)} \big|\cF_{\t_\d(\nu)}\big]
     \le  \underset{j \in \cI}{\esssup}\, \wt{\cE}_j\big[Z(\t_\d(\rho))
    +H^j_{\t_\d(\nu),\t_\d(\rho)} \big|\cF_{\t_\d(\nu)}\big]  \le Z(\t_\d(\nu))=
   Z^0_{\t_\d(\nu)},\q a.s.
 \eeas
 Then using Corollary \ref{cor_os2} and applying Proposition \ref{properties_3} (3) to $\wt{\cE}_{k_n}$, we obtain
 \beas
  && \hspace{-2cm} \wt{\cE}_{k_n}\big[Z^0_{\t_\d(\rho)} +H^{k_n}_{\nu,\t_\d(\rho)} \big|\cF_\nu \big]
   =  \wt{\cE}_{k_n}\big[\wt{\cE}_{k_n}\big[Z^0_{\t_\d(\rho)}
  +H^{k_n}_{\nu,\t_\d(\rho)}\big|\cF_{\t_\d(\nu)}\big]  \big|\cF_\nu\big]\\
   &=& \wt{\cE}_{k_n}\big[\wt{\cE}_{k_n}\big[Z^0_{\t_\d(\rho)}
  +H^{k_n}_{\t_\d(\nu),\t_\d(\rho)}\big|\cF_{\t_\d(\nu)}\big]+H^{k_n}_{\nu,\t_\d(\nu)}
  \big|\cF_\nu\big] \\
   &\le & \wt{\cE}_{k_n}\big[ Z^0_{\t_\d(\nu)}+H^{k_n}_{\nu,\t_\d(\nu)}  \big|\cF_\nu\big]
    \le  \underset{j \in \cI}{\esssup}\, \wt{\cE}_j\big[ Z^0_{\t_\d(\nu)}+H^j_{\nu,\t_\d(\nu)}\big|\cF_\nu\big]
       =J_\d(\nu), \q a.s.,
 \eeas
 which together with (\ref{eqn-j01}) shows that
 \beas
  \wt{\cE}_i\big[ \wt{\cE}_{j_n} \big[Z^0_{\t_\d(\rho)}
 +H^{j_n}_{\rho,\t_\d(\rho)}\big|\cF_\rho \big] +H^i_{\nu,\rho}\big|\cF_\nu  \big]
 \le J_\d(\nu), \q  a.s.
 \eeas
 For any $n \in \hN$, we see from (\ref{Z0_LB}), (\ref{H_LB}) and Proposition \ref{properties_3} (5) that
  \beas
 \wt{\cE}_{j_n}\big[Z^0_{\t_\d(\rho)}+ H^{j_n}_{\rho,\t_\d(\rho)} \big|\cF_\rho\big]
  + H^i_{\nu,\rho} \ge  \wt{\cE}_{j_n}\big[ C_Y+ 2C_H \big|\cF_\rho\big]
  + C_H = C_Y+3C_H,  \q a.s.
  \eeas
Then Fatou's Lemma (Proposition \ref{fatou2}) implies that
  \beas
 \wt{\cE}_i \big[ J_\d(\rho)+ H^i_{\nu,\rho} \big|\cF_\nu \big] \le \underset{n \to \infty}{\lim} \dneg \ua
  \wt{\cE}_i \big[\wt{\cE}_{j_n}\big[Z^0_{\t_\d(\rho)}+ H^{j_n}_{\rho,\t_\d(\rho)} \big|\cF_\rho\big]
  + H^i_{\nu,\rho} \big|\cF_\nu
  \big] \le  J_\d(\nu), \q a.s.
  \eeas
  For any
  $\si \in \cS_{0,T}$, Lemma \ref{lem_Jd}, (\ref{eqn-cxc01}) and
  (D2) show that
  $ 
   J^i_\d(\si) \dfnn J_\d(\si)+ H^i_\si \in Dom(\sE)
  $.
  A simple application of Proposition \ref{properties_3} (3) yields that
   \bea \label{eqn-q140}
  \wt{\cE}_i \big[ J^i_\d(\rho)\big|\cF_\nu \big] = \wt{\cE}_i \big[ J_\d(\rho)+ H^i_{\nu,\rho} \big|\cF_\nu \big]
  + H^i_\nu \le J_\d(\nu)+ H^i_\nu=J^i_\d(\nu), \q a.s.
   \eea
In particular, when $0 \le s < t \le T$, we have
 $ \wt{\cE}_i \big[ J^i_\d(t) \big|\cF_s \big]  \le  J^i_\d(s)$,
 a.s., which show that $\left\{J^i_\d(t)\right\}_{t \in [0,T]}$ is an $\wt{\cE}_i
 $-supermartingale. \\

 \noindent \emph{Proof of 2.} For any $i \in \cI$ and $\nu \in \cS_{0,T}$,
 \eqref{eqn-wxw01} and \eqref{H_LB} imply that
  \bea \label{eqn-new11}
  J^i_\d(\nu)=J_\d(\nu)+H^i_\nu \ge C_Y+3C_H
  , \q  a.s.
   \eea
    In particular, $J^i_\d(T) \ge C_Y+3C_H $, a.s.
 Proposition \ref{prop_tilde_E} and Theorem \ref{upcrossing} then
show that $ \Big\{ J^{\d,i,+}_t \dfnn \underset{n \to
\infty}{\liminf} J^i_\d\big(q^+_n(t)\big)\Big\}_{t \in [0,T]} $
  defines an RCLL process.
 \if{0}
  For any $t \in [0,T]$ and $\tilde{i} \in \cI $,
  (\ref{Z0_LB}), (\ref{H_LB}) as well as Proposition \ref{properties_3}
  (5) imply that
 \beas
 \wt{\cE}_{\tilde{i}}\big[Z^0_{\t_\d(t)}+ H^{\tilde{i}}_{t, \t_\d(t)}  \big|\cF_t\big]
 \ge \wt{\cE}_{\tilde{i}}\big[ C_Y+2C_H \big|\cF_t\big]=
 C_Y+2C_H, \q a.s.
 \eeas
 Taking essential supremum over $\tilde{i} \in \cI$, we can deduce from (\ref{eqn-exe03}) that
 \beas
  C_Y+2C_H &\le&   \underset{\tilde{i} \in \cI}{\esssup}\,
   \wt{\cE}_{\tilde{i}}\big[Z^0_{\t_\d(t)}+ H^{\tilde{i}}_{t, \t_\d(t)}  \big|\cF_t\big]
   =  J_\d(t) \le    Z(t) , \q a.s.,
 \eeas
Then (\ref{H_LB}) shows that
 \fi
 Then \eqref{eqn-new11} implies that
 \bea \label{eqn-dxd03}
   \underset{t \in [0,T]}{\essinf}\, J^i_\d(t) \ge C_Y+3C_H, \q  a.s.
 \eea

  If $\cE_j$ satisfies
(\ref{ass_fatou}) for some $j \in \cI$,
 Corollary \ref{cor_RCLL} and (\ref{eqn-dxd03}) imply that
 \bea
 J^{\d,j,+}_\nu \in Dom^\#(\cE_j)=Dom(\sE), \q \fa \nu \in \cS_{0,T}, \hspace{3.5cm} \label{eqn-lxl11}\\
 \hb{and that $J^{\d,j,+}$ is an RCLL $ \wt{\cE}_j$-supermartingale
 such that for any $t \in [0,T]$, $J^{\d,j,+}_t \le J^j_\d(t) $,\; a.s.} \label{eqn-lxl12}
 \eea
 Otherwise, if no member of $\sE$ satisfies (\ref{ass_fatou}),
  we suppose that
\eqref{ass_zi} holds for some $j \in \cI$. Then \eqref{eqn-new11},
\eqref{eqn-exe03} and \eqref{Z_UB} imply that for any $t \in \cD_T$,
 \beas
  C_Y+3C_H \le  J^j_\d(t) = J_\d(t)+H^j_t \le  Z(t)  + H^j_t  \le  \z_Y-C_H +\z^j   , \q  a.s.
 \eeas
 Taking essential supremum of $J^j_\d(t)$ over $t \in \cD_T$ yields that
 \beas
  C_Y+3C_H   \le \underset{t \in \cD_T}{\esssup}\, J^j_\d(t) \le  \z_Y-C_H+\z^j , \q  a.s.,
 \eeas
  where $ \z_Y - C_H+ \z^j \in  Dom(\sE)$ thanks to \eqref{eqn-m48},
\eqref{ass_zi} and (D2).  Hence Lemma \ref{lem_dom_sharp} implies
that $ \underset{t \in \cD_T}{\esssup}\, J^j_\d(t) \in
Dom(\sE)=Dom^\#(\cE_j) $. Applying Corollary \ref{cor_RCLL} and
(\ref{eqn-dxd03}) again yields \eqref{eqn-lxl11} and
\eqref{eqn-lxl12}.

  \ms To see that $ J^{\d,j,+}$ is a modification of $\big\{J^j_\d(t)\big\}_{t \in
  [0,T]}$, it suffices to show that for any $t \in [0,T]$, $  J^{\d,j,+}_t \ge J^j_\d(t) $, a.s.
   Fix $t \in [0,T]$. For any $i \in \cI  $, Definition
\ref{def_stable_class} assures that there exists a $k=k(j,i,t) \in
\cI$ such that $\wt{\cE}_k=\cE^t_{j,i}$. Moreover, (S1) and
(\ref{h_tau_A}) imply that
  \bea \label{eqn-q108b}
   H^k_t= H^k_{0,t} = H^j_{0,t} =H^j_t, \q \hb{and} \q
   H^k_{t, \t_\d(t)} = H^i_{t, \t_\d(t)} , \q  a.s.
  \eea
  For any $n \in \hN$, we set $t_n \dfnn q^+_n(t)$.
Let $m \ge n$, it is clear that $t_m \le t_n$, a.s. Then
(\ref{eqn-h42}), Corollary \ref{cor_os2}, Proposition
\ref{properties_3} (3) as well as (\ref{D_2}) imply that
 \beas 
     \wt{\cE}_k\big[Z^{k,0}_{\t_\d(t_n)}\big|\cF_{t_m}\big] & =&
     \wt{\cE}_k\big[Z^k(\t_\d(t_n))\big|\cF_{t_m}\big]=
      \wt{\cE}_k\big[\wt{\cE}_k\big[Z^k(\t_\d(t_n)) \big|\cF_{\t_\d(t_m)}\big]\big|\cF_{t_m}\big]\nonumber \\
      &=& \wt{\cE}_k\big[\wt{\cE}_k\big[Z(\t_\d(t_n))
      + H^k_{\t_\d(t_m),\t_\d(t_n)} \big|\cF_{\t_\d(t_m)}\big]
      +H^k_{t_m,\t_\d(t_m)}\big|\cF_{t_m}\big] + H^k_{t_m} \nonumber \\
     &\le & \wt{\cE}_k\big[ \underset{l \in \cI}{\esssup}\, \wt{\cE}_l\big[Z(\t_\d(t_n))
      + H^l_{\t_\d(t_m),\t_\d(t_n)}  \big|\cF_{\t_\d(t_m)}\big]
      +H^k_{t_m,\t_\d(t_m)}\big|\cF_{t_m}\big] + H^k_{t_m} \nonumber \\
     &\le & \wt{\cE}_k\big[  Z(\t_\d(t_m))
      +H^k_{t_m,\t_\d(t_m)}\big|\cF_{t_m}\big] + H^k_{t_m} 
       \le  \underset{l \in \cI}{\esssup}\, \wt{\cE}_l\big[
       Z^0_{\t_\d(t_m)}
      +H^l_{t_m,\t_\d(t_m)}\big|\cF_{t_m}\big] + H^k_{t_m} \nonumber \\
       &=& J_\d(t_m) + H^k_{t_m} = J^j_\d(t_m) + H^k_{t_m}-H^j_{t_m}  , \q a.s.
  \eeas
As $m \to \infty$, (\ref{eqn-q108b}) as well as the right-continuity
of the processes
$\wt{\cE}_k\big[Z^{k,0}_{\t_\d(t_n)}\big|\cF_\cd\big] $, $H^k$ and
$H^j$ imply that
 \beas
 \wt{\cE}_k\big[Z^{k,0}_{\t_\d(t_n)}\big|\cF_t\big]
 = \underset{m \to \infty}{\lim}\wt{\cE}_k\big[Z^{k,0}_{\t_\d(t_n)}\big|\cF_{t_m}\big]
 \le \underset{m \to \infty}{\liminf}J^j_\d(t_m)+H^k_t-H^j_t= \underset{m \to \infty}{\liminf}J^j_\d(t_m)
  =J^{\d,j,+}_t, \q a.s.
 \eeas
 Since  $\underset{n \to \infty}{\lim} \dneg \da \t_\d(t_n)=\t_\d(t)$
 a.s., the right-continuity of the process $Z^{k,0}$ implies
 that $Z^{k,0}_{\t_\d(t_n)}$ converges a.s. to $Z^{k,0}_{\t_\d(t)}$,
 which belongs to $Dom(\sE)$ thanks to
  Proposition \ref{Z0}. Then
 (\ref{Z0_LB}) and Fatou's Lemma (Theorem \ref{fatou}) imply that
  \beas
 \wt{\cE}_k\big[Z^{k,0}_{\t_\d(t)}\big|\cF_t\big]   \le  \underset{n \to \infty}{\liminf}
 \wt{\cE}_k\big[Z^{k,0}_{\t_\d(t_n)}\big|\cF_t\big] \le J^{\d,j,+}_t, \q a.s.
 \eeas
 Similar to (\ref{eqn-q122}), we can deduce from
(\ref{tau_ij2}) and (\ref{eqn-q108b}) that
 \beas
 J^{\d,j,+}_t \ge \wt{\cE}_i
\big[Z^0_{\t_\d(t)}+H^i_{t,\t_\d(t)} \big|\cF_t \big]+H^j_t, \q a.s.
 \eeas
Letting $i$ run throughout $\cI $ yields that
   \beas
   J^{\d,j,+}_t    \ge  \underset{i \in \cI }{\esssup}\,
   \wt{\cE}_i \big[Z^0_{\t_\d(t)}+H^i_{t,\t_\d(t)} \big|\cF_t \big]+H^j_t
  = J_\d(t) +H^j_t =J^j_\d(t) , \q a.s.,
    \eeas
 which implies that ${}J^{\d,j,+}$ is an RCLL modification of $\big\{J^j_\d(t)\big\}_{t \in
  [0,T]}$. 
Correspondingly, $  J^{\d,0}  \dfnn \big\{ J^{\d, j,+}_t -H^j_t
\big\}_{t\in [0,T]}$ is an RCLL modification of
$\big\{J_\d(t)\big\}_{t \in [0,T]}$. Moreover, for any $ i \in \cI$,
$ J^{\d,i,0} \dfnn \big\{ J^{\d,0}_t+H^i_t\big\}_{t \in [0,T]}$
defines an RCLL modification of
  $\big\{J^i_\d(t)\big\}_{t \in [0,T]}$, thus
  it is an $\wt{\cE}_i$-supermartingale. \\

\noindent \emph{Proof of 3.} Now let us show (3). Similar to
(\ref{eqn-q124}), we can deduce from Lemma \ref{lem_Jd} that for any
$\nu \in \cS^F_{0,T}$
  \bea \label{eqn-q124b}
     J^{\d,0}_\nu=J_\d(\nu), \q  a.s.
  \eea
For a general stopping time $\nu \in \cS_{0,T}$, we let
$\left\{\nu_n\right\}_{n \in \hN}$ 
 be a decreasing sequence
 in $\cS^F_{0,T}$ such that $\underset{n \to \infty}{\lim} \neg \da
\nu_n=\nu$, a.s. Thus for any $i \in \cI$,
 the right-continuity of the process $J^{\d, i,0}$ shows that
 \bea \label{eqn-dxd01b}
 J^{\d, i,0}_\nu = \underset{n \to \infty}{\lim} J^{\d,
 i,0}_{\nu_n}, \q a.s.
 \eea
 In light of (\ref{eqn-q124b}) and \eqref{eqn-new11}, it holds a.s.
 that
 \beas
  J^{\d,i,0}_t=  J^i_\d(t) \ge C_Y+3C_H, \q \fa t \in \cD_T.
   \eeas
Since $J^{\d,i,0}$ is an RCLL process, it holds except on a null set
$N$ that
   \bea   \label{eqn-q134b}
  J^{\d,i,0}_t  \ge C_Y+3C_H, \q \fa t \in [0,T], \q \hb{thus} \q J^{\d,i,0}_\si  \ge
  C_Y+3C_H, \q \fa \si \in \cS_{0,T}.
   \eea
 If $\cE_j$ satisfies
(\ref{ass_fatou}) for some $j \in \cI$, we can deduce from
(\ref{eqn-q140}), (\ref{eqn-exe03}) and
 (Y2) that
 \beas
     \wt{\cE}_j\big[ J^{\d,j,0}_{\nu_n} \big]=\wt{\cE}_j\big[ J^j_\d(\nu_n)
    \big] \le
   J^j_\d(0) =J_\d(0) \le  Z(0)   = \underset{(i, \rho ) \in \cI \times \cS_{0,T} }
 {\sup} \wt{\cE}_i\big[Y_\rho+ H^i_\rho \big]<\infty,
 \eeas
thus $\underset{n \to \infty}{\liminf} \wt{\cE}_j\big[
J^{\d,j,0}_{\nu_n} \big]<\infty $. Then Remark \ref{rem_fatou2} (1)
implies that $ J^{\d,j,0}_\nu \in  Dom(\sE)$.

 \ms On the other hand, if no member of $\sE$ satisfies (\ref{ass_fatou}),
we suppose that \eqref{ass_zi} holds for some $j \in \cI$. In light
of (\ref{eqn-q124b}),  (\ref{eqn-exe03}) and \eqref{Z_UB}, it holds
a.s. that
 \beas
    J^{\d,j,0}_t  = J^j_\d(t)=J_\d(t)+H^j_t  \le  Z(t)+H^j_t \le \z_Y-C_H +\z^j  , \q \fa t \in \cD_T,
 \eeas
 where $\z_Y-C_H +\z^j \in
Dom(\sE)$ thanks to \eqref{eqn-m48}, \eqref{ass_zi} and (D2). Since
$J^{\d,j,0}$ is an RCLL process, it holds except on a null set $N'$
that
 \bea \label{eqn-q734b}
        J^{\d,j,0}_t     \le \z_Y-C_H +\z^j
   , \q \fa t \in [0,T], \q
  \hb{thus} \q      J^{\d,j,0}_\nu   \le \z_Y-C_H
  +\z^j.
 \eea
  Then \eqref{eqn-q134b} and Lemma \ref{lem_dom_sharp} imply that
  $J^{\d,j,0}_\nu \in Dom(\sE)$.
  We have seen in both cases that $J^{\d,j,0}_\nu \in
Dom(\sE)$ for some $j \in \cI$.

\ms Similar to the arguments used in (\ref{eqn-lxl05}) through
 (\ref{eqn-k101}) \big(with (\ref{eqn-q124b})-\eqref{eqn-q134b}
   replacing (\ref{eqn-q124})-(\ref{eqn-q134}) respectively,
   and with \eqref{eqn-q140} replacing \eqref{eqn-xxx01}\big), we can deduce that
\bea \label{eqn-j20}
  J^{\d, j,0}_\nu =\underset{n \to \infty}{\lim }\dneg \ua \wt{\cE}_j\big[ J^{\d, j,0}_{\nu_n} \big|\cF_\nu\big]
     =\underset{n \to \infty}{\lim }\dneg \ua \wt{\cE}_j\big[ J^j_\d(\nu_n) \big|\cF_\nu\big]
     \le J^j_\d(\nu),  \q \hb{thus} \q
 J^{\d, 0}_\nu \le J_\d(\nu), \q a.s.
 \eea
 The right-continuity of the process $J^{\d,0}$, \eqref{eqn-q124b} and \eqref{eqn-wxw01} show that
 \beas
 J^{\d,0}_\nu = \underset{n \to \infty}{\lim} J^{\d,
 0}_{\nu_n}= \underset{n \to \infty}{\lim} J_\d  (\nu_n) \ge C_Y+2C_H , \q a.s.
 \eeas
 Lemma \ref{lem_Jd} and Lemma \ref{lem_dom_sharp} thus imply that
$J^{\d,0}_\nu \in Dom(\sE) $. For any $i \in \cI$, (\ref{eqn-cxc01})
and (D2) show that $J^{\d,i,0}_\nu = J^{\d,0}_\nu  + H^i_\nu \in
Dom(\sE)$.

 \ms On the other hand, for any $i \in \cI$ and $n \in \hN$, it is clear that $\nu \le \nu_n \le
\t_\d(\nu_n)$, a.s. Then Corollary \ref{cor_os2}, \eqref{eqn-j20}
and (\ref{eqn-q140}) imply that
  \beas
     \wt{\cE}_i\big[J^{\d,i,0}_{\t_\d(\nu_n)}\big|\cF_{\nu}\big]
    =   \wt{\cE}_i\big[\wt{\cE}_i\big[J^{\d,i,0}_{\t_\d(\nu_n)} \big|\cF_{\nu_n}\big]\big|\cF_{\nu}\big]
      \le \wt{\cE}_i\big[\wt{\cE}_i\big[J^i_\d \big(\t_\d(\nu_n)\big) \big|\cF_{\nu_n}\big]\big|\cF_{\nu}\big]
     \le \wt{\cE}_i \big[J^i_\d(\nu_n)  \big|\cF_{\nu} \big]
    , \q a.s.
  \eeas
  It is easy to see that $\underset{n \to \infty}{\lim} \dneg \da \t_\d(\nu_n)=\t_\d(\nu)$, a.s.
  Using the right continuity of the process $J^{\d,i,0}$, we
   can deduce from (\ref{eqn-q134b}), Fatou's Lemma (Proposition
   \ref{fatou2}) and \eqref{eqn-j20} that
  \beas
  \wt{\cE}_i\big[J^{\d,i,0}_{\t_\d(\nu)}\big|\cF_\nu\big] \le
  \underset{n \to \infty}{\liminf} \wt{\cE}_i\big[J^{\d,i,0}_{\t_\d(\nu_n)}\big|\cF_\nu\big]
  \le \underset{n \to \infty}{\lim} \dneg \ua \wt{\cE}_i \big[ J^i_\d(\nu_n) \big|\cF_\nu\big]= J^{\d,i,0}_\nu, \q a.s.
 \eeas
 Proposition \ref{properties_3} (3) further implies that
  \beas
    \wt{\cE}_i\big[J^{\d,0}_{\t_\d(\nu)}+H^i_{\nu,\t_\d(\nu)} \big|\cF_\nu\big]
   = \wt{\cE}_i\big[J^{\d,i,0}_{\t_\d(\nu)}\big|\cF_\nu\big] -H^i_\nu \le
    J^{\d,i,0}_\nu-H^i_\nu=J^{\d,0}_\nu, \q a.s.
  \eeas
Taking the essential supremum over $i \in \cI$ gives
  \beas
   J_\d(\nu)= \underset{i \in \cI}{\esssup}\,\wt{\cE}_i\big[J^{\d,0}_{\t_\d(\nu)}+H^i_{\nu,\t_\d(\nu)} \big|\cF_\nu\big]
   \le J^{\d,0}_\nu, \q a.s.,
   \eeas
   which together with (\ref{eqn-j20}) shows that $J^{\d,0}_\nu=J_\d(\nu)
 $, a.s.
 \qed

\ss \no {\bf Proof of Theorem \ref{SN_exist}:} We first show that
for any $\d \in (0,1)$ and $\nu \in
 \cS_{0,T}$
 \bea \label{eqn-k04}
 J_\d(\nu)= Z^0_\nu= Z(\nu) , \q a.s.
 \eea
 Fix $i \in \cI$. Lemma \ref{lem_pconvex} indicates that
 $\wt{\cE}_i$ is a convex $\bF$-expectation on $Dom(\sE)$.
  Since $Z^{i,0}$ and $J^{\d,i,0}$ are both
 $\wt{\cE}_i$-supermartingales, 
 we can deduce that for any $0\le s < t \le T$,
  \beas
  && \hspace{-2cm} \wt{\cE}_i\big[\d Z^0_t+(1-\d)J^{\d,0}_t+H^i_t \big|\cF_s\big]
 = \wt{\cE}_i\big[\d Z^{i,0}_t+(1-\d)J^{\d,i,0}_t \big|\cF_s\big]
  \le \d \wt{\cE}_i\big[  Z^{i,0}_t \big|\cF_s\big]+(1-\d) \wt{\cE}_i\big[ J^{\d,i,0}_t  \big|\cF_s\big] \\
 &\le& \d  Z^{i,0}_s  +(1-\d)  J^{\d,i,0}_s= \d  Z^0_s  +(1-\d)  J^{\d,0}_s + H^i_s , \q  a.s.,
  \eeas
 which shows that $\left\{\d Z^0_t+(1-\d)J^{\d,0}_t+H^i_t\right\}_{t \in [0,T]}$
 is an RCLL $\wt{\cE}_i$-supermartingale.

 \ms Now we fix $t \in [0,T]$ and define $A \dfnn \left\{\t_\d(t)=t\right\} \in \cF_t$.
  Using Proposition \ref{prop_Jd} (3), Lemma \ref{lem_ess} (2)
  as well as applying Proposition \ref{properties_3} (2) and (5) to each $\wt{\cE}_i$, we obtain
   \beas
  \b1_A J^{\d,0}_t
 &= &\b1_A J_\d(t) = \b1_A \underset{i \in
 \cI}{\esssup}\,\wt{\cE}_i\big[Z^0_{\t_\d(t)}+ H^i_{t,\t_\d(t)} \big|\cF_t\big]
 = \underset{i \in \cI}{\esssup}\,\wt{\cE}_i\big[\b1_A
 \big(Z^0_{\t_\d(t)}+ H^i_{t,\t_\d(t)}\big) \big|\cF_t\big] \nonumber \\
 &=& \underset{i \in \cI}{\esssup}\,\wt{\cE}_i\big[\b1_A
  Z^0_t  \big|\cF_t\big]=\b1_A  Z^0_t, \q a.s.
 \eeas
Then (\ref{eqn-h42}) and (\ref{eqn-q104x}) imply that
 \bea \label{eqn-k03}
 \b1_A \big(\d  Z^0_t+ (1-\d) J^{\d,0}_t\big)=\b1_A  Z^0_t = \b1_A  Z(t) \ge \b1_A Y_t, \q a.s.
  \eea
  Moreover, we see from the definition of $\t_\d(t)$ that for any $\o \in A^c$
  \bea \label{eqn-exe01}
  Y_s(\o)<  \d Z^0_s(\o)+(1-\d)(C_Y+2C_H), \qq \fa s \in \big[t, \t_\d\big(t\big)(\o)\big)
  \eea
 Since both $Z^0$ and $Y$ are right-continuous processes, (\ref{eqn-exe01}) and \eqref{eqn-wxw01} imply that
  \beas
   Y_t \le \d Z^0_t +(1-\d)(C_Y+2C_H) \le \d Z^0_t + (1-\d) J^{\d,0}_t \q  a.s. \hb{ on $A^c$},
   \eeas
 which in conjunction with (\ref{eqn-k03}) and Remark \ref{rem_domin} (2) shows that
  the RCLL process $\d Z^0 +(1-\d) J^{\d,0} $ dominates $Y$,
  thus dominates $Z^0$ thanks to Proposition \ref{Z0}. It follows that $J^{\d,0}$ also dominates
 $Z^0$. Then for any $\nu \in \cS_{0,T}$,
 Proposition \ref{prop_Jd} (3), Remark \ref{rem_domin} (1) and (\ref{eqn-h42}) imply that
 $  J_\d(\nu) = J^{\d,0}_\nu \ge Z^0_\nu =  Z(\nu) $, a.s.,
The reverse inequality comes from (\ref{eqn-exe03}). This proves
(\ref{eqn-k04}).

 \ms Next, we fix $\nu \in \cS_{0,T}$ and set $\d^n = \frac{n-1}{n}$, $n \in
 \hN$. It is clear that the sequence $ \big\{ \t_{\d^n}(\nu)\big\}_{n \in \hN}$ increasing a.s. to $\ol{\t}(\nu)$.
  Since the family of processes $ \{Y^i\}_{i \in
\cI }$ is ``$\sE$-uniformly-left-continuous", we can find a
subsequence $ \{ \d^{n_{\neg k}} \}_{k \in \hN }$ of $ \{ \d^n \}_{n
\in \hN }$ such that
  \bea \label{eqn-x23}
  \underset{k \to \infty}{\lim} \underset{i \in \cI}{\esssup}\,
\Big|\wt{\cE}_i\big[\hb{$\frac{n_{\neg k}}{n_{\neg k}-1}$}
Y_{\t_{\d^{n_{\neg k}}}(\nu) }+H^i_{\t_{\d^{n_{\neg k}}}(\nu)}
\big|\cF_\nu\big]- \wt{\cE}_i\big[ Y^i_{\ol{\t}(\nu)
}\big|\cF_\nu\big]\Big|=0, \q a.s.
  \eea
    For any $i \in \cI$ and $k \in \hN$, Remark \ref{delta_t} (1)
   implies that $ Y_{\t_{\d^{n_{\neg k}}}(\nu)} \ge \d^{n_{\neg k}} Z^0_{\t_{\d^{n_{\neg
   k}}}(\nu)}+\big(1-\d^{n_{\neg k}} \big) (C_Y+2C_H)
   $, a.s. Hence Proposition \ref{properties_3} (3) shows that
 \bea
 && \hspace{-1.5cm} \wt{\cE}_i \big[Z^0_{\t_{\d^{n_{\neg k}}}(\nu)}+H^i_{\nu,\t_{\d^{n_{\neg k}}}(\nu)} \big|\cF_\nu \big]
 +\hb{$\frac{1}{n_k-1}(C_Y+2C_H)$} \le
  \wt{\cE}_i\big[\hb{$\frac{n_{\neg k}}{n_{\neg k}-1}$}
Y_{\t_{\d^{n_{\neg k}}}(\nu) }+H^i_{\t_{\d^{n_{\neg k}}}(\nu)} \big|\cF_\nu\big]-H^i_\nu  \nonumber \\
 &=&\wt{\cE}_i\big[\hb{$\frac{n_{\neg k}}{n_{\neg k}-1}$}
Y_{\t_{\d^{n_{\neg k}}}(\nu) }+H^i_{\t_{\d^{n_{\neg k}}}(\nu)}
\big|\cF_\nu\big]- \wt{\cE}_i \big[ Y^i_{\ol{\t}(\nu) }
\big|\cF_\nu\big] + \wt{\cE}_i \big[ Y_{\ol{\t}(\nu) } +
H^i_{\nu,\ol{\t}(\nu)}  \big|\cF_\nu\big]  \nonumber \\
& \le & \underset{i \in \cI}{\esssup}\,
\Big|\wt{\cE}_i\big[\hb{$\frac{n_{\neg k}}{n_{\neg k}-1}$}
Y_{\t_{\d^{n_{\neg k}}}(\nu) }+H^i_{\t_{\d^{n_{\neg k}}}(\nu)}
\big|\cF_\nu\big]- \wt{\cE}_i\big[ Y^i_{\ol{\t}(\nu)
}\big|\cF_\nu\big]\Big|+
  \wt{\cE}_i \big[ Y_{\ol{\t}(\nu) } + H^i_{\nu,\ol{\t}(\nu)} \big|\cF_\nu\big] \label{eqn-q150} \\
& \le & \underset{i \in \cI}{\esssup}\,
\Big|\wt{\cE}_i\big[\hb{$\frac{n_{\neg k}}{n_{\neg k}-1}$}
Y_{\t_{\d^{n_{\neg k}}}(\nu) }+H^i_{\t_{\d^{n_{\neg k}}}(\nu)}
\big|\cF_\nu\big]- \wt{\cE}_i\big[ Y^i_{\ol{\t}(\nu)
}\big|\cF_\nu\big]\Big|+
 \underset{i \in \cI}{\esssup}\, \wt{\cE}_i
\big[ Y_{\ol{\t}(\nu) } + H^i_{\nu,\ol{\t}(\nu)} \big|\cF_\nu\big] ,
\q a.s. \nonumber
 \eea
Taking the $\esssup$ of the left-hand-side over $\cI$, we see from
(\ref{eqn-k04}) that
 \beas
 && \hspace{-2
 cm} \underset{i \in \cI}{\esssup}\,
\Big|\wt{\cE}_i\big[\hb{$\frac{n_{\neg k}}{n_{\neg k}-1}$}
Y_{\t_{\d^{n_{\neg k}}}(\nu) }+H^i_{\t_{\d^{n_{\neg k}}}(\nu)}
\big|\cF_\nu\big]- \wt{\cE}_i\big[ Y^i_{\ol{\t}(\nu)
}\big|\cF_\nu\big]\Big|+
 \underset{i \in \cI}{\esssup}\, \wt{\cE}_i
\big[ Y_{\ol{\t}(\nu) } + H^i_{\nu,\ol{\t}(\nu)}
\big|\cF_\nu\big]  \nonumber \\
 &\ge&  J_{\d^{n_{\neg k}}}(\nu)+\hb{$\frac{1}{n_k-1}(C_Y+2C_H)$}
 =Z(\nu)+\hb{$\frac{1}{n_k-1}(C_Y+2C_H)$}, \q a.s., 
 \eeas
As $k \to \infty$, (\ref{eqn-x23}), (\ref{eqn-q104x}) and
(\ref{D_2}) imply that
 \beas
 Z(\nu) \le \underset{i \in \cI}{\esssup}\, \wt{\cE}_i \big[ Y_{\ol{\t}(\nu) }
  + H^i_{\nu,\ol{\t}(\nu)} \big|\cF_\nu\big]
  \le  \underset{i \in \cI}{\esssup}\, \wt{\cE}_i \big[
  Z(\ol{\t}(\nu))
  + H^i_{\nu,\ol{\t}(\nu)} \big|\cF_\nu\big]\le Z(\nu) , \q a.s.,
 \eeas
  which shows that
 \bea \label{eqn-k22}
 Z(\nu)  =\underset{i \in \cI}{\esssup}\, \wt{\cE}_i \big[ Y_{\ol{\t}(\nu) }
  + H^i_{\nu,\ol{\t}(\nu)} \big|\cF_\nu\big]
  = \underset{i \in \cI}{\esssup}\, \wt{\cE}_i \big[ Z(\ol{\t}(\nu))
  + H^i_{\nu,\ol{\t}(\nu)} \big|\cF_\nu\big], \q a.s.
 \eea
 Now we fix $\rho \in \cS_{\nu, \ol{\t}(\nu)} $. For any $i \in \cI$,
  Corollary \ref{cor_os2} and (\ref{eqn-xxx01})  
 show that
 \beas
  \wt{\cE}_i\big[ Z^i(\ol{\t}(\nu)) | \cF_\nu\big]
  =  \wt{\cE}_i\big[\wt{\cE}_i\big[ Z^i (\ol{\t}(\nu)) \big|\cF_\rho\big] \big| \cF_\nu \big]
  \le \wt{\cE}_i\big[ Z^i(\rho) | \cF_\nu\big], \q a.s.
 \eeas
 Then Proposition \ref{properties_3} (3) implies that
 \beas
 \wt{\cE}_i\big[ Z (\ol{\t}(\nu))+H^i_{\nu,\ol{\t}(\nu)} | \cF_\nu\big]
  = \wt{\cE}_i\big[ Z^i(\ol{\t}(\nu)) | \cF_\nu\big]-H^i_\nu
  \le \wt{\cE}_i\big[ Z^i(\rho) | \cF_\nu\big]-H^i_\nu
   =\wt{\cE}_i\big[ Z(\rho)+H^i_{\nu,\rho} | \cF_\nu\big], \q a.s.
 \eeas
Taking the essential supremum of both sides over $\cI$, we can
deduce from (\ref{D_2}) that
 \beas
 \underset{i \in \cI}{\esssup}\, \wt{\cE}_i\big[ Z (\ol{\t}(\nu))+H^i_{\nu,\ol{\t}(\nu)} | \cF_\nu\big]
  \le \underset{i \in \cI}{\esssup}\, \wt{\cE}_i\big[ Z(\rho)+H^i_{\nu,\rho} | \cF_\nu\big] \le Z(\nu), \q a.s.,
 \eeas
 which together with (\ref{eqn-k22}) proves (\ref{eqn-a013}).

\ms Finally, we will prove that $\ol{\t}(\nu) = \tau_1(\nu) $. For
any $i \in \cI$ and $k \in \hN$, \eqref{eqn-h42}, (\ref{D_2}),
Proposition \ref{properties_3} (3) as well as
   Corollary \ref{cor_os2} imply that
 \beas
 && \hspace{-2.5cm} \wt{\cE}_i \big[Z^0_{\t_{\d^{n_{\neg k}}}(\nu)}
 +H^i_{\nu,\t_{\d^{n_{\neg k}}}(\nu)} \big|\cF_\nu\big]
  =    \wt{\cE}_i \big[Z(\t_{\d^{n_{\neg k}}}(\nu))+H^i_{\nu,\t_{\d^{n_{\neg k}}}(\nu)} \big|\cF_\nu \big]\\
  &\ge&   \wt{\cE}_i \big[  \underset{j \in \cI}{\esssup}\,
   \wt{\cE}_j \big[ Z(\ol{\t}(\nu)) +   H^j_{\t_{\d^{n_{\neg k}}}(\nu), \ol{\t}(\nu)}
    \big|\cF_{\t_{\d^{n_{\neg k}}}(\nu)} \big] +H^i_{\nu,\t_{\d^{n_{\neg k}}}(\nu)} \big| \cF_\nu \big] \\
 &\ge&  \wt{\cE}_i \big[ \wt{\cE}_i \big[ Z(\ol{\t}(\nu)) + H^i_{\t_{\d^{n_{\neg k}}}(\nu), \ol{\t}(\nu)}
  \big|\cF_{\t_{\d^{n_{\neg k}}}(\nu)} \big] +H^i_{\nu,\t_{\d^{n_{\neg k}}}(\nu)} \big| \cF_\nu   \big] \\
 &=& \wt{\cE}_i \big[ \wt{\cE}_i \big[ Z(\ol{\t}(\nu)) +  H^i_{\nu, \ol{\t}(\nu)}
  \big|\cF_{\t_{\d^{n_{\neg k}}}(\nu)} \big]  \big| \cF_\nu   \big]
  =  \wt{\cE}_i \big[ Z(\ol{\t}(\nu)) + H^i_{\nu, \ol{\t}(\nu)}    | \cF_\nu  \big] , \q
  a.s.,
 \eeas
which together with (\ref{eqn-q150}) shows that
 \beas
 && \hspace{-3cm}  \underset{i \in \cI}{\esssup}\,
\Big|\wt{\cE}_i\big[\hb{$\frac{n_{\neg k}}{n_{\neg k}-1}$}
Y_{\t_{\d^{n_{\neg k}}}(\nu) }\neg +\neg H^i_{\t_{\d^{n_{\neg
k}}}(\nu)} \big|\cF_\nu\big]\neg -\neg \wt{\cE}_i\big[
Y^i_{\ol{\t}(\nu) }\big|\cF_\nu\big]\Big|+
  \wt{\cE}_i \big[ Y_{\ol{\t}(\nu) }\neg +\neg H^i_{\nu,\ol{\t}(\nu)}
  \big|\cF_\nu\big] \neg \\
  &&\ge \wt{\cE}_i \big[ Z(\ol{\t}(\nu))\neg + \neg H^i_{\nu, \ol{\t}(\nu)}    | \cF_\nu  \big]
 +\hb{$\frac{1}{n_k-1}(C_Y+2C_H)$}, \q  a.s.
 \eeas
 As $k \to \infty$, (\ref{eqn-x23}) implies that
 \bea \label{eqn-q151}
 \wt{\cE}_i \big[
 Y_{\ol{\t}(\nu) } + H^i_{\nu, \ol{\t}(\nu)} \big|\cF_\nu\big]
  \ge \wt{\cE}_i \big[ Z(\ol{\t}(\nu)) +
H^i_{\nu, \ol{\t}(\nu)}  | \cF_\nu  \big]    , \q a.s.
 \eea
The reverse inequality follows easily from (\ref{eqn-q104x}), thus
(\ref{eqn-q151}) is in fact an equality. Then the second part of
Proposition \ref{properties_3} (1) and \eqref{eqn-h42} imply that
 \beas
  Y_{\ol{\t}(\nu)} =  Z(\ol{\t}(\nu))= Z^0_{\ol{\t}(\nu)},\q a.s.,
 \eeas
which shows that $ \inf\big\{t \in  [\nu,T ] :\, Z^0_t=Y_t\big\} \le
\ol{\t}(\nu)$, a.s. For any $\d \in (0,1)$, since $\big\{t \in
 [\nu,T ] :\,  Z^0_t=Y_t\big\} \subset \big\{t \in
 [\nu,T ] :\,   Y_t \ge \d Z^0_t  +(1-\d)(C_Y+2C_H) \big\}$, one can deduce that
 \beas
 \ol{\t}(\nu) \ge \inf\big\{t \in  [\nu,T ]:\, Z^0_t=Y_t\big\}
 \ge \inf\big\{t \in  [\nu,T ] :\, Y_t \ge \d Z^0_t  +(1-\d)(C_Y+2C_H) \big\} \land T =
 \t_\d(\nu),  \q a.s.
 \eeas
Letting $\d \to 1$ yields that
 \beas
 \ol{\t}(\nu) \ge \inf\big\{t \in  [\nu,T ] :\, Z^0_t=Y_t\big\} \ge
  \underset{\d \to 1}{\lim} \t_\d(\nu) = \ol{\t}(\nu),  \q a.s.,
 \eeas
 which implies that $\ol{\t}(\nu)  =   \inf\big\{t \in  [\nu,T ] :\, Z^0_t=Y_t\big\}$,
 a.s.  \qed

\subsection{Proofs of Section \ref{nonco_game}}

\begin{deff}
 A family $\{ \xi_i  \}_{i \in \cI} \subset L^0(\cF_T)$ is said to be {\it directed downwards} if  for any $i,j \in
\cI$, there exists a $k \in \cI$ such that $\xi_k  \le \xi_i
\land  \xi_j      $, a.s.
\end{deff}

\ss \no {\bf Proof of Lemma \ref{essinf_lim}:} In light of
\cite[Proposition VI-\b1-1]{Neveu_1975}, it suffices to show that
the family $\{ R^i(\nu) \}_{i \in \cI}$ is directed downwards.
To see this, we define the event
 $A \dfnn \left\{ R^i(\nu) \ge R^j(\nu) \right\}$
 and the stopping times
  \beas
  \rho \dfnn  \t^j(\nu) \b1_A +  \t^i(\nu) \b1_{A^c} \in
 \cS_{\nu,T} \q \hb{and} \q \nu(A)  \dfnn \nu \b1_A  +T \b1_{A^c}  \in \cS_{\nu,T} .
 \eeas
 By Definition \ref{def_stable_class}, there exists a $k=k\big(i,j,\nu(A)\big) \in \cI$
 such that $\wt{\cE}_k = \cE^{\nu(\neg A)}_{i,j}$. Similar to \eqref{eqn-fxf01} it holds
 for any $\xi \in Dom(\sE)$ that
 \bea \label{eqn-new101}
 \wt{\cE}_k[\xi |\cF_\nu]  =\b1_A \wt{\cE}_j [\xi|\cF_{\nu }] +\b1_{A^c} \wt{\cE}_i [  \xi  |\cF_\nu ], \q a.s.
 \eea
 Moreover, (\ref{h_tau_A}) implies that
  \beas
      && \q \;\; H^k_{\nu,\rho}
 \neg =\neg   H^i_{ \nu(A) \land \nu,  \nu(A)  \land \rho }
 \neg+\neg   H^j_{ \nu(A)  \vee \nu , \nu(A)  \vee \rho }
 \neg=\neg  \b1_{A^c} H^i_{\nu, \t^i(\nu)}    \neg +\neg  \b1_A H^j_{\nu, \t^j(\nu) }, \q
 a.s. \\
 \hb{and that} &&  H^k_{\nu, \t^k(\nu)}
 \neg =\neg  H^i_{ \nu(A) \land \nu, \nu(A)  \land {\t^k(\nu)} }
 \neg+\neg  H^j_{ \nu(A)  \vee \nu, \nu(A)  \vee {\t^k(\nu)} }
 \neg=\neg  \b1_{A^c} H^i_{\nu, {\t^k(\nu)}}
 \neg +\neg  \b1_A H^j_{\nu, {\t^k(\nu)}}  , ~\; a.s.
  \eeas

 Using (\ref{eqn-new101}) twice and applying Proposition \ref{properties_3} (2)
 to $\wt{\cE}_i$ and $\wt{\cE}_j$, we can deduce from
 (\ref{eqn-k25}) that
 \beas
 R^k(\nu) &\ge& \wt{\cE}_k\big[ Y_\rho+ H^k_{\nu,\rho} \big|\cF_\nu\big]
  =   \b1_A \wt{\cE}_j\big[ Y_\rho + H^k_{\nu,\rho}\big|\cF_\nu\big]
  + \b1_{A^c} \wt{\cE}_i\big[ Y_\rho + H^k_{\nu,\rho}\big|\cF_\nu\big]\\
 &=&   \wt{\cE}_j\big[ \b1_A Y_{\t^j(\nu)} + \b1_AH^j_{\nu,\t^j(\nu)}\big|\cF_\nu\big]
  + \wt{\cE}_i\big[ \b1_{A^c} Y_{\t^i(\nu)} + \b1_{A^c} H^i_{\nu,\t^i(\nu)}\big|\cF_\nu\big] \\
  &=&  \b1_A \wt{\cE}_j\big[  Y_{\t^j(\nu)} +  H^j_{\nu,\t^j(\nu)}\big|\cF_\nu\big]
  + \b1_{A^c} \wt{\cE}_i\big[  Y_{\t^i(\nu)} +  H^i_{\nu,\t^i(\nu)}\big|\cF_\nu\big]
  = \b1_A   R^j(\nu)  + \b1_{A^c} R^i(\nu)  \\
 &\ge & \b1_A  \wt{\cE}_j\big[ Y_{\t^k(\nu)} +H^j_{\nu,\t^k(\nu)}
 \big|\cF_\nu\big] + \b1_{A^c} \wt{\cE}_i\big[ Y_{\t^k(\nu)}+H^i_{\nu,\t^k(\nu)} \big|\cF_\nu\big] \\
  &=&   \wt{\cE}_j\big[ \b1_A Y_{\t^k(\nu)} +\b1_A H^k_{\nu,\t^k(\nu)}
 \big|\cF_\nu\big] +  \wt{\cE}_i\big[ \b1_{A^c} Y_{\t^k(\nu)}+\b1_{A^c} H^k_{\nu,\t^k(\nu)} \big|\cF_\nu\big]\\
  &=& \b1_A  \wt{\cE}_j\big[  Y_{\t^k(\nu)} +  H^k_{\nu,\t^k(\nu)}
 \big|\cF_\nu\big] + \b1_{A^c}  \wt{\cE}_i\big[ Y_{\t^k(\nu)} +  H^k_{\nu,\t^k(\nu)} \big|\cF_\nu\big]\\
 &=& 
   \wt{\cE}_k \big[Y_{\t^k(\nu)}+  H^k_{\nu,\t^k(\nu)} \big|\cF_\nu\big] = R^k(\nu), \q
   a.s.,
  \eeas
which shows that $R^k(\nu)= \b1_A R^j(\nu) + \b1_{A^c} R^i(\nu)=
R^i(\nu) \land R^j(\nu)$, a.s. In light of the basic properties of
the essential infimum (e.g., \cite[Proposition
VI-\b1-1]{Neveu_1975}), we can find a sequence
$\left\{i_n\right\}_{n \in \hN}$ in $\cI$ such that (\ref{eqn-p05})
holds.  \qed

\ss \no {\bf Proof of Lemma \ref{ul_nu}:} As in the proof of Lemma
\ref{essinf_lim}, it suffices to show that the family $ \{\t^i(\nu)
 \}_{i \in \cI}$ is directed downwards.
To see this, we define the stopping time $\si\dfnn \t^i(\nu) \land
\t^j(\nu) \in \cS_{\nu, T}$, the event $A \dfnn \{R^{i,0}_\si \ge
R^{j,0}_\si \} \in \cF_\si$ as well as the stopping time
 $  \si(A) \dfnn \si \b1_A + T \b1_{A^c} \in \cS_{\si, T}$.
 By Definition \ref{def_stable_class}, there exists a $k=k\big(i,j,\si(A)\big) \in
\cI$ such that $\wt{\cE}_k=\cE^{\si(\neg A)}_{i,j} $. Fix $t \in
[0,T]$, similar to \eqref{eqn-fxf01}, it holds for any $\xi \in
Dom(\sE)$ that
 \bea \label{eqn-fxf02}
 \wt{\cE}_k[\xi |\cF_{\si\vee t}]   =\b1_A \wt{\cE}_j [\xi|\cF_{\si\vee t}]
 +\b1_{A^c} \wt{\cE}_i [ \xi  |\cF_{\si\vee t} ], \q a.s.
 \eea
 Moreover, we can deduce from (\ref{h_tau_A}) that for any $\rho \in
\cS_{\si\vee t,T}$
  \beas
      H^k_{\si\vee t, \rho}
 \neg =\neg  H^i_{\si(A) \land (\si\vee t), \si(A)  \land \rho }
 \neg+\neg  H^j_{\si(A)  \vee (\si\vee t), \si(A)  \vee \rho }
 \neg=\neg  \b1_{A^c} H^i_{\si\vee t, \rho }
  \neg +\neg  \b1_A H^j_{\si\vee t, \rho} , \q a.s.,
  \eeas
 which together with \eqref{eqn-fxf02} and Proposition \ref{properties_3} (2)
 imply that
 \beas
 \wt{\cE}_k\big[Y_\rho \neg+\dneg H^k_{\si\vee t,\rho} \big|\cF_{\si\vee t}\big]
 & =& 
 \b1_A\wt{\cE}_j \big[Y_\rho \neg+\dneg H^k_{\si\vee t,\rho} \big|\cF_{\si\vee t}\big]
 + \b1_{A^c}\wt{\cE}_i\big[Y_\rho \neg+\dneg H^k_{\si\vee t,\rho} \big|\cF_{\si\vee t}\big]\\
& =&  \wt{\cE}_j \big[ \b1_A Y_\rho \neg +\b1_A   H^j_{\si\vee
t,\rho} \big|\cF_{\si\vee t}\big] +\wt{\cE}_i\big[ \b1_{A^c} Y_\rho
\neg
+\b1_{A^c} \neg H^i_{\si\vee t,\rho} \big|\cF_{\si\vee t}\big] \\
& =&    \b1_A \wt{\cE}_j \big[ Y_\rho \neg + \neg H^j_{\si\vee
t,\rho} \big|\cF_{\si\vee t}\big] +\b1_{A^c} \wt{\cE}_i\big[ Y_\rho
\neg
+  \neg H^i_{\si\vee t,\rho} \big|\cF_{\si\vee t}\big], \q a.s. 
 \eeas
 Then applying Proposition \ref{properties_3} (3), 
  Lemma \ref{lem_ess} (2) as well as (\ref{eqn-p04}), we obtain
 \beas
 R^{k,0}_{\si\vee t}   &=& R^k(\si\vee t)
  = \underset{\rho \in \cS_{\si\vee t,T}}{\esssup}\, \wt{\cE}_k \big[Y_\rho
 + H^k_{\si\vee t,\rho} \big|\cF_{\si\vee t}\big] \nonumber \\
 &=&    \b1_A \underset{\rho \in \cS_{\si\vee t,T}}{\esssup}\,
 \wt{\cE}_j \big[ Y_\rho \neg+ \neg H^j_{\si\vee t,\rho} \big|\cF_{\si\vee t}\big]
 +\b1_{A^c} \underset{\rho \in \cS_{\si\vee t,T}}{\esssup}\,
  \wt{\cE}_i\big[  Y_\rho \neg+  \neg H^i_{\si\vee t,\rho} \big|\cF_{\si\vee t}\big] \nonumber\\
 &=&   \b1_A  R^j(\si\vee t)  +\b1_{A^c} R^i(\si\vee t)
  =    \b1_A  R^{j,0}_{\si\vee t}  +\b1_{A^c} R^{i,0}_{\si\vee t},\q
a.s.    
 \eeas
 Since $R^{i,0}$, $R^{j,0}$ and $R^{k,0}$ are all RCLL processes, it holds except on a null set $N$ that
 \beas
  R^{k,0}_{\si\vee t} =  \b1_A  R^{j,0}_{\si\vee t}  +\b1_{A^c} R^{i,0}_{\si\vee t},\q
 \fa t \in [0,T],
 \eeas
   which 
   further implies that
 \bea \label{eqn-p02}
 \t^k(\nu) &\dneg =& \dneg \inf\left\{t \in [\nu,T] :\, R^{k,0}_t=Y_t\right\}  \le
   \inf\left\{t \in [\si,T] :\,  R^{k,0}_t=Y_t \right\}  \nonumber \\
 &\dneg =& \dneg \b1_A \inf\left\{t \in [\si,T] :\,  R^{j,0}_t=Y_t  \right\}
 +\b1_{A^c} \inf\left\{t \in [\si,T] :\, R^{i,0}_t=Y_t \right\}, \q a.s.
 \eea
Since $R^{i,0}_{\t^i(\nu)}=Y_{\t^i(\nu)}$,
$R^{j,0}_{\t^j(\nu)}=Y_{\t^j(\nu)}$, a.s. and since $\si= \t^i(\nu)
\land \t^j(\nu)$, it holds a.s. that $ Y_\si$ is equal either to
$R^{i,0}_\si$ or to $ R^{j,0}_\si$. Then the definition of the set
$A$ shows that $R^{j,0}_\si =Y_\si$ a.s. on $A$ and that
$R^{i,0}_\si =Y_\si$ a.s. on $A^c$, both of which further implies
that
 \beas
 \b1_A  \inf\left\{t \in [\si,T] :\, R^{j,0}_t =Y_t\right\} =\si\b1_A   \q \hb{and}
 \q \b1_{A^c} \inf\left\{t \in [\si, T] :\, R^{i,0}_t =Y_t\right\}=\si\b1_{A^c} , \q
 a.s.
 \eeas
Hence, we see from (\ref{eqn-p02}) that $\t^k(\nu)  \le   \si=
\t^i(\nu) \land \t^j(\nu) $, a.s.  Thanks
to the basic properties of the essential infimum (e.g.,
\cite[Proposition VI-\b1-1]{Neveu_1975}), we can find a sequence
$\left\{i_n\right\}_{n \in \hN}$ in $\cI$ such that
 \beas
 \ul{\t}(\nu) = \underset{i \in \cI}{\essinf}\, \t^i(\nu) = \underset{n \to \infty}{\lim} \dneg \da \t^{i_n}(\nu), \q a.s.
 \eeas
The limit $\underset{n \to \infty}{\lim} \dneg \da \t^{i_n}(\nu)$ is
also a stopping time, thus we have $\ul{\t}(\nu) \in \cS_{\nu,T}$.
\qed

  \ss \no {\bf Proof of Theorem \ref{V_process}:} In light of Lemma
\ref{ul_nu}, there exists a sequence $\left\{j_n\right\}_{n \in
\hN}$ in $\cI$ such that
  \beas
  \ul{\t}(\nu)  = \underset{n \to \infty}{\lim} \dneg \da \t^{j_n}(\nu), \q
  a.s.
  \eeas
   Since the family of processes $ \{Y^i \}_{i \in
\cI}$ is ``$\sE$-uniformly-right-continuous", we can find a
subsequence of $ \{j_n \}_{n \in \hN}$ (we still
  denote it by $ \{j_n \}_{n \in \hN}$) such that
  \bea \label{eqn-x25}
  \underset{n \to \infty}{\lim}  \underset{i \in \cI}{\esssup}\,
  \big|\wt{\cE}_i\big[Y^i_{\t^{j_n}(\nu)}\big|\cF_{\ul{\t}(\nu)}\big]-Y^i_{\ul{\t}(\nu)} \big|
 =0, \q a.s.
  \eea
 Fix $i \in \cI$ and $n \in \hN$, we know from Definition
 \ref{def_stable_class} that there exists a $k_n =k(i,j_n,\ul{\t}(\nu) )\in \cI$ such that
 $\wt{\cE}_{k_n}=\cE^{\ul{\t}(\nu)}_{i, j_n}$.
  \if{0}
 \bea \label{eqn-p08}
  \wt{\cE}_{k_n}=\cE^{\ul{\t}(\nu)}_{i, j_n}
  \q \hb{and}\q  h^{k_n}_t  =  \b1_{\left\{t < \ul{\t}(\nu) \right\}} h^k_t
 +  \b1_{ \left\{ t \ge \ul{\t}(\nu) \right\} } h^{i_n}_t, \q a.s.
  \eea
  \fi
 For any $t \in [0,T]$, Lemma \ref{lem_trancate} implies that
 $ R^{k_n,0}_{\ul{\t}(\nu) \vee t} =R^{j_n,0}_{\ul{\t}(\nu) \vee t}$,
 a.s. Since $R^{k_n,0} $ and $R^{j_n,0} $ are both RCLL processes,
 it holds except on a null set $N$ that
 \beas
 R^{k_n,0}_{\ul{\t}(\nu) \vee t} = R^{j_n,0}_{\ul{\t}(\nu) \vee t}, \q  \fa t \in [0,T],
 \eeas
which together with the fact that $\ul{\t}(\nu) \le  \t^{k_n}(\nu)
\land \t^{j_n}(\nu)$, a.s. implies that
 \bea  \label{eqn-k27}
 \t^{k_n}(\nu)& =& \inf\left\{t \in \big[\nu,T\big] :\, R^{k_n,0}_t=Y_t\right\}
 = \inf\left\{t \in \big[\ul{\t}(\nu),T\big] :\,
 R^{k_n,0}_t=Y_t\right\} \nonumber
 \\&=&\inf\left\{t \in \big[\ul{\t}(\nu),T\big] :\, R^{j_n,0}_t=Y_t\right\}
 = \inf\left\{t \in \big[\nu,T\big] :\,
 R^{j_n,0}_t=Y_t\right\}=\t^{j_n}(\nu), \q a.s.
 \eea
 Then \eqref{h_tau_A}, (\ref{eqn-k27}) and (\ref{tau_ij2}) show that
  \bea \label{eqn-new103}
  &&  \hspace{-1.5cm} R^{k_n}(\nu)+H^i_\nu=R^{k_n}(\nu)+H^{k_n}_\nu
  =\wt{\cE}_{k_n}\big[Y_{\t^{k_n}(\nu)}+H^{k_n}_{\nu, \t^{k_n}(\nu)}\big|\cF_\nu\big]+ H^{k_n}_\nu \nonumber\\
    &\tneg =& \tneg\wt{\cE}_{k_n}\big[Y^{k_n}_{\t^{k_n}(\nu)}\big|\cF_\nu\big]
  =\cE^{\ul{\t}(\nu)}_{i, j_n}\big[Y^{k_n}_{\t^{j_n}(\nu)}\big|\cF_\nu\big]
 = \wt{\cE}_i \big[ \wt{\cE}_{j_n}\big[ Y^{k_n}_{\t^{j_n}(\nu)}\big|\cF_{\ul{\t}(\nu)}\big]
 \big|\cF_\nu\big] \nonumber \\
 &\tneg = &\tneg  \wt{\cE}_i \big[ \wt{\cE}_{j_n}\neg\big[Y_{\t^{j_n}(\nu)}
 \neg +\neg   H^{j_n}_{\ul{\t}(\nu),\t^{j_n}(\nu)}\neg +\neg    H^i_{\ul{\t}(\nu)} \big|\cF_{\ul{\t}(\nu)}\big]
 \big|\cF_\nu\big]
  \neg =  \wt{\cE}_i \big[ \wt{\cE}_{j_n} \neg \big[Y^{j_n}_{\t^{j_n}(\nu)} \big|\cF_{\ul{\t}(\nu)}\big]
   -H^{j_n}_{\ul{\t}(\nu)}\neg + \neg H^i_{\ul{\t}(\nu)}
 \big|\cF_\nu\big] \nonumber \\
  &\tneg \le & \tneg  \wt{\cE}_i \Big[ \big|\wt{\cE}_{j_n}\big[Y^{j_n}_{\t^{j_n}(\nu)} \big|\cF_{\ul{\t}(\nu)}\big]
  \neg -\neg Y^{j_n}_{\ul{\t}(\nu)} \big| \neg+\neg  Y^i_{\ul{\t}(\nu)}
  \Big|\cF_\nu\Big]   \le   \wt{\cE}_i \Big[    \underset{l \in \cI}{\esssup}\,
  \big|\wt{\cE}_l\big[Y^l_{\t^{j_n}(\nu)}\big|\cF_{\ul{\t}(\nu)}\big]\neg-\neg Y^l_{\ul{\t}(\nu)} \big|
  \neg+\neg  Y^i_{\ul{\t}(\nu)} \Big|\cF_\nu\Big], ~\, a.s.\q
 \eea
 For any $l \in \cI$,  Proposition \ref{properties_3} (3),
 \eqref{Y_LB}, \eqref{H_LB} and (\ref{eqn-q104}) imply that
 \beas
     \big|\wt{\cE}_l\big[Y^l_{\t^{j_n}(\nu)}\big|\cF_{\ul{\t}(\nu)}\big]-Y^l_{\ul{\t}(\nu)}  \big|
   &=&  \big|\wt{\cE}_l\big[Y_{\t^{j_n}(\nu)} +H^l_{\ul{\t}(\nu),\t^{j_n}(\nu) }-C_*\big|\cF_{\ul{\t}(\nu)}\big]
  - (Y_{\ul{\t}(\nu)}- C_Y   )+C_H  \big|\\
   &\le & \big| \wt{\cE}_l\big[Y_{\t^{j_n}(\nu)} +H^l_{\ul{\t}(\nu),\t^{j_n}(\nu)
   }-C_*\big|\cF_{\ul{\t}(\nu)}\big]\big|
   + \big|Y_{\ul{\t}(\nu)}- C_Y\big|+|C_H|  \\
   &= &  \wt{\cE}_l\big[Y_{\t^{j_n}(\nu)} +H^l_{\ul{\t}(\nu),\t^{j_n}(\nu) }-C_*\big|\cF_{\ul{\t}(\nu)}\big]
  + (Y_{\ul{\t}(\nu)}- C_Y )-C_H \\
 &=& \wt{\cE}_l\big[Y^l_{\t^{j_n}(\nu)} +H^l_{\ul{\t}(\nu),\t^{j_n}(\nu) } \big|\cF_{\ul{\t}(\nu)}\big]+
 Y_{\ul{\t}(\nu)} -2C_*   \le 2 R^l\big(\ul{\t}(\nu)\big) -2C_*   , \q
 a.s. \qq
 \eeas
Taking the essential supremum over $l \in \cI$, we can deduce from
\eqref{Y_LB2} and (\ref{eqn-q104}) that
  \beas
 &&  C_* \le  \underset{l \in \cI}{\esssup}\,
  \big|\wt{\cE}_l\big[Y^l_{\t^{j_n}(\nu)}\big|\cF_{\ul{\t}(\nu)}\big]-Y^l_{\ul{\t}(\nu)} \big|
  +  Y^i_{\ul{\t}(\nu)}  \le  3 R^l \big(\ul{\t}(\nu)\big) -2C_* +H^i_{\ul{\t}(\nu)}, \q a.s.,
  \eeas
  where $3 R^l \big(\ul{\t}(\nu)\big) -2C_* + H^i_{\ul{\t}(\nu)} \in Dom(\sE)$
  thanks to Proposition \ref{prop_Ri} (1), (S1') and
(D2). Applying the Dominated Convergence Theorem (Proposition
\ref{DCT2}) and Proposition \ref{properties_3} (3), we can deduce
from \eqref{eqn-new103} and (\ref{eqn-x25}) that
 \beas
    \ol{V}(\nu) &=& \underset{j \in \cI}{\essinf}\, R^j(\nu) \le \underset{n \to \infty}{\liminf} R^{k_n}(\nu)
    \le  \neg \underset{n \to \infty}{\lim} \wt{\cE}_i \Big[\underset{l \in \cI}{\esssup}\,
  \big|\wt{\cE}_l\big[Y^l_{\t^{j_n}(\nu)}\neg\big|\cF_{\ul{\t}(\nu)}\big]\neg
  - \neg Y^l_{\ul{\t}(\nu)}\neg \big| \neg+\neg  Y^i_{\ul{\t}(\nu)} \neg
  \Big|\cF_\nu\Big]\neg - \neg H^i_\nu   \\
     &    =& \wt{\cE}_i \big[ Y^i_{\ul{\t}(\nu)}
  \neg \big|\cF_\nu\big] \neg- \neg H^i_\nu
  = \wt{\cE}_i \big[   Y_{\ul{\t}(\nu)}\neg + \neg H^i_{\nu,\ul{\t}(\nu)}  \big|\cF_\nu\big]
   ,~\; a.s.
 \eeas
Taking the essential infimum of the right-hand-side over $i \in \cI$
yields that
  \beas
   \qq 
  \ol{V}(\nu) \le \underset{i \in \cI}{\essinf}\, \wt{\cE}_i  \big[ Y_{\ul{\t}(\nu)}
   +\neg H^i_{\nu,\ul{\t}(\nu)}  \big|\cF_\nu \big]
   \le  \underset{\rho \in \cS_{\nu,T}}{\esssup}\,
 \Big(\, \underset{i \in \cI}{\essinf}\,\wt{\cE}_i \big[Y_\rho+ \neg H^i_{\nu,\rho}  \big|\cF_\nu\big]\Big)
 =\ul{V}(\nu) \le \ol{V}(\nu), \q a.s.
  \eeas
 Hence, we have
  \beas
    \ul{V}(\nu) = \underset{i \in \cI}{\essinf}\, \wt{\cE}_i  \big[ Y_{\ul{\t}(\nu)}
   +\neg H^i_{\nu,\ul{\t}(\nu)}  \big|\cF_\nu \big] = \ol{V}(\nu)
   =\underset{i \in \cI}{\essinf}\, R^i(\nu) \ge Y_\nu, \q a.s.,
  \eeas
  where the last inequality is due to \eqref{eqn-q104}. \qed

\ss \no {\bf Proof of Proposition \ref{V_Y_meet}:} By Lemma
\ref{ul_nu}, there exists a sequence $\{i_n\}_{n \in \hN}$ in $\cI$
such that
  \beas
 \si \dfnn \ul{\t}(\nu)  = \underset{n \to \infty}{\lim} \dneg \da \t^{i_n}(\nu), \q
  a.s.
  \eeas
 For any $n \in \hN$, since $\si \le \t^{i_n}(\nu)$, a.s., we have
 \beas
  \t^{i_n}(\nu)=\inf\{t \in [\nu,T] :\, R^{i_n,0}_t = Y_t \}  =\inf\{t \in [\si,T] :\, R^{i_n,0}_t = Y_t
  \}=\t^{i_n}(\si), \q a.s.
 \eeas
Then (\ref{eqn-k51}) and (\ref{eqn-k25}) imply that
 \beas
V(\si) &=& \ol{V}(\si) \le R^{i_n}(\si )
 =\wt{\cE}_{i_n}[Y_{\t^{i_n}(\si)}+H^{i_n}_{\si, \t^{i_n}(\si)}  |\cF_\si]
 =\wt{\cE}_{i_n}[Y_{\t^{i_n}(\nu)}+H^{i_n}_{\si, \t^{i_n}(\nu)}  |\cF_\si ]\\
 &=&  \wt{\cE}_{i_n}[Y^{i_n}_{\t^{i_n}(\nu)}|\cF_\si ]-H^{i_n}_\si
 = \wt{\cE}_{i_n}[Y^{i_n}_{\t^{i_n}(\nu)}|\cF_\si ]-Y^{i_n}_\si + Y_\si
 \le \underset{i \in \cI}{\esssup}\,\Big|\wt{\cE}_i[Y^i_{\t^{i_n}(\nu)}|\cF_\si]
 -Y^i_\si \Big|+ Y_\si, \q a.s.
 \eeas
As $n\to \infty$, the ``$\sE$-uniform-right-continuity" of
$\{Y^i\}_{ i \in \cI}$ implies that $  V(\si) \le Y_\si $, a.s.,
while the reverse inequality is obvious from (\ref{eqn-k51}).
\qed

\ss \no {\bf Proof of Proposition \ref{R_sub}:}
 In light of Lemma \ref{essinf_lim} and
(\ref{eqn-k51}), there exists a sequence $\left\{j_n\right\}_{n \in
\hN}$ in $\cI$ such that
  \beas
 V(\nu)= \ol{V}(\nu) = \underset{n \to \infty}{\lim} \dneg \da R^{j_n}(\nu), \q a.s.
  \eeas
For any $n \in \hN$, Definition \ref{def_stable_class} assures a
$k_n=k(i,j_n, \nu) \in \cI$ such that $\wt{\cE}_{k_n}=\cE^{\nu}_{i,
j_n}$.
 \if{0}
  \bea \label{eqn-k137}
 \wt{\cE}_{k_n}=\cE^{\nu}_{i, j_n} \q \hb{and}\q  h^{k_n}_t  =  \b1_{\left\{t < \nu \right\}} h^i_t
 +  \b1_{ \left\{ t \ge \nu \right\} } h^{j_n}_t, \q a.s.
 \eea
 \fi
 Applying Proposition \ref{properties_3} (5) to $\wt{\cE}_i$, 
  we can deduce from (\ref{tau_ij2}) and (\ref{eqn-k26}) that
 \beas
 \wt{\cE}_{k_n}\big[ V(\rho) +\neg H^{j_n}_{\nu,\rho} \big|\cF_\nu\big]
 &\le& \wt{\cE}_{k_n}\big[ R^{j_n}(\rho) +H^{j_n}_{\nu,\rho} \big|\cF_\nu\big]
    =    \cE^{\nu}_{i,j_n}\big[R^{j_n}(\rho) +H^{j_n}_{\nu,\rho}\big|\cF_\nu\big]
   = \wt{\cE}_i\big[  \wt{\cE}_{j_n} \big[ R^{j_n}(\rho) +H^{j_n}_{\nu,\rho} \big|\cF_\nu \big] \big|\cF_\nu\big]\\
  &\tneg =&\tneg \wt{\cE}_{j_n} \big[ R^{j_n}(\rho) +H^{j_n}_{\nu,\rho} \big|\cF_\nu \big]
    \le R^{j_n}(\nu), \q a.s.
 \eeas
 Then Proposition \ref{properties_3} (3) and \eqref{h_tau_A} imply that
 \beas 
  \underset{k \in \cI}{\essinf}\,\wt{\cE}_k\big[V^k(\rho)\big|\cF_\nu\big]
  &\tneg \le& \tneg \wt{\cE}_{k_n}\big[V^{k_n}(\rho)\big|\cF_\nu\big]
   =\wt{\cE}_{k_n}\big[ V(\rho) +\neg H^{j_n}_{\nu,\rho} \big|\cF_\nu\big]+\neg H^i_\nu \le R^{j_n}(\nu)+H^i_\nu
    , \q a.s.
   \eeas
As $n \to \infty$, (\ref{eqn-m31}) follows:
 \beas
 \underset{k \in \cI}{\essinf}\,\wt{\cE}_k\big[V^k(\rho)\big|\cF_\nu\big]
 \le  \underset{n \to \infty}{\lim} \dneg \da R^{j_n}(\nu)+H^i_\nu
 = V(\nu)+H^i_\nu =V^i(\nu) , \q  a.s.
 \eeas

 \ms Now we assume that $\nu \le \rho \le  \ul{\t}(\nu)$, a.s. Applying Lemma
\ref{essinf_lim} and (\ref{eqn-k51}) once again, we can find another
sequence $\left\{j'_n\right\}_{n \in \hN}$ in $\cI$ such that
  \beas
 V(\rho)= \ol{V}(\rho) = \underset{n \to \infty}{\lim} \dneg \da R^{j'_n}(\rho), \q a.s.
  \eeas
 For any $n \in \hN$, Definition \ref{def_stable_class} assures a
 $k'_n=k(i, j'_n, \rho) \in \cI$ such that $\wt{\cE}_{k'_n}=\cE^{\rho}_{i, j'_n}$.
Since $\rho \le \ul{\t}(\nu) \le \t^{k'_n}(\nu) $, a.s., using
\eqref{eqn-k25} with $i=k'_n$ and applying Proposition
\ref{properties_3} (5) to $\wt{\cE}_{j'_n}$, we can deduce from
(\ref{h_tau_A}), (\ref{tau_ij2})
 as well as Lemma \ref{lem_trancate} that
 \bea
  V^i(\nu)&\dneg=& \dneg V(\nu)+ H^i_\nu = V(\nu)+ H^{k'_n}_\nu \le
  R^{k'_n}(\nu)+ H^{k'_n}_\nu     =  \wt{\cE}_{k'_n}\big[
R^{k'_n}(\rho)+ H^{k'_n}_\rho \big|\cF_\nu\big]=\cE^{\rho}_{i,
j'_n}\big[
R^{k'_n}(\rho)+ H^{k'_n}_\rho  \big|\cF_\nu\big] \nonumber  \\
&\dneg =& \dneg \wt{\cE}_i\big[ \wt{\cE}_{j'_n} \big[R^{k'_n}(\rho)+
H^{k'_n}_\rho \big|\cF_\rho \big] \big|\cF_\nu\big]
 =\wt{\cE}_i \big[ R^{k'_n}(\rho)+ H^{k'_n}_\rho \big|\cF_\nu \big]
  =\wt{\cE}_i \big[ R^{j'_n}(\rho)+H^i_\rho  \big|\cF_\nu \big], \q
  a.s. \qq  \label{eqn-gxg03}
 \eea
Then \eqref{Y_LB2} and \eqref{eqn-q104} imply that
  \beas
  C_*  \le Y^i_\rho =  Y_\rho +H^i_\rho \le  R^{j'_n}(\rho)+H^i_\rho \le R^{j'_1}(\rho)+H^i_\rho, \q a.s.,
 \eeas
  where $ R^{j'_1}(\rho)+H^i_\rho \in  Dom(\sE)$ thanks to Proposition \ref{prop_Ri} (1), (S1') and  (D2).
  As $n \to \infty$ in \eqref{eqn-gxg03},
  the Dominated Convergence Theorem (Proposition \ref{DCT2}) imply that
 \beas
   V^i(\nu) \le  \underset{n \to \infty}{\lim}  \wt{\cE}_i \big[  R^{j'_n}(\rho)+H^i_\rho \big|\cF_\nu
 \big]
   =  \wt{\cE}_i \big[  V(\rho)+H^i_\rho \big|\cF_\nu \big] = \wt{\cE}_i \big[  V^i(\rho) \big|\cF_\nu \big] ,\q
   a.s.,
 \eeas
 which proves (\ref{eqn-l20}).

 \ms It remains to show that $\left\{V^i\big( \ul{\t}(0) \land
t\big)\right\}_{t \in [0,T] }$ is an $\wt{\cE}_i$-submartingale: To
see this, we fix $0 \le s<t \le T$ and set $\nu \dfnn \ul{\t}(0)
\land s $, $\rho \dfnn \ul{\t}(0) \land t $. It is clear that $ \nu
\le \rho \le  \ul{\t}(0) \le \ul{\t}(\nu) $,
 a.s., hence (\ref{eqn-l20}), Corollary \ref{cor_os2} and Proposition \ref{properties_3} (5) show that
 \beas
V^i\big( \ul{\t}(0) \land s \big) &=& V^i (\nu) \le \wt{\cE}_i \big[
V^i(\rho) \big|\cF_\nu \big]= \wt{\cE}_i\big[V^i\big(\ul{\t}(0)
\land t\big) \big|\cF_{\ul{\t}(0) \land s}\big]
 =\wt{\cE}_i \big[\wt{\cE}_i\big[V^i\big(\ul{\t}(0) \land t \big) \big|\cF_{\ul{\t}(0) } \big]
 \big|\cF_s\big]\\
 &=& \wt{\cE}_i \big[ V^i\big(\ul{\t}(0) \land t\big)  \big|\cF_s \big], \q a.s.,
 \eeas
 which implies that $\left\{V^i\big( \ul{\t}(0) \land t\big)\right\}_{t \in [0,T] }$
 is an  $\wt{\cE}_i$-submartingale.  \qed

\ss \no {\bf Proof of Theorem \ref{V_RC}: Proof of (1).}

\noindent {\bf Step 1:} For any $\rho, \nu \in \cS_{0,T}$, we define
 \beas
  \P^\rho(\nu)  \dfnn \underset{i \in \cI}{\essinf}\, \wt{\cE}_i \big[Y_\rho
  +H^i_{\rho \land \nu,\rho} \big|\cF_{ \rho \land \nu}\big] + H^{i'}_{\rho \land \nu} \in \cF_{\rho \land
  \nu}.
 \eeas
 It follows from \eqref{Y_LB}, \eqref{H_LB}, and Proposition
\ref{properties_3} (5) that
  \bea
  C_Y+2C_H &=& \underset{i \in \cI}{\essinf}\, \wt{\cE}_i[C_Y+C_H |\cF_{\rho \land
  \nu}]+C_H \nonumber \\
  &\le& \P^\rho(\nu) \le \wt{\cE}_{i'} \big[Y_\rho
  +H^{i'}_{\rho \land \nu,\rho} \big|\cF_{ \rho \land \nu}\big] + H^{i'}_{\rho \land \nu}
  \le R^{i'}(\rho \land  \nu) + H^{i'}_{\rho \land \nu}, \q a.s., \label{eqn-new107}
  \eea
  where $R^{i'}(\rho \land  \nu) + H^{i'}_{\rho \land \nu} \in
  Dom(\sE)$ thanks to Proposition
\ref{prop_Ri} (1), (S1') and (D2). Then Lemma \ref{lem_dom_sharp}
implies that $\P^\rho(\nu) \in Dom(\sE)$.
 Applying Proposition
\ref{properties_3} (2)-(3) and Lemma \ref{lem_ess}, we can
alternatively rewrite $\P^\rho(\nu)$ as follows:
 \beas
      \P^\rho(\nu)-\neg H^{i'}_{\rho \land
  \nu}
   & \tneg=&  \tneg \underset{i \in
\cI}{\essinf}\,
 \wt{\cE}_i \big[\b1_{\{\rho \le \nu \}} Y_{\rho \land \nu}\neg + \neg \b1_{\{\rho > \nu \}}\neg\left( Y_\rho
 \neg +\neg H^i_{ \nu,\rho}\right) \big|\cF_{ \rho \land \nu}\big] \\
  & \tneg =& \tneg \underset{i \in \cI}{\essinf}\,\neg
 \left(\b1_{\{\rho \le \nu \}} Y_{\rho \land \nu} \neg + \neg \b1_{\{\rho > \nu \}} \wt{\cE}_i \big[  Y_\rho
 \neg +\neg H^i_{ \nu,\rho} \big|\cF_{  \nu}\big] \right) \neg
  =  \b1_{\{\rho \le \nu \}} Y_{\rho  } \neg + \neg \b1_{\{\rho > \nu \}}
  \underset{i \in \cI}{\essinf}\, \wt{\cE}_i \big[  Y_\rho
  \neg+\neg H^i_{ \nu,\rho} \big|\cF_{  \nu}\big] , \q
  a.s.
 \eeas
Let $\si \in \cS_{0,T}$. Lemma \ref{lem_ess}
(2) and Proposition \ref{properties_3} (2) once again imply that
  \bea \label{eqn-k907}
   \b1_{\{ \nu =\si\}}\P^\rho(\nu) &=& \b1_{\{\rho \le \nu=\si \}} Y_{\rho  } \neg + \b1_{\{\rho > \nu =\si \}}
  \underset{i \in \cI}{\essinf}\,  \wt{\cE}_i \big[  Y_\rho
  + H^i_{ \nu,\rho} \big|\cF_{  \nu}\big]+ \b1_{\{ \nu =\si\}} H^{i'}_{\rho \land
  \nu} \nonumber \\
  &=& \b1_{\{\rho \le \nu=\si \}} Y_{\rho  } \neg + \b1_{\{\rho > \nu \}}
  \underset{i \in \cI}{\essinf}\,  \wt{\cE}_i \big[ \b1_{\{ \nu =\si\}} \big( Y_\rho
  + H^i_{ \si,\rho} \big) \big|\cF_{  \nu}\big]+ \b1_{\{ \nu =\si\}} H^{i'}_{\rho \land
  \si} \nonumber \\
  &=& \b1_{\{\rho \le \nu=\si \}} Y_{\rho  } \neg + \b1_{\{\rho > \nu \}}
  \underset{i \in \cI}{\essinf}\, \b1_{\{ \nu =\si\}} \wt{\cE}_i \big[  Y_\rho
  + H^i_{ \si,\rho}   \big|\cF_{  \si}\big]+ \b1_{\{ \nu =\si\}} H^{i'}_{\rho \land
  \si} \nonumber \\
  &=& \b1_{\{\rho \le \si = \nu  \}} Y_{\rho  } \neg + \b1_{\{\rho > \si = \nu \}}
  \underset{i \in \cI}{\essinf}\,  \wt{\cE}_i \big[  Y_\rho
  + H^i_{ \si,\rho}   \big|\cF_{  \si}\big]+ \b1_{\{ \nu =\si\}} H^{i'}_{\rho \land
  \si} = \b1_{\{ \nu =\si\}}\P^\rho(\si), \q a.s. \qq
  \eea

 \ms \no {\bf Step 2:}
 Fix $\rho \in \cS_{0,T}$. For any $\nu \in \cS_{0,T}$ and $\si  \in \cS_{\nu,T}$,
 letting 
 $(\nu, \cI',\cU) = (\rho \land \si, \cI, \left\{\rho\right\})$ and $X(\rho)=Y_\rho$
 in Lemma \ref{lem_02}, we can find a
sequence $\left\{ j_n  \right\}_{n \in \hN}$ in $ \cI $ such that
 \beas
 \underset{i \in \cI}{\essinf}\, \wt{\cE}_i \big[Y_\rho +H^i_{\rho \land \si,\rho} \big|\cF_{\rho \land \si}\big]
  = \underset{n \to \infty}{\lim} \dneg \da
   \wt{\cE}_{j_n}\big[Y_\rho+H^{j_n}_{\rho \land \si,\rho}\big|\cF_{\rho \land \si}\big],\q  a.s.
 \eeas
 Definition \ref{def_stable_class} assures the existence of a $k_n=k(i',j_n, \rho \land
\si ) \in \cI$  such that $ \wt{\cE}_{k_n}=\cE^{ \rho \land
\si}_{i'\neg,j_n} $. Applying Proposition \ref{properties_3} (3) to
$\wt{\cE}_{k_n}$, we can deduce from \eqref{h_tau_A} and
(\ref{tau_ij2}) that
 \beas
    \P^\rho(\nu) &\le& \wt{\cE}_{k_n}\big[Y_\rho +H^{k_n}_{\rho \land \nu,\rho}\big|\cF_{\rho \land
    \nu}\big]+H^{i'}_{\rho \land \nu}=\wt{\cE}_{k_n}\big[Y_\rho +H^{k_n}_{\rho \land \nu,\rho}\big|\cF_{\rho \land
    \nu}\big]+H^{k_n}_{\rho \land \nu}=\wt{\cE}_{k_n}\big[Y^{k_n}_\rho\big|\cF_{\rho \land
    \nu}\big] \\
  &=& \cE^{  \rho \land \si}_{i'\neg,j_n} \big[Y^{k_n}_\rho \big|\cF_{  \rho \land
  \nu}\big]
   =\wt{\cE}_{i'}\big[ \wt{\cE}_{j_n} \big[ Y^{k_n}_\rho \big|\cF_{  \rho \land \si} \big]\big|\cF_{ \rho \land
   \nu}\big] , \q a.s.
  \eeas
 For any $n \in \hN$, Proposition \ref{properties_3} (3) $\&$ (5), \eqref{Y_LB2} as well as \eqref{h_tau_A} imply that
 \beas
 \q  C_* &=&  \wt{\cE}_{j_n} \big[C_* \big|\cF_{  \rho \land \si} \big]
 \le \wt{\cE}_{j_n} \big[Y^{k_n}_\rho \big|\cF_{  \rho \land \si} \big]
    = \wt{\cE}_{j_n} \big[Y_\rho +H^{j_n}_{\rho \land \si,\rho}  \big|\cF_{  \rho \land \si} \big]
   +H^{i'}_{\rho \land \si}\\
   &\le& \wt{\cE}_{j_1}\big[Y_\rho+H^{j_1}_{\rho \land \si,\rho}\big|\cF_{\rho \land \si}\big] +H^{i'}_{\rho \land
   \si}  \le  R^{j_1} ( \rho \land \si  )+H^{i'}_{\rho \land
   \si} ,~\;  a.s.,
 \eeas
  where $R^{j_1} ( \rho \land \si  ) +H^{i'}_{\rho \land
   \si} \in Dom(\sE)$ thanks to Proposition \ref{prop_Ri} (1), (S1') and (D2).
Then the Dominated Convergence Theorem (Proposition \ref{DCT2}),
Corollary \ref{cor_os2} and Proposition \ref{properties_3} (5) show
that
 \bea
    \P^\rho(\nu)  &\le& \underset{n \to \infty}{\lim} \dneg \da
    \wt{\cE}_{i'}\big[ \wt{\cE}_{j_n} \big[ Y^{k_n}_\rho \big|\cF_{  \rho \land \si} \big]\big|\cF_{ \rho \land
   \nu}\big]= \underset{n \to \infty}{\lim} \dneg \da \wt{\cE}_{i'}\Big[ \wt{\cE}_{j_n} \big[Y_\rho
 +H^{j_n}_{\rho \land \si,\rho}  \big|\cF_{  \rho \land \si} \big]
   +H^{i'}_{\rho \land \si} \Big|\cF_{ \rho \land   \nu}\Big] \nonumber \\
  &=&  \wt{\cE}_{i'}\big[ \,\underset{n \to \infty}{\lim} \dneg \da \wt{\cE}_{j_n} \big[Y_\rho
+H^{j_n}_{\rho \land \si,\rho}  \big|\cF_{  \rho \land \si}
\big]+H^{i'}_{\rho \land \si}  \big|\cF_{ \rho \land \nu}\big]
       =\wt{\cE}_{i'}\big[ \underset{i \in \cI}{\essinf}\,
   \wt{\cE}_i \big[Y_\rho +H^i_{\rho \land \si,\rho} \big|\cF_{\rho \land \si}\big]
   +H^{i'}_{\rho \land \si}  \big|\cF_{ \rho \land \nu}\big] \nonumber \\
   & =& \wt{\cE}_{i'}\big[ \P^\rho(\si)   \big|\cF_{ \rho \land
\nu}\big]  =\wt{\cE}_{i'}\big[\wt{\cE}_{i'}\big[  \P^\rho(\si)
\big|\cF_\rho \big] \big|\cF_\nu \big]=\wt{\cE}_{i'}\big[
\P^\rho(\si) \big|\cF_\nu
 \big] , \q a.s., \label{eqn-m07}
 \eea
 which implies that $\{\P^\rho(t)\}_{t \in [0,T]}$ is an
 $\wt{\cE}_{i'}$-submartingale.  Hence, $\{-\P^\rho(t)\}_{t \in [0,T]} $ is an
 $\cE'$-supermartingale by assumption (\ref{ass_odd}).
 Since $\cE' $ satisfies (H0), (H1), \eqref{def_cE_concave} and since $Dom(\cE') \in \wt{\sD}_T$
 \big(which results from $Dom(\sE) \in \wt{\sD}_T$ and (\ref{ass_odd})\big), we know from
  Theorem \ref{upcrossing} that $ \P^{\rho,+}_t \dfnn
 \underset{n \to \infty}{\liminf} \P^\rho\big(q^+_n(t)\big) $, $t \in [0,T]$ is an RCLL
 process and that
  \bea \label{eqn-k909}
 P\Big(\P^{\rho,+}_t
 = \underset{n \to \infty}{\lim} \P^\rho\big(q^+_n(t)\big)   \hb{ for any }
 t \in [0,T] \Big)=1.
  \eea

  \ss \no {\bf Step 3:}  For any $\nu \in \cS_{0,T}$ and $n
\in \hN$, $ q^+_n(\nu)$ takes values
 in a finite set $\cD^n_T \dfnn \big([0,T)\cap \{k2^{-n}\}_{k \in \hZ}
\big) \cup \{T\}  $. Given an $\a \in \cD^n_T$, it holds for any $m
\ge n$ that $q^+_m (\a)=\a $ since $ \cD^n_T \subset \cD^m_T$. It
follows from \eqref{eqn-k909} that
 \beas
 \P^{\rho,+}_\a
 = \underset{m \to \infty}{\lim} \P^\rho\big(q^+_m(\a)\big)=\P^\rho ( \a
 ) , \q a.s.
 \eeas
 Then one can deduce from \eqref{eqn-k907} that
  \beas 
 \qq \P^{\rho,+}_{q^+_n(\nu)} =
   \sum_{\a \in \cD^n_T} \b1_{\{q^+_n(\nu) = \a
  \}}\P^{\rho,+}_{\a} = \sum_{\a \in \cD^n_T} \b1_{\{q^+_n(\nu) = \a
  \}}\P^\rho(\a) = \sum_{\a \in \cD^n_T} \b1_{\{q^+_n(\nu) = \a  \}}\P^\rho\big(q^+_n(\nu)\big)
  = \P^\rho\big(q^+_n(\nu)\big), \q a.s.
  \eeas
  Thus the right-continuity of the process $\P^{\rho,+}$ implies that
   \bea \label{eqn-k901}
   \P^{\rho,+}_{\nu}= \lmt{n \to \infty} \P^{\rho,+}_{q^+_n(\nu) }= \lmt{n \to \infty}
   \P^\rho\big(q^+_n(\nu)\big), \q a.s.
   \eea
 We have assumed that $\underset{ t  \in \cD_T} {\esssup}\, \wt{\cE}_j[Y^j_\rho |\cF_t] \in  Dom(\sE)$ for some
$j =j(\rho) \in \cI$. It holds a.s. that
  \beas
   \wt{\cE}_j \big[Y^j_\rho  \big|\cF_t\big]  \le
   \underset{s  \in \cD_T}{\esssup}\, \wt{\cE}_j \big[Y^j_\rho \big|\cF_s\big],
   \q \fa   t \in  \cD_T.
 \eeas
  Since $\wt{\cE}_j\big[Y^j_\rho \big|\cF_\cd\big]$ is an RCLL
process, it holds except on a null set $N$ that
 \beas
  \q \wt{\cE}_j \big[ Y^j_\rho  \big|\cF_t\big]  \le
   \underset{s  \in \cD_T}{\esssup}\, \wt{\cE}_j \big[Y^j_\rho \big|\cF_s\big],
    ~ \fa t \in [0,T],
 ~~  \hb{ thus} ~~\;  \wt{\cE}_j \big[Y^j_\rho  \big|\cF_{\rho \land q^+_n(\nu)}\big]  \le
   \underset{s  \in \cD_T}{\esssup}\, \wt{\cE}_j \big[Y^j_\rho
   \big|\cF_s\big],  ~ \fa n \in \hN.
 \eeas
 Then one can deduce from \eqref{eqn-new107}, \eqref{H_LB} and Proposition \ref{properties_3} (3) that
 \beas
C_Y+2C_H & \dneg \le & \dneg \P^\rho\big(q^+_n(\nu)\big)
 \le  \wt{\cE}_j \big[Y_\rho + H^j_{\rho \land q^+_n(\nu),\rho} \big|\cF_{ \rho \land q^+_n(\nu)}\big] + \neg
 \z^{i'} \neg =  \wt{\cE}_j \big[Y^j_\rho - H^j_{\rho \land q^+_n(\nu)} \big|\cF_{ \rho \land q^+_n(\nu)}\big] + \neg
   \z^{i'} \\
   & \dneg\le& \dneg \wt{\cE}_j \big[Y^j_\rho- C_H  \big|\cF_{ \rho \land q^+_n(\nu)}\big] + \neg
  \z^{i'} = \wt{\cE}_j \big[Y^j_\rho  \big|\cF_{ \rho \land q^+_n(\nu)}\big]
  - C_H + \neg  \z^{i'} \neg \le \underset{s  \in \cD_T}{\esssup}\,
  \wt{\cE}_j \big[Y^j_\rho \big|\cF_s\big]- C_H + \neg  \z^{i'} , \q
  a.s.,
 \eeas
 where the right hand side belongs to $Dom(\sE)$ thanks to (D2) and the assumption that $\z^{i'} \in Dom(\sE)$.
  Hence the Dominated Convergence Theorem (Proposition \ref{DCT2}), \eqref{eqn-k901}, (\ref{eqn-m07})
  as well as Proposition \ref{properties_3} (5) imply that $\P^{\rho,+}_\nu = \underset{n \to
\infty}{\lim} \P^\rho\big( q^+_n(\nu) \big)    \in Dom(\sE)$ and
that
 \bea \label{eqn-k902}
    \P^\rho(\nu) \le  \underset{n \to \infty }{\lim}
    \wt{\cE}_{i'}\big[
    \P^\rho(q^+_n(\nu)) \big|\cF_\nu\big]
   =\wt{\cE}_{i'}\big[   \P^{\rho,+}_\nu \big|\cF_\nu \big]
   =    \P^{\rho,+}_\nu , \q  a.s.,
 \eea
 where in the last equality we used the fact that $ \P^{\rho,+}_\nu = \underset{n
\to
\infty}{\lim} \P^\rho\big( q^+_n(\nu) \big) 
 \in \cF_\nu $, thanks to the right-continuity of the filtration $\bF$.

\ms \no {\bf Step 4:} Given $\nu \in \cS_{0, T}$, we set
   \beas
      \g \dfnn \ul{\t}(0) \land \nu, \q \g_n \dfnn \ul{\t}(0) \land
   q^+_n(\nu), \q \fa n \in \hN
   \eeas
  and let $\rho \in \cS_{\g,T}$. Since $\lmtu{n \to \infty} \b1_{\{\ul{\t}(0) > q^+_n (\nu)\}}
   = \b1_{\{\ul{\t}(0) >  \nu\}} $ and since
    \beas
    \{\ul{\t}(0)   > \nu  \} \subset
 \left\{ q^+_n (\nu) = q^+_n \big( \ul{\t}(0) \land \nu \big)
 \right\},\q \{\ul{\t}(0) > q^+_n (\nu)\} \subset \left\{ q^+_n(\nu) = \ul{\t}(0) \land
   q^+_n(\nu) \right\} , \q \fa n \in \hN,
 \eeas
  one can deduce
 from \eqref{eqn-k902}, \eqref{eqn-k901} and \eqref{eqn-k907}  that
 \bea \label{eqn-k920}
    \b1_{\{\ul{\t}(0)   > \nu\}} \P^\rho(\g) &\le & \b1_{\{\ul{\t}(0)   > \nu\}} \P^{\rho,+}_\g
    =\b1_{\{\ul{\t}(0)   > \nu\}}   \lmt{n \to \infty}  \P^\rho \big( q^+_n (\g) \big)
    = \lmt{n \to \infty} \b1_{\{\ul{\t}(0)   > \nu\}}
    \P^\rho \left(q^+_n \big( \ul{\t}(0) \land \nu \big) \right) \nonumber\\
    &=& \lmt{n \to \infty} \b1_{\{\ul{\t}(0)   > \nu\}}  \P^\rho \big(
    q^+_n(\nu)\big) =  \lmt{n \to \infty}\b1_{\{\ul{\t}(0) > q^+_n (\nu)\}}
   \P^\rho \big(  q^+_n(\nu)\big) \nonumber\\
   &=&  \lmt{n \to \infty}\b1_{\{\ul{\t}(0) > q^+_n (\nu)\}} \P^\rho \big( \ul{\t}(0) \land
   q^+_n(\nu) \big) = \b1_{\{\ul{\t}(0)   > \nu\}} \lmt{n \to \infty} \P^\rho \big(
     \g_n\big)
    , \q  a.s.
 \eea
For any $n \in \hN$, we see from (\ref{eqn-k51}) that
 \bea \label{eqn-k911}
 V(\g_n)&\tneg =& \tneg \ul{V}(\g_n) = \underset{\si \in \cS_{\g_n, T}}{\esssup}\,
 \Big(\, \underset{i \in \cI}{\essinf}\,\wt{\cE}_i
 \big[Y_\si\neg +\neg H^i_{\g_n,\si} \big|\cF_{\g_n}\big]\Big) \ge \underset{i \in
 \cI}{\essinf}\,\wt{\cE}_i
 \big[Y_{\rho \vee \g_n}\neg +\neg H^i_{\g_n, \rho \vee
 \g_n}\big|\cF_{\g_n}\big], \q a.s.
 \eea
  Since $ \{ \ul{\t}(0) \le \nu\} \subset \{ \g_n=\g=\ul{\t}(0)\} $,
 Proposition \ref{properties_3} (2) and (3) imply that for any $i \in \cI$
  \beas
  \q \b1_{\left\{  \ul{\t}(0) \le \nu \right\}}\wt{\cE}_i \big[  Y_{\rho }
  \neg+ \neg H^i_{\rho \land \g_n,\,\rho} \big|\cF_{\g_n}\big]
     =        \wt{\cE}_i \left[ \b1_{\{ \ul{\t}(0) \le \nu \}} \big( Y_{\rho }
  \neg+ \neg H^i_{\rho \land \g,\,\rho}\big) \big|\cF_{\g_n}\right]
   =    \b1_{\left\{ \ul{\t}(0) \le \nu \right\}}\wt{\cE}_i \left[  Y_{\rho }
  \neg+ \neg H^i_{\rho \land \g,\,\rho}  \big|\cF_{ \g }\right]    , \q
  a.s.,
 \eeas
and that
  \beas
  \qq  && \hspace{-2cm} \wt{\cE}_i
 \big[Y_{\rho \vee \g_n}\neg +\neg H^i_{\g_n, \, \rho \vee \g_n}\big|\cF_{\g_n}\big]
 =   \wt{\cE}_i  \big[\b1_{\left\{\rho \le \g_n \right\}} Y_{\g_n}
 \dneg +\neg  \b1_{\left\{\rho> \g_n \right\}} \big( Y_{\rho }
 \neg  + \neg H^i_{\rho \land \g_n,\, \rho}\big)  \big|\cF_{\g_n}\big] \\
 &=& \dneg \b1_{\left\{\rho \le \g_n \right\}}
  Y_{\g_n}\dneg + \neg \b1_{\left\{\rho> \g_n \right\}}\wt{\cE}_i \big[  Y_{\rho }
  \neg+ \neg H^i_{\rho \land \g_n,\, \rho} \big|\cF_{\g_n}\big]\\
 &=& \dneg \b1_{\left\{\rho \le \g_n \right\}}
  Y_{\g_n}\dneg + \neg \b1_{\left\{\rho> \g_n, \ul{\t}(0) > \nu \right\}}\wt{\cE}_i \big[  Y_{\rho }
  \neg+ \neg H^i_{\rho \land \g_n,\, \rho} \big|\cF_{\rho \land \g_n}\big]
  + \neg \b1_{\left\{\rho> \g_n, \ul{\t}(0) \le \nu \right\}}\wt{\cE}_i \big[  Y_{\rho }
  \neg+ \neg H^i_{\rho \land \g,\, \rho} \big|\cF_{\rho \land \g}\big], \q
  a.s.
  \eeas
  Then it follows from \eqref{eqn-k911} and Lemma \ref{lem_ess} that
 \beas 
 V(\g_n)&\tneg \ge & \tneg       \b1_{\left\{\rho \le \g_n \right\}}
  Y_{\g_n}\dneg + \neg \b1_{\left\{\rho> \g_n, \ul{\t}(0) > \nu \right\}}\underset{i \in
 \cI}{\essinf}\,\wt{\cE}_i \big[  Y_{\rho }
  \neg+ \neg H^i_{\rho \land \g_n,\, \rho} \big|\cF_{\rho \land \g_n}\big]
  + \neg \b1_{\left\{\rho> \g_n, \ul{\t}(0) \le \nu \right\}}\underset{i \in
 \cI}{\essinf}\,\wt{\cE}_i \big[  Y_{\rho }
  \neg+ \neg H^i_{\rho \land \g,\, \rho} \big|\cF_{\rho \land \g}\big]   \nonumber \\
 &\tneg =& \tneg   \b1_{\left\{\rho \le \g_n \right\}}  Y_{\g_n}
  \neg + \neg \b1_{\left\{\rho> \g_n, \ul{\t}(0) > \nu \right\}}
  \left( \P^\rho(\g_n)\neg -\neg H^{i'}_{\rho \land \g_n} \right)
  + \neg \b1_{\left\{\rho> \g_n, \ul{\t}(0) \le \nu \right\}}
  \left( \P^\rho(\g)\neg -\neg H^{i'}_{\rho \land \g} \right), \q
  a.s.\qq \q
 \eeas
 \if{0}
 Then \eqref{Y_LB}, \eqref{H_LB} and Proposition \ref{properties_3} (5) show that
 \beas
     \P^\rho(\g_n)=\underset{i \in \cI}{\essinf}\,
   \wt{\cE}_i \big[Y_\rho +H^i_{\rho \land \g_n,\rho} \big|\cF_{ \rho \land \g_n}\big] + H^{i'}_{\rho \land \g_n}
 \ge \underset{i \in \cI}{\essinf}\,
   \wt{\cE}_i \big[C_Y+ C_H \big|\cF_{ \rho \land \g_n}\big]
  +C_H = C_Y+2C_H ,\q a.s.
 \eeas
 \fi
As $n \to \infty$,  the right-continuity of processes $Y$ and
$H^{i'}$, \eqref{eqn-k920}, Lemma \ref{lem_ess} as well as
  Proposition \ref{properties_3} (2) \& (3) show that
 \beas 
   \underset{n \to \infty}{\liminf}V(\g_n)
   &\ge& \b1_{\left\{\rho=\g\right\}}  Y_{\g}
  \neg + \neg \b1_{\left\{\rho> \g, \ul{\t}(0) > \nu \right\}}
  \left( \lmt{n \to \infty} \P^\rho(\g_n)\neg -\neg H^{i'}_{\rho \land \g} \right)
  + \neg \b1_{\left\{\rho> \g, \ul{\t}(0) \le \nu \right\}}
  \left( \P^\rho(\g)\neg -\neg H^{i'}_{\rho \land \g} \right)\\
    &\ge& \b1_{\left\{\rho=\g\right\}}  Y_\g
  +\b1_{\left\{\rho> \g\right\}} \Big(  \P^\rho(\g) -\neg H^{i'}_{\rho \land  \g}
  \Big)
    =   \b1_{\left\{\rho = \g   \right\}}Y_\g + \b1_{\left\{\rho> \g \right\}}  \underset{i \in \cI}{\essinf}\,
  \wt{\cE}_i \big[   Y_\rho + \neg H^i_{\rho \land
  \g,\rho}   \big|\cF_{\rho \land
  \g}\big] \\
  &=&  \underset{i \in \cI}{\essinf}\,\Big(\b1_{\left\{\rho = \g   \right\}}Y_\g+\b1_{\left\{\rho> \g
  \right\}}\wt{\cE}_i\big[   Y_\rho + H^i_{\g,\rho}  \big|\cF_\g\big]  \Big)= \underset{i \in \cI}{\essinf}\,\wt{\cE}_i
  \big[\b1_{\left\{\rho = \g   \right\}}Y_\g +  \b1_{\left\{\rho> \g
  \right\}}\big( Y_\rho + H^i_{\g,\rho} \big) \big|\cF_\g\big] \\
   &=& \underset{i \in \cI}{\essinf}\,\wt{\cE}_i\big[Y_\rho+H^i_{\g,\rho}
  \big|\cF_\g\big], \q a.s.
  \eeas
Taking the essential supremum of the right-hand-side over $\rho \in
\cS_{\g, T}$, we obtain
  \bea  \label{eqn-m05}
 \underset{n \to \infty}{\liminf}V(\g_n) \ge \underset{\rho \in \cS_{\g, T}}{\esssup}\, \Big(
 \underset{i \in \cI}{\essinf}\,\wt{\cE}_i\big[Y_\rho+H^i_{\g,\rho} \big|\cF_\g\big]
 \Big) = \ul{V}(\g) = V(\g) , \q a.s.
  \eea

  On the other hand, for any $i \in \cI$ and $n \in \hN$ we have that
 $ V(\g_n)= \ol{V}(\g_n) = \underset{l \in \cI}{\essinf}\, R^l (\g_n) \le R^i (\g_n)$,
 a.s. Then (\ref{eqn-p04}) and the right continuity of the process $R^{i, 0}$
imply that
  \beas
      \underset{n \to \infty}{\limsup} V(\g_n) \le \underset{n \to \infty}{\lim}R^i (\g_n)
    = \underset{n \to \infty}{\lim}R^{i,0}_{\g_n}=R^{i,0}_\g=R^i(\g), \q a.s.
  \eeas
    Taking the essential infimum of $R^i(\g)$ over $i \in \cI$ yields that
    \beas
     \underset{n \to \infty}{\limsup} V(\g_n) \le  \underset{i \in \cI}{\essinf}\,
     R^i(\g) = \ol{V}(\g) =V(\g), \q a.s.
    \eeas
     This inequality together with (\ref{eqn-m05}) shows that
  $ \underset{n \to \infty}{\lim} V  ( \g_n )  = V (\g)$, a.s.,
  which further implies that for any $\nu \in \cS_{0, T}$ and $i \in \cI$
  \bea \label{eqn-k905}
       \lmt{n \to \infty} V^i \neg \left(\ul{\t}(0) \land q^+_n(\nu)\right)
    &=&   \lmt{n \to \infty} \dneg \left(V \big(\ul{\t}(0) \land
  q^+_n(\nu)\big) + H^i_{\ul{\t}(0) \land
  q^+_n(\nu)}\right)    \nonumber  \\
   & = &  V \big(\ul{\t}(0) \land  \nu\big) + H^i_{\ul{\t}(0) \land  \nu}
  =  V^i  \big( \ul{\t}(0) \land \nu \big) , \q  a.s.
    \eea

 \ss \no {\bf Step 5:} Proposition \ref{R_sub} shows that
 the stopped process $\big\{ V^{i'} \big(  \ul{\t}(0) \land t\big)\big\}_{t \in  [0,T ]} $
 is an $\wt{\cE}_{i'}$-submartingale, thus
 $\big\{ - V^{i'} \big(  \ul{\t}(0) \land t\big)\big\}_{t \in  [0,T ]}$ is an $\cE'$-supermartingale
 by (\ref{ass_odd}). Then Theorem \ref{upcrossing} implies that $ V^{i',+}_t \dfnn
 \underset{n \to \infty}{\liminf}  V^{i'} \big(  \ul{\t}(0) \land
q^+_n(t) \big)$, $t \in [0,T]$ is an RCLL
 process and that
 \beas
 P \Big( V^{i',+}_t =
\underset{n \to \infty}{\lim}  V^{i'} \big(  \ul{\t}(0) \land
q^+_n(t) \big) \hb{ for any }
 t \in [0,T] \Big)=1.
 \eeas
For any $\si,\z \in \cS_{0,T}$, Lemma \ref{lem_ess} and
\eqref{eqn-k777} show that
  \beas
  \q \b1_{\{\si = \z\}} V(\si) =\neg \b1_{\{\si = \z\}} \ol{V}(\si)
  = \underset{j \in \cI}{\essinf}
  \left(\b1_{\{\si = \z\}}  R^j(\si) \right)\neg = \underset{j \in
  \cI}{\essinf}
  \left(\b1_{\{\si = \z\}}  R^j(\z) \right)\neg
    =\neg \b1_{\{\si = \z\}} \ol{V}(\z)= \neg \b1_{\{\si = \z\}} V(\z), \q  a.s.,
  \eeas
  which implies that
   \bea \label{eqn-k801}
   \b1_{\{\si = \z\}} V^{i'}(\si) = \b1_{\{\si = \z\}}  V(\si)+\b1_{\{\si = \z\}}
   H^{i'}_{\si}  = \b1_{\{\si = \z\}}  V(\z)+\b1_{\{\si = \z\}} H^{i'}_{\z}
    = \b1_{\{\si = \z\}} V^{i'}(\z), \q a.s.
   \eea
  Let $\si \in \cS^F_{0,T}$ take values
 in a finite set $\{t_1 < \cds< t_m\}$. For any $\a \in \{1 \cds m\}$ and $n \in \hN$,
 \if{0}
 one can deduce from
 \eqref{eqn-k801} that
  \beas 
  \b1_{\{\{\ul{\t}(0) \land q^+_n(\si)=\ul{\t}(0) \land q^+_n(t_\a)  \} \}}
V^{i'}\left(\ul{\t}(0) \land q^+_n(\si)\right) = \b1_{\{\ul{\t}(0)
\land q^+_n(\si)=\ul{\t}(0) \land q^+_n(t_\a)  \}}
V^{i'}\left(\ul{\t}(0) \land q^+_n(t_\a)\right), \q a.s.
  \eeas
  \fi
 since  $  \{\si=t_\a  \} \subset \{\ul{\t}(0) \land
q^+_n(\si)=\ul{\t}(0) \land q^+_n(t_\a)  \}$,
 one can deduce from \eqref{eqn-k801} that
 \beas
 \b1_{\{\si=t_\a  \}}
V^{i'}\left(\ul{\t}(0) \land q^+_n(\si)\right) = \b1_{\{\si=t_\a \}}
V^{i'}\left(\ul{\t}(0) \land q^+_n(t_\a)\right), \q a.s.
 \eeas
As $n \to \infty$,  \eqref{eqn-k905} shows that
  \beas
   \b1_{ \{\si=t_\a  \}}V^{i',+}_\si &=& \b1_{ \{\si=t_\a
   \}}V^{i',+}_{t_\a}
   =\underset{n \to \infty}{\lim} \b1_{\{\si=t_\a\}} V^{i'}\left(\ul{\t}(0) \land q^+_n(t_\a)\right)
 = \underset{n \to \infty}{\lim} \b1_{\{\si=t_\a\}}
 V^{i'}\left(\ul{\t}(0) \land q^+_n(\si)\right) \\
   &=& \b1_{\{\si=t_\a\}}   V^{i'}  \left( \ul{\t}(0) \land \si \right)  , \q a.s.
  \eeas
 Summing the above expression over $\a$, we obtain
  $  V^{i',+}_\si =V^{i'}  \left( \ul{\t}(0) \land \si \right)$, a.s.
 Then the right-continuity of the
process $V^{i',+}$ and \eqref{eqn-k905} imply that
 \bea \label{eqn-k807}
  V^{i',+}_\nu= \lmt{n \to \infty}V^{i',+}_{q^+_n(\nu)}
  = \lmt{n \to \infty} V^{i'}  \left(\ul{\t}(0) \land q^+_n(\nu) \right)
   =    V^{i'}  \left( \ul{\t}(0) \land \nu \right) , \q  a.s.
 \eea
 In particular, $ V^{i',+} $ is an RCLL modification of
  the stopped process $\left\{ V^{i'}\big(\ul{\t}(0) \land t\big) \right\}_{t \in [0,T]}
  $. Therefore, \\$ V^0  \dfnn \left\{ V^{i',+}_t-H^{i'}_{\ul{\t}(0) \land t} \right\}_{t \in [0,T]}$
  is an RCLL modification of the stopped value process $\left\{ V\big(\ul{\t}(0) \land t\big) \right\}_{t \in [0,T]} $.
 For any $\nu \in \cS_{0,T}$, \eqref{eqn-k807} implies that
 \beas
 V^0_\nu = V^{i',+}_\nu-H^{i'}_{\ul{\t}(0) \land \nu}
 =V^{i'}\big(\ul{\t}(0)  \land  \nu \big)-H^{i'}_{\ul{\t}(0) \land \nu}
 = V\big(\ul{\t}(0)  \land  \nu \big),~\; a.s.,~\; \hb{ proving (\ref{V0}).}
 \eeas

\ss \noindent \textbf{Proof of (2).}
   \eqref{V0} and Proposition \ref{V_Y_meet} imply that
 $V^0_{\ul{\t}(0)}=  V\big(\ul{\t}(0)  \big) = Y_{\ul{\t}(0) }$, a.s.
 Hence, we can deduce from the right-continuity of processes $V^0$ and
 $Y$ that $\t_V$ in (\ref{eqn-m15}) is a stopping time belonging to $\cS_{0,
 \ul{\t}(0)}$ and that
  \beas
 Y_{\t_V} =V^0_{\t_V} =  V(\t_V)    , \q  a.s.,
  \eeas
 where the second equality is due to \eqref{V0}.
 Then it follows from \eqref{eqn-l20} that for any $i \in \cI$
 \beas
 V(0) = V^i(0) \le \wt{\cE}_i \big[V^i(\t_V)\big]
 = \wt{\cE}_i \big[Y^i_{\t_V}\big]
 =      \cE_i \big[Y^i_{\t_V}\big].
 \eeas
 \if{0}
In light of Lemma \ref{essinf_lim} and (\ref{eqn-k51}), we can find
a sequence $\left\{j_n\right\}_{n \in \hN}$ in $\cI$ such that
  \beas
 V(\t_V)= \ol{V}(\t_V) = \underset{n \to \infty}{\lim} \dneg \da R^{j_n}(\t_V), \q a.s.
  \eeas
 Fix $i \in \cI$. For any $n \in \hN$, Definition \ref{def_stable_class} assures that there exists a
 $k_n=k(i, j_n, \t_V) \in \cI$ such that $\wt{\cE}_{k_n}=\cE^{\t_V}_{i, j_n}$.
Since $\t_V  \le \ul{\t}(0) $, a.s., we can deduce from
 \eqref{eqn-k25}, Lemma \ref{lem_trancate}, \eqref{h_tau_A}, \eqref{tau_ij2} as well as Proposition \ref{properties_3}
  (5) that
 \bea \label{eqn-gxg11}
 V(0) &=& \ol{V}(0)= \underset{l \in \cI}{\inf} R^l(0) \le R^{k_n}(0) = \wt{\cE}_{k_n} \big[ R^{k_n}(\t_V)+ H^{k_n}_{\t_V}\big]
  = \wt{\cE}_{k_n} \big[ R^{j_n}(\t_V)+ H^i_{\t_V}\big] \nonumber \\
 &=& \cE^{\t_V}_{i, j_n}\big[ R^{j_n}(\t_V)+ H^i_{\t_V} \big]
 =\wt{\cE}_i\Big[ \cE_{j_n}\big[ R^{j_n}(\t_V)+ H^i_{\t_V} \big|\cF_{\t_V}\big]\Big]
 =\wt{\cE}_i\big[  R^{j_n}(\t_V)+ H^i_{\t_V}\big].
 \eea
 Moreover, \eqref{Y_LB2} and \eqref{eqn-q104} imply that
 \beas
   C_* \le Y^i_{\t_V} = Y_{\t_V}+ H^i_{\t_V} \le  R^{j_n}(\t_V)+H^i_{\t_V} \le R^{j_1}(\t_V)+ H^i_{\t_V}
     , \q a.s.,
 \eeas
where $R^{j_1}(\t_V)+ H^i_{\t_V} \in Dom(\sE)$ thanks to Proposition
\ref{prop_Ri} (1), (S1') and (D2). Letting $n \to \infty$ in
\eqref{eqn-gxg11} and using the Dominated Convergence Theorem
(Proposition \ref{DCT2}) and \eqref{eqn-kk01} we obtain that
 \beas
 V(0) &\le &  \underset{n \to \infty}{\lim} \dneg \da \wt{\cE}_i\big[  R^{j_n}(\t_V)+ H^i_{\t_V}\big]
  = \wt{\cE}_i \big[ \underset{n \to \infty}{\lim} \dneg \da R^{j_n}(\t_V) + H^i_{\t_V}\big]
    =\wt{\cE}_i \big[ V(\t_V)+ H^i_{\t_V}\big]  
   = \wt{\cE}_i \big[ Y_{\t_V} + H^i_{\t_V}\big] \\
  & =& \wt{\cE}_i \big[ Y^i_{\t_V}\big] =  \cE_i \big[ Y^i_{\t_V}\big].
 \eeas
 \fi
Taking the infimum of the right-hand-side over $i \in \cI$ yields
that
 \beas
 V(0) \le  \underset{i \in \cI}{\inf}  \cE_i \big[ Y^i_{\t_V}\big]
   \le  \underset{\rho \in \cS_{0,T}}{\sup} \Big( \underset{i \in \cI}{\inf}
    \cE_i \big[ Y^i_{\rho}\big]\Big)=\ul{V}(0)=V(0),
 \eeas
 which implies that
  $  \underset{i \in \cI}{\inf} \cE_i \big[ Y^i_{\t_V}\big]
   = \underset{\rho \in \cS_{0,T}}{\sup} \, \underset{i \in \cI}{\inf}
    \cE_i \big[ Y^i_{\rho}\big]$.
  \qed

\subsection{Proofs of Section \ref{section_remark}}

\ss \no {\bf Proof of Proposition \ref{prop_olsE}:} Fix $t \in
[0,T]$. For any $\xi \in Dom(\sE)$ and $i \in \cI$, the definition
of $Dom(\sE)$ assures that there exists a $c(\xi) \in \hR$ such that
$c(\xi) \le \xi$, a.s. Then Proposition \ref{properties_3} (5) shows
that
  \bea \label{eqn-vxv03}
 c(\xi) = \wt{\cE}_i[c(\xi)|\cF_t] \le \wt{\cE}_i[ \xi |\cF_t] , \q a.s.
  \eea
Taking the essential infimum of the right-hand-side over $i \in
\cI$, we obtain for an arbitrary $i' \in \cI$ that
  \beas
 c(\xi) \le \ul{\sE}[\xi|\cF_t] \le \wt{\cE}_{i'}[ \xi |\cF_t], \q a.s.
 \eeas
 Since $\wt{\cE}_{i'}[ \xi |\cF_t]
 \in    Dom^\#(\cE_{i'})=Dom(\sE)$, Lemma
 \ref{lem_dom_sharp} implies that $\ul{\sE}[\xi|\cF_t] \in Dom(\sE)$, thus $ \ul{\sE}[\cd|\cF_t] $
 is a mapping from $Dom(\sE)$ to $Dom_t(\sE)=Dom(\sE)\cap
 L^0(\cF_t)$.

 \ms A simple application of Lemma \ref{lem_ess} shows that $
\ul{\sE} $ satisfies (A3), (A4) and  \eqref{eq: mono_ulsE}. Hence,
it only remains to show (A2) for $ \ul{\sE} $. Fix $0 \le s <t \le
T$. Letting $(\nu, \cI', \cU)=(t, \cI, \{T\})$ and taking $X(T) =
\xi $ in Lemma \ref{lem_02}, we can find a sequence $\{ i_n \}_{n
\in \hN}$ in $ \cI  $ such that
 \bea \label{eqn-vxv09}
 \ul{\sE}[\xi|\cF_t] =    \underset{i \in \cI}{\essinf}\,  \wt{\cE}_i  [\xi|\cF_t]
           =  \underset{n \to \infty}{\lim} \dneg \da    \wt{\cE}_{i_n} [\xi|\cF_t] , \q a.s.
 \eea
Now fix $ j \in \cI$. For any $n \in \hN$, it follows from
Definition \ref{def_stable_class} that there exists
  a $k_n =k(j,i_n,t) \in \cI$ such that  $\wt{\cE}_{k_n}=\cE^t_{j,i_n}$.
Applying (\ref{tau_ij2}) yields that
 \bea \label{eqn-vxv05}
     \ul{\sE}[\xi|\cF_s]   \le \wt{\cE}_{k_n}[\xi|\cF_s]
   = \cE^t_{j,i_n} [\xi|\cF_s]
   = \wt{\cE}_j\big[ \wt{\cE}_{i_n} \big[\xi \big|\cF_t \big] \big|\cF_s\big] , \q  a.s.
  \eea
  For any $n \in \hN$, we see from \eqref{eqn-vxv03} and \eqref{eqn-vxv09} that
  \beas
  c(\xi) = \wt{\cE}_{i_n}[c(\xi)|\cF_t]  \le \wt{\cE}_{i_n} \big[\xi \big|\cF_t \big]
  \le \wt{\cE}_{i_1} \big[\xi \big|\cF_t \big], \q
  a.s.,
  \eeas
  where $ \wt{\cE}_{i_1} \big[\xi \big|\cF_t \big] \in Dom^\#(\cE_{i_1})= Dom(\sE)$.
  The Dominated Convergence Theorem (Proposition \ref{DCT2}) and (\ref{eqn-vxv05}) then imply that
\beas
   \wt{\cE}_j \big[ \ul{\sE}[\xi|\cF_t] \big|\cF_s \big] =  \underset{n \to \infty}{\lim}
    \wt{\cE}_j\big[ \wt{\cE}_{i_n} [\xi|\cF_t] \big|\cF_s \big]
   \ge \ul{\sE}[\xi|\cF_s]
     , \q a.s.
 \eeas
Taking the essential infimum of the left-hand-side over $j \in \cI$,
we obtain
 \bea \label{eqn-q102xx}
 \ul{\sE} \big[ \ul{\sE}[\xi|\cF_t] \big|\cF_s \big]
   \ge \ul{\sE}[\xi|\cF_s]    , \q a.s.
 \eea

  On the other hand, for any $i \in \cI$ and $\rho \in \cS_{t,
 T}$,
applying Corollary \ref{cor_os2}, we obtain
 \beas
  \wt{\cE}_i  [ \xi  |\cF_s  ]=\wt{\cE}_i\big[\wt{\cE}_i [ \xi  |\cF_t ]\big|\cF_s \big]
 \ge \wt{\cE}_i [ \ul{\sE}[\xi|\cF_t] |\cF_s ] \ge
  \ul{\sE} \big[ \ul{\sE}[\xi|\cF_t] \big|\cF_s \big], \q a.s.
 \eeas
 Taking the essential infimum of the left-hand-side over $i  \in \cI$ yields that
  $ \ul{\sE}  [ \xi  |\cF_s  ] \ge
  \ul{\sE} \big[ \ul{\sE}[\xi|\cF_t] \big|\cF_s \big]$, a.s.,
which together with (\ref{eqn-q102xx}) proves (A2) for $\ul{\sE}$.
\qed

 \subsection{Proofs of Section \ref{ch_app}}

\ss \no {\bf Proof of Proposition \ref{gexp}:} By (\ref{g-cond}),
it holds $\dtp$ that for any $z \in \hR^d$
 \beas
 |g(t,z)|= |g(t,z)-g(t,0)|\le K_g |z|, \q \hb{thus} \q  \tilde{g}(t,z) \dfnn  -K_g|z| \le
 g(t,z) .
 \eeas
Clearly, $\tilde{g}$ is a generator satisfying (\ref{g-cond}). It is
also positively homogeneous in $z$, i.e.
 \beas
 \tilde{g}(t,\a z)= -K_g|\a z|=-\a K_g|z|= \a \tilde{g}(t, z),\q \fa \a \ge 0, ~\fa
 z \in \hR^d.
 \eeas
 Then Example 10 of \cite{Peng-97} (or Proposition 8 of \cite{Eman}) and (\ref{eqn-y11}) imply that for any $n
\in \hN$ and any $A \in \cF_T$ with $P(A)>0$
  \bea \label{eqn-y12}
  n \cE_{\tilde{g}}[ \b1_A]= \cE_{\tilde{g}}[n \b1_A] \le  \cE_g[n \b1_A]  .
  \eea
  Since $\cE_{\tilde{g}}[ \b1_A
  ]>0$ \big(which follows from the second part of (A1)\,\big),
  letting $n \to \infty$ in (\ref{eqn-y12}) yields (H0).

\ms Next, we consider a sequence $\{\xi_n\}_{n \in \hN} \subset
L^2(\cF_T)$ with $ \underset{n \in \hN}{\sup} |\xi_n | \in
L^2(\cF_T)$. If $\xi_n$ converges a.s., it is clear that $\xi \dfnn
\underset{n \to \infty}{\lim} \xi_n \in L^2(\cF_T)$. Applying Lemma
\ref{g-inequ} with $\mu=K_g$, we obtain
 \beas
  \big| \cE_g[\xi_n]- \cE_g[\xi] \big|
  &\le&  \cE_{g_\mu}\big[|\xi_n-\xi | \big] = \big\| \cE_{g_\mu}\big[|\xi_n-\xi | \big]\big\|_{L^2(\cF_T)}
 \le \Big\| \underset{t \in [0,T]}{\sup} \cE_{g_\mu}\big[|\xi_n-\xi |\, \big|\cF_t \big] \Big\|_{L^2(\cF_T)} \\
     &\le& 
     C e^{(K_g +K^2_g)T}  \|\xi_n-\xi
     \|_{L^2(\cF_T)},
 \eeas
 where we used the fact that $K_{g_\mu}=\mu $ in the last inequality. As $n \to \infty$,
 thanks to the Dominated Convergence Theorem of the linear expectation
 $E$,
 we have that $\|\xi_n-\xi \|^2_{L^2(\cF_T)}=E |\xi_n-\xi |^2 \to 0 $;
thus  $  \underset{n \to \infty}{\lim} \cE_g[\xi_n] = \cE_g[\xi] $.
 Then (H1) and (H2) follow.

  \ms For any $\nu \in \cS_{0,T}$ and $\xi
\in L^{2,+}(\cF_T) \dfnn \{\xi \in L^2(\cF_T): \xi \ge0, ~a.s.\}$,
 Lemma \ref{g-inequ} (1) shows that $ \underset{t \in [0,T]}{\sup}\big|
\cE_g[\xi |\cF_t]\big| \in L^{2,+}(\cF_T)$, consequently $ \cE_g[\xi
|\cF_\nu] \in  L^{2,+}(\cF_T)$.
 Since $X^\xi\dfnn \cE_g[\xi |\cF_\cd]$ is a
 continuous process,
   $  X^{\xi,+}_\nu= X^\xi_\nu = \cE_g[\xi |\cF_\nu] \in L^{2,+}(\cF_T)$,
  which proves (H3).  \qed

\ss \no {\bf Proof of Proposition \ref{tau_A_g}:}
 Fix $\nu \in \cS_{0,T}$. It is easy to check that the generator
 $g^\nu$ satisfies (\ref{g-cond}) with Lipschitz
coefficient $K_1 \vee K_2$. For any $\xi \in L^2(\cF_T)$, we set
$\eta \dfnn    \G^{\xi, g_2}_\nu  \in \cF_\nu$ and define
  \beas 
  \widetilde{\Th}_t &\dfnn & \b1_{\{\nu\le t\}}  \Th^{\xi, g_2}_t
     +  \b1_{\{\nu > t\}}
\Th^{\eta, g_1}_t, \q \fa t \in [0,T].
 \eeas
It follows that
 \beas
 g^\nu(t, \widetilde{\Th}_t)
 =\b1_{\{\nu\le t\} }\, g_2(t, \widetilde{\Th}_t)
  + \b1_{\{\nu > t\}} g_1(t, \widetilde{\Th}_t)
=   \b1_{\{ \nu \le t \} }\, g_2(t,\Th^{\xi, g_2}_t)+ \b1_{\{ \nu
>t \}} g_1(t, \Th^{\eta, g_1}_t)  , \q \fa t \in [0,T].
\label{eqn-l03}
 \eeas
For any $t \in [0,T]$, since $\{\nu\le t\} \in \cF_t$, one can
deduce that
 \bea
&& \hspace{-2cm} \b1_{\{\nu\le t\}} \bigg(\xi+ \int_t^T \neg
g^\nu(s, \widetilde{\Th}_s)ds - \int_t^T  \neg \widetilde{\Th}_s
 dB_s\bigg)=\b1_{\{\nu\le t\}} \xi+ \int_t^T \neg \b1_{\{\nu\le t\}} g^\nu(s, \widetilde{\Th}_s)ds
 - \int_t^T \neg \b1_{\{\nu\le t\}}\widetilde{\Th}_s
 dB_s \nonumber \\
 &=&\b1_{\{\nu\le t\}}  \xi+ \int_t^T \neg \b1_{\{\nu\le t\}} g_2(s,\Th^{\xi, g_2}_s)\,ds
 - \int_t^T \neg \b1_{\{\nu\le t\}} \Th^{\xi, g_2}_s dB_s
 \nonumber \\ &=& \b1_{\{\nu\le t\}} \bigg(\xi+
\int_t^T \neg g_2(s,\Th^{\xi, g_2}_s)\,ds - \int_t^T \neg \Th^{\xi,
g_2}_s dB_s\bigg) =\b1_{\{\nu\le t\} } \G^{\xi, g_2}_t , \q a.s.
\label{eqn-l01}
 \eea
The continuity of processes $\int_\cd^T \neg g^\nu(s,
\widetilde{\Th}_s)ds$, $\int_\cd^T \neg \widetilde{\Th}_s
 dB_s$ and $\G^{\xi, g_2}_\cd$ then implies that except on a null set $N$
 \beas
 \b1_{\{\nu\le t\}} \bigg(\xi+ \int_t^T \neg
g^\nu(s, \widetilde{\Th}_s)ds - \int_t^T  \neg \widetilde{\Th}_s
 dB_s\bigg) =\b1_{\{\nu\le t\} }\, \G^{\xi, g_2}_t, \q \fa t \in
 [0,T].
 \eeas
 Taking $t=\nu(\o)$ for any $\o \in N^c$ yields that
 \bea \label{eqn-wxw11}
 \xi+ \int_\nu^T g^\nu(s, \widetilde{\Th}_s)ds
- \int_\nu^T \widetilde{\Th}_s dB_s =   \G^{\xi, g_2}_\nu =\eta , \q
a.s.
 \eea
Now fix $t \in [0,T]$. We can deduce from \eqref{eqn-wxw11} that
 \bea \label{eqn-l04}
 \b1_{\{\nu > t\}} \Big(\xi+ \int_t^T g^\nu(s, \widetilde{\Th}_s)ds - \int_t^T \widetilde{\Th}_s
 dB_s\Big)&\dneg=& \dneg\b1_{\{\nu > t\}} \Big( \eta  + \int_t^\nu g^\nu(s, \widetilde{\Th}_s) ds
 - \int_t^\nu \widetilde{\Th}_s dB_s\Big) \nonumber \\
 & \dneg=& \dneg \b1_{\{\nu > t\}} \Big( \eta  + \int_t^\nu g_1(s, \Th^{\eta, g_1}_s) ds
 - \int_t^\nu \Th^{\eta, g_1}_s dB_s\Big), \q a.s. \qq
 \eea
 Moreover, Proposition \ref{properties_3} (5) implies that
 \beas
 \cE_{g_1}[\eta |\cF_{t \land \nu}] &=&  \eta+ \int_{t \land \nu}^T g_1(s, \Th^{\eta, g_1}_s) ds
 - \int_{t \land \nu}^T \Th^{\eta, g_1}_s dB_s= \cE_{g_1}[\eta |\cF_\nu]
  + \int_{t \land \nu}^\nu g_1(s, \Th^{\eta, g_1}_s) ds
 - \int_{t \land \nu}^\nu \Th^{\eta, g_1}_s dB_s \q \\
  &=&  \eta + \int_{t \land \nu}^\nu g_1(s, \Th^{\eta, g_1}_s) ds
  - \int_{t \land \nu}^\nu \Th^{\eta, g_1}_s dB_s, \q a.s.
 \eeas
 Multiplying both sides with $\b1_{\{\nu > t\}}$ and using
 (\ref{eqn-l04}), we obtain
  \beas
\b1_{\{\nu > t\}} \Big(\xi+ \int_t^T g^\nu(s, \widetilde{\Th}_s)ds -
\int_t^T \widetilde{\Th}_s dB_s\Big) = \b1_{\{\nu > t\}}\cE_{g_1}[
\eta| \cF_t]= \b1_{\{\nu > t\}}\cE_{g_1}[   \G^{\xi, g_2}_\nu  \big|
\cF_t], \q a.s.,
 \eeas
 which in conjunction with (\ref{eqn-l01}) shows that for any $t\in [0,T]$
 \beas
 \xi+ \int_t^T g^\nu(s, \widetilde{\Th}_s)ds - \int_t^T
\widetilde{\Th}_s dB_s &=& \b1_{\{\nu\le t\}  } \G^{\xi, g_2}_t+
\b1_{\{\nu > t\}}\cE_{g_1}[  \G^{\xi, g_2}_\nu \big| \cF_t] \\
&=&\b1_{\{\nu\le t\}} \cE_{g_2}[\xi |\cF_t]  + \b1_{\{\nu>
t\}}\cE_{g_1}\big[  \cE_{g_2} [\xi |\cF_\nu] \big|\cF_t\big] =
\cE^{\nu}_{g_1,g_2}[\xi  |\cF_t], \q a.s. \qq
 \eeas
 Since $\int_\cd ^T g^\nu(s, \widetilde{\Th}_s)ds$, $\int_\cd^T
\widetilde{\Th}_s dB_s$ and $\cE^{\nu}_{g_1,g_2}[\xi  |\cF_\cd]$ are
all continuous processes, it holds except a null $N'$ that
 \beas
 \cE^{\nu}_{g_1,g_2}[\xi  |\cF_t]= \xi+ \int_t^T g^\nu(s, \widetilde{\Th}_s)ds - \int_t^T
\widetilde{\Th}_s dB_s  , \q  \fa t \in [0,T].
 \eeas
One can easily show that $\big(\cE^{\nu}_{g_1,g_2}[\xi  |\cF_\cd],
\widetilde{\Th}\big) \in \hC^2_\bF([0,T])\times
\cH^2_\bF([0,T];\hR^d)$. Thus the pair is the unique solution to the
BSDE$(\xi, g^\nu)$, namely $ \cE_{g^\nu}[\xi  |\cF_t]
=\cE^{\nu}_{g_1,g_2}[\xi  |\cF_t]$ for any $ t \in [0,T]$.   \qed

\ss \no {\bf Proof of Theorem \ref{prop_g_result}:} We first note
that for any $g \in \sG'$, (\ref{eqn-g-neg}) implies that for every
$ \cE_g$-submartingale $X$,\; $-X$ is an $\cE_{g^-}$-supermartingale
although $g^-$ is concave (which means that $\cE_{g^-} $ may not
belong to
 $\sE'$). Hence, condition \eqref{ass_odd} is satisfied.

 \ss Fix $g \in \sG'$. Clearly $H^g_0=0$. For any $s,t  \in \cD_T$ with $s<t$,
 we can deduce from ($\tilde{h}$1) and ($\tilde{h}$2) that
 \bea \label{eqn-ttt01}
  \q C_{\sH'} \dfnn c'T \le \neg \int_s^t  \neg c' ds  \le \neg \int_s^t \neg h^g_r dr
    =  H^g_{s,t}   \le \int_s^t  h'(r) dr \le \int_0^T  h'(r) dr, \q
    a.s.,
 \eea
 which implies that
   \beas
   C_{\sH'}  \le  \underset{s,t \in \cD_T;
 s<t}{\essinf}\, H^g_{s,t}   \le  \underset{s,t \in \cD_T;
 s<t}{\esssup}\, H^g_{s,t} \le   \int_0^T h'(r) dr, \q a.s.,
 \eeas
 thus (S2) holds. Since $\int_0^T h'(r) dr \in L^2(\cF_T)$, it follows that
  \beas
  \underset{s,t \in \cD_T;
s<t}{\esssup}\, H^g_{s,t} \in  L^{2,\#}(\cF_T)\dfnn \big\{\xi \in
L^2 (\cF_T): \xi \ge c,~a.s. \hb{ for some }c \in \hR
\big\}=Dom(\sE').
 \eeas
   We can also deduce from
 \eqref{eqn-ttt01} that except on a null set $N$
  \beas
  \q  C_{\sH'}   \le   H^g_{s,t}   \le   \int_0^T h'(r) dr , \q  \fa \, 0\le s
   < t \le T.
  \eeas
Hence, for any $ \nu, \rho \in \cS_{0, T}$ with $\nu \le \rho$,
a.s., we have
  \beas
   C_{\sH'}   \le   H^g_{\nu, \rho}  \le \int_0^T h'(r) dr, \q  a.s.,
   \eeas
    which implies that
    $  H^g_{\nu, \rho} \in L^{2,\#}(\cF_T)=Dom(\sE')$; so we got (S1').
    Moreover, (S3) directly follows from ($\tilde{h}$3).

  \ms   Next, we check that the process $Y$ satisfies (Y1) and \eqref{eqn-m48}. 
  By \eqref{ass_Y_gexp} and (Y3), it holds a.s. that
$ C_Y \le Y_t \le \z'_Y$ for any $ t \in \cD_T$. The
right-continuity of the process $Y$ then implies that except on a
null set $\tilde{N}$
  \bea \label{eqn-wxw21}
    C_Y \le  Y_t  \le \z'_Y , \q \fa t \in [0,T], \q \hb{thus} \q   C_Y    \le  Y_\rho \le   \z'_Y,
    \q \fa \rho \in  \cS_{0,T}.
  \eea
 Since $\z'_Y \in L^2(\cF_T)$, it follows that $Y_\rho \in
 L^{2,\#}(\cF_T)=Dom(\sE')$ for any $\rho \in  \cS_{0,T}$,
 thus (Y1) holds.
Moreover, for any $g \in \sG'$, $\rho \in \cS_{0,T}$ and $t \in
\cD_T$,  Proposition \ref{3addition} (2), \eqref{eqn-wxw21} and
Lemma \ref{g-inequ} (2) show that
 \beas
 C_Y + c'T &  =&   \cE_g[C_Y + c'T |\cF_t] \le \cE_g \bigg[Y_\rho +  \int_0^\rho c' ds  \bigg|\cF_t\bigg]
 \le  \cE_g[Y^g_\rho  |\cF_t] \le \big|\cE_g[Y^g_\rho  |\cF_t]\big|
    \\ &=&   \big|\cE_g[Y^g_\rho |\cF_t] - \cE_g[0  |\cF_t] \big|
     \le     \cE_{g_{{}_M}}\big[| Y^g_\rho|  \big|\cF_t\big]
  \le \cE_{g_{{}_M}} \bigg[ \big| Y_\rho\big|+  \int_0^\rho  |h^g_s| ds
  \bigg|\cF_t\bigg]\\
  &\le&  \underset{t \in [0,T]}{\sup}\cE_{g_{{}_M}}\bigg[ \z'_Y \vee (-C_Y)
  +  \int_0^T h'(s) \vee (-c')\, ds  \bigg|\cF_t\bigg], \qq  a.s.
 \eeas
 Taking essential supremum of $\cE_g[Y^g_\rho \big|\cF_t]$ over $(g, \rho, t)
 \in \sG' \times \cS_{0,T} \times \cD_T  $, we can deduce from (A4) that
  \bea \label{eqn-wxw03}
  C_Y + c'T \le \underset{(g, \rho, t) \in \sG' \times \cS_{0,T} \times \cD_T  }{\esssup}\,
  \cE_g[Y^g_\rho \big|\cF_t]
  \le \underset{t \in [0,T]}{\sup}\cE_{g_{{}_M}}\bigg[\z'_Y + \neg \int_0^T
  h'(s)  ds  \bigg|\cF_t\bigg] -C_Y -c'T , \q a.s.
  \eea
 Lemma \ref{g-inequ} (1) implies that
 \beas
        \bigg\|\underset{t \in [0,T]}{\sup} \cE_{g_{{}_M}}\bigg[\z'_Y  \neg
 + \neg  \int_0^T \neg h'(s) ds  \bigg|\cF_t\bigg]
 \bigg\|_{L^2(\cF_T)}
    \le   C e^{(M+M^2 )T}  \bigg\|\,\z'_Y   \neg
 + \neg  \int_0^T \neg h'(s) ds  \bigg\|_{L^2(\cF_T)} < \infty.
 \eeas
 Hence, we see from \eqref{eqn-wxw03} that
 $\underset{(g, \rho, t) \in \sG' \times \cS_{0,T} \times \cD_T  }{\esssup}\,
  \cE_g[Y^g_\rho  |\cF_t] \in L^{2,\#}(\cF_T)=Dom(\sE') $, which is exactly \eqref{eqn-m48}.

\ms Now we show that the family of processes $\big\{Y^g_t  ,~ t\in
 [0,T]\big\}_{g \in \sG'}$ is both ``$\sE'$-uniformly-left-continuous" and
``$\sE'$-uniformly-right-continuous". For any $\nu, \rho \in
\cS_{0,T}$ with $\nu\le \rho  $, a.s.,
 let $\{\rho_n\}_{n \in \hN}\subset \cS_{\nu, T}$ be a sequence increasing a.s. to
 $\rho$.
 For any $g \in \sG'$, Lemma \ref{g-inequ} (2) implies that
  \beas
 \qq \qq && \hspace{-2.5cm} \q \Big|\cE_g  \big[\hb{$\frac{n}{n-1}$}
 Y_{\rho_n} +  H^g_{\rho_n}   \big|\cF_\nu \big]  -
\cE_g  \big[ Y^g_\rho  \big|\cF_\nu \big]\Big|  \le
\cE_{g_{{}_M}}\bigg[\Big|\hb{$\frac{n}{n-1}$} Y_{\rho_n} - Y_\rho -
 \int_{\rho_n}^\rho  h^g(s) ds  \Big|\bigg|\,\cF_\nu\bigg] \\
&  & 
 \le \cE_{g_{{}_M}}\bigg[\Big|\hb{$\frac{n}{n-1}$} Y_{\rho_n} -
 Y_\rho\Big|+  \int_{\rho_n}^\rho \wt{h}'(s) ds  \bigg|\,\cF_\nu\bigg], \q a.s.,
 \eeas
where $g_{{}_M}(z) \dfnn M|z|$, $ z \in \hR^d$ and $\wt{h}'(t) \dfnn
 h'(t) -c'$, $t \in [0,T]$.
 Taking essential supremum of the left hand side over $g \in \sG'$ yields that
  \bea \label{eqn-x11}
  \underset{g \in \sG'}{\esssup}\, \Big|\cE_g
  \big[\hb{$\frac{n}{n-1}$} Y_{\rho_n}+ H^g_{\rho_n}  \big|\cF_\nu \big]-
\cE_g\big[ Y^g_\rho  \big|\cF_\nu \big]\Big| \le
\cE_{g_{{}_M}}\bigg[\Big|\hb{$\frac{n}{n-1}$} Y_{\rho_n} -
 Y_\rho\Big|+  \int_{\rho_n}^\rho  \wt{h}'(s) \, ds  \bigg|\,\cF_\nu\bigg], \q a.s.
 \eea
 Moreover, Lemma \ref{g-inequ} (1) implies that
  \bea
 \qq \qq  && \hspace{-3cm} \bigg\| \cE_{g_{{}_M}}\neg \bigg[\Big|\hb{$\frac{n}{n-1}$} Y_{\rho_n}\neg
 - \neg Y_\rho\Big| \neg + \neg \int_{\rho_n}^\rho \neg \wt{h}'(s) ds  \bigg|\,\cF_\nu\bigg] \bigg\|_{L^2(\cF_T)}
   \neg \le \bigg\| \underset{t \in [0,T]}{\sup}\; \cE_{g_{{}_M}} \neg \bigg[\Big|\hb{$\frac{n}{n-1}$} Y_{\rho_n} \neg -
  \neg Y_\rho\Big| \neg + \neg  \int_{\rho_n}^\rho \neg \wt{h}'(s)  ds  \bigg|\,\cF_t\bigg]
    \bigg\|_{L^2(\cF_T)}  \nonumber \\
     &\le&  C e^{(M+M^2)T}  \bigg\|\Big|\hb{$\frac{n}{n-1}$} Y_{\rho_n} \neg -
Y_\rho\Big|+  \int_{\rho_n}^\rho \wt{h}'(s)  ds
\bigg\|_{L^2(\cF_T)}. \label{eqn-x12}
  \eea
Since
 \beas
  \big|\hb{$\frac{n}{n-1}$} Y_{\rho_n} \neg - Y_\rho\big|
 \le \hb{$\frac{n}{n-1}$} \big| Y_{\rho_n} - Y_\rho \big| +\hb{$\frac{1}{n-1}$} \big|Y_\rho  \big|
 \le 2 \big| Y_{\rho_n} - Y_\rho \big| +\hb{$\frac{1}{n-1}$} \big|Y_\rho
 \big|, \q \fa n \ge 2,
 \eeas
 the continuity of $Y$ implies that $ \underset{n \to \infty}{\lim}\Big(
\big|\hb{$\frac{n}{n-1}$} Y_{\rho_n} - Y_\rho\big|+
\int_{\rho_n}^\rho \wt{h}'(s) ds \Big)=0$, a.s. It also holds for
any $n \ge 2$ that
 \beas
 \big|\hb{$\frac{n}{n-1}$} Y_{\rho_n} -
Y_\rho\big|+ \int_{\rho_n}^\rho \wt{h}'(s) ds  \le 3 \big( \z'_Y
-C_Y \big) + \int_0^T  h'(s) ds -c'T , \q a.s.,
 \eeas
where the right-hand sides belongs to $L^2(\cF_T)$. Thus the
Dominated Convergence Theorem implies that
 \beas
 \hb{ the sequence }
  \bigg\{\Big|\hb{$\frac{n}{n-1}$} Y_{\rho_n} -
Y_\rho\Big|+   \int_{\rho_n}^\rho \wt{h}'(s) ds   \bigg\}_{n \in
\hN} \hb{ converges to $0$ in $L^2(\cF_T)$,}
 \eeas
 which together with (\ref{eqn-x11}) and (\ref{eqn-x12}) implies that
 \beas
 \hb{ the sequence }
 \Big\{ \underset{g \in \sG'}{\esssup}\, \big|\cE_g
  \big[\hb{$\frac{n}{n-1}$} Y_{\rho_n}+ H^g_{\rho_n}  \big|\cF_\nu \big]-
\cE_g\big[ Y^g_\rho  \big|\cF_\nu \big]\big| \Big\}_{n \in \hN}
 \hb{ also converges to $0$ in $L^2(\cF_T)$.}
 \eeas
  Then we can find a
subsequence  $\{n_k\}_{k \in \hN}$ of $\hN$ such that
 \beas
   \underset{n \to \infty}{\lim} \underset{g \in \sG'}{\esssup}\,
\big|\cE_g
  [\hb{$\frac{n_k}{n_k-1}$} Y_{\rho_{n_k}}+ H^g_{\rho_{n_k}}  |\cF_\nu ]-
\cE_g[ Y^g_\rho  |\cF_\nu ]\big|=0, \q a.s.
 \eeas
Therefore, the family of process $ \{Y^g  \}_{g \in \sG'}$ is
``$\sE'$-uniformly-left-continuous". The
``$\sE'$-uniform-right-continuity" of $ \{Y^g  \}_{g \in \sG'}$ can
be shown similarly.
 \qed

 \ss \no {\bf Proof of Theorem \ref{prop_optim_control}:}
 For any $U \in \fU$, Theorem \ref{prop_g_result} and Proposition \ref{Z_RCLL}
imply that $Z^{U,0}=\left\{Z^0_t + \int_0^t h^U_s ds \right\}_{t \in
[0,T]}$ is an $\cE_{g_{{}_U}}$-supermartingale.
 In light of the Doob-Meyer Decomposition of $g$-expectation
 (see e.g. \cite[Theorem 3.3]{Peng-99}, or \cite[Theorem 3.9]{Pln}),
 there exists an RCLL increasing process $\D^U$ null at
$0$ and a process $\Th^U \in \cH^2_\bF([0,T];\hR^d)$ such that
 \bea \label{eqn-wxw51}
 Z^{U,0}_t =Z^{U,0}_T+\int_t^T
 g_{{}_U}(s, \Th^U_s)ds+\D^U_T-\D^U_t-\int_t^T \Th^U_s dB_s, \q
 t\in\left[0,T\right].
 \eea

 In what follows we will show that
 \beas 
  U^*(t, \o) \dfnn u^*\big(t, \o, \Th^{U^0}_t(\o) \big), \q  (t, \o)
  \in [0,T] \times \O
 \eeas
  is an optimal control desired, where $U^0 \equiv 0$ denotes the null control.

 \ms Recall that $\ol{\t}(0)  = \inf\big\{t \in [0,T]
\;|\; Z^0_t=Y_t\big\}$. Taking $t=\ol{\t}(0)$ and $t=\ol{\t}(0)
\land t $ respectively in \eqref{eqn-wxw51} and
 subtracting the former from the latter yields that
 \bea \label{eqn-n15}
  Z^{U,0}_{ \ol{\t}(0) \land t }  =  Z^{U,0}_{\ol{\t}(0)}+\int_{\ol{\t}(0) \land
  t}^{\ol{\t}(0)} g_{{}_U}(s, \Th^U_s)ds+\D^U_{\ol{\t}(0)}-\D^U_{\ol{\t}(0) \land t}
  -\int_{\ol{\t}(0) \land t}^{\ol{\t}(0)} \Th^U_s dB_s, \q t\in\left[0,T\right],  
 \eea
which is equivalent to
  \bea \label{eqn-n225}
 Z^0_{ \ol{\t}(0) \land t }
 =   Z^0_{\ol{\t}(0)}  +\int_{\ol{\t}(0) \land t}^{\ol{\t}(0)}
 H(s, \Th^U_s, U_s) ds +\D^U_{\ol{\t}(0)}-\D^U_{\ol{\t}(0) \land t}
 -\int_{\ol{\t}(0) \land t}^{\ol{\t}(0)} \Th^U_s dB_s, \q  t\in\left[0,T\right].
  \eea
In particular, taking $U=U^0$, we obtain
 \bea \label{eqn-n226}
 Z^0_{ \ol{\t}(0) \land t }
 =   Z^0_{\ol{\t}(0)}  +\int_{\ol{\t}(0) \land t}^{\ol{\t}(0)}
 H(s, \Th^{U^0}_s, U^0_s) ds +\D^{U^0}_{\ol{\t}(0)}-\D^{U^0}_{\ol{\t}(0) \land t}
 -\int_{\ol{\t}(0) \land t}^{\ol{\t}(0)} \Th^{U^0}_s dB_s, \q  t\in\left[0,T\right].
  \eea
Comparing the martingale parts of (\ref{eqn-n225}) and
(\ref{eqn-n226}), we see that for any $U \in \fU$,
 \bea \label{eqn-n229}
 \Th^U_t= \Th^{U^0}_t , \q   \dtp
 \eea
 on the stochastic interval $ \[0, \ol{\t}(0)\] \dfnn \left\{(t, \o) \in \left[0,T\right]
\times \O :\, 0 \le t \le \ol{\t}(0) \right\}$. Plugging this back
into (\ref{eqn-n225}) yields that
 \bea \label{eqn-n227}
 Z^0_{ \ol{\t}(0) \land t }  = Z^0_{\ol{\t}(0)}  +\int_{\ol{\t}(0) \land t}^{\ol{\t}(0)}
 H(s, \Th^{U^0}_s, U_s) ds +\D^U_{\ol{\t}(0)}-\D^U_{\ol{\t}(0) \land t}
 -\int_{\ol{\t}(0) \land t}^{\ol{\t}(0)} \Th^{U^0}_s dB_s, \q  t\in\left[0,T\right].
  \eea
   \if{0}
  which implies that
   \beas
  \int_{\ol{\t}(0) \land t}^{\ol{\t}(0)}
 H(t, \o, \Th^{U^0}_s, U_s) ds +\D^U_{\ol{\t}(0)}-\D^U_{\ol{\t}(0) \land t}
 = \int_{\ol{\t}(0) \land t}^{\ol{\t}(0)}
 H(t, \o, \Th^{U^0}_s, U^*_s) ds +\D^{U^*}_{\ol{\t}(0)}-\D^{U^*}_{\ol{\t}(0) \land t}
   \eeas
   \fi

\ss Let us define $g_{{}_{K_o}}(z) \dfnn K_o|z|$, $ z \in \hR^d$.
Note that it is not necessary that $g_{{}_{K_o}}=g_U$ for some $U
\in \fU$. For any $U \in \fU$,  we set $ \G_t \dfnn
\cE_{g_{{}_U}}\left[ Z^{U,0}_{\ol{\t}(0)} \Big|\cF_t
 \right]$ and $\hat{\G}_t \dfnn \cE_{g_{{}_{K_o}}}\left[ - \D^{U^*}_{\ol{\t}(0)} \Big|\cF_t
 \right]$, $t\in \left[0,T\right]$, which are the solutions to
the BSDE$\big(Z^{U,0}_{\ol{\t}(0)},g_{{}_U}\big)$ and BSDE$\big(-
\D^{U^*}_{\ol{\t}(0)},g_{{}_{K_o}}\big)$ respectively, i.e.,
 \beas
 \G_t&=&Z^{U,0}_{\ol{\t}(0)}   +\int_t^T g_{{}_U}(s, \Th_s) ds-\int_t^T  \Th_s
 dB_s, \q t \in \left[0,T\right], \\
 \hat{\G}_t&=&- \D^{U^*}_{\ol{\t}(0)}  +\int_t^T K_o \big|\hat{\Th}_s \big| ds-\int_t^T  \hat{\Th}_s dB_s
 , \q t \in \left[0,T\right],
 \eeas
where $ \Th, \hat{\Th} \in \cH^2_\bF([0,T]; \hR^d)$.
 Applying Proposition \ref{properties_3} (5) and Corollary
 \ref{cor_os2}, we obtain that for any $t \in \left[0,T\right]$
  \beas
  \G_{\ol{\t}(0)}-\G_{\ol{\t}(0) \land t}&=&\cE_{g_{{}_U}}\left[ Z^{U,0}_{\ol{\t}(0)}
  \Big|\cF_{\ol{\t}(0)} \right]-\cE_{g_{{}_U}}\left[ Z^{U,0}_{\ol{\t}(0)} \Big|\cF_{\ol{\t}(0) \land t}
  \right] = Z^{U,0}_{\ol{\t}(0)} -\cE_{g_{{}_U}}\left[\cE_{g_{{}_U}}
  \left[ Z^{U,0}_{\ol{\t}(0)} \Big|\cF_{\ol{\t}(0) }  \right]\Big|\cF_t  \right]\\
 & = & Z^{U,0}_{\ol{\t}(0)} - \cE_{g_{{}_U}}\left[ Z^{U,0}_{\ol{\t}(0)}  \Big|\cF_t  \right]
  =Z^{U,0}_{\ol{\t}(0)} -\G_t,  \q a.s.
  \eeas
Then the continuity of processes $\G_\cd$ and $Z^{U,0}_\cd$ imply
that
 \beas
 && \hspace{-1.2cm}\G_t-Z^{U,0}_{\ol{\t}(0) \land t }  =   Z^{U,0}_{\ol{\t}(0)}-Z^{U,0}_{\ol{\t}(0) \land t }
 +\G_{\ol{\t}(0) \land t}-\G_{\ol{\t}(0)}
  =Z^{U,0}_{\ol{\t}(0)}-Z^{U,0}_{\ol{\t}(0) \land t }
   +\int_{\ol{\t}(0) \land t}^{\ol{\t}(0)} g_{{}_U}(s, \Th_s) ds
   -\int_{\ol{\t}(0) \land t}^{\ol{\t}(0)}  \Th_s dB_s  \\
  &=& \tneg Z^0_{\ol{\t}(0)}-Z^0_{\ol{\t}(0) \land t }+\int_{\ol{\t}(0) \land t}^{\ol{\t}(0)} H(s, \Th_s, U_s) ds
  -\int_{\ol{\t}(0) \land t}^{\ol{\t}(0)}  \Th_s dB_s\\
       &=& \tneg  -\D^{U^*}_{\ol{\t}(0)}+\D^{U^*}_{\ol{\t}(0) \land t}
  +\int_{\ol{\t}(0) \land t}^{\ol{\t}(0)} \left[ H(s, \Th_s, U_s)-H(s, \Th^{U^0}_s, U^*_s)\right] ds
  -\int_{\ol{\t}(0) \land t}^{\ol{\t}(0)}  (\Th_s - \Th^{U^0}_s) dB_s, \q   t \in  [0,T ], 
 \eeas
where we used (\ref{eqn-n227}) with $U=U^*$ in the last inequality.
Since it holds $\dtp$ that
  \beas
    && \hspace{-2cm} H(t, \Th_t , U_t )-H(t, \Th^{U^0}_t , U^*_t )
  =  H(t, \Th_t , U_t )-H\big(t,  \Th^{U^0}_t, u^*(t, \Th^{U^0}_t)\big)
  \le  H(t,  \Th_t , U_t  )-H(t, \Th^{U^0}_t, U_t)\\
 &=&  g^o(t,  \Th_t , U_t )-g^o(t, \Th^{U^0}_t, U_t)
 \le   \big| g^o(t, \Th_t, U_t)-g^o(t, \Th^{U^0}_t, U_t) \big| \le
 K_o\big| \Th_t - \Th^{U^0}_t \big|,
  \eeas
  the comparison Theorem for BSDEs (see e.g. \cite[Theorem 35.3]{Peng-97})
 implies that
   \beas
 \hat{\G}_t \ge \G_t-Z^{U,0}_{\ol{\t}(0) \land t }-\D^{U^*}_{\ol{\t}(0) \land
 t}, \q  t \in \left[0,T\right].
   \eeas
 In particular, when $t=0$, we can deduce from \eqref{eqn-h42} that
 \beas 
     \cE_{g_{{}_{K_o}}}\Big[ - \D^{U^*}_{\ol{\t}(0)} \Big]
 \ge \cE_{g_{{}_U}}\left[ Z^{U}\big(\ol{\t}(0)\big) \right] -Z(0).
  \eeas
 Taking supremum of the right hand side over $U \in \fU$ and applying
Theorem \ref{SN_exist} with $\nu=0$, we obtain
 \beas
    0 \ge  \cE_{g_{{}_{K_o}}}\left[ - \D^{U^*}_{\ol{\t}(0)} \right]
 \ge  \underset{U \in \fU}{\sup}  \cE_{g_{{}_U}}\left[ Z^{U}\big(\ol{\t}(0)\big)
 \right]-Z(0)=0,
 \eeas
 thus $ \cE_{g_{{}_{K_o}}}\left[ - \D^{U^*}_{\ol{\t}(0)} \right]
 =0$. The strict monotonicity of $g$-expectation (see e.g. \cite[Proposition 2.2(iii)]{CHMP}) then implies that
$ \D^{U^*}_{\ol{\t}(0)}=0$, a.s. Plugging it back to (\ref{eqn-n15})
and using (\ref{eqn-n229}), we obtain
 \bea \label{eqn-n130}
 Z^{U^*\neg,\,0}_{ \ol{\t}(0) \land t }  &=&  Z^{U^*\neg,\,0}_{\ol{\t}(0)}+\int_{\ol{\t}(0) \land
  t}^{\ol{\t}(0)} g_{{}_{U^*}}(s, \Th^{U^0}_s)ds -\int_{\ol{\t}(0) \land t}^{\ol{\t}(0)} \Th^{U^0}_s
  dB_s  \nonumber \\
 &=& Z^{U^*\neg,\,0}_{\ol{\t}(0)}+\int_t^T g_{{}_{U^*}}(s, \b1_{\left\{ s \le \ol{\t}(0)\right\}}\Th^{U^0}_s)ds
 -\int_t^T  \b1_{\left\{ s \le \ol{\t}(0)\right\}}\Th^{U^0}_s dB_s, \q
 t\in\left[0,T\right],
 \eea
which implies that
 $   \cE_{g_{{}_{U^*}}} \left[Z^{U^*\neg,\,0}_{\ol{\t}(0)}\Big|\cF_t\right]=
Z^{U^*\neg,\,0}_{ \ol{\t}(0) \land t }$,  $  \fa
t\in\left[0,T\right]$.
  Namely, $\left\{Z^{U^*\neg,\,0}_{ \ol{\t}(0) \land t } \right\}_{t \in \left[0,T\right]}
 $ is a $g_{{}_{U^*}}$-martingale. Eventually, letting $t=0$ in (\ref{eqn-n130}),
 we can deduce from \eqref{eqn-h42} and Theorem \ref{SN_exist} that
 \beas
   \q  \qq Z(0) 
 =Z^{U^*\neg,\,0}_0=
  \cE_{g_{{}_{U^*}}} \left[Z^{U^*\neg,\,0}_{\ol{\t}(0)} \right]
 =\cE_{g_{{}_{U^*}}} \left[Z \big(\ol{\t}(0)\big) +\int_0^{\ol{\t}(0)} h^{U^*}_s ds \right]
  =\cE_{g_{{}_{U^*}}} \bigg[Y_{\ol{\t}(0)}+\int_0^{\ol{\t}(0)} h^{U^*}_s ds
  \bigg]. \qq \qq    \hb{\qed}
 \eeas

  \no {\bf Proof of Proposition \ref{prop_Karatzas}:}
 Because of its linearity in $z$, the primary generator
 \bea \label{defn_g0}
 g^o(t, \o, z, u)\dfnn \big\lan
\si^{-1}(t,X(\o))f\big(t,X(\o), u  \big),z \big\ran \q \fa (t, \o,
z, u) \in [0,T] \times \O \times \hR^d \times S
 \eea
  satisfies ($g^o2$) and ($g^o4$). Then ($g^o1$) follows from
 the continuity of the process $\{X(t)\}_{t \in [0,T]}$
 as well as the measurability of the volatility $\si$ and of the function
 $f$. Moreover, \eqref{eqn-ttt25} and \eqref{eqn-ttt51} imply that for a.e. $t \in [0,T]$
 \beas
  |g^o(t, \o, z_1, u)- \neg  g^o(t, \o, z_2, u)|
 & \tneg =&    \tneg  \big| \lan \si^{-1} \neg (t,X(\o))f\big(t,X(\o), u  \big), z\neg-\neg z'\ran
 \big|   \le    \big\|  \si^{-1} \neg (t,X(\o))\big\|\cd \big| f\big(t,X(\o), u  \big) \big| \cd | z \neg-\neg z' |
   \\  &  \tneg\le &   \tneg  K^2 | z-z'  |,  \qq  \fa   z_1, z_2 \in  \hR^d, ~\fa (\o, u) \in \O \times  S,
 \eeas
 which shows that $g^o$ satisfies ($g^o4$) with $ K_o=K^2$. Clearly, $\wt{\fU} = \cH^0_\bF([0,T];
S)$ is closed under the pasting in the sense of \eqref{def_U_nu}.
Hence, we know from last section that $\{\cE_{g_{{}_U}} \}_{ U \in
\wt{\fU}} $ is a stable class of $g$-expectations, where
${g_{{}_U}}$ is defined in \eqref{defn_gU}.

\ms Fix $U \in \wt{\fU}$. For any $\xi \in
 L^2(\cF)$, we see from \eqref{def_g_exp} that
  \beas
   \cE_{g_{{}_U}}[\xi|\cF_t] &=&  \xi+\int_t^T g_U(s, \Th_s)ds-\int_t^T  \Th_s dB_s\\
      &=&  \xi+\int_t^T \big\lan
\si^{-1}(s,X)f\big(s,X,U_s \big), \Th_s \big\ran ds-\int_t^T \Th_s
dB_s   =  \xi -\int_t^T \Th_s dB^U_s,  \q t \in [0,T],
  \eeas
where $  B^U_t \dfnn B_t-\int_0^t \si^{-1}(s,X)f(s,X,U_s)\,ds$, $t
\in  [0,T]$ is a Brownian Motion with respect to $P_U$.
 For any $t \in [0,T]$, taking $E_U[\cd|\cF_t]$ on both sides above yields that
  \bea \label{eqn-ttt63}
  \cE_{g_{{}_U}}[\xi|\cF_t] = E_U\big[ \cE_{g_{{}_U}}[\xi|\cF_t] \big|\cF_t\big] =
     E_U[\xi|\cF_t]-E_U\bigg[\int_t^T \Th_s
     dB^U_s\bigg|\cF_t\bigg]=E_U[\xi|\cF_t], \q a.s.
  \eea
Hence the $g$-expectation $\cE_{g_{{}_U}}$ coincides with the linear
expectation $E_U$ on $L^2(\cF_T)$.

 \ms Clearly, the process $Y \dfnn \big\{ \vf(X(t)) \big\}_{t \in
 [0,T]}$ satisfies $(Y3)$ since $\vf $ is bounded from
 below by $-K$. We see from \eqref{eqn-ttt55} that for any $t \in
[0,T]$
     \beas
     Y_t = \vf(X(t)) \le K| X(t ) | \le K \|X \|^*_T.
    \eeas
 Taking essential supremum of $Y_t$ over $t \in \cD_T$ yields that
 \bea \label{eqn-ttt61}
    \z'_Y \dfnn \Big( \underset{t  \in \cD_T }{\esssup}\, Y_t \Big)^+ \le  K \|X
\|^*_T, \q a.s.
 \eea
For any $t \in [0,T]$,   the Burkholder-Davis-Gundy inequality, ($\si$1),
\eqref{eqn-ttt25} as well as Fubini Theorem imply that
 \beas
   E\Big[\big( \| X \|^*_t\big)^2 \Big]&=& E\Big[  \underset{s \in [0,t]}{\sup}  | X(s)|^2 \Big]
 \le 2x^2+2 E \bigg\{  \underset{s\in[0,t]}{\sup} \Big|
  \int_0^s \si(r,X) dB_r \Big|^2 \bigg\} \le   2x^2+  2C  E\int_0^t |\si(s,X)|^2 ds \\
 &\le &   2x^2+  4C   \int_0^t |\si(s,\vec{0})|^2 ds+4C  E\int_0^t |\si(s,X)-\si(s,\vec{0})|^2 ds    \\
 &\le&  2x^2+  4C   \int_0^T |\si(s,\vec{0})|^2 ds+4C n^2K^2  \int_0^t E\Big[\big(\|X\|^*_s\big)^2\Big]
ds.
 \eeas
 Then applying Gronwall's inequality yields that
 \bea \label{eqn-new120}
  E\Big[\big( \| X \|^*_T\big)^2 \Big] \le \bigg(2x^2+  4C   \int_0^T |\si(s,\vec{0})|^2
ds \bigg)e^{4C  n^2K^2 T} < \infty,
 \eea
 which together with
\eqref{eqn-ttt61} shows that $\z'_Y \in L^2(\cF_T)$, proving
\eqref{ass_Y_gexp}.

 \ms Next, we define a function $ h^o(t, \o, u) \dfnn  h(t, X(\o) , u)
  $, $\fa (t, \o, u) \in [0,T] \times \O \times S$.
 The continuity of the process $\{X(t)\}_{t \in [0,T]}$
 and the measurability of the function $h$ imply that $h^o $
 is $\sP  \otimes \fS/ \sB(\hR)$-measurable. We see from
\eqref{eqn-ttt51} that $h^o$ satisfies ($\hat{h}$1).
 It also follows from \eqref{eqn-ttt51} that for a.e. $t \in [0,T]$ and for any $\o \in \O$,
  \beas
    h^U_t(\o) \dfnn  h^o(t, \o, U_t(\o))=  h(t, X(\o) , U_t(\o)) \le K\| X(\o)
\|^*_T, \q  \fa U \in \wt{\fU}.
  \eeas
 Taking essential supremum of $h^U_t(\o)$ over $U \in \wt{\fU} $ with respect to the product measure
 space $([0,T]\times \O,\, \sP, 
 \l \times P )$ yields that
 \beas
    \hat{h}(t,\o) \dfnn \Big( \underset{U \in \wt{\fU}}{\esssup}\, h^U_t(\o)\Big)^+ \le K\| X(\o) \|^*_T, \q
    \dtp,
 \eeas
 which leads to that
  $  \int_0^T  \hat{h}(t,\o) dt \le K T    \| X(\o) \|^*_T  $ ,
a.s. Hence, \eqref{eqn-new120} implies that $ \int_0^T
 \hat{h}(t,\o) dt \in L^2(\cF_T) $, proving ($\hat{h}$2) for $h^o$.

\ms We can apply the optimal stopping theory developed in Section
\ref{co_game} to the triple $\big(\{\cE_{g_{{}_U}} \}_{ U \in
\wt{\fU}}, \{ h^U \}_{U \in \wt{\fU}}, Y \big)$ and use
\eqref{eqn-ttt63} to obtain \eqref{eqn-ttt43}.
 In addition, if there exists a measurable mapping $u^*: [0,T]
\times \O \times \hR^d \mapsto S$ satisfying \eqref{eqn-ttt45}, then
\eqref{defn_g0} indicates that for any $(t, \o, z) \in [0,T] \times
\O \times \hR^d$
 \beas
 \underset{u \in S}{\sup} \Big(g^o(t,\o, z, u)+h^o(t,\o,u) \Big)
&=&\underset{u \in S}{\sup} \, \wt{H}(t, X(\o), z, u)
 = \wt{H}\big(t, X(\o), z, u^*(t, X(\o), z) \big) \\&=& \tneg g^o (t,\o, z, u^*(t, X(\o), z))
 +h^o (t, \o, u^*(t, X(\o),
z)),
 \eeas
which shows that \eqref{eqn-ttt71} holds for the mapping
$\wt{u}^*(t, \o, z)= u^*(t, X(\o), z)$, $(t, \o, z) \in [0,T] \times
\O \times
 \hR^d$. Therefore, an application of Theorem
\ref{prop_optim_control} yields \eqref{eqn-ttt47} for some $U^* \in
\wt{\fU}$. \qed

\ss \no {\bf Proof of Proposition \ref{thm_strict_com}:}
\eqref{scomp1} directly follows from \cite[Theorem 5]{BH-07}.
   To see (\ref{scomp2}), we set $\D \G \dfnn \G^{\xi_1, \hat{g}}-\G^{\xi_2, \hat{g}}  $ and $\D \Th
 \dfnn \Th^{\xi_1, \hat{g}}-\Th^{\xi_2, \hat{g}} $, then \eqref{g2-cond}(i)
 implies that
 \beas
  d \D \G_t &=&  -\big(\hat{g}(t, \Th^{^{\xi_1, \hat{g}}}_t)-\hat{g}(t, \Th^{^{\xi_2, \hat{g}}}_t)\big)dt +\D \Th_t d B_t
  =  -\int_0^1   \frac{\pa \hat{g}}{\pa z}(t, \l\D \Th_t + \Th_t^{^{\xi_2, \hat{g}}})\D \Th_t   d\l \,d t  +\D \Th_t d
  B_t \\
  &=&  \D \Th_t \big(-a_tdt+ d B_t\big),  \q \;\, t \in [0,T],
 \eeas
 where $  a_t \dfnn \int_0^1 \frac{\pa \hat{g}}{\pa z}(\l\D \Th_t +
 \Th^{\xi_2, \hat{g}}_t)\,
d \l $, $t \in [0,T]$. Since $M_\bF([0,T];\hR^d) \subset
 M^2_\bF([0,T];\hR^d)=\cH^2_\bF([0,T];\hR^d)$, one can deduce from \eqref{g2-cond}(ii)
that
 \beas
E \int_0^T |a_t|^2 dt &\leq& E \int_0^T \neg \int_0^1
\bigg|\frac{\pa \hat{g}}{\pa z}(\l\D \Th_t + \Th^{\xi_2,
\hat{g}}_t)\bigg|^2 d \l \, dt
\leq  2 \k^2 T + 2 \k^2 E \int_0^T \neg \int_0^1    |\l \Th^{\xi_1, \hat{g}}_t + (1-\l) \Th^{\xi_2, \hat{g}}_t|^2  d \l \, dt \\
&\leq & 2 \k^2 T + \frac{4}{3} \k^2 E  \int_0^T \Big(\big|
\Th^{\xi_1, \hat{g}}_t\big|^2+ \big| \Th^{\xi_2, \hat{g}}_t
\big|^2\Big)  dt <\infty. \eeas Moreover, Doob's martingale inequality shows
that
 \bea
 \label{bddab}
 E \bigg[  \underset{t\in[0,T]}{\sup} \Big|
 \int_0^t a_sd B_s \Big|^2 \bigg] \le  4 E \bigg[  \Big|
 \int_0^T a_s d B_s \Big|^2 \bigg]\neg =  4  E\int_0^T |a_t|^2 dt  < \infty.
 \eea
  Thus, we can define process $Q_t\dfnn \exp \Big\{ \neg
-\frac{1}{2}\int_0^t |a_s|^2 d s +\int_0^t a_sd B_s\Big\}$, $t \in
[0,T]$ as well as stopping times
 \beas
 \nu_n \dfnn \inf\big\{ t \in [\nu,T]: Q_t\vee |\D \G_t |>n
\big\}\land T, \q ~\;  \fa  n \in \hN.
 \eeas
 It is clear that $\underset{n \to \infty}{\lim} \neg \ua \nu_n= T$, a.s., and (\ref{bddab}) assures that there exists
 a null set $ N$ such that for any $\o\in N^c$, $T=\nu_m (\o)$ for
some $m =m(\o)\in \hN$.

 \ms For any $n \in \hN$, integrating by parts on $[\nu,\nu_n]$ yields
that
 \beas
 Q_{\nu_n}\D \G_{\nu_n} &\tneg=\neg\dneg& Q_\nu \D \G_\nu
- \int_\nu^{\nu_n}  Q_t \D \Th_t a_t dt   + \neg \int_\nu^{\nu_n}
\neg Q_t \D \Th_t d B_t  +\neg \int_\nu^{\nu_n}\neg \D \G_t Q_t a_t
d B_t
  +\neg \int_\nu^{\nu_n}\neg   Q_t \D \Th_t a_t d t \\
 &\tneg=\tneg&   \int_\nu^{\nu_n} \dneg \big(Q_t \D \Th_t   +   \D \G_t Q_t a_t \big) d B_t.
 \eeas
which implies that $ E\big[Q_{\nu_n}\D \G_{\nu_n} \big] =0$. Thus we
can find a null set $N_n$ such that       $     \D \G_{\nu_n(\o)}(\o)=0$, $ \fa  \o \in N^c_n $.
Eventually, for any $\o \in \Big\{ N \cup \big( \underset{n \in
\hN}{\cup} N_n \big)\Big\}^c$, we have
 \beas
  \hspace{3.2cm} \xi^1 (\o)= \G^{\xi_1, \hat{g}}_T (\o)
  =  \underset{n \to \infty}{\lim} \G^{\xi_1, \hat{g}}_{\nu_n(\o)}(\o)
  =  \underset{n \to \infty}{\lim} \G^{\xi_2, \hat{g}}_{\nu_n(\o)}(\o) =
  \G^{\xi_2, \hat{g}}_T (\o)=\xi^2 (\o). \hspace{3.1cm} \hb{\qed}
\eeas

\ss \no {\bf Proof of Proposition \ref{g2exp}: } 
 Let $\{A_n\}_{n \in \hN}$ be any sequence in $\cF_T$ such that
  $\underset{n \to \infty}{\lim} \neg \da \b1_{A_n}=0 $,
 a.s. For any $\xi, \eta \in L^{\neg e,+}(\cF_T)
 \dfnn \{\xi \in L^{\neg e}(\cF_T) : \xi \ge 0, ~a.s.\}$, since
 \beas
    E\big[e^{\l |\xi|  }\big]< \infty \q \hb{and since} \q
 \underset{n \in \hN}{\sup}E\Big[e^{\l |\xi+\b1_{A_n}\eta|
 }\Big] \le   E\Big[e^{\l |\xi|}e^{\l |\eta|
 }\Big]\le \frac12 E\big[ e^{2\l | \xi | }\big]+\frac12 E\big[ e^{2\l | \eta |
 }\big] < \infty
 \eeas
 holds for each $\l > 0$, Lemma \ref{lem_g2_stable} implies that
  \beas
 0= \underset{n \to \infty}{\lim}  E\bigg[   \underset{t \in [0,T]}{\sup}
 \Big|\cE_{\hat{g}}[\xi+\b1_{A_n}\eta|\cF_t]-\cE_{\hat{g}}[\xi|\cF_t]\Big|
 \bigg]\ge
 \underset{n \to \infty}{\lim}
 \big|\cE_{\hat{g}}[\xi+\b1_{A_n}\eta ]-\cE_{\hat{g}}[\xi ]\big|  \ge 0,
  \eeas
 thus $\cE_{\hat{g}}$ satisfies (H2). Similarly, we can show that (H1) also holds for $\cE_{\hat{g}}$.

  \ms Moreover, for any $\nu \in \cS_{0,T}$ and $\xi
\in  L^{\neg e,+}(\cF_T)$,
 since the process $\G^{\xi,\hat{g}}$ belongs to $ \hC^e_\bF([0,T])$, one can deduce that
  $ \cE_{\hat{g}}[\xi |\cF_\nu]=\G^{\xi,\hat{g} }_\nu \in L^{\neg
e,+}(\cF_T)$.
 Then the continuity of the process $X^\xi\dfnn \cE_{\hat{g}}[\xi
 |\cF_\cd]$ implies that
  $  X^{\xi,+}_\nu= X^\xi_\nu= \cE_{\hat{g}}[\xi |\cF_\nu] \in L^{\neg e,+}(\cF_T)$,
  which proves (H3).  \qed

\ss \no {\bf Proof of Theorem \ref{prop_g2exp_result}:} This
  proof is just an application of the optimal stopping theory
developed in Section \ref{co_game} to the singleton
$\{\cE_{\hat{g}}\}$. 
Hence, it suffices to check that $Y$ satisfies (Y1), (Y2) and
\eqref{eqn-q170}.

 \ms  Similar to \eqref{eqn-wxw21}, it holds except on a null set $N$ that
  \bea \label{eqn-wxw21_b}
   C_Y \le Y_t \le \hat{\z}_Y , \q \fa t \in [0,T], \q \hb{thus} \q   C_Y    \le  Y_\rho \le  \hat{\z}_Y,
    \q \fa \rho \in  \cS_{0,T}.
  \eea
 Since $\hat{\z}_Y \in L^{\neg e}(\cF_T)$, it holds for any $\rho \in
 \cS_{0,T}$ that
  \bea \label{eqn-ttt20}
   E\big[e^{\l |Y_\rho| }\big] \le E\Big[e^{\l  (\hat{\z}_Y -C_Y )
 }\Big] =  e^{-\l C_Y }E\big[e^{\l  \hat{\z}_Y    }\big] <\infty , \q \fa  \l >
 0,
 \eea
  which implies that $Y_\rho \in
 L^{\neg e,\#}(\cF_T)=Dom\big( \{\cE_{\hat{g}}\}\big)$. Hence (Y1) holds.

 \ms Next, for any $\rho \in \cS_{0,T}$ and $t \in \cD_T$, Proposition
\ref{3addition} (2), \eqref{eqn-wxw21_b} show that
 \beas
 C_Y  =   \cE_{\hat{g}}[C_Y   |\cF_t] \le \cE_{\hat{g}}  [Y_\rho  |\cF_t ]
 \le  \cE_{\hat{g}} [ \hat{\z}_Y  |\cF_t ]=\G^{ \hat{\z}_Y , \hat{g} }_t
 \le  \underset{t \in [0,T]}{\sup}\big|\G^{ \hat{\z}_Y , \hat{g}
 }_t\big| , \q  a.s.
 \eeas
 Taking essential supremum of $\cE_{\hat{g}}[Y_\rho \big|\cF_t]$ over $(\rho, t)
 \in  \cS_{0,T} \times \cD_T  $ yields that
  \beas
  C_Y  \le \underset{( \rho, t) \in \cS_{0,T} \times \cD_T  }{\esssup}\,
  \cE_{\hat{g}}[Y_\rho \big|\cF_t]
  \le \underset{t \in [0,T]}{\sup}\big|\G^{ \hat{\z}_Y , \hat{g} }_t\big| , \q a.s.
  \eeas
  Since $\G^{ \hat{\z}_Y , \hat{g} } \in \hC^e_\bF([0,T])$,
  or equivalently $\underset{t \in [0,T]}{\sup}\big|\G^{
\hat{\z}_Y , \hat{g} }_t\big| \in L^{\neg e}(\cF_T) $, we can deduce
that $\underset{( \rho, t) \in \cS_{0,T} \times \cD_T  }{\esssup}\,
  \cE_{\hat{g}}[Y_\rho \big|\cF_t] \in L^{\neg e,\#}(\cF_T)=Dom\big( \{\cE_{\hat{g}}\}\big)
  $, which together with Remark \ref{rem_fatou2} (2) proves (Y2).

 \ms Moreover, for any $\nu, \rho \in \cS_{0,T}$ with $\nu \le \rho$,
 a.s. and any sequence $\{\rho_n\}_{n \in \hN} \subset \cS_{\nu, T}$
 increasing a.s. to $\rho$, the continuity of the process
 $Y$ implies that $\frac{n}{n-1} Y_{\rho_n} $ converges to $Y_\rho $ a.s.
 By \eqref{eqn-wxw21_b}, one can deduce that
 \beas
 \underset{n \in \hN}{\sup}E\Big[\exp\big\{\l |\hb{$\frac{n}{n-1}$} Y_{\rho_n}
  \big|\big\} \Big]  \le  \underset{n \in \hN}{\sup}  E\Big[e^{2\l  |  Y_{\rho_n} | }
  \Big]\le  E\Big[e^{2\l  ( \hat{\z}_Y -C_Y )
 }\Big] =  e^{-2\l C_Y }E\Big[e^{2\l  \hat{\z}_Y    }\Big]   < \infty, \q   \fa \l > 0 ,
 \eeas
 which together with \eqref{eqn-ttt20} allows us to apply Lemma
 \ref{lem_g2_stable}:
  \beas
 0= \underset{n \to \infty}{\lim}  E\bigg[  \underset{t \in [0,T]}{\sup}
 \Big|\cE_{\hat{g}}[\hb{$\frac{n}{n-1}$} Y_{\rho_n}
    |\cF_t]-\cE_{\hat{g}}[Y_\rho |\cF_t]\Big| \bigg]\ge
 \underset{n \to \infty}{\lim}
 E\bigg[
 \Big|\cE_{\hat{g}}[\hb{$\frac{n}{n-1}$} Y_{\rho_n}
    |\cF_\nu]-\cE_{\hat{g}}[Y_\rho |\cF_\nu]\Big|
 \bigg]  \ge 0,
  \eeas
 thus $ \underset{n \to \infty}{\lim}
 E\bigg[  \Big|\cE_{\hat{g}}[\hb{$\frac{n}{n-1}$} Y_{\rho_n}
    |\cF_\nu]-\cE_{\hat{g}}[Y_\rho |\cF_\nu]\Big|
 \bigg]=0$. Then we can find a
subsequence $\{n_k\}_{k \in \hN}$ of $\hN$ such that
  \beas
   \underset{n \to \infty}{\lim}
   \Big|\cE_{\hat{g}}[\hb{$\frac{n_k}{n_k-1}$} Y_{\rho_{n_k}}
    |\cF_\nu]-\cE_{\hat{g}}[Y_\rho |\cF_\nu]\Big|  =0, \q  a.s.,
    \eeas
  proving \eqref{eqn-q170} for $Y$. \qed

\bibliographystyle{abbrvnat}
\bibliography{erhan_song}

\end{document}